\documentclass[a4paper, 10pt]{article}

\usepackage{header2}

\usepackage{permutation}
\usepackage{multirow}

\usepackage[font={small}]{caption}

\begin{document}

\title
{\makebox[0pt][c]{A combinatorial approach to Rauzy-type dynamics}\\
\makebox[0pt][c]{I:\ permutations and the
  Kontsevich--Zorich--Boissy}\\
classification theorem}

\author{Quentin De Mourgues\\
  \multicolumn{1}{p{.7\textwidth}}{\centering\emph{\small 
\rule{0pt}{14pt}%
Institut Fourier, Universit\'e Grenoble Alpes\\
100, rue des Maths, 38610 Gières, France\\
 and \\
LIPN, Universit\'e Paris 13\\
99, av.~J.-B.~Cl\'ement, 93430 Villetaneuse, France\\
{\tt quentin.de-mourgues@univ-grenoble-alpes.fr}}}
\\
\rule{0pt}{16pt}%
Andrea Sportiello\\
  \multicolumn{1}{p{.7\textwidth}}{\centering\emph{\small 
\rule{0pt}{14pt}%
CNRS, and LIPN,
Universit\'e Paris 13\\
99, av.~J.-B.~Cl\'ement, 93430 Villetaneuse, France\\
{\tt andrea.sportiello@lipn.univ-paris13.fr}}}
}

\date{\small 3rd May 2017}

\maketitle


\begin{center}
\begin{minipage}{.9\textwidth}
\noindent
\small
{\bf Abstract.}
Rauzy-type dynamics are group actions on a collection of combinatorial
objects. The first and best known example concerns an action on
permutations, associated to interval exchange transformations (IET)
for the Poincar\'e map on compact orientable translation surfaces. The
equivalence classes on the objects induced by the group action are
related to components of the moduli spaces of Abelian differentials
with prescribed singularities, and, in two variants of the problem,
have been classified by Kontsevich and Zorich, and by Boissy, through
methods involving both combinatorics and algebraic geometry.

\qquad We provide here a purely combinatorial proof of both
classification theorems, and in passing establish a few previously
unnoticed features.  As will be shown elsewhere, our methods extend
also to other Rauzy-type dynamics, both on labeled and unlabeled
structures. Some of these dynamics have a geometrical interpretation
(e.g., matchings, related to IET on non-orientable surfaces), while
some others do not have one so far.
\end{minipage}
\end{center}

\newpage

\tableofcontents

\newpage

\section{Algebraic setting}
\label{sec.algsett}

\subsection{Permutational diagram monoids and groups}

Let $X_I$ be a set of labeled combinatorial objects, with elements
labeled from the set $I$.\footnote{In the most general context $I$ may
  be a \emph{multi-set}, i.e.\ allow for repetitions. However, at the
  purposes of the present paper, we will mainly concentrate on the
  case in which each label is repeated exactly once.}  We use the
shortcut $[n]=\{1,2,\ldots,n\}$, and
$X_n \equiv X_{[n]}$. The symmetric group $\kS_n$ acts naturally on
$X_n$, by producing the object with permuted labels.\footnote{For $I$
  a multiset, the natural action would be the one of $\kS_{|I|}/\Aut(I)$.}

Vertex-labeled graphs (or digraphs, or hypergraphs) are a typical
example.  Extra structure may be added, e.g.\ in hypergraphs we can
complement hyperedges with a cyclic ordering of the incident vertices,
and in both graphs and hypergraphs we can specify a cyclic ordering of
the incident edges at each vertex. This provides an embedding of the
abstract graph on a surface.

Set partitions are a special case of hypergraphs (all vertices have
degree 1).
Matchings are a special case of partitions, in which all blocks have
size 2.
Permutations $\s$ are a special case of matchings, in which each block
$\{i,j\}$ has $i \in \{1,2,\ldots,n\}$ and $j = \s(i)+n \in
\{n+1,n+2,\ldots,2n\}$.

We will consider dynamics over spaces of this type, generated by
operators of a special form that we now introduce:
\begin{definition}[Monoid and group operators]
\label{def.monoidoperator}
We say that $A$ is a \emph{monoid operator} on set $X_n$, if, for the
datum of a finite set $Y_n$, a map $a: X_n \to Y_n$, and a map
$\aa: Y_n \to \kS_n$, it consists of the map on $X_n$ defined by
\be
A(x) = \aa_{a(x)} x
\ef,
\ee
where the action $\aa x$ is in the sense of the symmetric-group
action over $X_n$.  We say that $A$ is a \emph{group operator} if,
furthermore, $a(A(x))=a(x)$.
\end{definition}
\noindent
Said informally, the function $a$ ``poses a question'' to the
structure $x$. The possible answers are listed in the set $Y$. For
each answer, there is a different permutation, by which we act on $x$.
Actually, as anticipated, we only use $a$ and $\a$ in the combination
$A=\a_{a(\cdot)}$, so that the use of two symbols for the single
function $A$ is redundant. This choice is done for clarity in our
applications, where the notation allows to stress that $Y_n$ has a
much smaller cardinality than $X_n$ and $\kS_n$, i.e.\ very few
`answers' are possible. In our main application, $|X_n|=|\kS_n|=n!$
while $|Y_n|=n$. The asymptotic behaviour is similar ($|X_n|$ at least
exponential, $|Y_n|$ at most linear) in all of our applications.

Clearly we have:


\begin{proposition}
\label{prop.Groupoperatorsareinvertible}
Group operators are invertible.
\end{proposition}

\pf The property $a(A(x))=a(x)$ implies that, for all $k \in \bN$,
$A^k(x) = (\aa_{a(x)})^k x$. Thus, for all $x$ there exists an integer
$d_A(x) \in \bN^+$ such that $A^{d_A(x)}(x)=x$. More precisely,
$d_A(x)$ is the l.c.m.\ of the cycle-lengths of $\aa_{a(x)}$. Call
$d_A = \lcm_{x \in X_n} d_A(x)$. Then 
$d_A$ is a finite integer,
and we can pose 
$A^{-1} = A^{d_A-1}$. The reasonings above shows that $A$ is a
bijection on $X_n$, and $A^{-1}$ is its inverse.
\qed

\begin{definition}[monoid and group dynamics]
We call a \emph{monoid dynamics} the datum of a family of spaces
$\{X_n\}_{n \in \bN}$ as above, and a finite collection $\cA=\{A_{i}\}$ of
monoid operators.
We call a \emph{group dynamics} the analogous structure, in which all
$A_i$'s are group operators.

For a monoid dynamics on the datum $S_n = (X_n, \cA)$, we say that 
$x, x' \in X_n$ are \emph{strongly connected}, $x \sim x'$, if there
exist words $w, w' \in \cA^*$ such that $w x = x'$ and $w' x' = x$.

For a group dynamics on the datum $S_n = (X_n, \cA)$, we say that $x,
x' \in X_n$ are \emph{connected}, $x \sim x'$, if there exists a word
$w \in \cA^*$ such that $w x = x'$.
\end{definition}
\noindent
Here the action $wx$ is in the sense of monoid action.  Being
connected is clearly an equivalence relation, and coincides with the
relation of being graph-connected on the Cayley Graph associated to
the dynamics, i.e.\ the digraph with vertices in $X_n$, and edges $x
\longleftrightarrow_{i} x'$ if $A_i^{\pm 1} x = x'$. An analogous
statement holds for strong-connectivity, and the associated Cayley
Digraph.


\begin{definition}[classes of configurations]
Given a dynamics as above, and $x \in X_n$, we define 
$C(x) \subseteq X_n$, the \emph{class of $x$}, as the set of
configurations connected to $x$, $C(x)=\{ x' \,:\, x\sim x'\}$.
\end{definition}

\noindent
Two natural distances on the classes of $X_n$ can be associated to a
group dynamics with generators $\cA=\{A_i\}$.
\begin{definition}[distance and alternation distance]
Let $\s$, $\t$ configurations of $X_n$ in
the same class. The \emph{graph distance} $d_G(\s,\t)$ is the
ordinary graph distance in the associated Cayley Graph,
i.e.\ $d_G(\s,\t)$ is the minimum $\ell \in \bN$ such that there
exists a word $w \in \{A_i,A_i^{-1}\}^*$,
of length $\ell$, such that $\t = w \s$. A \emph{graph geodesic} is a
$w$ realising the minimum.  We also define the
\emph{alternation distance} $d(\s,\t)$, and 
\emph{alternation geodesics}, as the analogous quantities, for words
in the infinite alphabet $\{(A_i)^j\}_{j \geq 1}$.
\end{definition}

\begin{example}
If $w = a^3\, b\, a^{-2}\, b^{-3}\, a$, and $w x = x'$, we know that the
distance $d_G(x,x')$ is at most $3+1+2+3+1=10$, and the alternation
distance $d(x,x')$ is at most 5.
\end{example}

\noindent
Contrarily to graph distance, the alternation distance is stable
w.r.t.\ a number of combinatorial operations that we will perform on
our configurations, and which are defined later on in the text or in
future work (restriction to `primitive classes', study of `reduced
dynamics', permutations in $\kS_{\infty}$,~\ldots). These facts
suggests that alternation distance is a more natural notion in this
family of problems.

The goal of the series of papers, of which this is the first one, is
to provide \emph{classification theorems} for dynamics of this kind.

In this and a companion paper we will be concerned with the
classification of classes, and further combinatorial study of their
structure, for three special group dynamics. These dynamics are two
version of the \emph{Rauzy dynamics}, first introduced by Rauzy
\cite{Rau79}, and whose study has been pioneered by Veech \cite{Vee82}
and Masur \cite{Mas82} (the connected classes are called \emph{Rauzy
  classes}). In these cases the group action is related to the
interval exchange map on translation surfaces. Thus, on one side, its
study is motivated by questions in dynamical systems. On the other
side, the previously obtained classification of classes relies on
notions and known results in algebraic geometry.  Section
\ref{ssec.3families} describes these three dynamics, and Section
\ref{sec.geome} gives a short account of these connections.

The classification theorems associated to these families have been
provided (with some \emph{caveat} discussed later on) in a series of
papers, starting with the seminal work of Kontsevich and Zorich
\cite{KZ03}, and followed by \cite{Boi12}, so our approach provides
just an alternative derivation of these results. Nonetheless, it has
two points of interest. A first point is that, as our approach is
quite different from the previous ones, along the way we happen to
extract some extra information on the combinatorial structure of the
Rauzy classes. An example, which is mentioned in Section
\ref{ssec.surgery_geo_combi} and that will be illustrated in more
detail in future work, is the construction of `many' representatives
for each Rauzy class, extending previous results of Zorich
\cite{Zor08}.  The second, more methodological point of interest is
that our approach is completely combinatorial, self-contained, and in
particular it makes no use of any facts from algebraic geometry. In
this sense, it answers a question posed by Kontsevich and Zorich
in~\cite{KZ03}:
\begin{quote}
{\it The extended Rauzy classes can be defined in purely combinatorial
  terms [...] thus the problem of the description of the extended
  Rauzy classes, and hence, of the description of connected components
  of the strata of Abelian differentials, is purely
  combinatorial. However, it seems very hard to solve it
  directly. [...] we give a classification of extended Rauzy classes
  using not only combinatorics but also tools of algebraic geometry,
  topology and of dynamical systems.}
\end{quote}
Note that the task of providing a purely combinatorial proof of these
classification theorems has been also carried over, with
different methods than ours, by~\cite{Fic16}.

In further papers, we will address also other variants of these group
actions on combinatorial structures, that, by analogy, we name
\emph{Rauzy-type dynamics}. A number of these other versions do not
have at the moment any algebraic or geometric interpretation,
nonetheless reasonings on similar grounds of the ones presented here
provide a classification theorem also in these other cases.

\subsection{Three basic examples of dynamics}
\label{ssec.3families}

\begin{figure}[tb!]
\[
\begin{array}{cc}
m=\big( (16)(24)(37)(58) \big) \in \matchs_8
&
\s=[41583627] \in \kS_8 \subseteq \matchs_{16}
\\
&
\textrm{matching diagram representation}
\\
\raisebox{7pt}{
\setlength{\unitlength}{10pt}
\begin{picture}(10,4)(0,-1)
\put(0,0.2){\includegraphics[scale=2.5]{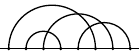}}
\put(.5,-.25){\rotatebox{30}{\makebox[0pt][r]{\scriptsize{$\bm{1}$}}}}
\put(6.8,-.25){\rotatebox{30}{\makebox[0pt][r]{\scriptsize{$m(1)=\bm{6}$}}}}
\end{picture}
}
&
\rule{0pt}{72pt}%
\raisebox{7pt}{
\setlength{\unitlength}{10pt}
\begin{picture}(20,4)(0,-1)
\put(0,0){\includegraphics[scale=2.5]{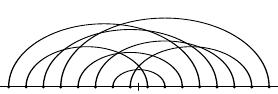}}
\put(.5,-.25){\rotatebox{30}{\makebox[0pt][r]{\scriptsize{$\bm{1}$}}}}
\put(14.5,-.25){\rotatebox{30}{\makebox[0pt][r]{\scriptsize{$n+\s(1)=\bm{12}$}}}}
\end{picture}
}
\\
\rule{0pt}{22pt}%
\s=[41583627] \in \kS_8
&
\s=[41583627] \in \kS_8
\\
\textrm{diagram representation}
&
\textrm{matrix representation}
\\
\raisebox{12.5pt}{%
\setlength{\unitlength}{10pt}
\begin{picture}(11,3)
\put(0,0){\includegraphics[scale=2.5]{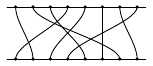}}
\put(.9,-.15){\rotatebox{-30}{\scriptsize{$\bm{1}$}}}
\put(3.6,5.4){\rotatebox{-30}{\makebox[0pt][c]{\scriptsize{$\s(1)=\bm{4}$}}}}
\end{picture}
}
&
\rule{0pt}{90pt}\raisebox{0pt}{%
\setlength{\unitlength}{10pt}
\begin{picture}(8,8)
\put(0,0){\includegraphics[scale=1]{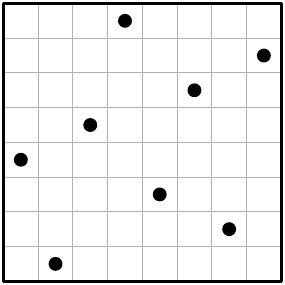}}
\put(0.4,-.6){\scriptsize{$\bm{1}$}}
\put(-.3,3.4){\makebox[0pt][r]{\scriptsize{$\s(1)=\bm{4}$}}}
\end{picture}
}
\end{array}
\]
\caption{\label{fig_match_representation} Diagram representations of
  matchings and permutations, and matrix representation of permutations.}
\end{figure}
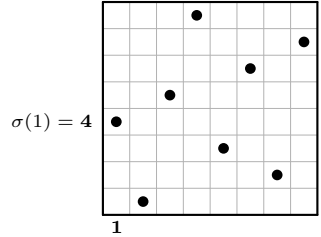

Let $\kS_n$ denote the set of permutations of size $n$, and $\kM_n$
the set of matchings over $[2n]$, thus with $n$ arcs.  Let us call
$\coxe$ the permutation $\coxe(i)=n+1-i$.

A permutation $\s \in \kS_n$ can be seen as a special case of a
matching over $[2n]$, in which the first $n$ elements are paired to
the last $n$ ones, i.e.\ the matching $m_{\s}\in \kM_n$ associated to
$\s$ is $m_{\s}=\{ (i,\s(i)+n) \,|\, i\in [n]\}.$

We say that $\s \in \kS_n$ is \emph{irreducible} if $\coxe \s$ doesn't
leave stable any interval $\{1,\ldots,k\}$, for $1 \leq k < n$,
i.e.\ if $\{\s(1), \ldots, \s(k)\} \neq \{n-k+1,\ldots,n\}$ for any
$k=1,\ldots,n-1$.
We also say that $m \in \kM_n$ is \emph{irreducible} if
it does not match an interval
$\{1,\ldots,k\}$ to an interval $\{2n-k+1,\ldots,2n\}$.
Let us call $\kSirr_n$ and $\kMirr_n$
the corresponding sets of irreducible configurations.


We represent matchings over $[2n]$ as arcs in the upper half plane,
connecting pairwise $2n$ points on the real line (see figure
\ref{fig_match_representation}, top left). Permutations, being a
special case of matching, can also be represented in this way (see
figure \ref{fig_match_representation}, top right), however, in order
to save space and improve readibility, we rather represent them as
arcs in a horizontal strip, connecting $n$ points at the bottom
boundary to $n$ points on the top boundary (as in
Figure~\ref{fig_match_representation}, bottom left).  Both sets of
points are indicised from left to right. We use the name of
\emph{diagram representation} for such representations.

We will also often represent configurations as grids filled with one
bullet per row and per column (and call this \emph{matrix
  representation} of a permutation). We choose here to conform to the
customary notation in the field of Permutation Patterns, by adopting
the algebraically weird notation, of putting a bullet at the
\emph{Cartesian} coordinate $(i,j)$ if $\s(i)=j$, so that the identity
is a grid filled with bullets on the \emph{anti-diagonal}, instead
that on the diagonal. An example is given in figure
\ref{fig_match_representation}, bottom right.

Let us define a special set of permutations (in cycle notation)
\begin{subequations}
\label{eqs.opecycdef}
\begin{align}
\gamma_{L,n}(i)
&=
(i-1\;i-2\;\cdots\;1)(i)(i+1)\cdots(n)
\ef;
\\
\gamma_{R,n}(i)
&=
(1)(2)\cdots(i)(i+1\;i+2\;\cdots\;n)
\ef;
\end{align}
\end{subequations}
i.e., in a picture in which the action is diagrammatic, and acting on
structures $x \in X_n$ from below,
\begin{align*}
\gamma_{L,n}(i):
&
\quad
\setlength{\unitlength}{8.75pt}
\begin{picture}(8,2)(-.1,.8)
\put(-.3,0.3){\includegraphics[scale=1.75]{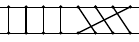}}
\put(0,0){$\scriptstyle{1}$}
\put(3,0){$\scriptstyle{i}$}
\put(7,0){$\scriptstyle{n}$}
\end{picture}
&
\gamma_{R,n}(i):
&
\quad
\setlength{\unitlength}{8.75pt}
\begin{picture}(8,2)(-.1,.8)
\put(-.3,0.3){\includegraphics[scale=1.75]{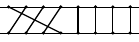}}
\put(0,0){$\scriptstyle{1}$}
\put(4,0){$\scriptstyle{i}$}
\put(7,0){$\scriptstyle{n}$}
\end{picture}
\end{align*}
Of course, $\coxe \gamma_{L,n}(i) \coxe = \gamma_{R,n}(n+1-i)$.

\begin{figure}[tb!]
\begin{align*}
L\;
\Big(
\;
\raisebox{-15pt}{\reflectbox{\includegraphics[scale=1.5]{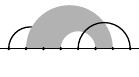}}}
\;
\Big)
&=
\raisebox{-15pt}{\reflectbox{\includegraphics[scale=1.5]{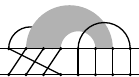}}}
&
R\;
\Big(
\;
\raisebox{-15pt}{\includegraphics[scale=1.5]{Figure1_fig_Mlr2bul.pdf}}
\;
\Big)
&=
\raisebox{-15pt}{\includegraphics[scale=1.5]{Figure1_fig_Mlr1bul.pdf}}
\\
L\;\Big(
\;
\raisebox{14.2pt}{\includegraphics[scale=1.75, angle=180]{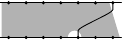}}
\;
\Big)
&=
\raisebox{27.2pt}{\includegraphics[scale=1.75, angle=180]{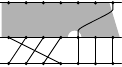}}
\rule{0pt}{38pt}
&
R\;\Big(
\;
\raisebox{-9pt}{\includegraphics[scale=1.75]{Figure1_fig_PPlr1_simp.pdf}}
\;
\Big)
&=
\raisebox{-24pt}{\includegraphics[scale=1.75]{Figure1_fig_PPlr2_simp.pdf}}
\\
L'\;\Big(
\;
\raisebox{-9pt}{\reflectbox{\includegraphics[scale=1.75]{Figure1_fig_PPlr1_simp.pdf}}}
\;
\Big)
&=
\raisebox{-24pt}{\reflectbox{\includegraphics[scale=1.75]{Figure1_fig_PPlr2_simp.pdf}}}
&
R'\;\Big(
\;
\raisebox{14.2pt}{\reflectbox{\includegraphics[scale=1.75, angle=180]{Figure1_fig_PPlr1_simp.pdf}}}
\;
\Big)
&=
\raisebox{27.2pt}{\reflectbox{\includegraphics[scale=1.75, angle=180]{Figure1_fig_PPlr2_simp.pdf}}}
\rule{0pt}{38pt}
\end{align*}
\caption{\label{fig.defDyn3}Our three main examples of dynamics, in
  diagram representation. Top: the $\matchs_n$ case. Middle: the
  $\perms_n$ case. Bottom: in the $\permsex_n$ case we have the same
  operators $L$ and $R$ as in the $\perms_n$ case, plus the two
  operators $L'$ and $R'$.}
\end{figure}

\begin{figure}[b!!]
\begin{align*}
L\;\left(
\;
\raisebox{-25pt}{\includegraphics[scale=.499]{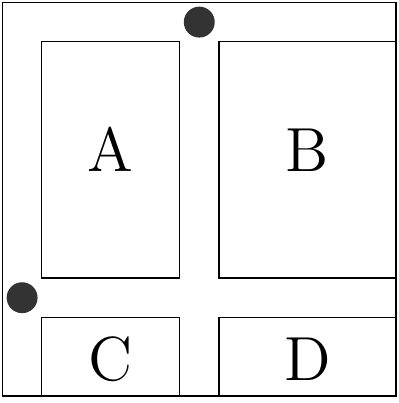}}
\;
\right)
&=
\raisebox{-25pt}{\includegraphics[scale=.499]{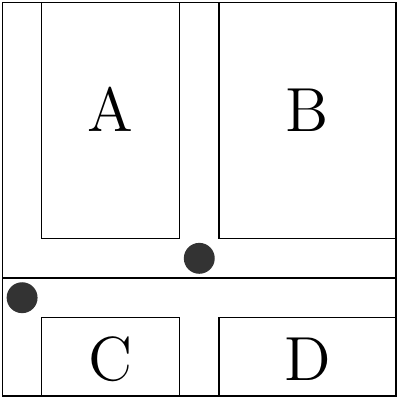}}
\rule{0pt}{38pt}
&
R\;\left(
\;
\raisebox{-25pt}{\includegraphics[scale=.499]{Figure1_fig_matrix_dyna_1.pdf}}
\;
\right)
&=
\raisebox{-25pt}{\includegraphics[scale=.499]{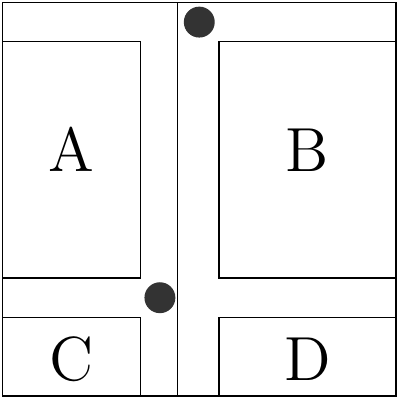}}
\\
L'\;\left(
\;
\raisebox{-25pt}{\includegraphics[scale=.499]{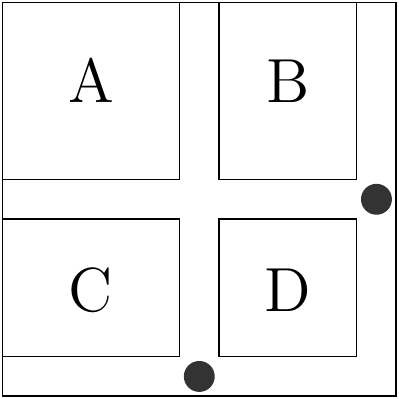}}
\;
\right)
&=
\raisebox{-25pt}{\includegraphics[scale=.499]{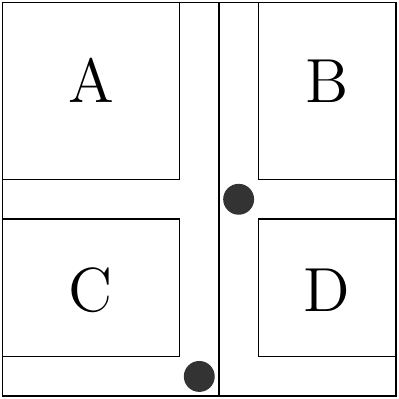}}
\rule{0pt}{38pt}
&
R'\;\left(
\;
\raisebox{-25pt}{\includegraphics[scale=.499]{Figure1_fig_matrix_dyna_3.pdf}}
\;
\right)
&=
\raisebox{-25pt}{\includegraphics[scale=.499]{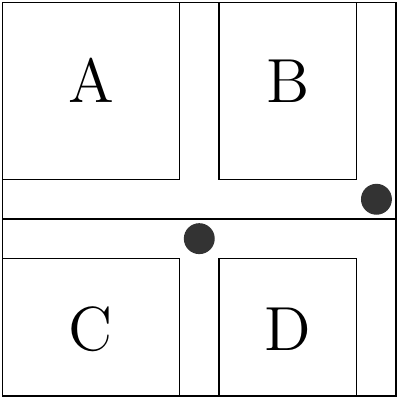}}
\end{align*}
\caption{\label{fig.defDyn3mat}Our two examples of dynamics concerning
  permutations, in matrix representation. Top: the $\perms_n$ case.
  Bottom: in the $\permsex_n$ case we have the same operators $L$ and
  $R$ as in the $\perms_n$ case, plus the two operators $L'$ and
  $R'$. See also Figure \ref{fig.defDyn3}.}
\end{figure}


The three group dynamics we mainly analyse in this and
companion papers~are
\begin{description}
\item[$\matchs_n$\;:]\ The space of configuration is $\kMirr_n$,
  irreducible $n$-arc matchings. There are two generators, $L$ and
  $R$, with $a_L(m)$ being the index paired to 1, $a_R(m)$ the index
  paired to $2n$, $\a_{L,i}=\gamma_{L,i,2n}$, and
  $\a_{R,i}=\gamma_{R,i,2n}$ ($a$ and $\a$ are as in Definition
  \ref{def.monoidoperator}).  See Figure~\ref{fig.defDyn3}, top.
\item[\phantom{$\matchs_n$}\gostrR{$\perms_n$}\;:]\ The space of
  configuration is $\kSirr_n$, irreducible permutations of size
  $n$. Again, there are two generators, $L$ and $R$. If permutations
  are seen as matchings such that indices in $\{1,\ldots,n\}$ are
  paired to indices in $\{n+1,\ldots,2n\}$, the dynamics coincide with
  the one given above. See Figure~\ref{fig.defDyn3}, middle.
\item[\phantom{$\matchs_n$}\gostrR{$\permsex_n$}\;:]\ The space of
  configuration is $\kSirr_n$. Now we have four generators, $L$ and
  $R$ are as above, and $L'$ and $R'$ act as $L' \s = S ( L (S \s))$,
  and $R' \s = S ( R (S \s )$ where $S$ is the anti-diagonal axial
  symmetry in the matrix representation.
  See Figure~\ref{fig.defDyn3}, middle and bottom.
\end{description}
More precisely, we study $\perms_n$ and $\permsex_n$ here, and delay
the study of $\matchs_n$ to future work.

\begin{figure}
\begin{align*}
&
\includegraphics{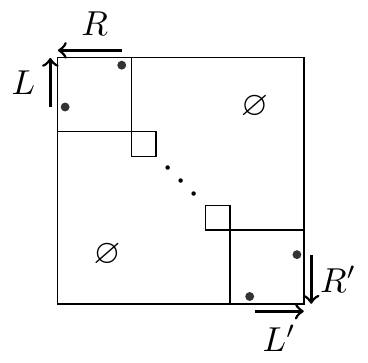}
&&
\includegraphics{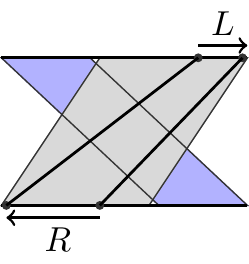}
\end{align*}
\caption{Left: A reducible permutation in matrix representation; the
  $\permsex$ dynamic acts on the first block with $L$ and $R$ and on
  the last block with $L'$ and $R'$. An analogous statement holds for
  the dynamics $\perms$.  Right: A reducible permutation in diagram
  representation; the $\perms$ dynamics acts on the gray part with $L$
  and $R$, while leaving the blue part unchanged. An analogous
  statement holds for the dynamics $\permsex$.
\label{fig.ex_reducibility}}
\end{figure}

The motivation for restricting to irreducible permutations and
matchings shall be clear at this point: a non-irreducible permutation
is a grid with a non-trivial block-decomposition. The operators $L$
and $R$ only act on the first block (say, of size $k$), while, in
$\permsex_n$, the operators $L'$ and $R'$ only act on the last block
(say, of size $k'$), so that the study of the dynamics trivially
reduces to the study of the $\perms_k$ dynamics, or of the direct
product $\perms_k \times \perms_{k'}$, on these blocks (see figure
\ref{fig.ex_reducibility}).

This simple observation, however, comes with a discaimer: in our
combinatorial operations on configurations, which produce the required
induction steps in the classification theorem, we shall always
guarantee that the outcome of our manipulations on irreducible
configurations is still irreducible. Proving such a condition will
occasionally be a subtle task.

\subsection{What kind of results?}
\label{ssec.whichkind}

In the series of papers initiated by the present one, we study
collections $X_n$ of discrete combinatorial objects, where $n$ is a
size parameter, and the cardinality of $X_n$ is exponential or
super-exponential
(i.e., $\ln X_n = \Theta(n (\ln n)^{\gamma}))$ for some $\gamma$).  We
introduce two or more bijections, called `operators', like the
operators $L$, $R$, \ldots\ of the previous section, and define
classes as the connected components of the associated Cayley Graph. It
will come out that there is a super-polynomial number of classes (in
the examples investigated here, roughly $\exp(\alpha \sqrt{n})$ for a
certain $\alpha$), the majority of them having super-exponential
cardinalities, which are not `round' numbers\footnote{I.e., numbers
  which most likely do not have simple factorised formulas, contrarily
  e.g.\ to the existence of hook formulas for the number of standard
  and semi-standard Young tableaux, or MacMahon formulas for the
  number of plane partitions, and various other examples in Algebraic
  Combinatorics.}.

Thus, a natural question is: in which form can we expect to `solve'
such a problem, if the structure of the problem is apparently so wild,
and any possible complete answer is most likely not encompassed by a
compact formula?

We may expect, and in fact give here, results of the following forms:
\begin{itemize}
\item A \emph{classification} of the classes, i.e., the identification
  of a natural labeling of the classes, and a criterium that, for a
  configuration of size $n$, gives the label of its class through an
  algorithm which is polynomial in $n$ (and possibly linear). As a
  corollary, this would give an algorithm that, for any two
  configurations, determines in polynomial time if they are in the
  same class or not.\footnote{We insist on the algorithm complexity
    aspect as, if we do not pose a complexity bound, the mere greedy
    exploration of the Cayley Graph trivialises the question at hand.}
\item The exact characterisation of the Cayley Graph of the few
  \emph{exceptional classes} that may exist, which in fact \emph{do}
  have round formulas, and an intelligible structure. Interestingly,
  in the case of the Rauzy dynamics, we found two exceptional classes:
  the hyperelliptic one that has been understood since Rauzy
  \cite[sec.~4]{Rau79}, and a second one, which is primitive only in
  the $\perms$ dynamics (not in the $\permsex$ one), and, to the best
  of our knowledge, was never described before. These two classes play
  an important role in the forementioned classification theorem, and
  their detailed study is performed in Appendix~\ref{sec.excp_class}.
\item Upper/lower bounds on various interesting quantities, e.g.\ the
  cardinalities of the classes, or their diameter.
\item Finally, what we consider the most important and original
  contribution of our work, the elucidation of a combinatorial
  structure which is recursive in $n$: even if the structure of the
  classes at size $n$ is intrinsically too complex for being described
  by a simple formula (see \cite{Del13} for the most compact formulas
  known so far), nothing prevents in principle from the existence of a
  reasonably explicit description of this structure at size $n$, in
  terms of the structures at sizes $n' < n$. This interplay of
  structures at different sizes, obtained through peculiar ``surgery''
  operations at the level of the corresponding Cayley Graphs
  (partially analogous to the well-established surgeries at the level
  of Riemann surfaces of \cite{EMZ03,KZ03}, see the following
  Section~\ref{ssec.surgery_geo_combi}), appears at several stages in
  our proofs.
\end{itemize}

\subsection{Definition of the invariants}
\label{ssec.invardefs}

The main purpose of this paper, and its first companion, is to
characterise the classes appearing in the dynamics introduced above:
of $\perms_n$ and $\permsex_n$, in this paper, and of $\matchs_n$, in
the companion one. The characterisation is based roughly on three
steps: (1)~we identify some ``exceptional classes'', for which the
structure of the configurations is completely elucidated; (2)~we
identify a collection of data structures which are invariant under the
dynamics; (3)~we prove, through a complex induction based on surgery,
that these invariants are complete, i.e.\ we characterise the
admissible invariant structures, and prove that configurations with
the same invariant, and not in an exceptional class, are also in the
same class.

In this section we sketch the step (2), i.e.\ we present the
invariants. Some of the proofs are postponed.

\subsubsection{Cycle invariant}
\label{sssec.cycinv}

\noindent
Let $\s$ be a permutation, identified with its diagram.  An edge of
$\s$ is a pair $(i^-,j^+)$, for $j=\s(i)$, where $-$ and $+$ denote
positioning at the bottom and top boundary of the diagram.  Perform
the following manipulations on the diagram: (1)~replace each edge with
a pair of crossing edges; more precisely, replace each edge endpoint,
say $i^-$, by a black and a white endpoint, $i_b^-$ and $i_w^-$ (the
black on the left), then introduce the edges $(i_b^-,j_w^+)$ and
$(i_w^-, j_b^+)$.  (2)~connect by an arc the points $i_w^\pm$ and
$(i+1)_b^\pm$, for $i=1,\ldots,n-1$, both on the bottom and the top of
the diagram; (3)~connect by an arc the top-right and bottom-left
endpoints, $n_w^+$ and $1_b^-$. Call this arc the ``$-1$ mark''.

\begin{figure}[tb!]
\begin{center}
\begin{tabular}{cp{0cm}c}
\includegraphics[scale=5]{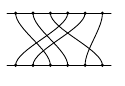}
&&
\includegraphics[scale=5]{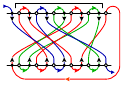}
\end{tabular}
\end{center}
\caption{\label{fig.exLamMatch}Left: an irreducible permutation,
  $\s=[\,451263\,]$. Right: the construction of the cycle
  structure. Different cycles are in different colour. The length of a
  cycle or path, defined as the number of top (or bottom) arcs, is
  thus 2 for red and violet, and 1 for blue.  The green arrows are the
  endpoints of the rank path, which is in blue.  As a result, in this
  example $\lam(\s)=(2,2)$ (for the cycles of color red and violet),
  $r(\s)=1$ (corresponding to the rank path of color blue), and
  $\ell(\s)=2$.}
\end{figure}

The resulting structure is composed of a number of closed cycles, and
one open path connecting the top-left and bottom-right endpoints, that
we call the \emph{rank path}. If it is a cycle that goes through the
$-1$ mark (and not the rank path), we call it the \emph{principal
  cycle}.  Define the length of an (open or closed) path as the number
of top (or bottom) arcs (connecting a white endpoint to a black
endpoint) in the path. These numbers are always positive integers (for
$n>1$ and irreducible permutations).  The length $r$ of the rank path
will be called the \emph{rank} of $\s$, and $\lam = \{ \lam_i \}$, the
collection of lengths of the cycles, will be called the \emph{cycle
  structure} of $\s$. Define $\ell(\s)$ as the number of cycles in
$\s$ (this does not include the rank).  See
Figure~\ref{fig.exLamMatch}, for an example.

Note that this quantity does \emph{not} coincide with the ordinary
path-length of the corresponding paths. The path-length of a cycle of
length $k$ is $2k$, unless it goes through the $-1$ mark, in which
case it is $2k+1$. Analogously, if the rank is $r$, the path-length of
the rank path is $2r+1$, unless it goes through the $-1$ mark, in
which case it is $2r+2$. (This somewhat justifies the name of ``$-1$
mark'' for the corresponding arc in the construction of the cycle
invariant.)

In the interpretation within the geometry of translation surfaces, the
cycle invariant is exactly the collection of conical singularities in
the surface (we have a singularity of $2k\pi$ on the surface, for
every cycle of length $\lam_i=k$ in the cycle invariant, and the rank
corresponds to the `marked singularity' of \cite{Boi12}, see
Section~\ref{sec.geome}).

It is easily seen that
\be
r + \sum_i \lam_i = n-1
\ef,
\label{eq.size_inv_cycle}
\ee
this formula is called the 
\emph{dimension formula}. 
Moreover, in
the list $\{r,\lam_1,\ldots,\lam_{\ell}\}$, there is an even number of
even entries (This is part of Lemma \ref{lem.induSign} in Section
\ref{ssec.signindu}, but it could also be proven easily already at
this stage\footnote{This is a simple induction. Adding an edge at any
  given position changes the cycle invariant either by 
  $\lam_i \to a+1 \oplus \lam_i - a$ for some $a$, or by 
  $\lam_i \oplus \lam_j \to \lam_i + \lam_j +1$. A case analysis on
  the parity of $\lam_i$, $a$ and $\lam_j$ leads to our claim.}).

We have
\begin{proposition}
\label{prop.cycinv1}
The pair $(\lam,r)$ is invariant in the $\perms$ dynamics.
\end{proposition}
\noindent
Now connect also the endpoints $n_w^-$ and $1_b^+$ of the rank
path. Call $\lam' = \lam \cup \{r\}$ the resulting collection of cycle
lengths.  It is easily seen that
\be
\sum_i \lam'_i = n-1
\ef.
\label{eq.size_inv_cycle_extended}
\ee
We have
\begin{proposition}
\label{prop.cycinv2}
The quantity $\lam'$ is invariant in the $\permsex$ dynamics.
\end{proposition}
\noindent
These two propositions are proven in Section~\ref{ssec.arcbased}.

\subsubsection{Sign invariant}
\label{ssec.arf_inv_intro}

\begin{figure}[tb!]
\[
\begin{array}{cp{4mm}c}
\includegraphics[scale=4]{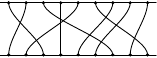}
&&
\includegraphics[scale=4]{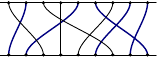}
\end{array}
\]
\caption{\label{fig.exarf}Left: an example of permutation,
  $\s=[\,251478396\,]$. Right: an example of subset $I=\{1,2,6,8,9\}$
  (labels are for the bottom endpoints, edges in $I$ are in
  blue). There are two crossings, out of the maximal number
  $\binom{|I|}{2}=10$, thus $\chi_I=8$ in this case, and this set
  contributes $(-1)^{|I|+\chi_I}=(-1)^{5+8}=-1$ to $A(\s)$.}
\end{figure}

\noindent
For $\s$ a permutation, let $[n]$ be identified to the set of edges
(e.g., by labeling the edges w.r.t.\ the bottom endpoints, left to
right). For $I \subseteq [n]$ a set of edges, define $\chi(I)$ as the
number of pairs $\{i',i''\} \subseteq I$ of non-crossing edges.  Call
\be
\Abar(\s) := \sum_{I \subseteq [n]} (-1)^{|I| + \chi(I)}
\label{eq.arf1stDef}
\ee
the 
\emph{Arf invariant} of $\s$ (see Figure \ref{fig.exarf} for an
example). Call $s(\s)= \mathrm{Sign}(\Abar(\s)) \in \{-1,0,+1\}$ the
\emph{sign} of $\s$. We have

\begin{proposition}\label{pro_sign_inv}
The sign of $\s$ can be written as
$s(\s)=2^{-\frac{n+\ell}{2}}\Abar(\s)$, where $\ell$ is the number of
cycles of $\s$. The quantity $s(\s)$ is invariant both in the $\perms$
and $\permsex$ dynamics.
\end{proposition}

\noindent
The proof of invariance claimed in this proposition is given in
Section~\ref{ssec.arfcalcseasy}, while the proof that
$s(\s)=2^{-\frac{n+\ell}{2}}\Abar(\s)$ is given in
Section~\ref{ssec.signindu}, namely in Lemma \ref{lem.induSign}.  In
section \ref{sec.geome} we motivate the name of Arf invariant by
making explicit the connection to the associated quantity in the
theory of translation surfaces (and, more generally, of Riemann
surfaces).  

Here we give the main idea in the proof of the invariance.  The four
operators of $\permsex$ are related by a dihedral symmetry of the
diagram, and the expression for $\Abar(\s)$ is manifestly invariant
under these symmetries, so it suffices to consider just one operator,
say $L$. This operator uses an edge $e$ of the diagram of $\s$ to
determine the permutation, then it moves a number of edge endpoints by
a single position on the left, and the endpoint of a second edge $f$
on the right.

Any function of the form
$F(\s) = \sum_{I \subseteq [n]} f_{\s}(I)$
can be decomposed in the form
\[
F(\s) = \sum_{I \subseteq [n]\setminus \{e,f\}} \big( 
f_{\s}(I) +
f_{\s}(I \cup \{e\}) + f_{\s}(I \cup \{f\}) + f_{\s}(I \cup \{e,f\})
\big)
\ef.
\]
Call $\t=L(\s)$. It can be verified that, for the function $\Abar(\s)$,
\begin{align*}
f_{\s}(I) &= f_{\t}(I)
\ef,
&
f_{\s}(I \cup \{e\}) &= f_{\t}(I \cup \{e\})
\ef,
\\
f_{\s}(I \cup \{f\}) &= f_{\t}(I \cup \{e,f\})
\ef,
&
f_{\s}(I \cup \{e,f\}) &= f_{\t}(I \cup \{f\})
\ef,
\end{align*}
so that the sum of the four terms is invariant for any 
given $I \subseteq [n] \setminus \{e,f\}$.

\subsubsection{Cycles of length 1 and primitivity}
\label{ssec.primi}

\noindent
We have stated in Section \ref{sssec.cycinv} that a certain graphical
construction leads to the definition of a ``cycle invariant''. In this
section we discuss how cycles of length 1 have an especially simple
behaviour.

Establishing this simplifying feature is helpful in our proof of
the main classification theorem. It goes in the direction
of understanding as many as possible relevant combinatorial features
of these classes, and allows us to rule out a large part of the massive numerics associated to
this problem. 

\begin{figure}[bt!]
\begin{center}
\begin{tabular}{ccc}
\includegraphics[scale=.8]{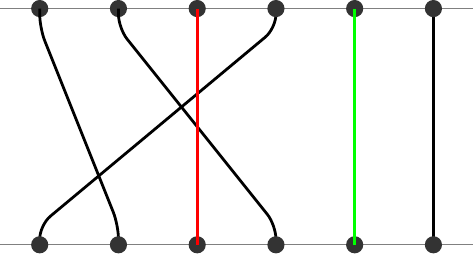}
&
&
\includegraphics[scale=.8]{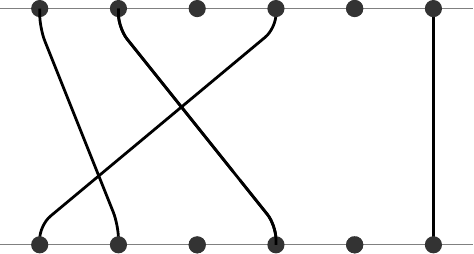}\\
\includegraphics[scale=.8]{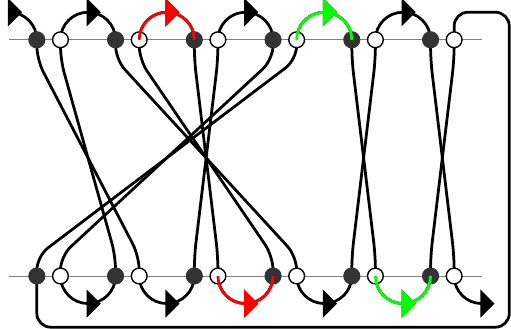}
&&
\includegraphics[scale=.8]{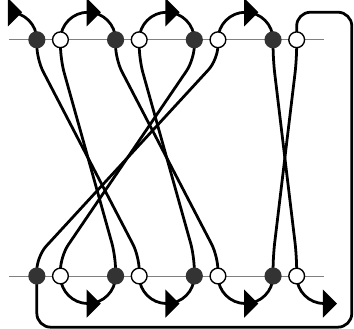}
\end{tabular}
\end{center}
\caption{Top left: A non-primitive permutation $\s \in \perms_{6}$
  with a descent in red and the special descent in green. Its cycle
  invariant is $(\{1,1\},r=3)$ as can be seen explicitly in the cycle
  construction shown on the bottom left of the figure. Top right: the
  permutation $\prim{\s}$, the primitive of $\s$. Its cycle invariant
  is $(\emptyset,r=3)$, as can be seen on the bottom right of the
  figure.\label{fig.ex_primitive_non_primitive}}
\end{figure}

\begin{definition}[descent and special descent]
For a permutation $\s$ in the dynamics $\perms_n$, we say that the
edge $(i,\s(i))$ is a \emph{descent} if $\s(i+1)=\s(i)-1$, and it is a
\emph{special descent} if $\s(1)=\s(i)-1$ and $\s(i+1)=n$.
For $\s$ in the dynamic $\permsex_n$, we say that
$(i,\s(i))$ is a \emph{special descent} also
if $\s^{-1}(i-1)=1$ and $\s^{-1}(n)=\s(i)+1$.
\end{definition}
\noindent
Note that the descents of a permutation (special or not) are
associated to cycles of length~1 (see figure
\ref{fig.ex_primitive_non_primitive} left). In particular, the number
of descents is preserved by the dynamics.
Say that $\s$ is \emph{primitive} if it has no descents.
Hence in
a given class $C$ either all permutations are primitive (in this case
we say that $C$ is a \emph{primitive class}), or none is.  We define
$\prim{\s}$, the \emph{primitive of $\s$}, as the configuration obtained by
removing descent- and special-descent-edges
from $\s$ (see figure \ref{fig.ex_primitive_non_primitive} for an
example in the $\perms$ dynamics).  We have
\begin{proposition}
\label{prop.primhomoPre}
$\s \sim \tau$ iff $\prim{\s} \sim \prim{\tau}$ and $|\s|=|\t|$.
\end{proposition}
\noindent
In other words, within a class we can `move the descents freely'. In
particular, this gives
\begin{corollary}
\label{cor.primhomo}
The map $\prim{\,\cdot\,}$ is a homomorphism for the dynamics.
\end{corollary}
\noindent
I.e., if $\s \sim \tau$, then $\prim{\s} \sim \prim{\tau}$.  
In fact, much more is true
(various other properties are discussed in
Appendix \ref{app.primitheo}).


\subsection{Exceptional classes}
\label{sssec.ididp}

As we have outlined, a simple preliminary analysis allows to restrict to
`irreducible' and `primitive' classes.  Then, as anticipated, the
invariants described above allow to characterise all primitive classes
for the dynamics on irreducible configurations, \emph{with two
  exceptions}. These two exceptional classes, for the $\perms_n$
dynamics, are called $\Id_n$ and $\tree_n$. In the $\permsex_n$ dynamics,
only $\Id_n$ remains exceptional and primitive (in the literature on
Rauzy dynamics, $\Id_n$ is called the \emph{`hyperelliptic class'},
because the Riemann surface associated to $\Id_n$ is hyperelliptic. Similarly $\Id'_n$ is often referred to as the \emph{`hyperelliptic class' with a marked point}).

The most compact definition of $\Id_n$ and $\tree_n$ is as being the
classes containing $\id_n$ and $\id'_n$, respectively, where $\id_n$
is the identity of size $n$, and, for $n \geq 3$, $\id'_n$ is the
permutation $\s$ of size $n$ such that $\s(1)=1$, $\s(2)=n-1$,
$\s(n)=n$ and $\s(i)=i-1$ elsewhere, for example
\begin{align*}
\id_6 &=\ 
\raisebox{-32pt}{\begin{tikzpicture}[scale=0.4]
\permutation{1,2,3,4,5,6}
\end{tikzpicture}}
&
\id'_6 &=\ 
\raisebox{-32pt}{\begin{tikzpicture}[scale=0.4]
\permutation{1,5,2,3,4,6}
\end{tikzpicture}}
\end{align*}
%

\begin{table}[b!!]
\[
\begin{array}{|c|c|c|}
\cline{2-3}
\multicolumn{1}{c|}{}
&
\textrm{$n$ even}
&
\textrm{$n$ odd}
\\
\hline
(\lam,r) \text{ of } \Id_n & 
(\emptyset,n-1)&
\raisebox{-5pt}{\rule{0pt}{15pt}}%
(\{\frac{n-1}{2}\}, \frac{n-1}{2})
\\ 
\hline
(\lam,r) \text{ of } \tree_n &
\raisebox{-5pt}{\rule{0pt}{15pt}}%
(\{ \frac{n-2}{2}, \frac{n-2}{2} \},1)& 
(\{n-2\},1)
\\ \hline
\end{array}
\]
\[
\begin{array}{|c|cccccccc|}
\cline{2-9}
\multicolumn{1}{c|}{n \mod 8}
&
0&1&2&3&4&5&6&7
\\ \hline
s \text{ of } \Id_n & 
+&0&-&-&-&0&+&+
\\ 
\hline
s \text{ of } \tree_n &
+&+&0&-&-&-&0&+\\
\hline
\end{array}
\]
\caption{Cycle, rank and sign invariants of the exceptional classes. The
  sign \hbox{$s \in \{-1,0,+1\}$} is shortened into $\{-,0,+\}$.
\label{table_invariant_id_tree}}
\end{table}

The cycle and sign invariants of these classes depend from their size
mod~4, and are described in Table~\ref{table_invariant_id_tree}.  The
classes $\Id_n$ and $\tree_n$ are atypical w.r.t.\ other classes, in
various respects:
\begin{itemize}
\item They have `small' cardinality. More precisely, at size $n$,
  $|\Id_n| = 2^{n-1}-1$ and $|\tree_n| = (2^{n-2}+n-2)$ (compare this
  to the fact that all other classes have size $\geq \exp(c\, n \ln
  n)$, with $c \geq \frac{2}{3}$, i.e.\ at least size $\sim
  {n!}^{2/3}$).
\item Contrarily to all other classes, the Cayley Graph associated to
  the dynamics can be described in a compact way, in terms of one or
  more complete binary trees of a certain height, in the case of
  $\tree_n$ complemented by `few' other nodes and transitions (a
  linear number).
\item The configurations of these classes have a simple structure,
  labeled by a certain array of integers, and related to the position
  of the configuration in the Cayley graph. This structure makes easy
  to verify, in linear time, if a configuration is in $\Id_n$, in
  $\tree_n$, or in none of the above, and thus allows to restrict the
  quest for a classification to non-special classes.
\end{itemize}
These results are illustrated at length in Appendix~\ref{app.ididp}.

\subsection{The classification theorems}
\label{ssec.theointro}

For the case of the $\perms$ dynamics, we have a classification
involving the cycle structure $\lambda(\s)$, the rank $r(\s)$ and sign
$s(\s)$ described in Section \ref{ssec.invardefs}

\begin{theorem}\label{thm.Main_theorem}
Two permutations $\s$ and $\s'$ are in the same class iff they have
the same number of descents, and $\prim{\s}$, $\prim{\s'}$, 
are in the same class.

A permutation $\s$ is primitive iff
$\lam(\s)$ has no parts of size 1 (note that the rank may be 1).
Besides the exceptional classes $\Id$ and $\tree$, which have cycle
and sign invariants described in Table~\ref{table_invariant_id_tree},
the number of primitive classes with cycle invariant $(\lam,r)$ (no
$\lambda_i=1$) depends on the number of even elements in the list
$\{\lambda_i\} \cup \{r\}$, and is, for $n \geq 9$,
\begin{description}
\item[zero,] if there is an odd number of even elements;
\item[\phantom{zero,}\gostrR{one,}] if there is a positive even number of even elements; the class then has sign 0.
\item[\phantom{zero,}\gostrR{two,}] if there are no even elements at all. The two classes then have
non-zero opposite sign invariant.
\end{description}
For $n \leq 8$ the number of primitive classes with given cycle
invariant may be smaller than the one given above, and the list in
Table \ref{tab.smallsizeThm1} gives a complete account.

As a consequence, two primitive permutations $\s$ and $\s'$, not of
$\Id$ or $\tree$ type, are in the same class iff they have the same
cycle and sign invariant.
\end{theorem}

\begin{table}[t!!]
\[
\begin{array}{r||c|c|cccccc}
n & \Id & \Id' & \multicolumn{6}{c}{\textrm{non-exceptional classes}}
\\
\hline
4 & \emptyset|3- &  \\
5 & 2|2 & 3|1- \\
6 & \emptyset|5+ & 22|1 & \emptyset|5- \\
7 & 3|3+ & 5|1+ & 2|4 & 4|2 & 3|3- & 5|1- \\
8 & \emptyset|7+ & 33|1+ & 
22|3 & 32|2 & 42|1 & 33|1- & \emptyset|7+ & \emptyset|7-
\end{array}
\]
\caption{\label{tab.smallsizeThm1}List of invariants $(\lam,r,s)$ for
  $n \leq 8$, for which the corresponding class in the $\perms_n$
  dynamics exists. We shorten $s$ to $\{-,+\}$ if valued $\{-1,+1\}$,
  and omit it if valued~0.}
\end{table}

\begin{table}[b!!]
\[
\begin{array}{r||c|cccccc}
n & \Id & \multicolumn{6}{c}{\textrm{non-exceptional classes}}
\\
\hline
4 & 3- & \\
5 & 22 \\
6 & 5+ & 5- \\
7 & 33+ & 24 & 33-
\end{array}
\]
\caption{\label{tab.smallsizeThm2}List of invariants $(\lam',s)$ for
  $n \leq 7$, for which the corresponding class in the $\permsex_n$
  dynamics exists. We shorten $s$ to $\{-,+\}$ if valued $\{-1,+1\}$,
  and omit it if valued~0.}
\end{table}

\noindent
For the case of the $\permsex$ dynamics, we have only two invariants
left, as defined in Section \ref{ssec.invardefs}: the cycle structure
$\lam'(\s)$ and the sign $s(\s)$. Recall that $\lam'(\s) = \lam(\s)
\cup \{r(\s)\}$, when the same configuration is considered under the
$\perms$ dynamics. Note however that a permutation that is primitive
and of rank 1 in $\perms$, is not primitive in $\permsex$. In
particular, the class $\tree_n$, which has rank $1$ for all $n$, is
non-primitive in $\permsex$.  We have

\begin{theorem}\label{thm.Main_theorem_2}
Two permutations $\s$ and $\s'$ are in the same class iff they have
the same number of descents, and $\prim{\s}$, $\prim{\s'}$, 
are in the same class.

A permutation $\s$ is primitive iff $\lam'(\s)$ has no parts of size
1.  Besides class $\Id$, which has invariant as in the previous
theorem (the rank is just added to $\lam'$), the number of primitive
classes with cycle invariant $\lam'$ (no $\lam'_i=1$) depends on the
number of even elements in the list $\{\lambda'_i\}$, and is, for 
$n \geq 8$,
\begin{description}
\item[zero,] if
there is an odd number of even
elements;
\item[\phantom{zero,}\gostrR{one,}] if there is a positive even number of even elements;
\item[\phantom{zero,}\gostrR{two,}] if there are no even elements at all. The two classes have
non-zero opposite sign invariant.
\end{description}
For $n \leq 7$ the number of primitive classes with given cycle
invariant may be smaller than the one given above, and the list in
Table \ref{tab.smallsizeThm2} gives a complete account.

As a consequence, two primitive permutations $\s$ and $\s'$, not of
$\Id$ type, are in the same class iff they have the same
cycle and sign invariant.
\end{theorem}

\noindent
We will first obtain Theorem \ref{thm.Main_theorem} and then Theorem
\ref{thm.Main_theorem_2} as an almost-straighforward corollary.
Curiously enough, historically, the first and fundamental article
\cite{KZ03} proved Theorem \ref{thm.Main_theorem_2}, and it was only a
few years later that the proof technique was adapted from the
$\permsex$ case to $\perms$, and the article \cite{Boi12} proved
Theorem~\ref{thm.Main_theorem}.

On the contrary, within our techniques, it is both (a bit) easier to
prove Theorem \ref{thm.Main_theorem} than Theorem
\ref{thm.Main_theorem_2}, if we had to do both of them from scratch,
and is considerably easier to prove Theorem \ref{thm.Main_theorem_2}
as corollary of \ref{thm.Main_theorem}, while we are not aware of a
simple derivation of Theorem \ref{thm.Main_theorem} from
Theorem~\ref{thm.Main_theorem_2}.


\subsection{Surgery operators}
\label{ssec.TQSintro}

The hardest part in our proof of Theorem \ref{thm.Main_theorem},
concerning the dynamics $\perms_n$, is the proof by induction that all
the classes with a given admissible set of invariants are non-empty
and connected. Classes of different rank (in the three cases $r=1$,
$r=2$ and $r>2$) are treated differently, so we do this with 3
operators, that we call $q_1$, $q_2$ and $T$, for these three cases,
respectively. These operators produce (representants within) classes
of given rank at size $n$ from (representants within) classes at
smaller sizes, by performing local manipulations of the
configurations, which, in the language of translation surfaces (see
section \ref{ssec.surgery_geo_combi}), correspond in a hidden way to
the insertion of a `handle' or of a `cylinder' in a Riemann
surface. Thus we call them
\emph{surgery operators}.
In the next paragraphs we make a preliminary discussion on which
combinatorial properties shall such operators have, in order to
provide the appropriate tool in our proof's scheme.

First of all, for each operator $X$, we will require that it is a
homomorphism, i.e., that for all pairs of configurations in the same
class, $\s\sim \t$, we have $X(\s)\sim X(\t)$.  As a consequence,
$\bar{X}(C)$, the class obtained by applying $X$ to any configuration
in $C$, is well-defined.

Then, we have a simple yet crucial constraint on the action the
operators may possibly have, given by the 
dimension
formula,
equation (\ref{eq.size_inv_cycle}): for a permutation $\s$ of size $n$
with cycle invariant $(\lambda,r)$, we have $r+\sum_{i} \lam_i=n-1$.
We know in retrospective, from the statement of Theorem
\ref{thm.Main_theorem} (see Section \ref{ssec.signindu}, Lemma
\ref{lem.induSign} for a proof), that the number of even elements in
$\lam \cup \{r\}$ must be even, so that the cycle invariant shall
change consistently under the action of the surgery operator.

For example, suppose that we aim to construct an operator $X$ such
that, for every $\s$ with invariant $(\lambda,r,s)$, $X(\s)$ has
invariant $(\lambda,r+1,s)$. This operator would be a viable candidate
for constructing an induction.  However, the number of even elements
in the two lists $\{\lambda_i\} \cup \{r\}$ and $\{\lambda_i\} \cup
\{r+1\}$ have different parity, so that we know in advance that such
an operator $X$ cannot exist.  The best hope is to identify instead an
operator -- it will be our operator $T$ -- such that, for every $\s$
with invariant $(\lambda,r,s)$, $T(\s)$ has invariant
$(\lambda,r+2,s)$. Such an operator increases the size of a
configuration by 2 (by the 
dimension
formula), while preserving the
forementioned parity constraint.

The classes in the image of $T$ have rank $r>2$, thus for the cases of
rank 1 and rank 2 we need to define two other operators, $q_1$ and
$q_2$. For the first one, tentatively, we may hope to have an operator
such that, if $\s$ has invariant $(\lambda,r,s)$ and size $n$, then
$q_1(\s)$ has invariant $(\lambda',1,s)$ and size $n+1$\footnote{We
  may consider also a larger increase of size, as was for $T$, but let
  us consider the simplest scenario first.}  Again, by equation
(\ref{eq.size_inv_cycle}) we have that
$1+\sum \lambda'_i=r+\sum \lambda_i +1$. What should we hope as a
possible natural choice of $\lambda'$? The answer is deceivingly
simple, namely $\lambda':=\lambda\cup\{r\}$, i.e.\ we add a new cycle,
of length $r$, in order to satisfy equation~(\ref{eq.size_inv_cycle}).

Finally, it remains the case of $q_2$. This operator works on a
similar basis as $q_1$, with the following difference: if $\s$ has
invariant $(\lambda,r,s)$ and size $n$, then $q_2(\s)$ has invariant
$(\lambda'',2,0)$ and size $n+1$ with
$\lambda'':=\lambda\cup\{r-1\}$. Note that the sign invariant of
the image class is set to zero, regardless of the sign in the
pre-image.  This is consistent with the statement of the theorem: the
new rank is $2$, so the number of even entries in $\lambda \cup \{r\}$
becomes positive (while keeping its parity due to our definition of $q_2$).

Which properties of these operators shall we establish?

As we will prove the theorem by induction, we can assume that at size
$n-1$ and $n-2$ the classes are fully charaterized by the triple
$(\lambda,r,s)$ (where $s = s(\s) \in \{-1,0,1\}$), as described in
the statement of the theorem, and then concentrate on size~$n$.

For the operator $T$, we need the following to hold: let $C$ be some
class with invariant $(\lambda,r,s)$, with $r>2$, and let $B$ be the
(unique by induction) class with invariant $(\lambda,r-2,s)$, then $\bar{T}(B)=C$.
The mere existence of this operator, once established that it is a
homomorphism, would provide the inductive step, from size $n-2$ to
size $n$, for every class with rank more than 2.  Indeed the operator
$\bar{T}$ is a bijection between classes of invariant $(\lambda,r,s)$
and classes of invariant $(\lambda,r+2,s)$.

For $q_1$, we have a slightly more complex requirement.  Let $C$ be a
class with invariant $(\lambda,1,s)$. For every $i$ an integer
appearing in $\lambda$ with positive multiplicity, let $\lambda(i)=
\lambda\setminus\{i\}$.  Define $\Delta(C,i)$ as the (unique by
induction) class with invariant $(\lambda(i),i,s)$.  Clearly, for all
such $i$, $\bar{q}_1(\Delta(C,i))$ has invariant $(\lambda,1,s)$.
Moreover, it turns out that no other classes $B$ are such that
$\bar{q}_1 B$ has this invariant.  Now, establishing that
$q_1$ is a homomorphism is not enough. We need to further establish
that $\bar{q}_1(\Delta(C,i))=C$ for every $i$, ruling out the
possibility that, e.g., if $\s \in \Delta(C,i)$ and $\t \in
\Delta(C,j)$, then
$q_1(\s) \not\sim q_1(\t)$ despite the fact that they have the same
invariants.  A crucial simplifying property of the operator $q_1$ we
will introduce is that $q_1$ is a injective map that is almost
surjective, because, when $q_1^{-1}(\tau)$ is not defined, then
$q_1^{-1}(R^{-1}\tau)$ is.

Finally, it remains the case of $q_2$.  This operator works on a
similar basis as $q_1$, but the loss of memory of the sign invariant
introduces a further subtlety.  Now we need the following.  Let $C$
have invariant $(\lambda,2,0)$. The set of $(\lambda',r',s')$ such
that it may be $\bar{q}_2(C')=C$ for some $C'$ with invariant
$(\lambda',r',s')$ is a bit complicated. First of all, its elements
have the form $(\lambda(i-1),i,s)$ for $i\in \lambda$ and $s \in
\{0,\pm1\}$.  Then, if there are no even cycles in $\lambda(i-1)$, we
have $s\in \{\pm 1\}$, otherwise $s=0$.  For $i$ and $s$ in the set as
above (depending on $\lam$), define $\mathcal{E}(C,i,s)$ as the class
with invariant $(\lambda(i-1),i,s)$ (by induction, there is at most
one class with such invariant). Thus $\bar{q}_2(\mathcal{E}(C,i,s))$
has invariant $(\lambda,2,0)$.  We ask that
$\bar{q}_2(\mathcal{E}(C,i,s))=C$ for every such $i$ and $s$.
Again, besides the list we have just described, there are no other
classes $B$ such that $\bar{q}_2 B$ has invariant $(\lam,2,0)$.
Similarly as was the case for $q_1$, $q_2$ is injective and almost
surjective.

The construction of these operators would produce the core of the
induction steps. Some further subtleties are in order, though. In
particular, the existence of exceptional classes requires to verify
the bevaviour of classes $T(\Id_n), q_1(\Id_n), q_2(\id_n)$ and
$T(\tree_n), q_1(\tree_n), q_2(\tree_n)$, at all sizes $n$, as these
could create a proliferation of new classes, even if the properties
outlined above are determined for non-exceptional classes.  For
example, in such a case 
$\{T^{2k}(\Id_4),T^{2(k-1)}(\Id_8),\ldots T^2(\Id_{4k}),\Id_{4k+4}\}$
could be $k+1$ different classes, according to Table
\ref{table_invariant_id_tree}, all with invariant $(\emptyset,k+1,1)$.
Luckily enough, what really happens (and we manage to prove) is that
for each of the three operators $X\in\{\bar{T},\bar{q}_1,\bar{q}_2\}$,
and each of the two families of exceptional classes $C\in \{
\Id_n,\tree_n\}$, for $n$ above a finite (small) value, there exists a
non-exceptional class $D$ such that $X(C)=X(D)$.

In our proof, it is convenient to establish separately this last fact
(this is done in Section \ref{ssec.tech}), and the behaviour of the
operators on non-exceptional classes (this is done in Sections
\ref{ssec.T_op} and \ref{ssec.q_op}).

Table \ref{tab.exinduct}, which provides the worked-out procedure of
the induction at size 11, may be illuminating.

\section{Connection with the geometry of translation surfaces}
\label{sec.geome}

In this section we illustrate how the combinatorial operators and invariants
introduced in Section \ref{sec.algsett} are related to certain
operations, called \emph{Rauzy--Veech induction}, acting on interval
exchange transformations associated to the Poincar\'e map of an
interval on a translation surface. This correspondance implies that the
classification of equivalence classes w.r.t.\ the combinatorial
operators is equivalent to the classification of strata of translation
surfaces with a given set of conical singularities. 

A first part reviews the geometrical background (and explains all the
terms used in the paragraph above), while a second part defines the
geometrical version of the invariants (illustrated in Section
\ref{ssec.invardefs}) that turned out to provide a complete set in the
classification theorem. Finally a last part describes the geometrical
interpretation of our combinatorial surgery operators (illustrated in
Section~\ref{ssec.TQSintro}).

Note that this section introduce the minimal amount of terminology
required to make sense of the invariants and the surgery operators in
the geometrical context. The interested reader may find a more
extensive account in the three surveys \cite{Zor06}, \cite{YocPisa}
and \cite{Forni2014271}, which develop respectively the theory of
(quasi-)flat surfaces, of interval exchange transformations and of
Teichm\"uller spaces (the three texts have a large amount of overlap
between them).

\subsection{Strata of translation surfaces}

A \emph{translation surface} is a connected compact oriented surface
$M$ of genus $g$, equipped with a Riemannian metric which is flat
except for a discrete set $S$ of conical singualities, and has trivial
$\mathrm{SO}(2)$-holonomy. In this context, trivial
$\mathrm{SO}(2)$-holonomy means that the parallel transport of a
tangent vector along any closed curve (avoiding the singularities)
brings the vector back to itself (instead of having it rotated). A
condition for this to happen is that the angle of each sigularity is a
multiple of $2\pi$. Indeed, for `small' closed curves, which encircle
one conical singularity of angle $\alpha$ while staying far from all
the others, the parallel transport of a vector would rotate the vector
by the angle $\alpha$, which is defined modulo~$2 \pi$.

\begin{figure}[tb!]
    \begin{center}
     \includegraphics[scale=1]{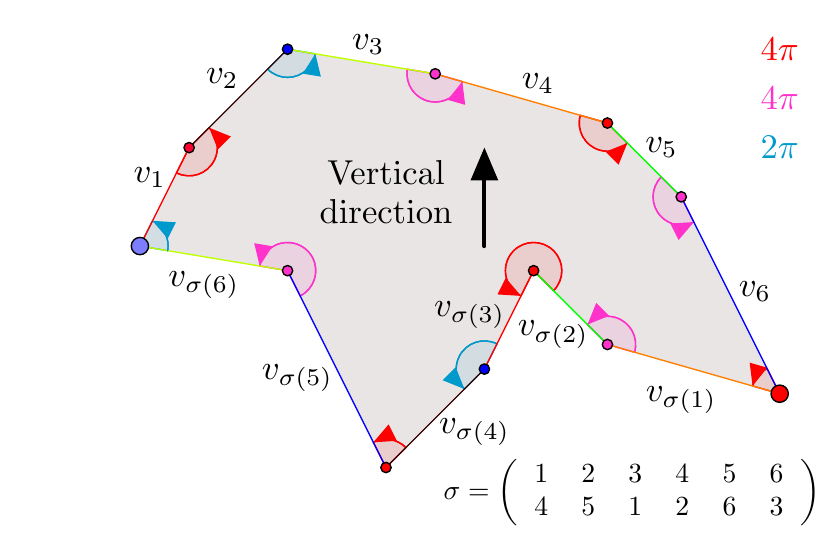}
\rule{30pt}{0pt}
    \end{center}
\caption[caption]{\label{fig_cons_surf}From the datum of $n$ vectors
  $v_1,\ldots,v_n$ in $\mathbb{C}$ and a permutation 
  $\sigma\in \mathfrak{S}_n$, we construct a $2n$-gon. Then every pair
  of parallel sides of the polygon are identified by translation.  The
  conical singularities of the associated translation surface are
  represented directly on the figure, by small angles, and the stratum
  $H(d_1^{m_1},\ldots,d_k^{m_k})$ to which the surface belongs can be
  evinced from this representation: the vertices of the polygon are
  partitioned into equivalence classes (due to the identification of
  the sides of the polygon by translation) represented by the coloured
  arcs in the figure. Each class is associated to a singularity, and
  the angle of each singularity is the sum of the angles around the
  vertices in its equivalence class. We can also directly compute the
  angle of a singularity, divided by $2\pi$, by counting the number of
  arcs of the given colour on any of the two broken lines $P_{\pm}$ of
  the polygon (the two arcs around the vertices $0$ and
  $\sum_{i=1}^n v_i$, there where $P_+$ and $P_-$ meet, are not
  counted).  In this example there are three equivalence classes,
  denoted in blue, red and magenta. The angle of the blue, red and
  magenta classes are $2\pi$, $4\pi$ and $4\pi$, respectively. Indeed,
  there are 1 blue, 2 red and 2 magenta arcs both on the top and on
  the bottom broken lines.}
\end{figure}

Additionally, the trivial holonomy allows us to choose a vertical
direction. Choose one (non-singular) point on the surface, and one
vector $x$ in the tangent space at that point. Then, one can define a
vector field via the parallel transport of~$x$. The resulting parallel
vector field is defined only outside of the singularities, but it can
be extended to the singularities in a multivalued way: in a
singularity of angle $2\pi(d_i+1)$ it would take $d_i+1$ distinct
values.


Equivalently, a translation surface with the choice of a vertical
direction can be described as a triplet $(M,U,S)$, where $M$ is a
topological connected surface, $S=\{s_1,\ldots,s_k\}$ is the discrete
subset of singularities of $M$, and $U= \{U_i,z_i\}$ is an atlas on
$M\setminus S$ such that every transition function is a translation:
\[
\forall i,j \qquad    
\begin{array}{ccccc}
    z_i\circ z_j^{-1} & : & z_j(U_i\cap U_j)  & \to & z_i(U_i\cap U_j) \\
     & & z & \mapsto & z+c \\
\end{array}
\]
and for each $s_i\in S$ there exists a neighborhood that is isometric
to a cone of angle $2\pi(d_i+1)$. Since $z_i=z_j +c$ for all pairs of
intersecting $U_i$ and $U_j$, we can define on each chart a holomorphic
complex 1-form $\omega_{z_i}= dz_i$, and thus the abelian differential
$\omega$ on $M\setminus S$. This abelian differential $\omega$ extends
to the points $s_i\in S$, where it has zeroes of degree $d_i$.

Conversely, given a non-zero abelian differential $\omega$ on a
compact connected Riemann surface $M$, with a finite set of zeroes
$S=\{s_1,\ldots s_k\}$, of degree $d_1,\ldots,d_k$, we define a
\emph{translation atlas} as follows. Let $(U,\zeta)$ be a chart
containing $p\in M\setminus S$ such that
$\omega_\zeta=\Phi(\zeta)d\zeta$. From this, one can obtain a chart
$(U,\xi)$ of adapted local coordinates at $p$ with $\omega_\xi=d\xi$ on
$U$ by defining $\xi=\int_p^z \Phi(w)dw$. Likewise, we can show that a neighbourhood of the zero $s_i$ is
isometric to an Euclidean cone of angle $2\pi(d_i+1)$, by finding
local coordinates $\eta$ for which $\omega_\eta=\eta^{d_i}d\eta$.
    
A geometric construction of a translation surface can be done as
follows. Choose $n$ distinct vectors $v_1,\ldots,v_n$ of $\mathbb{C}$ and a
permutation $\sigma\in \mathfrak{S}_n$ and construct the two broken
lines (i.e., polygonal chains) $P_+$ and $P_-$ in $\mathbb{C}$, with
vertices
\[
P_+ = (0,v_1,v_1+v_2,\ldots, v_1+\cdots+v_n)
\]
and
\[
P_- = (0,v_{\sigma(n)},v_{\sigma(n)}+v_{\sigma(n-1)},\ldots,
v_{\sigma(n)}+\cdots + v_{\sigma(1)} )
\ef.
\] 
Occasionally, we will use $u_1$, $u_2$, \ldots, $u_n$ as synonym of
$v_{\sigma(n)}$, $v_{\sigma(n-1)}$, \ldots, $v_{\sigma(1)}$.

If the vectors $v_i$ satisfy certain inequalities (see
Definition~\ref{def.suspdata} below), $P_+$ and $P_-$ intersect only
at their endpoints, so the concatenation of $P_+$ and the reverse of
$P_-$ defines a polygon embedded in the plane, with $2n$ sides, which
are naturally paired: if $\sigma(i)=j$, the $j$-th side of $P_+$ has
the same direction and orientation as the $(n+1-i)$-th side of
$P_-$. These pairs of sides can thus be identified, so that the
resulting surface has no border. By this procedure, illustrated in
figure \ref{fig_cons_surf}, we obtain a translation surface.
            
It is known that any translation surface can be represented in such a
way, and under the canonical choice of the vertical direction for the
vector field, the corresponding abelian differential is the canonical
$dz$ of the complex plane.  As explained on figure
\ref{fig_cons_surf}, the set of conical singularities can be directly
read on the polygon.

\begin{figure}[t!!]
\begin{center}
 \includegraphics[scale=1]{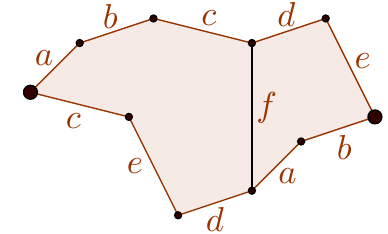}
 \includegraphics[scale=1]{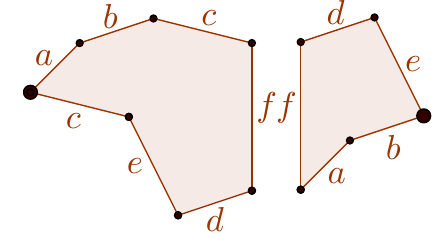}
 \includegraphics[scale=1]{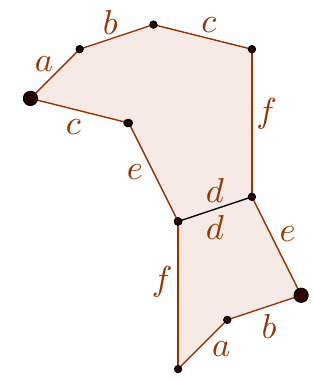}
\end{center}
\caption[caption]{\label{fig_cut_and_paste}A cut and paste move.}
\end{figure}

For $g \geq 2$, the moduli space of abelian differentials $H_g$ is the
space of pairs $(M,\omega)$, where $M$ a compact connected Riemann
surface of genus $g$ and $\omega$ is a non-zero abelian differential,
with the following identification: the points $(M,\omega)$ and
$(M',\omega')$ are equivalent if there is an analytic isomorphism $f:
M \to M'$ such that $f^*\omega'=\omega$. $H_g$ is a complex orbifold
of dimension~$4g-3$.

Let $d_1 \leq \cdots \leq d_k$ and $m_1,\ldots,m_k$ be sequences of
nonnegative integers.
We denote by $H(d_1^{m_1},\ldots,d_k^{m_k})$ the subspace of $H_g$
such that for all points $(M,\omega)$ the abelian differential
$\omega$ has $m_i$ zeroes of degrees $d_i$, for $i=1,\ldots,k$. It was
shown by Masur and Veech that the stratum
$H(d_1^{m_1},\ldots,d_k^{m_k})$ is a complex orbifold of dimension
$n=2g-1+\sum_{i} m_i$.  By the Gauss--Bonnet formula, if $M$ has genus
$g$ then the degrees of the conical singularities must verify
\be
\label{eq.GauBon}
\sum_{i=1}^k m_i d_i=2g-2
\ef.
\ee
In fact, $H_g$ can be partitioned in terms of the strata
$H(d_1^{m_1},\ldots,d_k^{m_k})$ for $d_i$'s and $m_i$'s such that
(\ref{eq.GauBon}) is satisfied.

\begin{figure}[b!!]
\begin{center}
\setlength{\unitlength}{10pt}
\begin{picture}(19,17)(-1.5,.5)
\put(0,0){\includegraphics[scale=1]{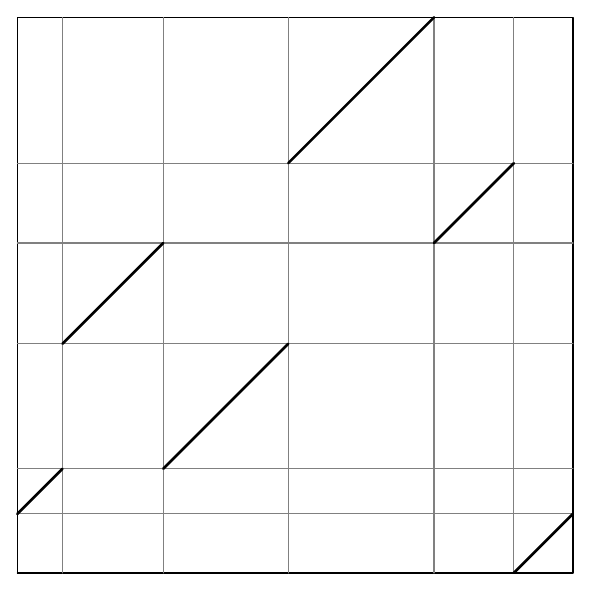}}
\put(.9,-.7){$I_1$}
\put(2.9,-.7){$I_2$}
\put(5.9,-.7){$I_3$}
\put(9.9,-.7){$I_4$}
\put(13.4,-.7){$I_5$}
\put(15.4,-.7){$I_6$}
\put(-.2, 1.4){\makebox[0pt][r]{$I'_{6\equiv \tau(6)}$}}
\put(-.2, 2.9){\makebox[0pt][r]{$I'_{1\equiv \tau(5)}$}}
\put(-.2, 5.4){\makebox[0pt][r]{$I'_{3\equiv \tau(4)}$}}
\put(-.2, 8.4){\makebox[0pt][r]{$I'_{2\equiv \tau(3)}$}}
\put(-.2,10.9){\makebox[0pt][r]{$I'_{5\equiv \tau(2)}$}}
\put(-.2,14.4){\makebox[0pt][r]{$I'_{4\equiv \tau(1)}$}}
\end{picture}
\qquad
\setlength{\unitlength}{14pt}
\begin{picture}(5,7)(0,-1)
\put(0,-1){1}
\put(1,-1){2}
\put(2,-1){3}
\put(3,-1){4}
\put(4,-1){5}
\put(5,-1){6}
\put(-1,0){1}
\put(-1,1){2}
\put(-1,2){3}
\put(-1,3){4}
\put(-1,4){5}
\put(-1,5){6}
\put(-.7,-.3){
\begin{tikzpicture}[scale=.5]
\permutation{5,3,4,1,2,6}
\end{tikzpicture}
}
\end{picture}
\end{center}
\caption[caption]{\label{fig_IETfun_1}Left: an IET $\phi$. Right: the
  associated permutation $\tau=\tau(\phi)$. We have a bullet at
  coordinate $(i,j)$ if $\tau(i)=j$. Note that there is a reflexion
  along the horizontal axis in the correspondence between $\tau$ and
  the obvious diagram interpretation of the graphics of $\phi$,
  i.e.\ $\tau$ describes the graphics of $-\phi(x)$.}
\end{figure}
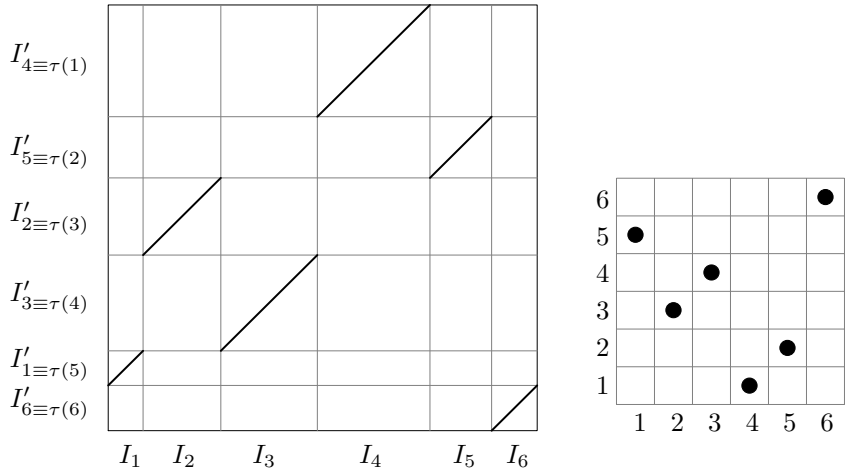

In term of polygons, two translation surfaces $S_1$ and $S_2$ with the
same area are in the same stratum $H(d_1^{m_1},\ldots,d_k^{m_k})$ if
and only if there exists a sequence of cut and paste moves from $S_1$
to $S_2$ (see figure~\ref{fig_cut_and_paste}). The `if' part is
obvious, while the `only if' part is a result of~\cite{Mas82}.

Let $M$ be a translation surface with a vertical direction, and let
$I$ be a horizontal open segment on $M$. The \emph{first return map}
$T$ (or \emph{Poincar\'e map}) of the translation flow\footnote{I.e.,
  the flow associated to the vertical foliation with singularities at
  the conical points defined by taking the integral curve of the
  parallel vector field i.e. defined by $\phi_t(x)=x+it$} from $I$ is
an \emph{interval exchange transformation} (IET).  An IET is a
one-to-one map $\phi$ from $I\setminus\{x_1,\ldots,x_k\}$ to
$I\setminus\{x'_1,\ldots,x'_k\}$ that is piecewise of the form
$\phi(z)=z+c_i$, i.e., there exists some minimal set of points
$\{x_1,\ldots,x_k\} \subset I$, and $\{x'_1,\ldots,x'_k\} \subset I'
\equiv \phi(I)$, such that, if we subdivide
$I\setminus\{x_1,\ldots,x_k\}$ into $k+1$ maximal open subintervals
$\{I_1,\ldots,I_{k+1}\}$, of lengths $\lambda_1,\ldots,\lambda_{k+1}$,
and $I'\setminus\{x_1',\ldots,x_k'\}$ into subintervals
$\{I'_1,\ldots,I'_{k+1}\}$ of the same lengths, the map $\phi$ is
identified by the datum of $\tau \in \mathfrak{S}_{k+1}$, and of
$(\lambda_1,\ldots,\lambda_{k+1}) \in (\mathbb{R}_+)^{k+1}$, with a
convention on the labelings as shown in the example of
Figure~\ref{fig_IETfun_1} (i.e., the intervals of
$I'\setminus\{x_1',\ldots,x_k'\}$, in their natural order from left to
right, have lengths $\lam_{\tau(k+1)}$, $\lam_{\tau(k)}$, \ldots,
$\lam_{\tau(1)}$).



\begin{figure}[t!!]
\begin{center}
 \includegraphics[scale=1]{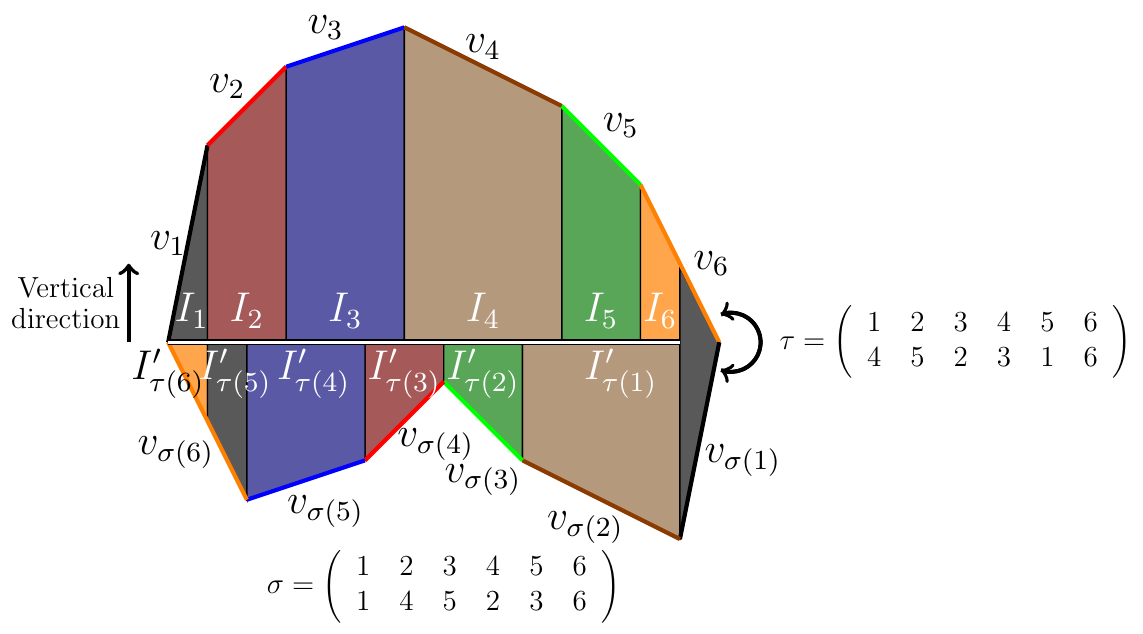}
\end{center}
\caption[caption]{\label{fig_IET_1}Let $\sigma$ and $(v_i)_{1\leq
    i\leq n}$ define a translation surface and choose a horizontal
  open interval $I$. The interval exchange map $T$ is obtained by
  following the vertical flow from almost each point of $I$, until it
  returns for the first time to $I$. The image by $T$ of a point of
  $I$ is not defined if its vertical flow meets a singularity before
  reaching $I$, and these are the points $x_i$ mentioned in the text.
  In the figure, the top sub-intervals $(I_i)$ of $I$ are numbered
  from left to right, and the bottom ones $(I'_{\tau(j)})$
  are numbered with $j$ increasing from right to left.}
\end{figure}

Consequently, also a first return map $T$ is parametrised by the same
combinatorial datum, of the permutation $\tau\in \mathfrak{S}_{k+1}$,
and a ``vector of lengths'' $\lambda=(\lambda_1,\ldots,\lambda_{k+1})$
(see figure \ref{fig_IET_1}). Moreover, the permutation $\tau$ is
irreducible (i.e., it has no non-trivial blocks) since the vertical
foliation is uniquely ergodic (see \cite{Vee82} and \cite{Mas82} for a
proof of the unique ergodicity).

The construction of a polygon with identified sides for a given datum
of $(\tau,\lam)$, with some useful geometric properties, involves the
notion of \emph{suspension}:
\begin{definition}[Suspension]
\label{def.suspdata}
Let $\s$ be an irreducible permutation of size $n$ and $\lambda$ a
vector of lengths. We say that the vectors $(v_i)_{i=1,\ldots,n} \in
\mathbb{C}^n$ are a \emph{suspension data} of $(\s,\lambda)$ if
\begin{itemize}
\item $\forall \;i$, $\reof(v_i)=\lambda_i$.
\item $\forall \;i<n$, $\sum_{j=1}^i \imof(v_j)>0$ and
  $\sum_{j=1}^i \imof(u_j)<0$.
\end{itemize}
\end{definition}

\begin{figure}[b!!]
\begin{center}
\makebox[0pt][l]{\raisebox{70pt}{$\eps=0$}}%
 \includegraphics[scale=.95]{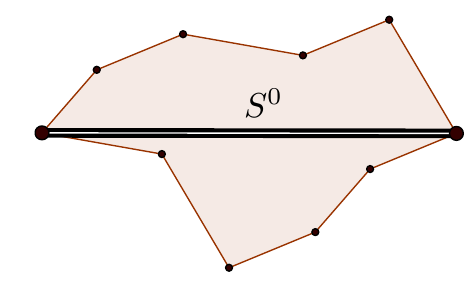}
\makebox[0pt][l]{\raisebox{70pt}{$\eps=+$}}%
 \includegraphics[scale=.95]{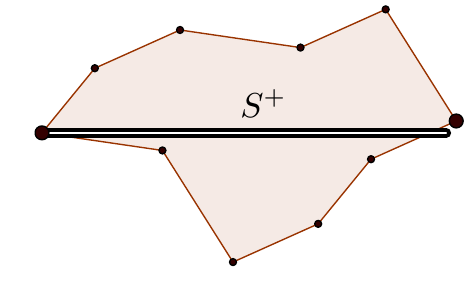}
\makebox[0pt][l]{\raisebox{70pt}{$\eps=-$}}%
 \includegraphics[scale=.95]{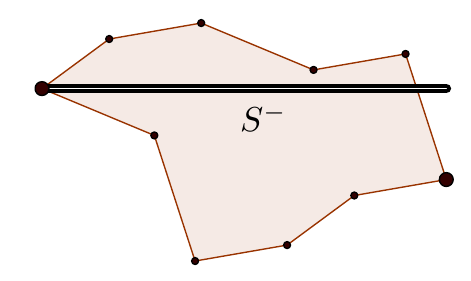}
\end{center}
\caption[caption]{\label{fig_suspension_1}Let $(\s,\lambda)$ describe
  an IET. The three figures describe three suspesions of
  $(\s,\lambda)$, $(v^{\eps}_i)_{1\leq i\leq n}$ for $\eps=0,+,-$. The
  second condition in the definition of a suspension means that for
  the translation surface $S^\eps=(\s,(v^\eps_i)_i)$ the top and
  bottom broken lines remain always above and below the real line,
  except possibly for their last segment. In the figure, the endpoint
  of the broken lines $P_+$ and $P-$ are on, above or below the real
  axis, for $\eps=0$, $+$ and $-$, respectively.}
\end{figure}

\begin{figure}[t!!]
\begin{center}
 \includegraphics[scale=.95]{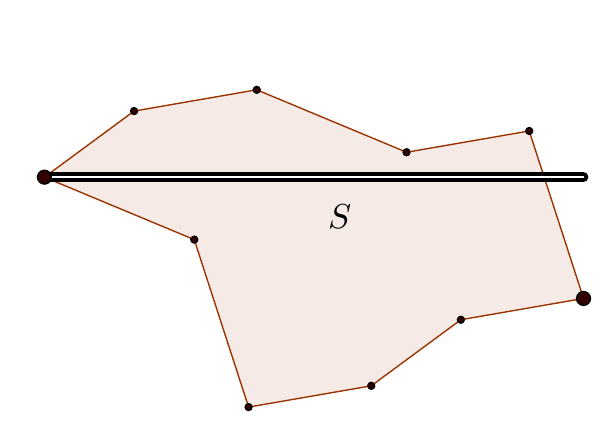}
 \includegraphics[scale=.95]{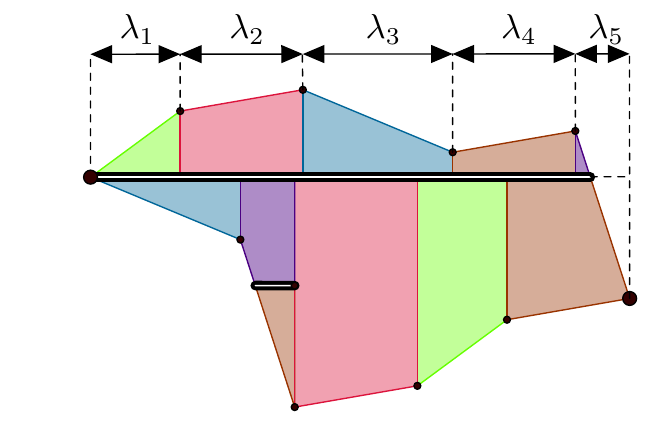}
\end{center}
\caption[caption]{\label{fig_suspension_2}As illustrated by Figure
  \ref{fig_suspension_1}, and here on the left picture, in the
  construction of the translation surface $S=(\s,(v_i)_i)$ from the
  suspension $(v_i)_i$ of $(\s,\lambda)$, the interval $I$ can be in
  part outside the polygon.  In this case, the construction of the
  interval exchange map involves a small subtlety, detailed on the
  right picture.  We cut the portion of $I$ outside the polygon, and
  translate it according to the identification of opposite sides.
  When this is done, the interval exchange map has a clear behaviour
  also on the right-most interval $I_n$.}
\end{figure}

\noindent
(Recall, $u_j=v_{\s(n+1-j)}$).  Given a suspension datum $(v_i)_i$ of
$(\s,\lambda)$ we can build the translation surface
$S=(\sigma,(v_i)_i)$ and choose the interval 
$]0,\sum_{i=1}^n \lambda_i[$ on the real line
so that the interval exchange map $T$ on $I$ is encoded exactly by
$(\sigma,\lambda)$. (See figures \ref{fig_suspension_1} and
\ref{fig_suspension_2}). In this case we say that $S$ is a suspension
of $(\sigma,\lambda)$. Veech proved that almost any translation
surface could be obtained as a suspension.

Let $(M,\omega) \in H(d_1^{m_1},\ldots,d_k^{m_k})$, and introduce the
shortcut $n=2g-1+\sum_i m_i$ for the dimension of the stratum.
In analogy with the previous section, we state that the Rauzy class of
$(M,\omega)$ is the subset of $\mathfrak{S}_{n}$ 
containing all the permutations which are the
combinatorial datum of some interval exchange map $T$ over an interval
$I$ on $(M,\omega)$. 

If $(\s,(v_i)_i)$ is a polygon representation of $S=(M,\omega)$, then
the Rauzy class of $S$ contains $\s$, since $S$ can be obtained as a
suspension of $(\s,\lambda)$ by Veech theorem, and is thus the class
$C(\s)$, in the notations of the previous section.

More generally, on
$H(d_1^{m_1},\ldots,d_k^{m_k})$ we have a notion of
connectivity due to the orbifold structure of this space.
Veech proved the two following things : \begin{itemize}
\item Two points $(M,\omega)$ and $(M',\omega')$ of the stratum
  $H(d_1^{m_1},\ldots,d_k^{m_k})$ are in the same connected component
  if and only if they are in the same extended Rauzy class.
\item Two permutations $\tau,\tau' \in \mathfrak{S}_{n}$ are in the
  same extended Rauzy class if and only if they are connected with
  respect to the extended Rauzy dynamics $\permsex_n$.
\end{itemize}
Thus, if $(M,\omega)$ and $(M',\omega') \in
H(d_1^{m_1},\ldots,d_k^{m_k})$ are suspensions of $(\tau,\lambda)$ and
$(\tau',\lambda')$, respectively, we have $(M,\omega) \sim
(M',\omega')$ iff $\tau \sim \tau'$ (the two $\sim$ symbols are
different, the first one is the connectivity within the orbifold, the
second one is the connectivity on the Cayley graph of the Rauzy
extended dynamics), and consequently the extended Rauzy classes
completely determine the connected components of the strata of~$H_g$.


Let us now interpret the combinatorial definition of the extended
Rauzy classes, given in Section \ref{ssec.3families}, in terms of the
\emph{Rauzy--Veech induction}. Let $M$ be a translation surface with a
vertical direction, $I$ an open interval and $T=(\tau,\lambda)$ the
associated interval exchange map. We label the top and bottom
sub-intervals induced by $T$ as in Figure \ref{fig_IET_1} above.
The Rauzy--Veech induction concerns the study of the possible
configurational changes in the structure of the IET, when the surface
$(M,\omega)$ is kept fixed, while the interval $I$ is modified.  A
\emph{right step} of the Rauzy--Veech induction consists in shortening
the interval $I$, from its right-end, by removing the shortest
sub-interval among $I_n$ and $I'_{\tau(1)}$. If we call $J$ the new
(total) interval, we have
\[
\left\{
\begin{array}{ll}
   J= I\setminus I'_{\tau(1)} & \text{if } |I_n|>|I'_{\tau(1)}|\\
   J= I\setminus I_n & \text{if } |I_n|<|I'_{\tau(1)}|
\end{array}
\right.
\]
The case $|I_n|=|I'_{\tau(1)}|$ never occurs if the $\lambda_i$'s are
generic in $(\mathbb{R}_+)^n$ (and thus they do not satisfy any linear
relation with coefficients in $\mathbb{Q}$), so we omit any discussion
of this possibility.  Conforming to the customary language, in the
first case of the analysis above we say that $T$ is
\emph{of right type 0}, and in the second case we say that $T$ is
\emph{of right type 1}. We also shorten right type with R-type.

\begin{figure}[b!!]
\begin{center}
 \includegraphics[scale=1]{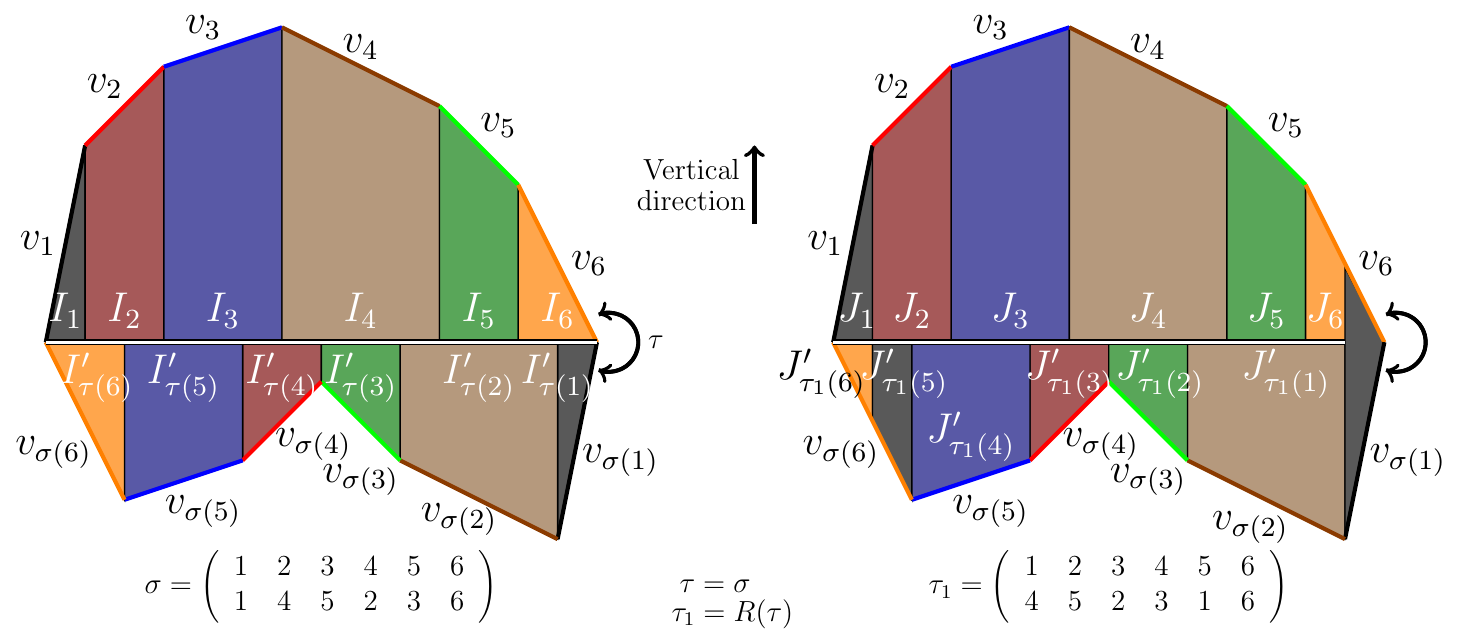}
\end{center}
\caption[caption]{\label{fig_IET_2}We apply the right step of the
  Rauzy--Veech induction on the surface shown on the left. In this
  case $T=(\tau,\lambda)$ is of right type 0, since $|I_6|>
  |I'_{\tau(1)}|$, thus $J=I\setminus I'_{\tau(1)}$. The new interval
  exchange map $T_1=(\tau_1,\lambda_1)$ on $J$ verifies
  $\tau_1=R(\tau)$, namely, in our example,
\\
\rule{.4\textwidth}{0pt}%
\makebox[0pt][c]{$
\displaystyle{\tau=\raisebox{-11pt}{\begin{tikzpicture}[scale=.15]\permutation{1,4,5,2,3,6}\end{tikzpicture}}
  \qquad \textrm{and} \qquad
  \tau_1=R(\tau)=\raisebox{-11pt}{\begin{tikzpicture}[scale=.15]\permutation{4,5,2,3,1,6}\end{tikzpicture}}\;.}$}
}
\end{figure}

Let us call $R_{\textrm{RV-ind}}$ the operation that sends the `old' IET
$T$ to the `new' one, $T_1=(\tau_1,\lambda_1)=R_{\textrm{RV-ind}}(T)$.  We
use this verbose notation in order to make a distinction with
operators $L$ and $R$ defined in Section \ref{ssec.3families}, whose role
here is as follows: if $T$ is of R-type 0, then $\tau_1=R(\tau)$,
while if $T$ is of R-type 1, then $\tau_1=L(\tau)$ (see figure
\ref{fig_IET_2}).

Likewise, a \emph{left step} of the Rauzy--Veech induction consists
in shortening the interval $I$ from the left-end. We have
\[
\left\{
\begin{array}{ll}
   J= I\setminus I'_{\tau(n)} & \text{if } |I_1|>|I'_{\tau(n)}|\\
   J= I\setminus I_1 & \text{if } |I_1|<|I'_{\tau(n)}| 
\end{array}
\right.
\]
and we say that $T$ is
\emph{of left type 0} in the first case and \emph{of left type 1} in
the second case.  Let us call 
$T_2=(\tau_2,\lambda_2) = L_{\textrm{RV-ind}}(T)$ the resulting new IET.
Then, the map $L_{\textrm{RV-ind}}$ has the property that $\tau_2=R'(\tau)$
or $L'(\tau)$ depending on the left-type of $T$, still with $L'$ and
$R'$ defined as in Section~\ref{ssec.3families}.
 
The extended Rauzy dynamic $\permsex_n$ acts on a permutation
$\sigma$, while the Rauzy--Veech induction acts on IET
$(\s,\lambda)$. As we have seen, the two operations are strictly
related.  Finally, we can define a third version of the operation, an
induction based on cut and paste operations (which we also call the
Rauzy--Veech induction), that acts on a suspension $(\sigma,(v_i)_i)$.

\begin{figure}[tb!]
\begin{center}
 \includegraphics[scale=1]{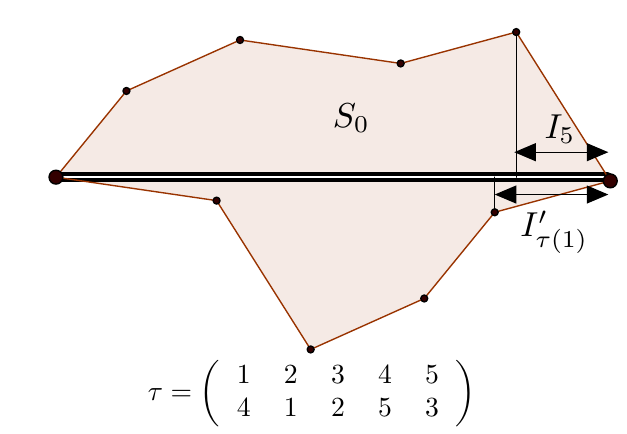}
 \includegraphics[scale=1]{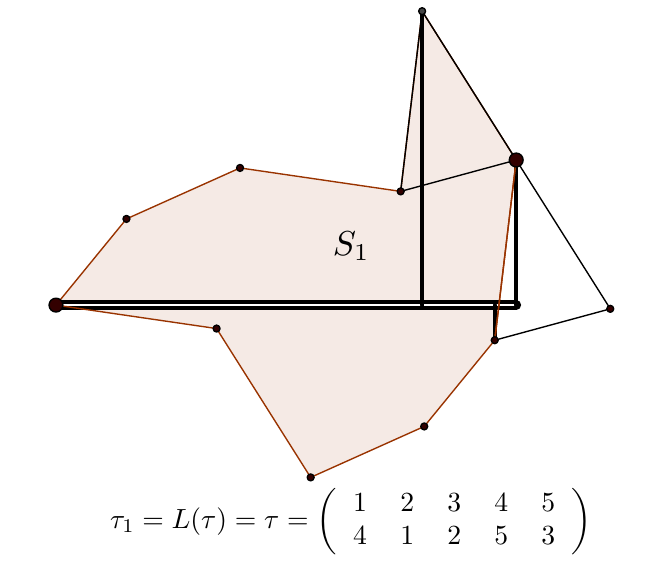}
 \includegraphics[scale=1]{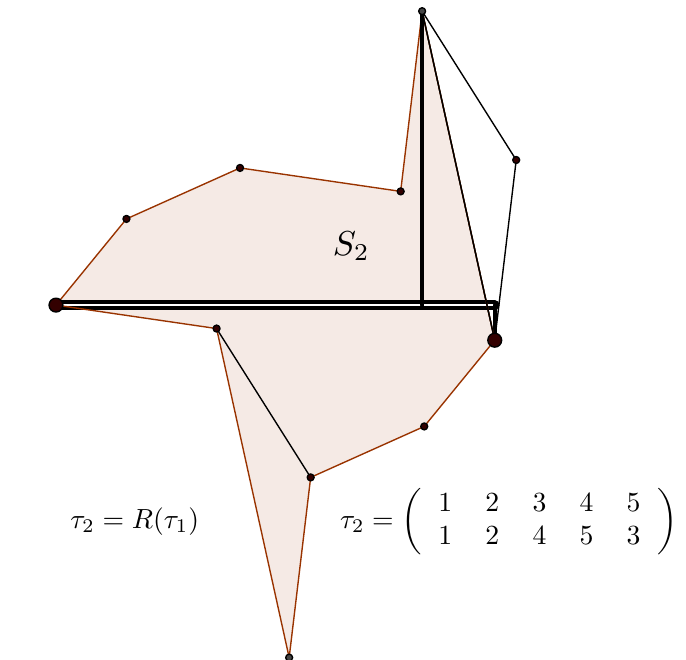}
\end{center}
\caption[caption]{\label{fig_suspension_cut_and_paste}We apply two
  times the right version of the Rauzy--Veech induction on the left
  figure. In the first case $T=(\tau,\lambda)$ is of right type 1,
  since $|I_5|< |I'_{\tau(1)}|$, thus $J=I\setminus I_{5}$. The new
  interval exchange map $T_1=(\tau_1,\lambda_1)$ on $J$ verifies
  $\tau_1=L(\tau)$.  We cut the triangle on the right part of $S_0$
  and paste it on the top broken line, since we are in type 1. In this
  way we obtain the suspension $S_1$ of $(\tau_1,\lambda_1)$.  Then,
  with a second right-step of the Rauzy--Veech induction for
  suspensions, we obtain the suspension $S_3$ from $S_2$.  In this
  case the step is of type 0, thus $\tau_2=R(\tau_1)$.}
\end{figure}

Given a suspension $(\sigma,(v_i)_i)$ of $(\sigma,\lambda)$, a right
step of the Rauzy--Veech induction for suspensions shortens the
interval $I$ into an interval $J$ as it would for (the right step of)
the Rauzy--Veech induction for IETs, then identifies a triangle within
the polygon, namely, the one with vertices
$p=\sum_{i=1}^{n} v_i$, $p-v_n$ and $p-v_{\sigma(1)}$,
and pastes it back on the polygon, accordingly to the identification
of the sides of the polygon and of the triangle, either on top of the
broken line $P_+$, or below the broken line $P_-$ in case of a Rauzy
induction of right type 1 or 0, respectively.  (see
Figure~\ref{fig_suspension_cut_and_paste}).

Likewise, a left step cuts the left-end of the interval, and deplaces
a triangle on the far left, with vertices $0$, $v_1$ and
$v_{\sigma(n)}$, in a way that depends on the type of the left step of
the Rauzy induction.

The relation of this version of the induction with the previous ones
is based on the following fact.  If $(\sigma,(v_i)_i)$ is a suspension
of $(\sigma,\lambda)$, $(\sigma_1,(v_i^1)_i)$ is the surface obtained
by the Rauzy--Veech induction for suspensions, and
$(\sigma_1,\lambda_1)$ is the IET obtained by Rauzy--Veech induction
for IETs.  Then, $(\sigma_1,(v_i^1)_i)$ is a suspension of
$(\sigma_1,\lambda_1)$.

As we have mentioned above, in \cite{Vee82} Veech showed that the
extended Rauzy classes were generated by the four operators
$L,R,L',R'$ that we have introduced, and, more generally, that for
every finite word $w$ in the alphabet
$\{L_{\textrm{RV-ind}},R_{\textrm{RV-ind}} \}$, and every word
$\omega$ of the same length in the alphabet $\{0,1\}$, there exists an
IET $T$ such that the induction performed according to the word $w$
visits IET's $T_i$ of type $\omega_i$ (namely, of left-type or
right-type $\omega_i$, according to $w_i$), and this crucial result is
at the basis of our combinatorial study of the classes determined by
the action of $L$, $R$, $L'$ and $R'$ on permutations.
In other words, if we call $\Pi_1$ the projection from a suspension
$S$ to the associated IET $T=(\tau,\lambda)$, and $\Pi_2$ the
projection from $T$ to the associated permutation $\tau$, we have the
commuting diagram
\be
\rule{80pt}{0pt}
\begin{CD}
\makebox[0pt][r]{\textrm{
\begin{minipage}{3.8cm}
\begin{center}
Rauzy--Veech induction\\
for suspensions:
\end{center}
\end{minipage}
}\qquad}
\makebox[50pt][c]{$S=(\tau,(v_i)_i)$}
@>\sim \textrm{~w.r.t.~}L_{\textrm{RV-ind}};\  R_{\textrm{RV-ind}}>
>
S'=(\tau',(v'_i)_i) 
\\
@VV{\Pi_1}V @VV{\Pi_1}V\\
\makebox[0pt][r]{\textrm{
\begin{minipage}{3.8cm}
\begin{center}
Rauzy--Veech induction\\
for IETs:
\end{center}
\end{minipage}
}\qquad}
\makebox[50pt][c]{$T=(\tau,\lambda)$}
@>\sim \textrm{~w.r.t.~}L_{\textrm{RV-ind}};\  R_{\textrm{RV-ind}}>
>
T'=(\tau',\lambda') 
\\
@VV{\Pi_2}V @VV{\Pi_2}V\\
\makebox[0pt][r]{\textrm{
\begin{minipage}{3.8cm}
\begin{center}
Combinatorial\\
operators:
\end{center}
\end{minipage}
}\qquad}
\makebox[50pt][c]{$\tau$}
@>\sim \textrm{~w.r.t.~}L, R;\  L', R'>> \tau'
\end{CD}
\label{eq.commuRVindu}
\ee

\bigskip

The (standard) Rauzy classes defined in Section \ref{ssec.3families}
also have a geometric interpretation. Given a stratum
$H(d_1^{m_1},\ldots,d_k^{m_k})$, we denote by
$H(d_1^{m_1},\ldots,\bar{d_i}^{m_i},\ldots,d_k^{m_k})$ the space with
a marked singularity of degree $d_i$. A point $(M,\omega)$ is in
$H(d_1^{m_1},\ldots,\bar{d_i}^{m_i},\ldots,d_k^{m_k})$ if, for
$(M,\omega)$ represented by a polygon as in Figure
\ref{fig_cons_surf}, the angle of the singularity of the equivalence
class containing the vertex $0$ is $2\pi(d_i+1)$. Thus
$H(d_1^{m_1},\ldots,d_k^{m_k})$ is partitioned into $k$ parts, one for
every possible value for the angle of the marked singularity.  For
example, in the figure \ref{fig_cons_surf} the angle around the
singularity passing through $(0,0)$ is $2\pi$, so the corresponding
translation surface $(M,\omega)$ in the figure is in
$H(\bar{0}^1,1^2)$.

Finally the \emph{Rauzy--Veech induction} is also defined and it
corresponds to the induction where only right steps are allowed. It
should be clear from figure \ref{fig_suspension_cut_and_paste} that
the degree of the singularity around the vertex 0 does not change
after a right move of the Rauzy--Veech induction for suspensions,
since the cut and paste at the level of the triangle can never change
the colour class of the left-most vertex, thus the colour class of the
left-most vertex is an invariant of the (non-extended) Rauzy dynamics,
and this corresponds to the marking.

In analogy to the forementioned result of Veech, Boissy proved in
\cite{Boi12} that the standard Rauzy classes are in one-to-one
correspondence with the connected components of the strata of the
moduli space of abelian differentials with a marked singularity, so
that we have a commuting diagram analogous to the one in
(\ref{eq.commuRVindu}), justifying the study of the combinatorial
dynamics~$\perms$.

To conclude this section, let us describe some more recent works
closely connected to the classification of the connected components of
the strata of the moduli space of abelian differentials of Kontsevich
and Zorich and of Boissy \cite{KZ03,Boi12}.

In \cite{Lan08}, Lanneau classified the connected components of the
strata of the moduli space of \emph{quadratic differentials}, Then,
Lanneau together with Boissy \cite{BL09} have formulated a
combinatorial definition of extended Rauzy classes for quadratic
differentials. These are defined as IETs on half-translation surfaces
and are once again in one-to-one correspondence with the connected
components of the strata of quadratic differentials. In this context
the combinatorial datum is no longer a permutation but a linear
involution (i.e.\ a matching with a one marked point between two given
arcs), and are related, to a certain extent, to the dynamics $\matchs$
of Section~\ref{ssec.3families}.

More recently, Boissy studied the \emph{labelled Rauzy classes} in
\cite{Boi13}, an object that we will also consider in forthcoming
work, because this notion will lead us to a second and independent
proof of the classification of the Rauzy classes, and it will emerge
as a general method to classify Rauzy-type dynamics in more general
circumstances. Boissy also classified the connected components of the
strata of the moduli space of \emph{meromorphic differentials} in
\cite{Boi15}. Delecroix studied the cardinality of Rauzy classes in
\cite{Del13}, and Zorich provided representatives of every Rauzy class
in the form of Jenkins--Strebel differentials \cite{Zor08}.  Finally,
in \cite{Fic16} Fickenscher has presented a combinatorial proof,
independent of the one presented here, of the Kontsevich-Zorich-Boissy
classification of Rauzy classes.

\subsection{Geometric interpretation of the invariants}
\label{ssec.geo_inv}

In this section we define the geometric version of two invariants of
the strata, and show that they correspond to our combinatorial
definitions of the \emph{cycle-} and the \emph{sign-}, or
\emph{arf-invariant} presented in Section~\ref{ssec.invardefs}.
We start by working at the level of the connected components of a
stratum $H(d_1^{m_1},\ldots,\bar{d_i}^{m_i},\ldots,d_k^{m_k})$ of
abelian differentials with a marked singularity.  

The first invariant that we introduce consists of an integer partition
with one marked part, and can be represented by the list
\[
(\underbrace{d_1+1,\ldots,d_1+1}_{m_1\text{~times}},
\ldots,
\underbrace{d_i+1,\ldots,d_i+1}_{m_i-1\text{~times}},
\ldots,
\underbrace{d_k+1,\ldots,d_k+1}_{m_k\text{~times}};
d_i+1)
\]
Each entry $d_i+1$ is associated to a conical singularity, and
corresponds to its angle, divided by $2\pi$. As said above, in the
standard (non-extended) dynamics there is a singularity passing
through the vertex $(0,0)$ of the polygon representation, and its
degree does not change in the dynamics (because the Rauzy--Veech
induction is acting on the other endpoint of the interval i.e.\ only
the right step of the induction is allowed for the standard dynamics),
and this accounts for the fact that one entry, the one after the
semicolon, is singled out. This invariant corresponds to the cycle
invariant
$(\lambda,r)$ of any permutation of the associated Rauzy class, as is
apparent from what we said above.
More precisely, let 
$(M,\omega) \in H(d_1^{m_1},\ldots,\bar{d_i}^{m_i},\ldots,d_k^{m_k})$,
let $(\sigma,(v_i)_{1\leq i \leq n})$ be a suspension data of
$(M,\omega)$ and let $C$ be the (standard) Rauzy class associated to
the connected component containing $(M,\omega)$. As we said before,
$\sigma\in C$, since the translation surface is a suspension of
$(\sigma,\mu)$ for some $\mu \in (\mathbb{R}_+)^n$.

\begin{figure}
\begin{center}
 \includegraphics[scale=.95]{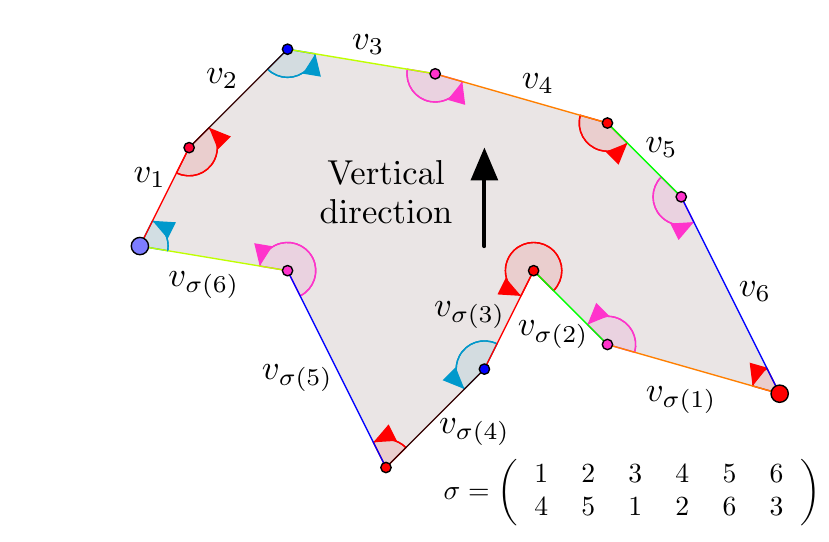}
 \includegraphics[scale=4.5]{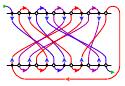}
 \caption[caption]{\label{fig_correspondance_cycle_angle}On the left,
   the polygon construction associated to the permutation
   $\s=(4,5,1,2,6,3)$, and a certain set of $v_i$'s. The arcs around
   the singularities are represented, with colours denoting their
   classes. On the right, the cycle invariant construction associated
   to the same permutation.  One can see that the arcs (from left to
   right) of the top broken line of the polygon correspond to the top
   arcs (from left to right) of the cycle invariant construction,
   while the arcs (now from right to left) of the bottom broken line
   of the polygon correspond to the bottom arcs (still from left to
   right) of the cycle invariant construction.  As a result, the angle
   of a singularity divided by $2\pi$ corresponds to the length of a
   cycle.}
 \end{center}
 \end{figure}

The results now follows from
Figure~\ref{fig_correspondance_cycle_angle}, which indicates the
correspondance between the angles around the vertices in the polygon
and the counting of top (or bottom) arcs of the cycle invariant in the
permutation.
 
The second invariant is more subtle, and corresponds to the parity of
the spin structure of the surface, that we now describe. Let $\gamma:
S^1 \to M$ be a smooth closed curve on $M$, avoiding the conical
singularities. Using as a reference our fixed parallel vector field
(i.e., the vertical direction), we can define the \emph{Gauss map} $G$
from $\gamma$ to $[0,2\pi[ \simeq S^1$
as follows: to every point of the curve, $\gamma(x)$, we associate the
angle $\theta(x) \in S^1$ corresponding to the angle between the
tangent vector $\gamma'(x)$ and the vector field in $x$. The index
$\ind(\gamma)$, called the \emph{degree} of the Gauss map, is the
integer counting (with sign) the number of times the curve
$G \cdot \gamma$ is wrapped around~$S^1$.

This allows us to define:
\begin{definition}[Parity of the spin structure] Let
  $(\alpha_i,\beta_i)_{1\leq i\leq g}$ be a collection of smooth
  closed curves avoiding the singularities, and representing a
  symplectic basis for the homology $H_1(M,\mathbb{Z})$. The
  \emph{parity of the spin structure} of $M$ is the quantity
\be
\Phi(M)= (-1)^{\sum_{i=1}^g (\ind(\alpha_i)+1)(\ind(\beta_i)+1)}
\label{eq.arfxx1}
\ee
\end{definition}
\noindent
It follows from \cite{Ati71} that this is an invariant of the
connected components of the stratum, and from \cite{John80} that it is
well-defined and independent of the choice of symplectic basis. More
precisely, Johnson in \cite{John80} shows that the quantity is the
\emph{Arf invariant} of a certain quadratic form
$\Phi: H_1(M,\mathbb{Z}_2) \to \mathbb{Z}_2$. We will describe this
second point of view in a few paragraphs, and we now proceed to
construct $\Phi$ concretely in the special case of translation
surfaces.

Let $c$ be a cycle of $H_1(M,\mathbb{Z}_2)$, and let $\gamma$ be a
smooth simple closed curve avoiding the singularities, and
representing $c$.  We define the function $\Phi: H_1(S,\mathbb{Z}_2)
\to \mathbb{Z}_2$ as
\[
\Phi(c)=\ind(\gamma)+1 \quad (\textrm{mod}~2).
\]
This function is only well-defined when all the zeroes of $\omega$
have even degree (i.e.~%
$(M,\omega)\in H(2d_1^{m_1},\ldots,\bar{2d_i}^{m_i},\ldots,2d_k^{m_k})$),
otherwise two curves representing the same cycle $c$ might have
opposite index.

Suppose we are in the even-degree case. For $c$ and $c'$ distinct
cycles, and $\gamma$ and $\gamma'$ representing $c$ and $c'$
respectively, and in generic mutual position (i.e., they may cross,
but only in generic way), we define the bilinear
\emph{intersection form} $\Omega(c,c')$ as the number of intersections
between $\gamma$ and $\gamma'$, mod~2.

The following theorem (obtained as a corollary of a theorem of \cite{John80})
certifies that our function is a quadratic form:
\begin{theorem}
\label{def_form_quadra}
The function $\Phi$ is well-defined on $H_1(M,\mathbb{Z}_2)$. It is a
quadratic form associated to the bilinear intersection
form $\Omega$
in the
following sense: 
\[
\forall\ c,c'\qquad  
\Phi(c+c')=\Phi(c)+\Phi(c')+\Omega(c,c')
\ef.
\]
\end{theorem}
\noindent
For our quadratic form $\Phi: H_1(S,\mathbb{Z}_2) \to \mathbb{Z}_2$
the Arf invariant is defined as
\be
\label{arf_definition_geo}
\arf(\Phi)=
\frac{1}{\sqrt{|H_1(S,\mathbb{Z}_2)|}}
\sum_{x\in H_1(S,\mathbb{Z}_2)} (-1)^{\Phi(x)}
=2^{-g}
\sum_{x\in H_1(M,\mathbb{Z}_2)} (-1)^{\Phi(x)} 
\ef.
\ee
Let $(a_i,b_i)_{1\leq i\leq g}$ be a sympletic basis of
$H_1(S,\mathbb{Z}_2)$, then it can be shown
that
\be
\arf(\Phi)=(-1)^{\sum_{i=1}^g \Phi(a_i)\Phi(b_i)}
\ef.
\label{eq.arfxx2}
\ee
Let $(\alpha_i, \beta_i)_{1\leq i \leq g}$ be a family of smooth
closed curves (avoiding the singularities) representing the family of
cycles $(a_i,b_i)_{1\leq i\leq g}$. Since we defined our quadratic
form by $\Phi(c)= \ind(\gamma)+1 \mod 2$ for any cycle $c$ and curve
$\gamma$ representing $c$, by comparing (\ref{eq.arfxx1}) and
(\ref{eq.arfxx2}) we have
\[
\arf(\Phi)=\Phi(S)
\ef,
\]
i.e., the parity of the spin structure is the Arf invariant of $\Phi$.
 
Let us now describe how to obtain our definition of the Arf invariant,
given in section \ref{ssec.arf_inv}, starting from the formula in
(\ref{arf_definition_geo}).

\begin{figure}[tb!]
\begin{center}
 \includegraphics[scale=1]{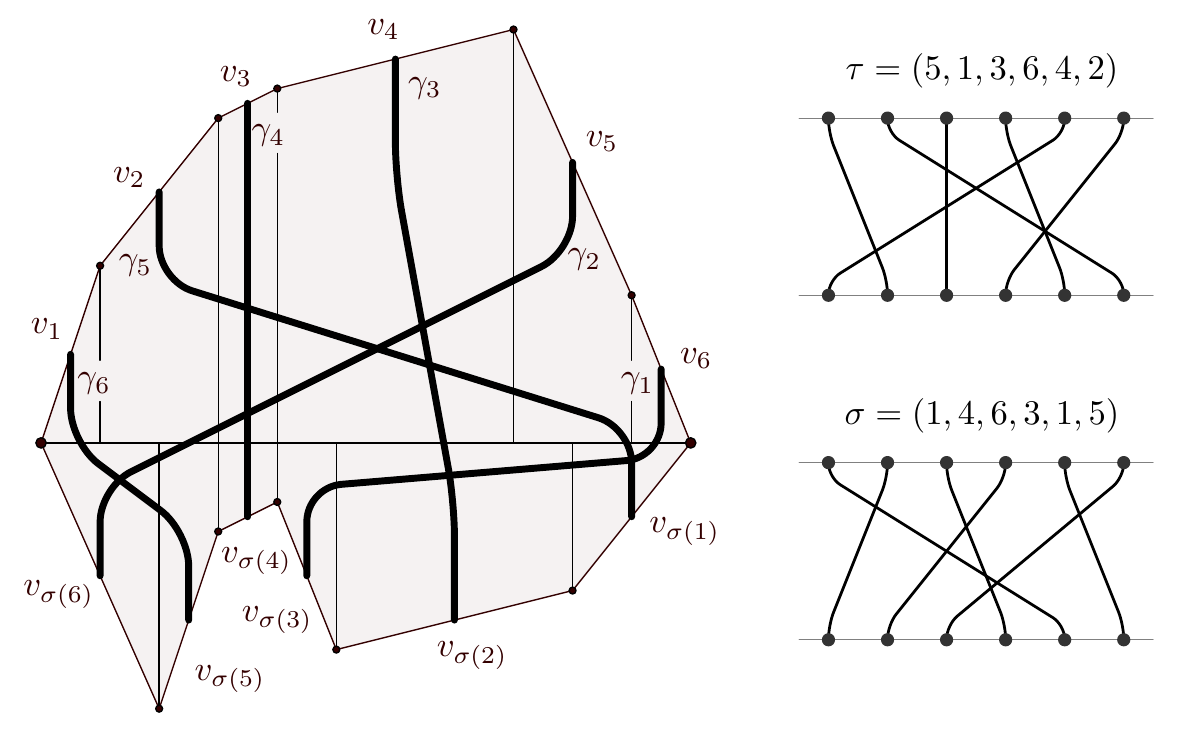}
\caption[caption]{\label{surface_with_curve} On the left, the
  polygon construction associated to the permutation
  $\s=(1,4,6,3,1,5)$ with a collection of simple smooth closed curves
  $(\gamma_i)_{1\leq i\leq n}$.  On the right, the permutation
  $\tau=(5,1,3,6,4,2)$ represents the collection of curves
  $(\gamma_i)_{1\leq i\leq n}$: two closed curves $\gamma_i$ and
  $\gamma_j$ intersect if and only if the edges $(i,\tau(i))$ and
  $(j,\tau(j))$ of $\tau$ intersect, since the set of curves
  $(\gamma_i)_{1\leq i \leq n}$ and the permutation $\tau$ have the
  same topology.  The permutation $\s$ is the reverse of $\tau$,
  i.e.\ it is obtained from $\tau$ by multiplying with the longest
  permutation $\coxe=(n,n-1,\ldots,1)$, thus two closed curves
  $\gamma_i$ and $\gamma_j$ intersect if and only if the edges
  $(n-i+1,\s(n-i+1))$ and $(n-j+1,\s(n-j+1))$ of $\s$ do \emph{not}
  intersect.}
 \end{center}
\end{figure}

Let $(M,\omega)$ be a translation surface defined as the suspension
$(\s,(v_i)_i)$ of an IET $T$. Say that $\sigma(i)=j$.  For $x \in [0,1]$
sufficiently small, there exists a curve $\gamma_{n+1-j}$ from the
point $p_i^- = v_{\sigma(n)}+v_{\sigma(n-1)}+\cdots+x\, v_{\sigma(i)}$ on
$P_-$ to the point $p_j^+ = v_{1}+v_{2}+\cdots+x\, v_{j}$ on $P_+$, such
that the Gauss map $G$ discussed above remains within the interval
$[-\pi/2,\pi/2]$. Indeed, this fact is true for all $x$, whenever
$j<n$, but only for $x$ small enough for $j=n$ and when the suspension
has $\eps \neq 0$ w.r.t.\ the notations of Figure
\ref{fig_suspension_1}. This can be seen from the existence of a
positive $\delta$ such that a strip of width $\delta$ around the
interval $I$ leaves the points $p_i^-$ and the points $p_i^+$ on
opposite half-planes.  Figure \ref{surface_with_curve} provides an
example of this construction.

Because of the identification of the sides of the polygon, the
$\gamma_i$'s are closed curves.  As $G \cdot \gamma_i(x) \in
[-\pi/2,\pi/2]$, it is easily seen that they all have index
zero. 
Let $(c_i)_i$ be the family of cycles represented by the curves
$(\gamma_i)_i$, then it can be shown that 
\begin{lemma}
The family of cycles $\{c_1,\ldots,c_n\}$ is a generating family of
$H_1(M,\mathbb{Z}_2)$.
\end{lemma}
 
\noindent
For a permutation $\pi$, define the symmetric matrix $K(\pi)$ to be
zero on the diagonal, and elsewhere
\[
\big( K(\pi)_{ij} \big)_{1\leq i,j\leq n}=
\left\{\begin{array}{ll}
1 & \text{if $(\pi(i)-\pi(j) )(i-j) < 0$} \\
0 & \text{if $(\pi(i)-\pi(j) )(i-j) > 0$}
\end{array}
\right.
\]
Then, for all $c_i$ and $c_j$ ($i \neq j$),
$\Omega(c_i,c_j)=K(\tau)_{i,j}$ (see figure \ref{surface_with_curve}
for this correspondence).  Defining $\s$ as the (right-)reverse of
$\tau$, i.e.\ $\s=\tau \coxe$ where $\coxe=(n,n-1,\ldots,1)$ is the
longest permutation, we also have $\Omega(c_i,c_j)=1-K(\s)_{i,j}$.
From the fact, proven above, that $\ind(\gamma_j)=0$, we get
$\Phi(c_j)=\ind(\gamma_j)+1 \mod 2 =1$.

For a set $I\subseteq \{1,\ldots,n\}$ we define
$\Phi(I)=\Phi(\sum_{i\in I} c_i)$ and 
$\chi'_I(\pi)= \sum_{\{i,j\}\subset I} K(\pi)_{ij}$ as the number of
pairs $i<j$ of edge labels, both in $I$, which are crossing in $\pi$.
In Section \ref{ssec.arf_inv_intro} we already have defined the
analogous quantity
$\chi_I(\pi)= \sum_{(i,j)\subset I} \big( 1-K(\pi)_{ij} \big)$, as the
number of pairs $i<j$ in $I$ which do not cross in $\pi$. Then we have
\begin{proposition}
If the surface $S$ is in the Rauzy class $C(\s)$,
for all $I\subset \{1,\ldots,n\}$ we have 
\be
\Phi(I)=|I|+\chi'(I)(\s \coxe)=|I|+\chi(I)(\s)
\ef.
\ee
\end{proposition}
\proof
Call again $\tau= \s \coxe$.
We have, by iterated application of Theorem \ref{def_form_quadra},
\be
\Phi(I)
:=\Phi \Big( \sum_{i\in I} c_i \Big)
= \sum_{i\in I} \Phi(c_i) +\sum_{(i,j)\subseteq I} K(\tau)_{ij}
\ef.
\ee
Recall that, on one side, $\Phi(c_i)=1$ for all $i$, and, on the other
side, $\sum_{(i,j)\subseteq I} K(\tau)_{ij}$ is the definition of
$\chi'(I)(\tau)$. Thus
\be
\Phi(I)
=
|I|+\chi'(I)(\tau)
\ef.
\ee
Similarly for $\s$, 
as $K(\tau)_{ij} = 1-K(\s)_{ij}$
and
$\sum_{(i,j)\subseteq I} \big(1-K(\s)_{ij}\big)$ corresponds to the
definition of $\chi(I)(\s)$,
\be
\Phi(I)
=
|I|+\chi(I)(\s)
\ef.
\ee
\qed

\noindent
As a result, we can rewrite the definition (\ref{eq.arf1stDef}) of
$\Abar(\s)$ as
\be
\bar{A}(\sigma)
:= \sum_{I \subseteq \{1,\ldots,n\}} (-1)^{|I| + \chi(I)}
=\sum_{I\subseteq \{1,\ldots,n\}} (-1)^{\Phi(I)}
\ee 
which, identifying sets $I$ with cycles $\sum_{i \in I} c_i$, is our
definition (\ref{eq.arfxx2}) of the Arf invariant, up to an overall
factor that we shall now discuss.
Let us start by observing:
  
\begin{lemma}\label{deux_Gauss_B}
Let $(M,\omega)\in
H(d_1^{m_1},\ldots,\bar{d_i}^{m_i},\ldots,d_k^{m_k})$ be a translation
surface of genus g represented by a permutation $\sigma$ of size
$n$. Finally, let $k'=\sum_{j=1}^l m_j$ be the total number of
singularities, then
\be
\label{genus_size}
n=2g+k'-1
\ef.
\ee 
\end{lemma}

\proof
This follows from the Gauss--Bonnet and the dimension formulas
\begin{align}
\sum_{j=1}^k m_j d_j &=
2g-2
\ef;
&
r+\sum_{j=1}^k m_j \lambda_j &=
n-1
\ef;
\end{align}
with $r=\bar{d}_i+1$ and $\lambda_j=d_j+1$ for the non-marked
singularities $z_j$, since
\[
n-1=r+\sum_{j=1}^k m_j d_j+1= 2g-2 +\sum_{j=1}^l m_j=2g-2+k'
\]
\qed

\noindent
Note that the integer $n \equiv 2g+k'-1$ on the two sides of
(\ref{genus_size}) is exactly the dimension of the complex orbifold
$H(d_1^{m_1},\ldots,\bar{d_i}^{m_i},\ldots,d_k^{m_k})$. We have
\begin{proposition}
\[
\bar{A}(\s)=2^{k'-1+g}
\, \arf(\Phi)
\ef.
\]
\end{proposition}
\proof
We already know that $\{c_i,\ldots,c_n\}$ is a generating set. Since
the basis of $H_1(S,\mathbb{Z}_2)$ has
size $2g$, and $n=2g+k'-1$, whenever $k'>1$ this family
is linearly dependent on $\mathbb{Z}_2$, thus there exists an index
$j$ and a set $I$ not containing $j$
such that $c_j=\sum_{i\in I} c_i$. 
Call $X=\{1,\ldots,\hat{j},\ldots,n\}$.
Denoting by
$\Delta$ the symmetric difference of two sets,
we can rearrange the summands in
$\bar{A}(\s)$ as
\[
\begin{array}{lcl}
\bar{A}(\s) & = & 
\sum_{J\subseteq X}
(-1)^{\Phi(J)}
+\sum_{J\subseteq X} (-1)^{\Phi(J\cup\{j\})}
\\  
& = & \sum_{J\subseteq X}
(-1)^{\Phi(J)}
+\sum_{J\subseteq X} (-1)^{\Phi(J\Delta I)}
\\
& = & 
2\sum_{J\subseteq X} (-1)^{\Phi(J)}
\ef,
\end{array}
\]
where in the last passage we use the fact that, as $J$ ranges over
subsets of $X$ and $I \subseteq X$, then also 
$J'=J\Delta I$ ranges over
subsets of $X$.

Hence by induction, performing $k'-1$ steps,
and up to relabeling the indices of the cycles so that the remaining
family $(c_1,\ldots,c_{2g})$ forms a basis,
\[
\bar{A}(\s)=2^{k'-1} \sum_{J\subseteq \{1,\ldots,2g\}}
(-1)^{\Phi(J)}=2^{k'-1+g}2^{-g}\sum_{x\in H_1(M,\mathbb{Z}_2)}
(-1)^{\Phi(x)}=2^{k'-1+g}\arf(\Phi)
\ef.
\]
\qed

\noindent
Recall that in Proposition \ref{pro_sign_inv} we claimed that
$s(\s)=2^{-\frac{n+k'-1}{2}}\bar{A}(\s)$ (this fact is proven later
on), thus, when $\Phi$ is defined, that is, for 
$(M,\omega) \in H(2d_1^{m_1},\ldots,\bar{2d_i}^{m_i},2d_k^{m_k})$,
this quantity is equal to $\arf(\Phi)$, since
$\frac{n+k'-1}{2}=\frac{2g-2+2k'}{2}=k'-1+g$ by
Lemma~\ref{deux_Gauss_B}.

When $(M,\omega) \notin
H(2d_1^{m_1},\ldots,\bar{2d_i}^{m_i},2d_k^{m_k})$, i.e.\ when $\omega$
has some zeroes of odd degree, the function $\Phi$ is not
defined. Nonetheless the function $s(\s)$ remains well-defined, and
its value is $0$ in this case (see the following Lemma
\ref{lem.induSign} for a proof).

In summary, we have
\[
\begin{array}{lcl}
s(\s)\in\{\pm 1\} &\Longleftrightarrow & s(\s)=\arf(\Phi)
\textrm{~and~} (M,\omega) \in
H(2d_1^{m_1},\ldots,\bar{2d_i}^{m_i},2d_k^{m_k})\ef;
\\
 s(\s)\in\{0\} &\Longleftrightarrow & \text{An even number of conical singularities
   with angle multiple of $4\pi$. }
\end{array}
\]

\subsection{Surgery operators}
\label{ssec.surgery_geo_combi}

\begin{figure}[tb!]
\begin{center}
 \includegraphics[scale=.8]{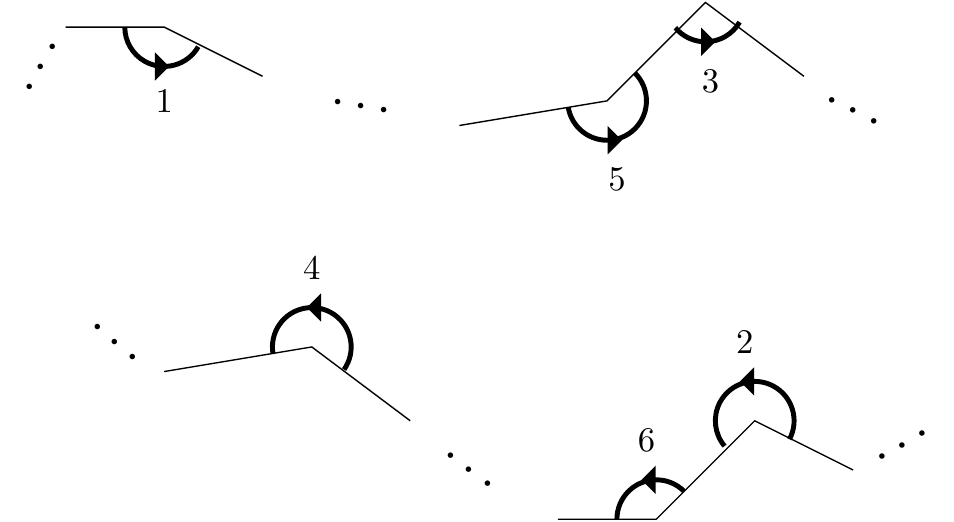}\\[2mm]
\raisebox{1.48ex}{ \includegraphics[scale=.8]{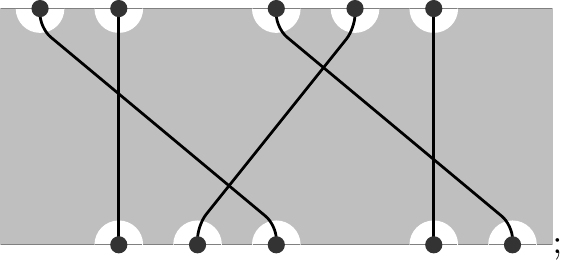}}
\includegraphics[scale=.8]{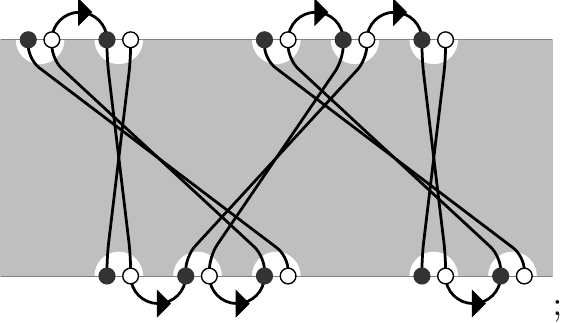}
\end{center}
\caption{\label{fig_6pi}Top: a translation surface with a conical
  singularity of angle $6\pi$. Bottom: the associated permutation.}
\end{figure}

In this section we describe two families of `geometric surgery
operations' which can be performed on translation surfaces, and more
notably on suspensions. Some special cases of these will correspond to
the combinatorial surgery operators, introduced in Section
\ref{ssec.TQSintro}, and crucially used in the main body of the
paper. 

\begin{figure}[b!!]
\begin{center}
 \includegraphics[scale=.8]{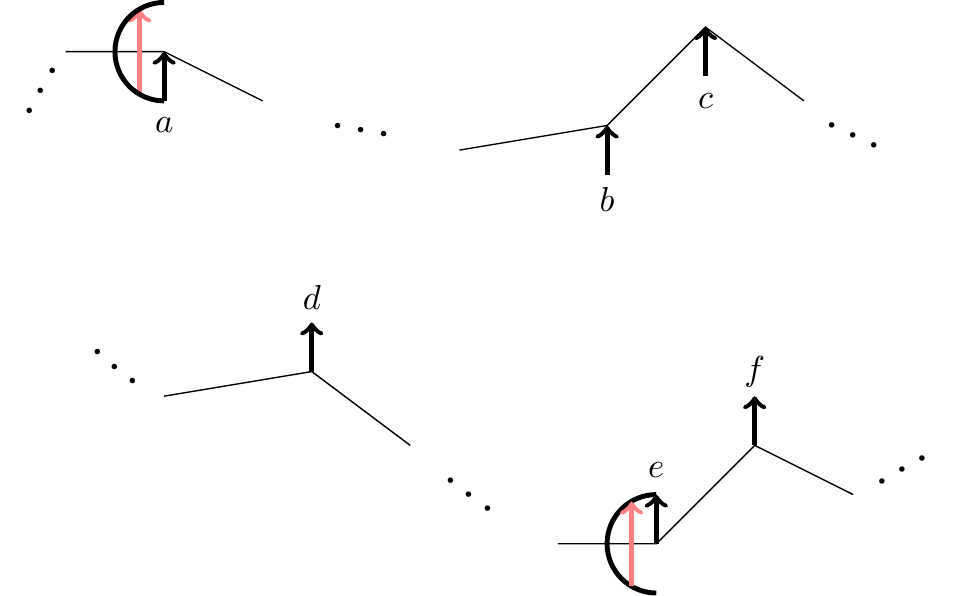}  $\quad$
\raisebox{6ex}{  \includegraphics[scale=1]{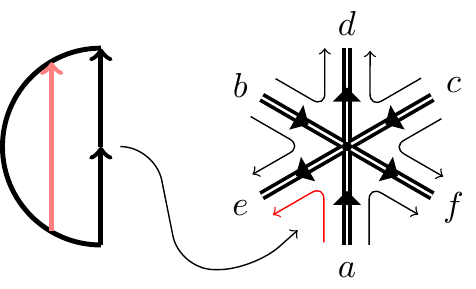}}
\end{center}
\caption{\label{fig_6pi_nei}On the left picture, the black arrows
  represent the leaves of the foliation on the critical point
  $z_i$. Since $z_i$ has conical angle $6\pi$, there are 6 half-discs
  (of radius $\epsilon$) forming together the neighbourhood of
  $z_i$. One such half-disc, adjacent to the leaves $a$ and $e$, is
  represented. On the right picture, we have a representation of the
  leaves of the vertical foliation in the neighbourhood of the
  singularity $z_i$. Each angular sector in the picture is of angle
  $\pi/3$, although it corresponds to a half-disc of radius
  $\epsilon$. The red trajectory, associated to a leaf passing in a
  neighbourhood of the singularity, helps in visualising the glueing
  procedure.}
\end{figure}

For this purpose it will be useful to introduce a convention on the
visualisation of conical singularities.  We have seen that the degree
$d_i$ of a zero $z_i$ can be computed directly on the translation
surface, either by adding the angles of the polygon which are
associated to the singularity (then dividing by $2 \pi$, and
subtracting 1) or also, more conveniently, just by counting the number
of top (or bottom) arcs around the vertices of the singularity, minus
1 (cf.\ figure \ref{fig_cons_surf} and figure~\ref{fig_6pi}).  In
order to visualise a singularity, say at $z_i$, it is useful to
consider a small neighbourhood, say of radius $\epsilon \ll |v_i|$.
Such a neighbourhood is isometric to a cone of angle $2\pi(d_i+1)$,
i.e.\ the glueing of $2(d_i+1)$ half-disks of radius $\epsilon$.  If,
in a picture, we take care of drawing some trajectories of the
vertical foliation, we can trade half-disks for angular sectors of
circle with smaller angles $\alpha_j$, so that, in total, we can
squeeze the $2\pi(d_i+1)$ cone into an ordinary disk (by choosing
$\alpha_j$'s such that $\sum_{j=1}^{2(d_i+1)} \alpha_j = 2 \pi$), and
still provide a picture in which the foliation is reconstructed
unambiguously.  Such a representation is illustrated in
figure~\ref{fig_6pi_nei}.

\begin{figure}[tb!]
\[
\includegraphics[scale=.8]{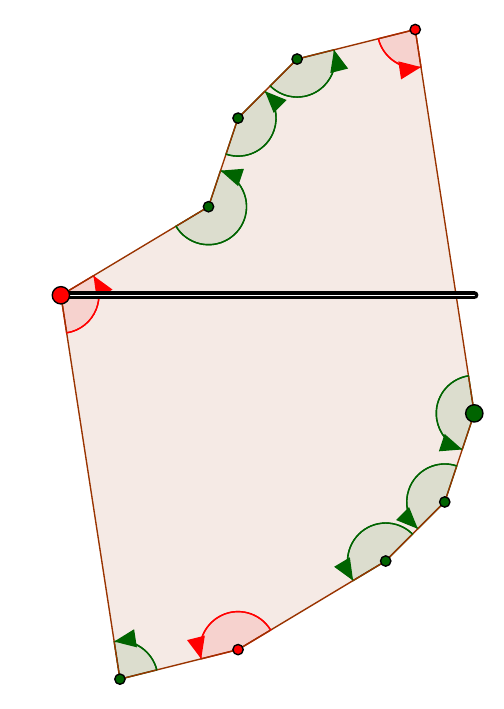}
\hspace*{-5pt}\raisebox{80pt}{$\left. \rule{0pt}{16pt} \right\} y$}
\]
\caption{A regular suspension for $\tau=(23145)$. This permutation has
  rank 1, and one cycle of length 3, i.e.\ the marked singularity has
  degree zero (red angles) and there is a single further singularity,
  of degree two (green angles).\label{fig.exRegSusp}}
\end{figure}

We now introduce a family of suspensions, with certain special
properties guaranteeing that the behaviour of the geometric surgery
operations is more simple than what is the case on a generic
suspension: (recall, $u_j$ are synonims for $v_{\tau(n+1-j)}$)
\begin{definition}[Regular suspension]
\label{def.regsusp}
We say that a suspension $S=(\tau;v_1,\ldots,v_n)$ is \emph{regular},
and \emph{of shift $y$} if
\begin{itemize}
\item $-\frac{\pi}{2} < \arg(u_1) < 0 < \arg(u_2) < \arg(u_3) < \cdots
  < \arg(u_n) < \frac{\pi}{2}$,
\item $\imof (v_1+v_2+\cdots+v_n) = \imof (u_1+u_2+\cdots+u_n) = -y$, with $y\geq 0$,
\item $\tau(n)=n$.
\end{itemize}
\end{definition}
\noindent
An example is given in Figure~\ref{fig.exRegSusp}.
In particular, a regular suspension, with reference to the notation of
Figure \ref{fig_suspension_1}, has $\eps \in \{0,-\}$.

\subsubsection{Adding a cylinder}

The first family of geometric surgery operators consists in
\emph{breaking up a zero} (of the abelian differential $\omega$), as
discussed, e.g., in \cite[sect.\ 4.2]{KZ03} and \cite{EMZ03}.
Geometrically, this consists in adding a cylinder on one given
singularity of the translation surface.

Let $S$ be a translation surface, represented as a polygon with
identified sides, and let $z_i$ be a zero of degree $d_i$. A
\emph{saddle connection} from $z_i$ to $z_i$ is a geodesic (i.e., in
our case, a straight line segment which may wrap along the identified
sides of the polygon) connecting a vertex of the equivalence class of
$z_i$ on the top broken line to a vertex of the equivalence class of
$z_i$ on the bottom broken line. So, for a zero of degree $\ell=d+1$,
there are $\ell^2$ families of saddle connections, one per top/bottom
choice of angle, with segments in the same family differing for the
choice of the winding.  

In order to break up a zero $z_i$, we choose one such saddle
connection,
in one such family, and ``add a cylinder''.  This operator may be
denoted as $\mathcal{Q}_{i,h,k,s}(v)$, with $i$ an index associated to
the singularity, $h, k \in \{1,\ldots,d_i+1\}$ associated to the
choice of top/bottom angle, $s$ associated to the choice of saddle
connection within the family, and $v$ a two-dimensional vector with
positive scalar product with the versor of the saddle connection.

\begin{figure}[tb!]
\begin{center}
 \includegraphics[scale=.8]{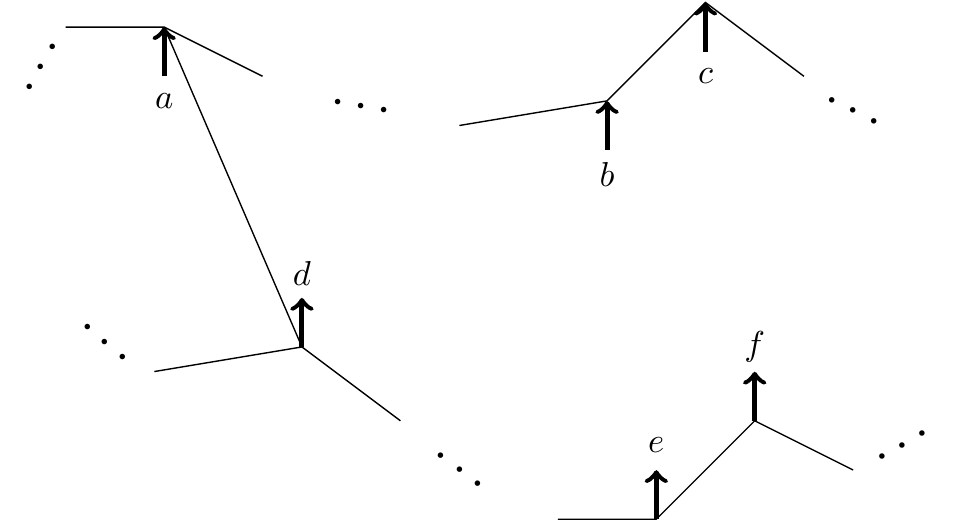} $\quad$
 \raisebox{2.5ex}{\includegraphics[scale=1.2]{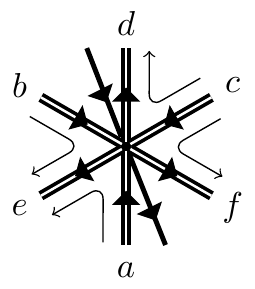}}
 \includegraphics[scale=.8]{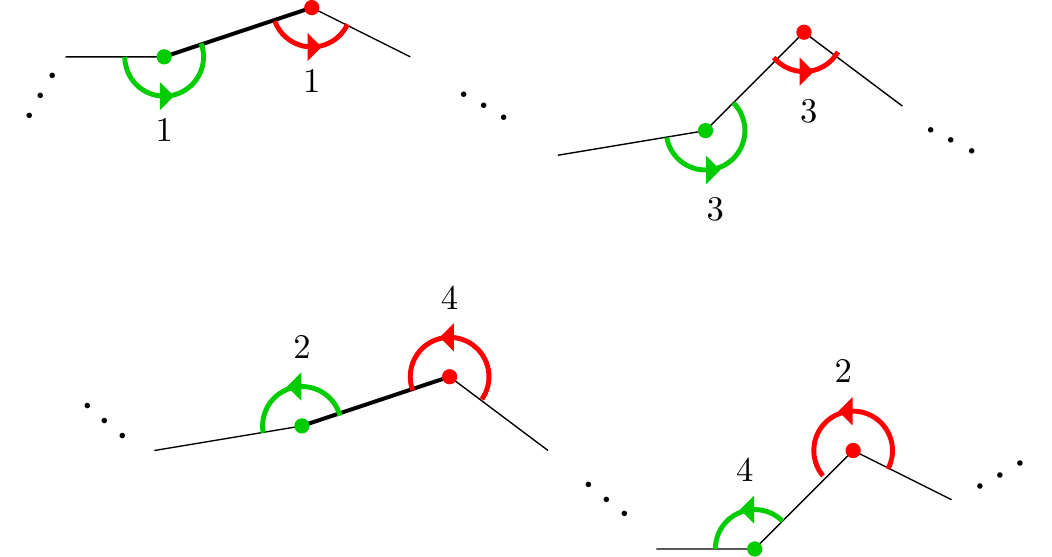}\\[2mm]
\raisebox{1.48ex}{ \includegraphics[scale=.8]{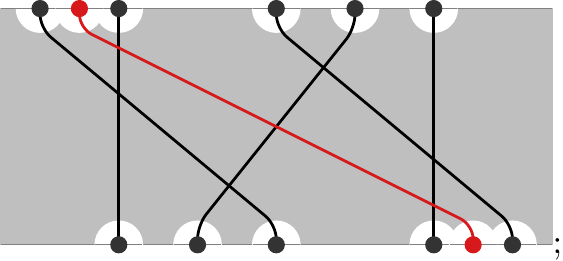}}
\includegraphics[scale=.8]{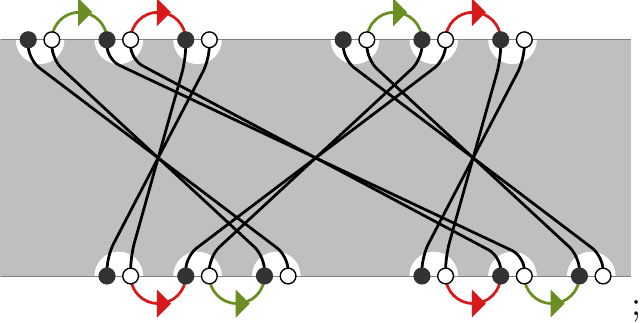}
\end{center}
\caption{\label{fig_add_cylinder_1}Top left: a saddle connection
  between two vertices of the conical singularity of angle $6\pi$ on
  the translation surface. Top right: the same geodesic in a
  neighbourhood of $z_i$. Middle: we cut the surface along the
  geodesic and glue the two boundary of a cylinder to the two newly
  formed boundaries of the surface. Adding the cylinder breaks up the
  singularity of angle $6\pi$ into two singularities of angle $4\pi$
  (in green and red respectively). Bottom: the projection on the
  permutation of this procedure consists in adding an edge (in red in
  the figure) in between two arcs belonging to the same cycle (cf
  figure \ref{fig_6pi} for the drawing of the cycle invariant).}
\end{figure}

This construction is especially simple if, for the given polygon $S$
and triple $(i,h,k)$ as above,
there exists one saddle connection $s_0$ in the family which consists
just in a straight segment, and doesn't go through the sides of the
polygon (if it exists, it is unique). I.e., if the segment in the
plane between the two angles $h$ and $k$ of the polygon is contained
within the polygon.  In this case adding a cylinder corresponds, in
the polygon representation, to adding a parallelogram, with one pair
of opposite sides associated to the saddle connection, and the other
pair associated to a vector $v$, which is thus inserted (at the
appropriate position) in the list of $v_j$'s describing the polygon
$S=(\tau;v_1,\ldots,v_n)$. This also adds one edge to the permutation
$\tau$, or, in other words, it adds the same vector $v$ both at a
position in the list $(v_j)_j$ and in the list $(u_j)_j$.  The
procedure in this simplified situation is described in
Figure~\ref{fig_add_cylinder_1}. In this situation we denote the
corresponding operator as $\mathcal{Q}_{i,h,k}(v)$, i.e.\ we omit the
index $s=s_0$ for the choice of saddle connection.

When, instead, we choose a saddle connection with a non-trivial
winding, in the polygon representation we shall add the vector $v$ at
several places along the list of $v_j$'s. The other copies are
associated to the introduction of false singularities (i.e.\ zeroes of
$\omega$ of degree zero, i.e., in the permutation, edges associated to
descents). This does not change the topology of the surface, but gives
non-primitive representants, or forces to study further manipulations
which allow to remove false singularities. We will not discuss this
analysis here.

Let us describe in more detail why the construction depicted above is
legitimate, now at the level of the translation surface as a whole
(and not of the polygon representation).  The saddle connection on the
surface is a closed curve, since its endpoints are identified in the
polygonal construction. Neglecting the identification, it consists of
some vector $w$, if the trajectory is followed from vertex $h$ to $k$,
in the local flat metric.  We cut the surface along this curve,
obtaining a new surface $S'$ with two boundary components homeomorphic
to circles.  Finally, we construct a cylinder by taking a
parallelogram with sides $(+v,+w,-v,-w)$, gluing its two $\pm w$ sides
to the two boundaries of $S'$, and its two $\pm v$ sides among
themselves.  This gives the desired surface
$S_{\rm new} = \mathcal{Q}_{i,h,k,s}(v) S$.

This operation breaks up the singularity $z_i$ of degree $d_i$ into
two singularities $z$ and $z'$ of degrees $d$ and $d'$ such that
$d+d'=d_i$. The value of $d$ and $d'$ depends on which pair of
vertices (in the equivalence class of $z_i$) are connected by the
geodesic. In particular, for each $h$ there exists exactly one value
$k$ such that the pair $(d,d')$ is attained, and similarly for $k$
(see figure~\ref{fig_add_cylinder_1}).

\begin{figure}[tb!]
\begin{center}
 \includegraphics[scale=.8]{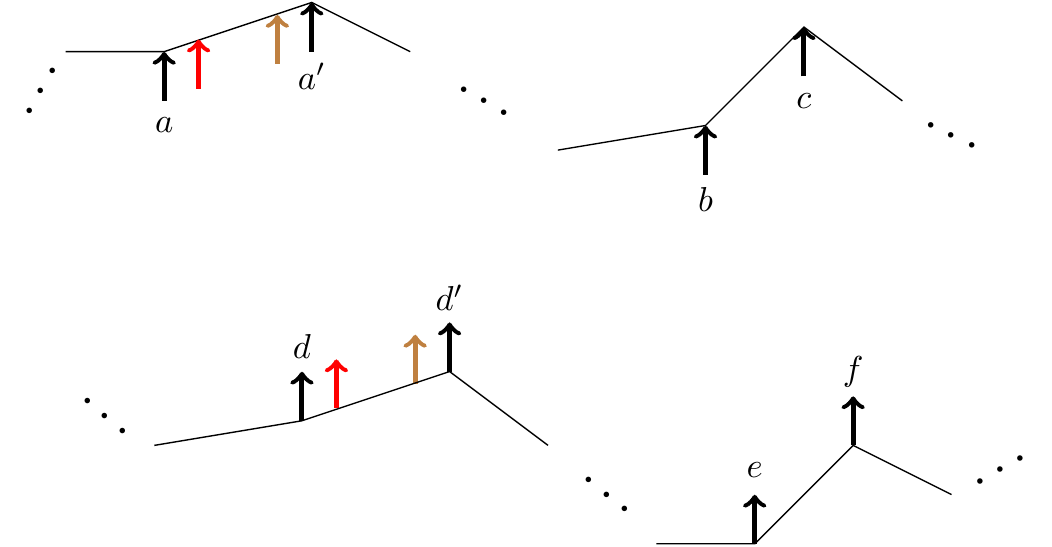}$\quad$
 \raisebox{2.5ex}{\includegraphics[scale=1.2]{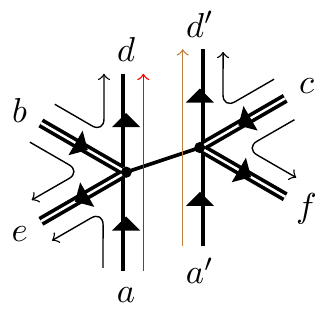}}\\[2mm]
 \includegraphics[scale=1.2]{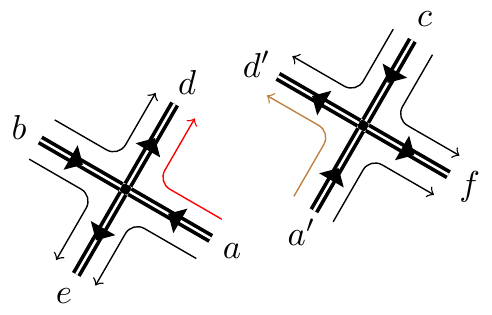}
\end{center}
\caption{\label{fig_add_cylinder_2} Top right: adding the cylinder
  breaks up the singularity of angle $6\pi$ in two, as is shown
  globally on the polygon, on the left, and in a neighbourhood of the
  saddle connection, on the right.  Bottom: after this operation, we
  see that there are two singularities in this case each of total
  angle $4\pi$ (neighbourhoods centered around each of them are
  shown). In red and brown, two leaves of a foliation parallel to the
  saddle connection, passing next to the singularities, help in
  visualising the local modification. In this case, in order to
  simplify the drawing, we have a description of the local
  neighbourhoods which corresponds to the case in which the saddle
  connection is itself vertical.}
\end{figure}

This operation is especially clear in the neighbourhood of $z_i$: the
introduction of the cylinder separates the singularity into two new
singularities, each one with its own neighbourhood represented by a
concatenation of half-planes (cf.\ Figure~\ref{fig_add_cylinder_2}).

\begin{figure}[tb!]
\[
\includegraphics[scale=.8]{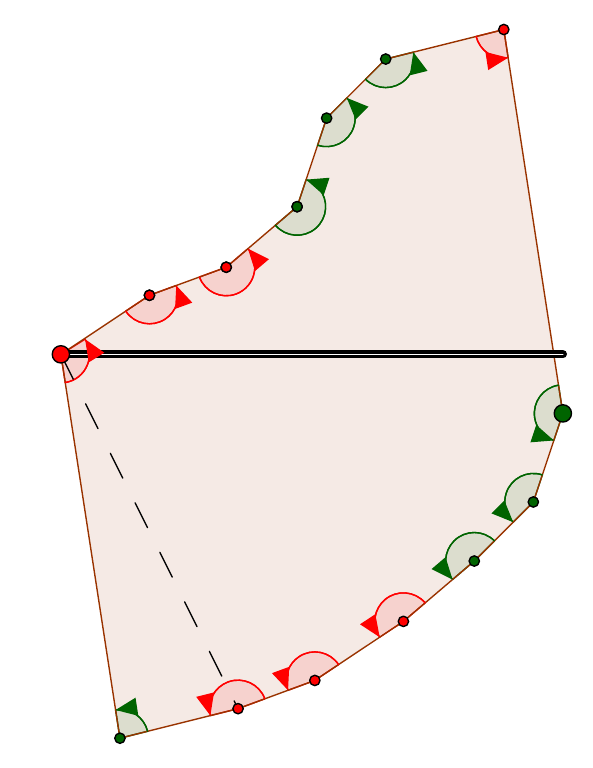}
\quad
\includegraphics[scale=.8]{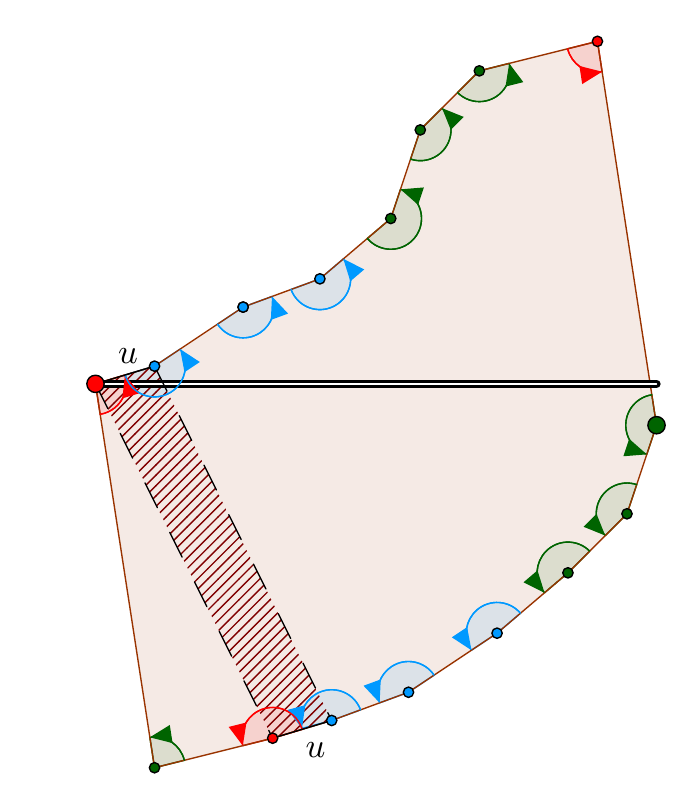}
\]
\caption{Left: a regular suspension for $\tau=(413256)$. This permutation has
  rank 3 (red angles), thus the marked singularity is of degree 2, and
  one cycle of length 3 (green), i.e.\ a singularity of degree 2. Right: the
  result of applying the operator $\mathcal{Q}$, which breaks the
  marked singularity into a marked singularity of degree zero and a
  new singularity of degree~2.\label{fig.exRegSuspQi}}
\end{figure}

In order to avoid the complicancy on the polygon representation
deriving by the presence of winding geodesics (and the resulting
introduction of false singularities), we will remark the crucial
property:
\begin{proposition}
\label{prop.cQgeo}
If $S$ is a regular suspension of shift $y$ for a permutation $\tau$
of rank $r$, $\bar{r}$ is the singularity associated to the rank,
$r_0$ is the top angle at the beginning of the rank cycle, and $k<r$
is the $k$-th bottom angle of the rank cycle, there exists one saddle
connection in the family $(r_0,k)$ which is contained within the
polygon, and, for all $v$ with argument in a non-empty range and
$\imof(v) < y$, $\mathcal{Q}_{\bar{r},r_0,k}(v) S$ is a regular
suspension of shift $y-\imof(v)$.

The range for $\arg(v)$ is as follows. If the new edge of the
permutation is added to the list $(u_j)_j$ at $j$, we shall have
\begin{align}
\label{eq.565865}
\theta_- 
&< \arg(v) <
\theta_+
&
\theta_- 
&=
\left\{
\begin{array}{ll}
\arg(u_j) & j>1 \\
0         & j=1
\end{array}
\right.
&
\theta_+
&=
\left\{
\begin{array}{ll}
\arg(u_{j+1})  & j<n \\
\frac{\pi}{2}  & j=n
\end{array}
\right.
\end{align}
\end{proposition}
\noindent
This construction is illustrated in Figure~\ref{fig.exRegSuspQi}.

The two surgery operators $q_1$ and $q_2$ correspond to
$\mathcal{Q}_{\bar{r},r_0,k}(v)$ as in Proposition \ref{prop.cQgeo},
for $k$ being such that the degree $d_i$ of the marked singularity is
divided into $d$ and $d_i-d$, with the degree of the new marked
singularity being $d=0$ or $1$, for $q_1$ and $q_2$, respectively.

We note in passing that the construction described so far, of breaking
up a zero, can be made local, i.e.\ can be realised in a way such that
the flat metric does not change outside of some neighbourhood of $z_i$
(see \cite[sect.\ 4.2]{KZ03} for details). This, however, and at
difference with Proposition \ref{prop.cQgeo}, yet again comes at the
price of introducing false singularities in the polygon (and
permutation) representation.

\subsubsection{Adding a handle}

\begin{figure}[tb!]
\begin{center}
 \includegraphics[scale=1.2]{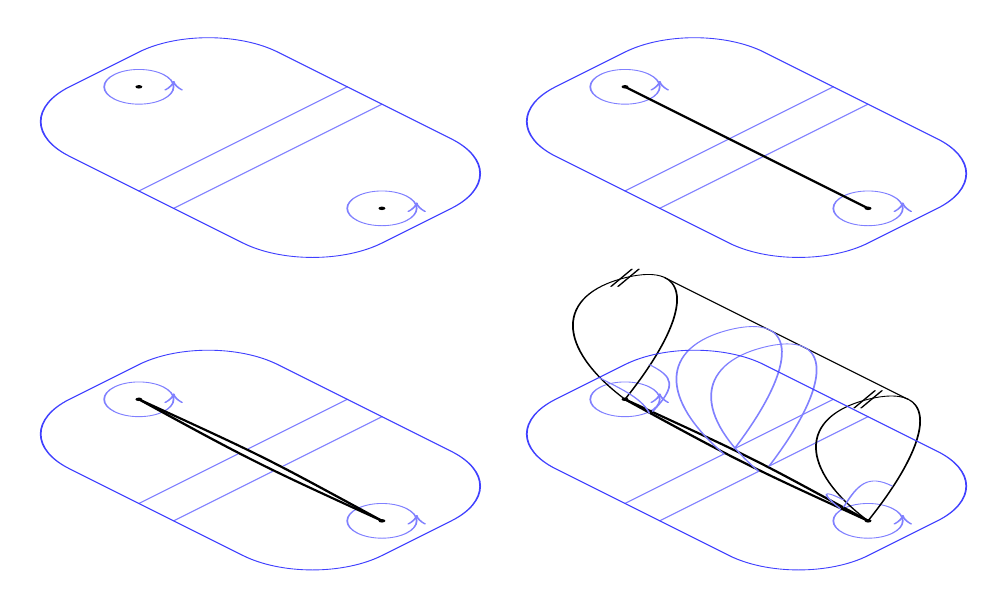}
\end{center}
\caption{\label{fig_add_handle_with_cyl}Top left: two singularities on
  a portion of a translation surface. In this case, in order to
  simplify the drawing, they are both false-singularities.  A small
  circle around each of them, and a couple of typical leaves of the
  vertical foliation, are shown. Top right: a saddle connection
  between the two singularities is added.  Bottom left: the surface is
  cut along this line. Bottom right: a handle is added along these
  boundaries. Note how both the circles around the singularities and
  the leaves are modified. In particular, now we have a unique
  singularity, in this case of angle $6 \pi$, as the arcs on the
  cylinder are connected in such a way to concatenate the two
  circles.}
\end{figure}

Our second family of geometric surgery operations $\mathcal{T}$
corresponds to adding a handle, and will be related to the operator
$T$.  Let us start by describing this operation on the translation
surface. Let $z_i$ be a singularity.
We choose a direction $\theta \in [0,2(d_i+1)\pi]$ and a length
$\rho$, and consider the arc of geodesic $\gamma$ starting from $z_i$
with direction $\theta$ and length $\rho$.  Assume that this arc does
not contain any singularity except for its starting point.  We call
$w$ the corresponding vector.
We cut the surface $S$ along the vector $w$, and call $S'$ the
resulting surface. Choose a vector $v$ with positive scalar product
with the versor of $w$ (, i.e., with an angle within $[0,\pi]$
w.r.t.\ the local metric at the starting point), and construct a
cylinder, by folding a parallelogram $(+v,+w,-v,-w)$ analogously to
how was the case for $\mathcal{Q}$, i.e.\ gluing its two $\pm w$ sides
to the two boundaries of $S'$, and its two $\pm v$ sides among
themselves.  This gives the desired surface
$S_{\rm new} = \mathcal{T}_{i,w}(v) S$.  This procedure is
illustrated, in a somewhat simplified picture, in
Figure~\ref{fig_add_handle_with_cyl}.

Let us describe this now at the level of the polygon representation.
Say that the direction $\theta$ is within the angle $x$ on the top
boundary of the polygon, and call $P$ the endpoint of $w$.  Similarly
to the case of $\mathcal{Q}$, if $P$ is on the boundary of the
polygon, and the vector $w$ is completely contained within the polygon
(i.e., it goes from $x$ to $P$ without passing through the identified
sides), then the construction is especially simple also at the level of
the polygon: it does not require the introduction of false
singularities, and consists in the introduction of a parallelogram.
More precisely, if $P$ is on the bottom side $u_j$ of the polygon, the
operation consists of two parts: first, it cuts the side $u_j$ into
two consecutive and collinear sides $u'_j$ and $u''_j$, thus adding
(temporarily) a false singularity of angle $2\pi$ on $P$, and then it
inserts a cylinder following $w$, which is now a saddle connection
between two angles of the polygon (just as in $\mathcal{Q}$, but with
the two angles associated to distinct singularities).

\begin{figure}[tb]
\begin{center}
 \includegraphics[scale=.8]{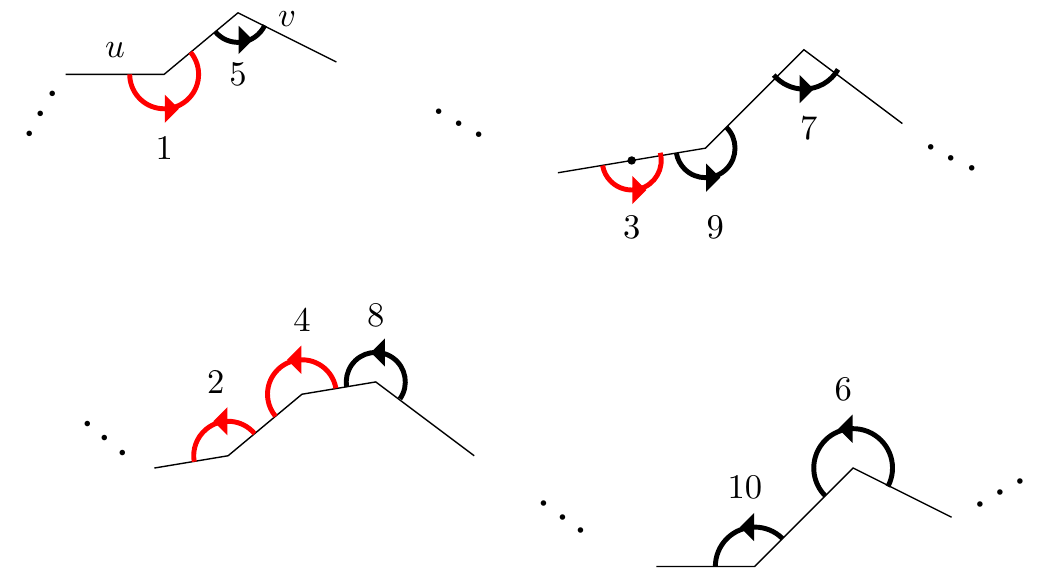}\\[-2mm]
\raisebox{1.58ex}{ \includegraphics[scale=.8]{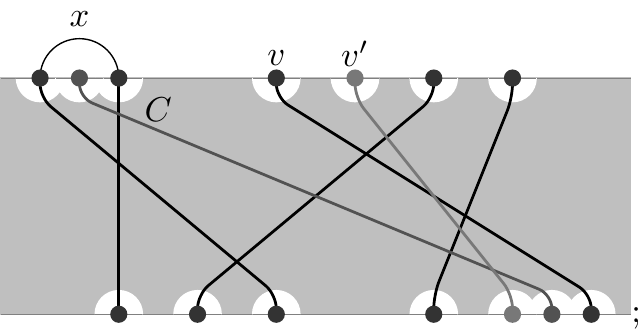}}
\includegraphics[scale=.8]{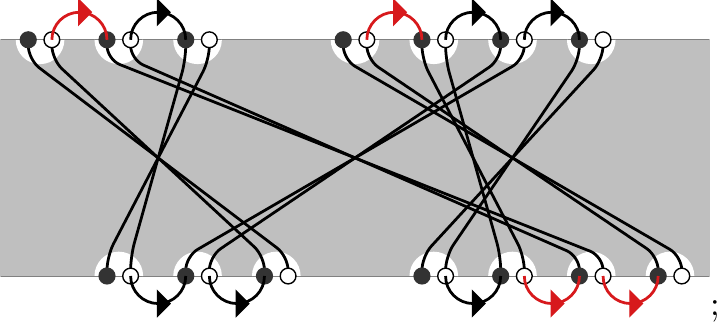}
\end{center}
\caption{\label{fig_add_handle_with_perm} Top : In A portion of a
  translation surface $S'$ corresponding to an operator $\mathcal{T}$
  acting on the example surface $S$ of Figure \ref{fig_6pi}. For what
  concerns the construction of small circles around the singularities,
  while in $S$ we have an arc from $u$ to $v$, here instead the arcs
  labeled with indices from 1 to 5 ultimately connect $u$ to $v$,
  while also encircling the newly introduced false sigularity, and the
  rest remains unchanged. The four arcs in red (two on the top broken
  line, and two on the bottom) are added to the same cycle, and thus
  increase by 2 the degree of the corresponding singularity, in this
  example passing from angle $6\pi$ to angle $10\pi$. Bottom: the
  procedure at the level of permutations. Let us describe the surgery
  procedure in details, comparing the language of permutations with
  the geometric counterpart:} 
\end{figure}

\begin{figure}[tb]
\begin{center}
\includegraphics[scale=.8]{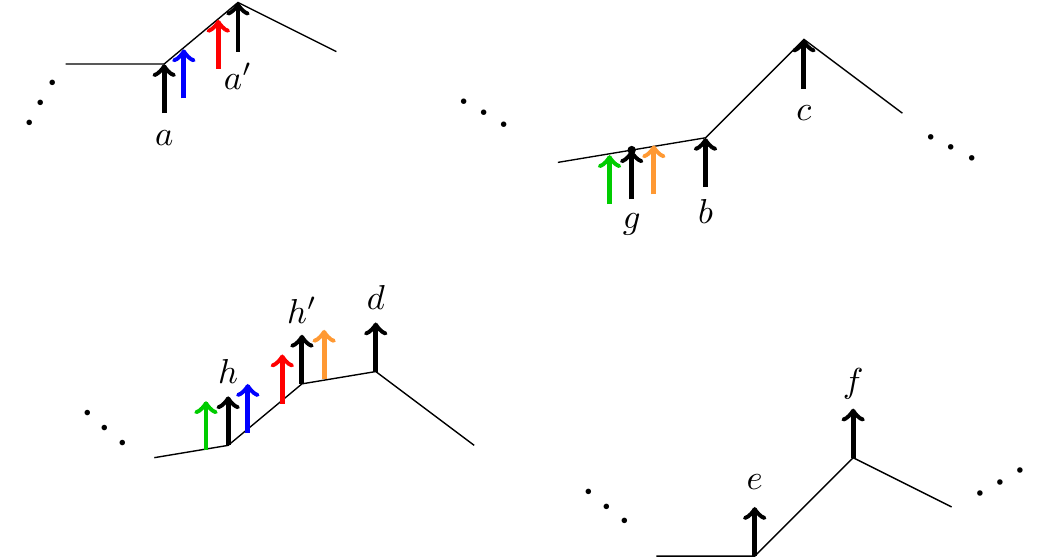}
\includegraphics[scale=1]{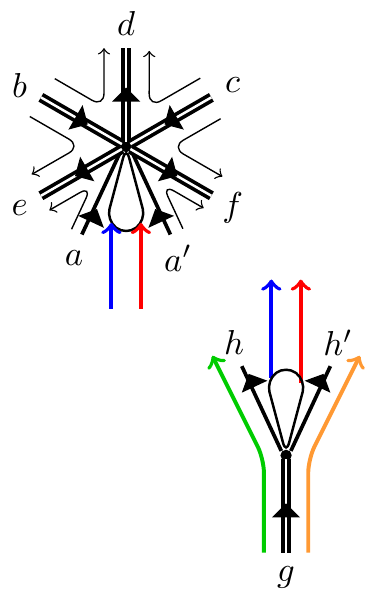}
\includegraphics[scale=1.2]{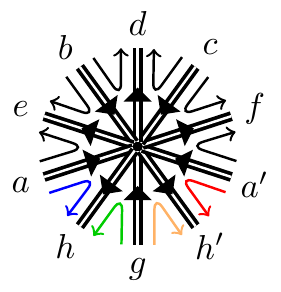}
\end{center}
\caption{\label{fig_add_handle_neigh} Adding the point $P$ on the
  boundary $v$ of the polygon divides this side in two, thus $P$
  becomes a false singularity of angle $2\pi$ (i.e. a zero of
  degree~0). The two half-disks are denoted here by the presence of a
  green and an orange leaf in the vertical foliation.  Adding the
  cylinder glues together the two singularities, resulting in this
  case in a $10\pi$ singularity, whose neighbourhood, obtained by the
  insertion of the parallelogram with blue and red leaves on the two
  singularities (see top-right), is isomorphic to the neighbourhood
  illustrated at the bottom of the image.}
\end{figure}

By this operation we have added a handle to the surface $S$ on the
singularity $z_i$.  If the degree of $z_i$ on $S$ is $d_i$, then its
degree on
$S_{\rm new}$ is $d_i+2$, as explained in figure
\ref{fig_add_handle_with_perm}.  Seen in a neighbourhood of the
singularity, this operation corresponds to adding four angular sectors
between two consecutive angular sectors (see figure
\ref{fig_add_handle_neigh}), and, in particular, the false singularity
that was temporarily added is now merged to the original
singularity~$z_i$.

If the point $P$ is not on a side of the polygon, we can proceed by a
slight generalisation of the first part of the operation: choose a
bottom side $u_j$, call $P'$ and $P''$ its endpoints, and suppose that
the triangle $(P,P',P'')$ can be added/removed on the top/bottom
copies of the edge (i.e., $u_j$ on the bottom and
$v_{\tau^{-1}(n+1-j)}$ on top), so to obtain an equivalent
representation of the polygon, with the introduction of a unique false
singularity (this happens when the segments $(P,P')$ and $(P,P'')$ on
the plane do not cross the boundary of the polygon). If we do so, we
are in the same situation as before, i.e.\ we have decomposed
$u_j = u'_j + u''_j$, with $u'_j = (P',P)$ and
$u''_j = (P,P'')$, the only difference being the fact that $u'_j$ and
$u''_j$ are not collinear. 

In order to describe the procedure at the level of the polygon, and
with a generic point $P$ as above (not necessarily on the boundary),
we need to supplement the choice of $x$ and $j$ to the notation, so we
will write $S_{\rm new} = \mathcal{T}_{x,j,w}(v) S$.

Similarly as was the case with adding a cylinder, if we follow this
procedure on a regular suspension, and the point $P$ is chosen
appropriately, we can have a regular suspension as outcome, as
illustrated in the following:
\begin{proposition}
\label{prop.zhxjgsd}
Given a side $u_j = (P',P'')$ on the bottom side of a regular
suspension, excluded the leftmost one, define $\theta_{\pm}$ as in
(\ref{eq.565865}).
Consider the triangle $\Delta$ with one side $u_j$, and the other two
sides with slopes $\theta_{\pm}$. For each $P$ in the interior of
$\Delta$, if we add the triangle $(P,P',P'')$ to the polygon, below
$u_j$, and remove the corresponding copy of the triangle below
$v_{\tau^{-1}(n+1-j)}$, we obtain a new suspension which is regular,
has the same shift parameter, and has a false singularity in~$P$.
\end{proposition}

\begin{figure}[tb!]
\[
\includegraphics[scale=.8]{Figure2b_fig_reg_susp_1.pdf}
\quad
\includegraphics[scale=.8]{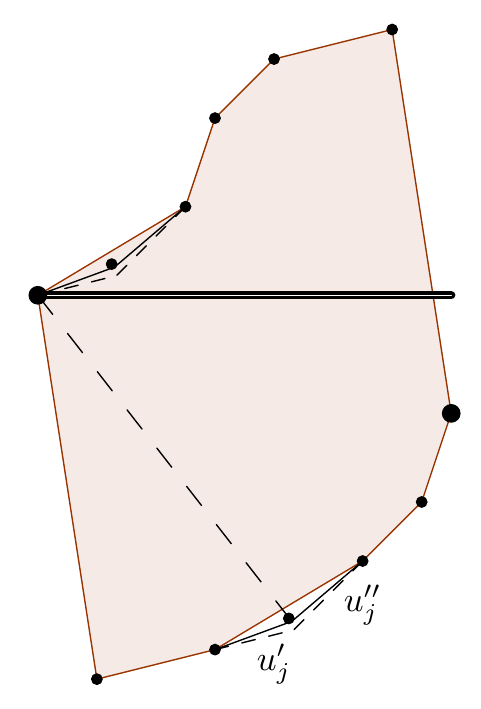}
\quad
\includegraphics[scale=.8]{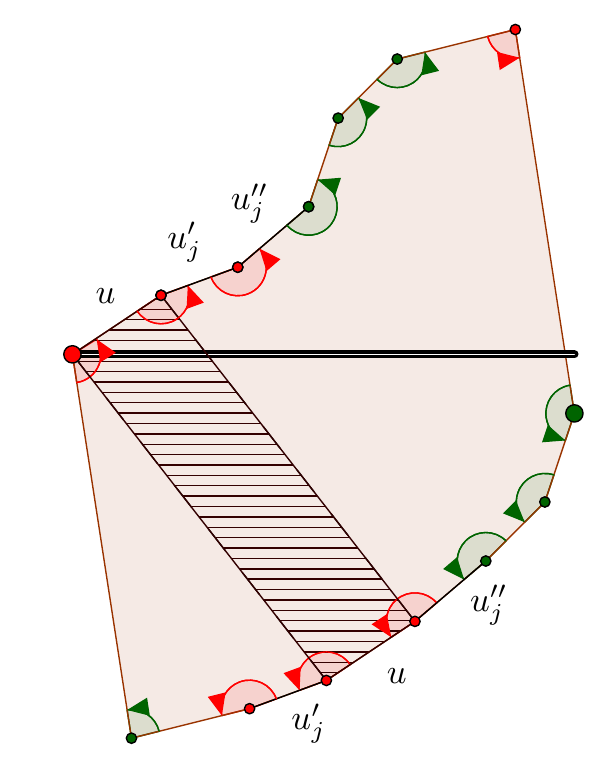}
\]
\caption{Left: a regular suspension $S$ for $\tau=(23145)$ (the same
  as in figure \ref{fig.exRegSusp}). Middle: the construction of $S'$
  described by Proposition \ref{prop.zhxjgsd}. Right: the operation
  $\mathcal{T}_{r_0,j,w}(v) S$ described by Proposition
  \ref{prop.cTgeo} (this suspension is as in the left part of figure
  \ref{fig.exRegSuspQi}).  The resulting permutation, $\tau=(413256)$,
  has the same cycle structure as the starting one, except for the
  rank, which has increased by~2.
\label{fig.exRegSuspTi}}
\end{figure}

\begin{proposition}
\label{prop.cTgeo}
Let $S$ be a regular suspension of shift $y$, for a permutation $\tau$
of rank $r$, $\bar{r}$ is the singularity associated to the rank,
$r_0$ is the top angle at the beginning of the rank cycle, $1 \leq j
\leq n-1$ is the index of a vector on the bottom side, excluded the
left-most one.

Let notations be as in Proposition \ref{prop.zhxjgsd}, with $S'$ being
the outcome polygon, and $w$ be the saddle connection connecting $r_0$
to $P$, by passing through no sides of the polygon $S'$.
Let $v$ be a vector with $\arg(u'_j) < \arg(v) < \arg(u''_j)$, and
$\imof(v)<y$. Then the suspension $\mathcal{T}_{r_0,j,w}(v) S$ is
regular with shift $y-\imof(v)$.
\end{proposition}
\noindent
This construction is illustrated in Figure~\ref{fig.exRegSuspTi}.

Our $T$ operator is a special case of this operation, corresponding to
$j$ being the right-most vector of the polygon, on the bottom side. We
can define, more generally, operators $\{T_j\}_{1 \leq j \leq n-1}$,
of which $T\equiv T_1$ is the special case above. For what concerns
the invariants, all of these $T_j$ behave in the same way: the degree
of the marked singularity increases by 2, and all the rest of the
invariant stays put.


\subsubsection{Regular suspensions in a given Rauzy class}

Our classification theorem, based on the analysis of surgery operators
$q_{1,2}$ and $T$, and the results of the present section, in
particular Propositions \ref{prop.cQgeo} and \ref{prop.cTgeo}, allow
to describe large families of regular suspensions (with no false
singularities) which are certified to be within one given (primitive)
Rauzy class.

For a class $C$, let us call $\cR(C)$ the set of such regular
suspensions. Thus we have $\cR(C) = \{ S \}$, with
$S=(\tau;v_1,\ldots,v_n)$ satisfying certain properties, which we now
explicitate. We describe explicitly $\cR(\Id_n)$ and $\cR(\Id'_n)$,
and, for $C$ a non-exceptional primitive class with invariant
$(r,\lam,s)$, we describe $\cR(C) \equiv \cR(r,\lam,s)$ recursively in
terms of the invariant.

We say that $(v_1,\ldots,v_n)$ is a \emph{regular vector} if it
satisfies the first two conditions of Definition \ref{def.regsusp}. We
can then anticipate a corollary of the results presented in this
section and in the body of the paper:
\begin{corollary}
\label{cor.represclass}
The surfaces in the sets $\cR(C)$ below are all regular suspensions:
\begin{align}
\cR(\Id_n)
&=
\Big\{
\Big(
\raisebox{-6pt}{\includegraphics[scale=1]{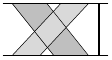}};
v_1,\ldots,v_n
\Big)
\textrm{~with $(v_1,\ldots,v_n)$ regular}
\Big\}
\\
\cR(\Id'_n)
&=
\Big\{
\Big(
\raisebox{-6pt}{\includegraphics[scale=1]{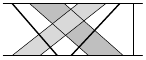}};
v_1,\ldots,v_n
\Big)
\textrm{~with $(v_1,\ldots,v_n)$ regular}
\Big\}
\\
\cR(1,\lam,s)
&=
\bigcup_{\ell \;|\; \exists \lam_i=\ell}
\{ \mathcal{Q}_{\bar{r},r_0,1}(v) S\}_{S \in \cR(\ell,\lam \setminus
  \ell,s); \textrm{~$v$ as in Prop.\ \ref{prop.cQgeo}} }
\\
\label{eq.345545}
\cR(2,\lam,0)
&=
\bigcup_{\substack{\ell \;|\; \exists \lam_i=\ell \\ s=0,\pm 1}}
\{ \mathcal{Q}_{\bar{r},r_0,2}(v) S\}_{S \in \cR(\ell+1,\lam \setminus
  \ell,s); \textrm{~$v$ as in Prop.\ \ref{prop.cQgeo}} }
\\
\cR(r,\lam,s)
&=
\bigcup_{j=1}^{n-1}
\{ \mathcal{T}_{r_0,j,w}(v) S\}_{S \in \cR(r-2,\lam,s); \textrm{~$w$ and
    $v$ as in Props.\ \ref{prop.zhxjgsd} and \ref{prop.cTgeo}} }
&
r \geq 3
\end{align}
where, in (\ref{eq.345545}), the set of invariants which contribute to
the set-union
is as follows: if $\lambda \setminus \ell$ has a positive number of
even cycles, then $s=0$, otherwise, $s=\pm 1$.
\end{corollary}
\noindent
(The domains in the sums are explained in detail in the following
Corollaries~\ref{cor.goback123a} and \ref{cor.goback123b}).

This result is complementary to the work present in an article of
Zorich \cite{Zor08}, in which he constructs representatives of every
connected component of every strata, consisting of Jenkins--Strebels
differentials. One advantage of our approach is the fact that the
representatives given in Corollary \ref{cor.represclass} are a large
number (the smallest asymptotics for a non-exceptional class,
corresponding to the iterated application of $\cT \circ
\mathcal{Q}_2$, is of the order of $n!!$). A second advantage is that
each representative has a simple interpretation in terms of nested
insertions of handles and cylinders.

\subsection{A summary of terminology}
\label{ssec.tabglossary}

We end this section by collecting a list of notions which appear both
in our approach and in the geometric construction (but, sometimes,
under different names). First, we recall in words some notational
shortcuts for the strata which were used in \cite{KZ03} and
\cite{Boi12}:
\begin{itemize}
\item They call $H^{\rm hyp}(2g-2)$ and $H^{\rm hyp}(g-1,g-1)$ the
  \emph{hyperelliptic classes}, which, by Lemma~\ref{deux_Gauss_B},
  correspond to $\Id_n$ with $n=2g$
  and $\Id_n$ with $n=2g+1$,
  respectively.

\item They call $H^{\rm even}(2d_1^{m_1},\ldots,2d_k^{m_k})$ and
  $H^{\rm odd}(2d_1^{m_1},\ldots,2d_k^{m_k})$ the two connected
  components of the stratum $H(2d_1^{m_1},\ldots,2d_k^{m_k})$ with Arf
  invariant $+1$ and $-1$ respectively.

\item We shall call $H^{\rm
  even}(2d_1^{m_1},\ldots,\bar{2d_i}^{m_i},\ldots,2d_k^{m_k})$ and
  $H^{\rm odd}(2d_1^{m_1},\ldots,\bar{2d_i}^{m_i},\ldots,2d_k^{m_k})$
  the two connected components of the marked stratum
  $H(2d_1^{m_1},\ldots,\bar{2d_i}^{m_i},\ldots,2d_k^{m_k})$ with Arf
  invariant $1$ and $-1$ respectively (in analogy with their
  notation).

\item Finally, we use symbols $H(d_1^{m_1},\ldots,d_k^{m_k})$
  (respectively
  $H(d_1^{m_1},\ldots,\bar{d_i}^{m_i},\ldots,d_k^{m_k})$) only when a
  positive number of the $d_i$'s are odd. In this case, these symbols
  denote a stratum (respectively, a marked stratum) with some zeroes
  of odd degree, so that the Arf invariant is not defined, and the
  stratum is connected.
\end{itemize}

\noindent
Then, Table \ref{table_correspondance} gives a list of correspondences
between the terminology adopted in \cite{KZ03} and \cite{Boi12},
related to the geometric notions associated to the strata, and the
one, issued from the combinatorial approach, used in this article

\begin{table}[tb!]
\begin{tabular}{|p{\dimexpr .55\textwidth -2\tabcolsep}|p{\dimexpr .45\textwidth -2\tabcolsep}|}\hline
\rule{0pt}{13pt}%
\bf{Geometric and topological objects} & 
\bf{Combinatorial objects}
\raisebox{-6pt}{\rule{0pt}{13pt}}%
\\ \hline \hline
\rule{0pt}{13pt}%
Translation surface with vertical direction 
\raisebox{-6pt}{\rule{0pt}{13pt}}%
&
\multirow{2}{*}{Suspension data ($\s$,$(v_i)_i)$}\\ 
Riemann surface with abelian differential $(M,\omega)$\raisebox{-6pt}{\rule{0pt}{13pt}}%
&\\ \hline
\rule{0pt}{13pt}%
Conical singularity of angle $2\pi(d_i+1)$ 
\raisebox{-6pt}{\rule{0pt}{13pt}}%
& 
\multirow{2}{*}{Cycle of length $\lambda_i=d_i+1$}\\ 
Zero of $\omega$ of degree $d_i$ 
\raisebox{-6pt}{\rule{0pt}{13pt}}%
& \\ 
False singularity
(zero of $\omega$ of degree 0)
&
descent (cycle of length 1)
\\
\hline
\rule{0pt}{13pt}%
Marked conical singularity of angle $2\pi(d+1)$ 
\raisebox{-6pt}{\rule{0pt}{13pt}}%
& 
\multirow{2}{*}{Rank path of length $r=d+1$}\\ 
Marked zero of $\omega$ of degree $d$ & 
\\ \hline
\rule{0pt}{13pt}%
Arf invariant $\mathrm{arf}(\Phi)$ & 
Arf (or sign) invariant $s(\s)$
\raisebox{-6pt}{\rule{0pt}{13pt}}%
\\\hline
\rule{0pt}{13pt}%
$H^{\rm hyp}(2g-2)$& $\Id_n$ with $n=2g$ even
\raisebox{-6pt}{\rule{0pt}{13pt}}%
\\ 
\rule{0pt}{13pt}%
$H^{\rm hyp}(g-1,g-1)$ & $\Id_n$ with $n=2g+1$ odd
\raisebox{-6pt}{\rule{0pt}{13pt}}%
\\ \hline
\rule{0pt}{13pt}%
$H^{\rm hyp}(2g-2,\bar{1})$& $\Id'_n$ with $n=2g+1$ odd
\raisebox{-6pt}{\rule{0pt}{13pt}}%
\\ 
\rule{0pt}{13pt}%
$H^{\rm hyp}(g-1,g-1,\bar{1})$ & $\Id'_n$ with $n=2g+2$ even
\raisebox{-6pt}{\rule{0pt}{13pt}}%
\\ \hline
\begin{minipage}{\textwidth}
\rule{0pt}{13pt}%
Connected component in $\perms$\\
$H^{\rm even}(2 d_1^{m_1},\ldots,2 \bar{d_i}^{m_i},\ldots,2
d_k^{m_k})$
\raisebox{-6pt}{\rule{0pt}{13pt}}%
\end{minipage}
 &
\begin{minipage}{\textwidth}
\rule{0pt}{13pt}%
Rauzy class with invariant\\
$(\lambda,r=2d_i+1,s=+1)$\\
$\lambda = \big( (2d_1+1)^{m_1} \cdots (2d_i+1)^{m_i-1}$
\\ \hspace*{100pt}$ \cdots (2d_k+1)^{m_k} \big)$
\raisebox{-6pt}{\rule{0pt}{13pt}}%
\end{minipage}
\\
\begin{minipage}{\textwidth}
\rule{0pt}{13pt}%
Connected component in $\perms$\\
$H^{\rm odd}(2 d_1^{m_1},\ldots,2 \bar{d_i}^{m_i},\ldots,2 d_k^{m_k})$
\raisebox{-6pt}{\rule{0pt}{13pt}}%
\end{minipage}
&
as above, with $s=-1$
\\
\begin{minipage}{\textwidth}
\rule{0pt}{13pt}%
Connected component in $\perms$\\
$H(d_1^{m_1},\ldots,\bar{d_i}^{m_i},\ldots,d_k^{m_k})$
\raisebox{-6pt}{\rule{0pt}{13pt}}%
\end{minipage}
 &
\begin{minipage}{\textwidth}
\rule{0pt}{13pt}%
Rauzy class with invariant\\
$(\lambda,r=d_i+1,s=0)$\\
$\lambda = \big( (d_1+1)^{m_1} \cdots (d_i+1)^{m_i-1}$
\\ \hspace*{100pt}$ \cdots (d_k+1)^{m_k} \big)$
\raisebox{-6pt}{\rule{0pt}{13pt}}%
\end{minipage}
\\ \hline
\begin{minipage}{\textwidth}
\rule{0pt}{13pt}%
Connected component in $\permsex$\\
$H^{\rm even}(2 d_1^{m_1},\ldots,2 d_k^{m_k})$
\raisebox{-6pt}{\rule{0pt}{13pt}}%
\end{minipage}
 &
\begin{minipage}{\textwidth}
\rule{0pt}{13pt}%
Extended Rauzy class with invariant\\
 $(\lambda,s=+1)$\\
$\lambda = \big( (2d_1+1)^{m_1} \cdots (2d_k+1)^{m_k} \big)$
\raisebox{-6pt}{\rule{0pt}{13pt}}%
\end{minipage}
\\
\begin{minipage}{\textwidth}
\rule{0pt}{13pt}%
Connected component in $\permsex$\\
$H^{\rm odd}(2 d_1^{m_1},\ldots,2 d_k^{m_k})$
\raisebox{-6pt}{\rule{0pt}{13pt}}%
\end{minipage}
&
as above, with $s=-1$
\\
\begin{minipage}{\textwidth}
\rule{0pt}{13pt}%
Connected component in $\permsex$\\
$H(d_1^{m_1},\ldots,d_k^{m_k})$
\raisebox{-6pt}{\rule{0pt}{13pt}}%
\end{minipage}
 &
\begin{minipage}{\textwidth}
\rule{0pt}{13pt}%
Extended Rauzy class with invariant\\
 $(\lambda,s=0)$\\
$\lambda = \big( (d_1+1)^{m_1} \cdots (d_k+1)^{m_k} \big)$
\raisebox{-6pt}{\rule{0pt}{13pt}}%
\end{minipage}
\\ \hline
\rule{0pt}{13pt}%
    Bubbling a handle & $T$ operator
\raisebox{-6pt}{\rule{0pt}{13pt}}%
\\ 
\rule{0pt}{13pt}%
    Adding a cylinder & $q_1$ and $q_2$ operators
\raisebox{-6pt}{\rule{0pt}{13pt}}%
\\ \hline
 \end{tabular}
 \caption{\label{table_correspondance}Correspondence between the
   terminology adopted in \cite{KZ03} and \cite{Boi12} and the one
   used in this article. $Id'_n$ (i.e the hyperelliptic class with a marked point) was studied in \cite{BL12}.}
\end{table}

\section{Some basic facts}
\label{sec.basic_facts}

\subsection{Properties of the cycle invariant}
\label{ssec.arcbased}

\begin{figure}[b!!]
\begin{center}
\begin{tabular}{ccc}
\raisebox{-28pt}{\includegraphics[scale=2.3]{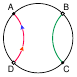}}
\ 
\raisebox{-40pt}{\includegraphics[scale=2.3]{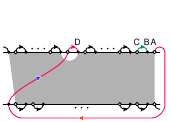}}
&
$\xrightarrow{\  L \ }$
&
\raisebox{-40pt}{\includegraphics[scale=2.3]{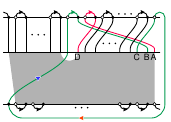}}
\ 
\raisebox{-28pt}{\includegraphics[scale=2.3]{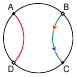}}
\\
\rule{0pt}{55pt}%
\raisebox{-28pt}{\includegraphics[scale=2.3]{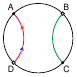}}
\ 
\raisebox{-40pt}{\includegraphics[scale=2.3]{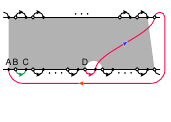}}
&
$\xrightarrow{\ R \ }$
&
\raisebox{-40pt}{\includegraphics[scale=2.3]{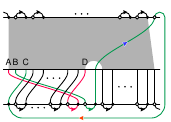}}
\ 
\raisebox{-28pt}{\includegraphics[scale=2.3]{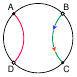}}
\end{tabular}
\end{center}
\caption{\label{fig.cycleinv} Invariance of $\lam(\s) \cup \{r\}$
  in $\perms_n$, and of $\lam'(\s)$ in $\permsex_n$, illustrated for
  the operators $L$ (top) and $R$ (bottom).  For $\permsex_n$, the cases of operators
  $L'$ and $R'$ are deduced analogously.}
\end{figure}

In this section we come back to the cycle invariant introduced in
Section \ref{sssec.cycinv}. In particular, we shall prove Propositions
\ref{prop.cycinv1} and \ref{prop.cycinv2}, concerning the fact that
$(\lam,r)$ is invariant in the $\perms$ dynamics, and $\lam$ is
invariant in the $\permsex$ dynamics.  These facts are especially
evident in the diagrammatic representation, which we illustrate in
Figure~\ref{fig.cycleinv}.

The idea is that the operators of the dynamics, i.e.\ the permutations
$\gamma_{L,n}(i)$ and $\gamma_{R,n}(i)$ defined in
(\ref{eqs.opecycdef}), perform local modifications on portions of
paths, without changing their lengths.  For all cycles (or the rank
path) except for at most two of them, this is completely evident: each
arc of the cycles is deformed by (say) the permutation
$\gamma_{L,n}(i)$ in a way which can be retracted without changing the
topology of the connections, this in turns certifying that the length
of the involved cycle does not change. (We recall that the length is
defined as the number of top or bottom arcs). Furthermore, it also
preserves the path-lengths of these cycles.

For the paths passing from fours special endpoints (marked as $A$, $B$,
$C$ and $D$ in the figure), the invariance of the cycle structure
holds through a more subtle mechanism, involving the exchange
of two arcs, and, for what concerns path-lengths, by a crucial use of
the $-1$ mark. The way in which the topology of connections among
these four points is modified by the permutation $\gamma_{L,n}(i)$ is
hardly explained in words, but is evident from
Figure~\ref{fig.cycleinv}.  This proves that the list
$\lam'=\lam \cup \{r\}$ is invariant. We shall prove also that, in the
$\perms_n$ dynamics, the rank length is also preserved (and thus, by
set difference, $\lam$ is).  This comes from the fact that the two
endpoints of the rank path cannot move under the action of $L$ and $R$
(while they \emph{would} move under $L'$ and $R'$).
Figure~\ref{fig.cycleinv} shows that path lengths are preserved also
if the paths are coloured (e.g., in the figure, if the lengths of the
green and purple paths in $\s$ are $\ell_1$ and $\ell_2$,
respectively, they are still $\ell_1$ and $\ell_2$ in $L \s$).  As we
perform local modifications on portions of paths adjacent to arcs,
these do not affect the rank-path endpoints, and thus do not affect
the rank-path `colour'.  From this we deduce that also the length $r$
of the rank path is invariant. This completes the proof.

Note that the length of the principal cycle (i.e.\ the cycle going
through the $-1$ mark) is \emph{not} preserved. Also, is not preserved
the boolean value [does the rank-path pass through the $-1$ mark].
Indeed, in the figure, the purple path goes through the mark in $\s$,
while the green path does, in $L \s$.  Nonetheless, these quantities
evolve in a rather predictable way.  If we annotate which cycles/paths
go through the arcs in the top-right portion of the diagram (namely,
at the right of $\sigma(1)$), we can observe that iterated powers of
$L$ make a circular shift on these labels (yet again, represented by
the green and purple colours on the arcs), and this immediately
reflects on which cycle (or the rank path) goes through the $-1$
mark. This fact is analysed in more detail in Section
\ref{ssec.typeXH} below.

\subsection{Standard permutations}
\label{ssec.std_perm}

As we have seen in the introduction, irreducibility is a property of
classes: either all permutations of a class are irreducible or none
is.  It is a trivial and graphically evident property, as the block
structure is clear both from the matrix representation of $\s$, and
from the diagram representation of $\s$. There is a further, less
obvious, characterisation of irreducibility, which leads us to the
definition of `standard permutations'. Both notions
are already present in the literature. For example, standard
permutations are introduced by Rauzy \cite{Rau79}, they play a crucial
role in \cite{KZ03} (and various other papers), and their enumerations
are investigated in~\cite{Del13}.

\begin{definition}[standard permutation]
The permutation $\s$ is \emph{standard} if $\s(1)=1$.
\end{definition}

\begin{figure}[b!!]
\begin{center}
\includegraphics{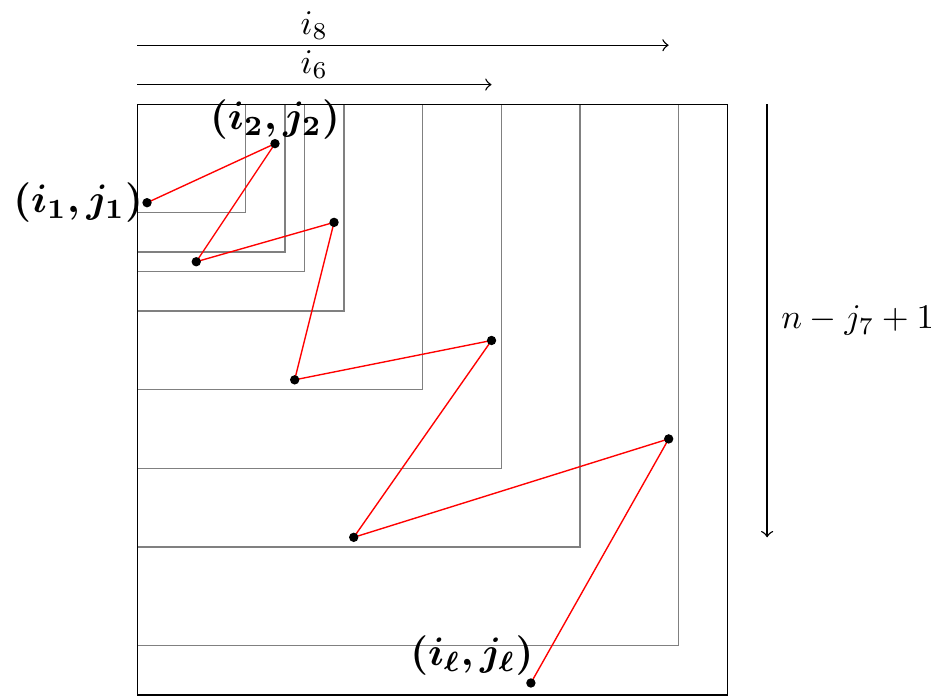}
\caption{\label{figure.zigzag}%
  Structure of a zig-zag path, in matrix representation.  It is
  constructed as follows: for each pair $(i_a,j_a)$, put a bullet at
  the corresponding position, and draw the top-left square of side
  $\max(i_a,j_a)$. Draw the segments between bullets $(i_a,j_a)$ and
  $(i_{a+1},j_{a+1})$ (here in red). The resulting path must connect
  the top-left boundary of the matrix to the bottom-right boundary.
  If the bullets are entries of $\s$, then a top-left $k \times k$
  square cannot be a block of the matrix-representation of $\s$,
  because, in light of the inequalities (\ref{eq.6537564}), one of the
  two neighbouring rectangular blocks (the $k \times (n-k)$ block to
  the right, or the $(n-k) \times k$ block below) must contain a
  bullet of the path, and thus be non-empty.}
\end{center}
\end{figure}

\noindent
(This definition is slightly different from the common one, where $\s$
is standard if $\s(1)=1$ and $\s(n)=n$. As we will see below, in
Proposition \ref{trivialbutimp}, this is a minor problem, because the
two notions are easily related.)  As mentioned above, we have:
\begin{lemma}
\label{lem.stdPerms}
A class $C \subseteq \perms_n$ is irreducible if and only if it
contains a standard permutation.
\end{lemma}
\noindent
This lemma is proven towards the end of this subsection. Before doing
this, we need to introduce \emph{zig-zag paths} (see also
Figure~\ref{figure.zigzag}).

\begin{definition}
A set of edges
$\big((i_1,j_1),(i_2,j_2),\ldots,(i_{\ell},j_{\ell})\big)$ is a
\emph{$L$ zig-zag path} if $i_1=1$, the indices satisfy the pattern of
inequalities
\begin{align}
j_{2b} &> j_{2b-1}
\ef,
&
i_{2b} &> n-j_{2b-1}+1
\ef,
&
i_{2b} &> i_{2b+1}
\ef,
&
n-j_{2b+1}+1 >i_{2b}
\ef,
\label{eq.6537564}
\end{align}
and either $i_\ell=n$ or $j_\ell=1$.  The analogous structure
starting with $j_1=n$ is called
\emph{$R$ zig-zag path}.
\end{definition}

\noindent
It is easy to see that, if $\s$ has a $L$ zig-zag path, no set
$[s]\subsetneq [n]$ can have $\s [s] = [n,\ldots,n-s+1]$, as, from the
existence of the path, there must exist either a pair
$(u,v)=(i_{2b+1},j_{2b+1})$ such that $u\leq s$, $n-s+1>v$ and
$v=\s(u)$, or a pair $(u,v)=(j_{2b},i_{2b})$ such that $u \geq n-s+1$,
$ s< v$, and $v=\s^{-1}(u)$ (this is seen in more detail below). A
similar argument holds for $R$ zig-zag paths.  Thus zig-zag paths
provide a concise certificate for irreducibility.

We now illustrate an algorithm that, for a given $\s \in
\mathfrak{S}_n$, produces in linear time either a certificate of
reducibility, by exhibiting the top-left block of the permutation, or
a certificate of irreducibility in the form of a canonically-chosen
$L$ zig-zag path, that we shall call the 
\emph{greedy $L$ zig-zag path} of the permutation. A completely
analogous algorithm gives the greedy $R$ zig-zag path. This algorithm
will imply the following:

\begin{lemma}
A permutation $\s$ is irreducible if and only if it has a zig-zag path.
\end{lemma}

\pf
For $X =\{x_a\} \subseteq [n]$, we use $\s(X)$
as a shortcut for $\{ \s(x_a)\}$.
Set $J_1 := \s(\{1\})$ and
$i_1=1$.  Now, either $\s$ is reducible, and has top-left block of size 1, or
$j_1 := \max(J_1) < n$.
If $j_1=1$, then $\big((i_1,j_1)\big)$ is a $L$ zig-zag path.  If
$j_1>1$, consider the set $I_1 = \s^{-1}(\{n,n-1,\ldots,n-j_1+1\})$.
Either the largest element of this set is $n-j_1+1$, in which case we
certify reducibility of $\s$ with a first block of size $n-j_1+1$, or
it is some index $i_2 > n-j_1+1$.  Say $\s(i_2)=j_2<j_1$.  If $i_2=n$,
then $\big((i_1,j_1),(i_2,j_2)\big)$ is a $L$ zig-zag path.  If
$i_2<n$, consider the set $J_2 = \s\{i_1+1=2,\ldots,i_2\}$.  Either
the smallest element of this set is $n-i_2+1$, in which case we
certify reducibility of $\s$, and with top-left block of size $i_2$, or it
is some index $j_3$ such that $n-j_3+1 > i_2$.  Say
$i_3 := \s^{-1}(j_3)<i_2$. Continue in this way.
At each round we query the image or pre-image of values in some
interval $(i_a+1,i_a+2,\ldots,i_{a+1})$, or same with $j$, so the
total amount of queries is bounded by $2n$.  At the various rounds at
which we do not halt we have $i_{b+2}>i_{b}$ and $j_{b+2}<j_{b}$,
As the indices $i_a$, $j_a$ are in the range $[n]$, the algorithm must
terminate in at most $\sim n$ rounds, and can only terminate either
because we have found a reducible block, or because $i_{2b}=n$ for
some $b$, or because $j_{2b-1}=1$ for some~$b$.  This proves both the
lemma and the statement on the linearity of the algorithm.  \qed

\begin{figure}[t!!]
\begin{center}
\includegraphics{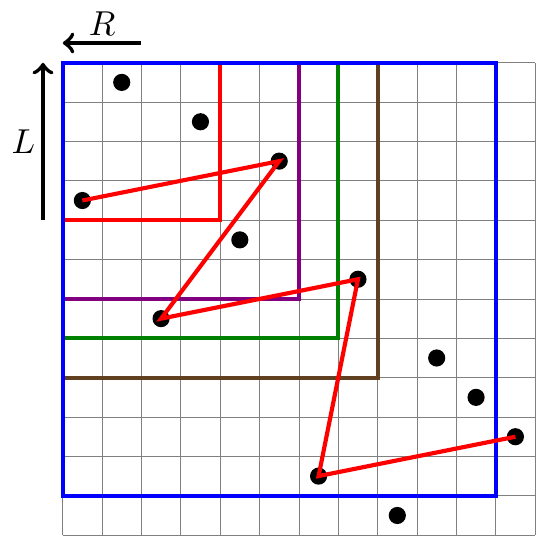}
\end{center}
\caption{\label{fig.irred_greedy}We construct the shortest $L$ zig-zag
  path of $\s$ in a greedy way: we take $(i_1,j_1)=(1,\s(1))$, then
  draw the square (here in red) and take $(i_2,j_2)$ the rightmost point above
  row $j_1$.  If it were inside the square, the permutation would have been
  reducible.  It is not the case here, so we draw a larger square
  (here in violet), and take $(i_3,j_3)$ as the lowest point to the left
  of column $i_2$. We continue this way, up to when the point
  $(i_6,j_6)$ is reached. As $i_6=n$, we terminate by certifying the
  irreducibility of~$\s$.
}
\end{figure}

Figure \ref{fig.irred_greedy} illustrates an example of the
construction of the greedy $L$ zig-zag path, in matrix representation,
and the caption illustrates the run of the algorithm.

Let us analyse a bit more carefully the complexity.  One can see that
the worst case is for irreducibility, and is attained on the
permutation $\s$ such that the vectors $(i_{a+1},j_{a+1})-(i_a,j_a)$
of the greedy $L$ or $R$ zig-zag path are 
\ldots $(3,1)$, $(-1,-3)$, $(3,1)$, $(-1,-3)$\ldots, with an exception at the
first and last elements (either a $(3,1)$ replaced by $(2,1)$, or a
$(-1,-3)$ replaced by $(-1,-2)$, depending on parities).
In this case the algorithm takes $n-1$ rounds, and overall performs
all of the $2n$ queries.

An analogous construction can be done by starting with $j_1=n$. By
this, in the reducible case we find the same certificating block, and
in the irreducible case we construct the greedy $R$ zig-zag path.

One can see that at least one among the two greedy paths ($L$ and $R$)
has minimal length among all zig-zag paths, and that the length of the
two greedy paths differ by at most 1. We call this value the
\emph{level} of~$\s$, and denote it by the symbol~$\lev(\s)$.

Of course, an irreducible permutation $\s$ of level 1 is a standard
permutation. More is true: 
\begin{lemma}
\label{lem.zigzagred}
An irreducible permutation $\s$ is at alternation distance at most
$\lev(\s)-1$ from a standard permutation.
\end{lemma}
\pf
Say that the level of $\s$ is $\ell>1$, and is realised by its greedy
$L$ zig-zag path (the argument for $R$ is symmetric). So we have
indices $j_2 < j_1$ and $i_1 = 1 < i_2$.

The action of $L^{n-j_2}$ moves the bullet $(i_2,j_2)$ in position
$(i_2,n)$ without affecting any other edge of the path.
Indeed, because of the inequalities satisfied by the indices of the
greedy $L$ zig-zag path, $j_a < j_1$ for all $a \geq 3$.  Then, the
path obtained by dropping the first edge is a $R$ zig-zag path of
level $\ell-1$ for the new configuration $\s'=L^{n-j_2} \s$, and the
two configurations $\s$, $\s'$ are nearest neighbour for alternation
distance.  The reasoning in the other case is analogous, with $R
\leftrightarrow L$ and $i_a \leftrightarrow n-j_a+1$.~\qed

Now Lemma \ref{lem.stdPerms} follows as an immediate corollary.



\subsection{Standard families}
\label{ssec.typeXH}

Here we introduce a property of permutations which is not invariant
under the dynamics, nonetheless it is useful in a combinatorial
decomposition of the classes, both for classification and enumeration
purposes.

\begin{definition}
\label{def.XHtype}
A permutation $\s$ is \emph{of type $H$} if the rank path goes through
the $-1$ mark, and \emph{of type $X$} otherwise. In the case of a
type-$X$ permutation, we call \emph{principal cycle} the cycle going
through the $-1$ mark.
\end{definition}
\noindent
The choice of the name is done for mnemonic purposes. Imagine to cut
the cycle invariant at the $-1$ mark. Then we have two open paths. In
a type-$H$ permutation, these paths have a $\asymp$ connectivity
pattern (which reminds of a $H$), and in a type-$X$ permutation they
have a $\times$ pattern (which reminds of a $X$).

\begin{figure}[tb!]
\begin{center}
\includegraphics{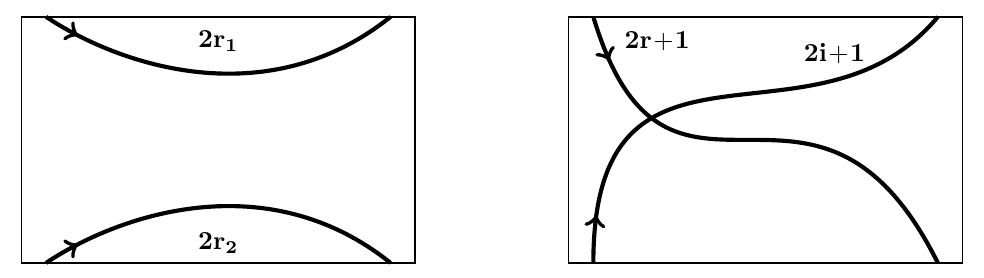}
\caption{\label{fig.typeXH}Left: a schematic representation of a
  permutation of type $H(r_1,r_2)$. Right: a representation of a
  permutation of type $X(r,i)$. These configurations have rank
  $r_1+r_2-1$ and $r$, respectively.}
  \end{center}
\end{figure}

A more refined definition is as follows
\begin{definition}
A permutation $\s$ is \emph{of type $H(r_1,r_2)$} if it is of type
$H$, and the portions of the rank path before and after the $-1$ mark
have path-length $2r_1$ and $2r_2$, respectively (so the rank length
is $r=r_1+r_2-1$). It is
\emph{of type $X(r,i)$} if it is of type $X$, has rank $r$, and the
principal cycle has path-length $2i+1$ (i.e., it contributes a $i$ entry to
the cycle structure $\lambda$).
\end{definition}
\noindent
See Figure \ref{fig.typeXH} for a schematic illustration.

A standard permutation is a notion with a simple definition. More
subtle is the associated notion:
\begin{definition}[standard family]
Let $\s$ be a standard permutation. The collection of $n-1$
permutations $\{ \s^{(i)} := L^{i} \s \}_{0 \leq i \leq n-2}$ is
called the \emph{standard family} of $\s$.
\end{definition}

\noindent
The properties established in Section \ref{ssec.arcbased} ultimately
imply the following statement:

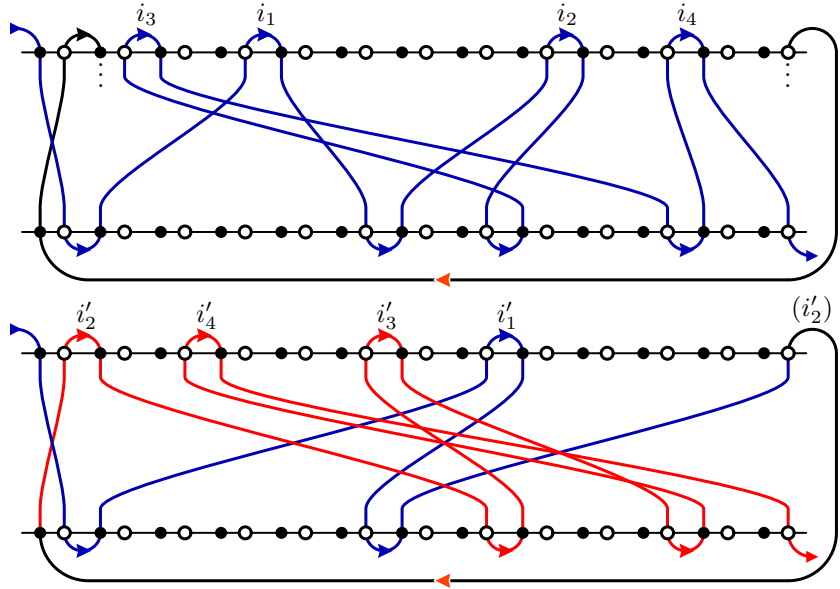
\begin{figure}[b!!]
\begin{center}
\setlength{\unitlength}{22.5pt}
\begin{picture}(14,5)
\put(0,0){\includegraphics[scale=4.5]{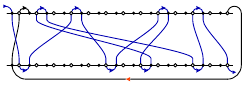}}
\put(2.3,4.7){$i_3$}
\put(4.3,4.7){$i_1$}
\put(9.3,4.7){$i_2$}
\put(11.3,4.7){$i_4$}
\put(1.65,3.6){$\vdots$}
\put(13.08,3.6){$\vdots$}
\end{picture}
\\
\setlength{\unitlength}{22.5pt}
\begin{picture}(14,5)
\put(0,0){\includegraphics[scale=4.5]{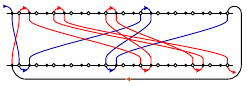}}
\put(6.3,4.7){$i'_3$}
\put(8.3,4.7){$i'_1$}
\put(13.2,4.8){$(i'_2)$}
\put(1.3,4.7){$i'_2$}
\put(3.3,4.7){$i'_4$}
\end{picture}
\caption{\label{fig.giraHstd}Illustration of the proof of an aspect of
  the fourth property in Proposition~\ref{trivialbutimp}.  Top:
  portions of a configuration $\s$, with rank 4, of type
  $X(4,\ast)$. The top arcs in the rank are at positions $i_1$,
  \ldots, $i_4$, and numbered according to their order along the
  path. Bottom: the result of applying $L^{n-i_2}$ to $\s$. The new
  configuration is of type $H(2,3)$, as the blue and red paths have
  path-length 4 and 6, respectively.}
\end{center}
\end{figure}

\begin{proposition}[Properties of the standard family]
\label{trivialbutimp}
Let $\s$ be a standard permutation, and 
$S=\{\s^{(i)}\}_i = \{L^i(\s)\}_i$ its standard family. The latter has
the following properties:
\begin{enumerate}
  \item Every $\tau \in S$ has $\tau(1)=1$;
  \item The $n-1$ elements of $S$ are all distinct;
  \item There is a unique $\tau \in S$ such that $\tau(n)=n$;
  \item Let $m_i$ be the multiplicity of the integer $i$ in $\lambda$
    (i.e.\ the number of cycles of length $i$), and $r$ the rank.
    There are $i\, m_i$ permutations of $S$ which are of type
    $X(r,i)$, and $1$ permutation of type $H(r-j+1,j)$, for each $1
    \leq j\leq r$.\footnote{Note that, as $\sum_i i\, m_i+r=n-1$ by
      the dimension formula (\ref{eq.size_inv_cycle}), this list
      exhausts all the permutations of the family.}
  \item Among the permutations of type $X(r,i)$ there is at least one
    $\tau$ with $\tau^{-1}(2) <\tau^{-1}(n)$.
 \end{enumerate}
\end{proposition}

\pf
The first three statements are obvious.

The fourth statement is based on the definitions of Section
\ref{ssec.arcbased}.  Let us analyse the top arcs: the cycle or path
going through the $-1$ mark is also the one going through the
rightmost arc, thus, as the family $L^i(\s)$, for $1\leq i \leq n-1$,
corresponds to a cyclic shift of the `colours' on these arcs, we
deduce that each cycle/path goes through the $-1$ mark a number of
times identical to its length value, contributing to the cycle
invariant. I.e., there are 
$i\, m_i$ permutations in $S$ of type $X(r,i)$, and $r$ permutations
of type $H$, as claimed.  In order to prove that, among the latter, we
have exactly one permutation per type $H(r-j+1,j)$, for $1 \leq j\leq
r$, we need to add some more structure. Let us number the top arcs of
the rank path in $\s$ from $1$ to $r$, following the order by which
they are visited by the path, starting from the top-left corner. Then
it becomes clear that $L^i(\s)$ has type $H(r-j+1,j)$ exactly when the
arc $n-i$ is the $(r-j+1)$-th arc of the rank path (this is also
illustrated in Figure~\ref{fig.giraHstd})

Let us now pass to the fifth statement.
As we are in type $X$, we have at least one cycle.  Let $a_1,\ldots,
a_k$ be the arcs associated to a cycle of length $k$. For
$j=1,\ldots,k$ we consider the edge $(\ell_{j,1},\s(\ell_{j,1}))$
incident to the left endpoint of the arc $a_j$, and the edge
$(\ell_{j,2},\s(\ell_{j,2})=\s(\ell_{j,1})+1)$ incident to the right
endpoint. Since we have a cycle (and not the rank path), at least one
of these arcs (say $a_j$) must have $\ell_{j,1}>\ell_{j,2}$, because
the $\ell_j$'s form a cyclic sequence of distinct integers, thus must
have at least one descent in this list.

This reasoning may not work if we deal with the principal cycle, as we
may have only one descent, in correspondence of the $-1$ mark. There
are two possible workarounds to this apparent problem. The first one
is that we can start by looking at a configuration $\tau_0$ of type
$H$ (there must be at least one such configuration in the family,
because the rank is at least~$1$). In this configuration no cycle is
the principal cycle, and we can choose one arc per cycle with the
required property. The second workaround is the fact that the
reasoning doesn't apply only when the principal cycle has only one
descent, and in correspondence of the $-1$ mark, but in fact in this
very case we just have $\tau^{-1}(2) <\tau^{-1}(n)$ with no need of
further analysis. 

\begin{figure}[t!!]
\begin{center}
\begin{tabular}{cp{0cm}c}
\includegraphics[scale=4.5]{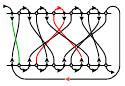}
&&
\includegraphics[scale=4.5]{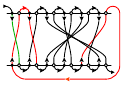}
\end{tabular}
\caption{\label{fig.lemtrivbutimp}Illustration of the proof of the
  fifth property in Proposition~\ref{trivialbutimp}.  Left: a
  configuration $\s$, with cycle invariant $(3,1;1)$, of type $X(1,1)$
  (in order to produce a small but generic example, we have taken a
  non-primitive configuration). The cycle of length $3$ is not the
  principal cycle, thus it must have an arc such that
  $\ell_{j,1}>\ell_{j,2}$. One such arc, and the two incident edges,
  are denoted in red.  Right: the configuration $\s'$ obtained by
  applying $L$ the appropriate number of times, so to make this arc go
  through the $-1$ mark. At this point, if we remove the left-most
  edge, the resulting configuration is irreducible.}
\end{center}
\end{figure}

Let us continue the analysis of the generic case, of an arc $a_j$ with
$\ell_{j,1}>\ell_{j,2}$.  Then $\s' = L^{n-a_j}(\s)$ has a principal
cycle of length $k$, and these two edges, in $\s'$, become
$(\ell_{j,1},n)$ and $(\ell_{j,2},2)$ (see figure
\ref{fig.lemtrivbutimp}). This allows to establish the fifth property.
\qed

\begin{corollary}
\label{trivialbutimpCor}
As a consequence of the property (5), let $\tau'$ be the permutation
resulting from removing the edge $(1,1)$ from a permutation $\tau$ as
in (5).  Then $\tau'$ is irreducible, and there exists $\ell$ (namely,
$\ell=\tau^{-1}(2)-2$) such that $R^{\ell}(\tau')$ is standard.
\end{corollary}
\noindent
This corollary will be of crucial importance in establishing the
appropriate conditions of the surgery operators $q_1$ and $q_2$,
outlined in Section~\ref{ssec.TQSintro}. (More precisely it will be a
key element of the proof of Theorems \ref{thm.q1} and \ref{thm.q2}.)

\pf For $\tau$ as in the proof of the property (5), note that removing
the edge $(1,1)$ provides an irreducible permutation $\tau'$.
Moreover it is clear that $\tau'$ has level $\lev(\tau')=2$, as in
fact $R^{\ell_{j,2}-2}(\tau')$ is standard.\footnote{We have an
  exponent $\ell_{j,2}-2$, instead of $\ell_{j,2}-1$ because we have
  removed the edge $(1,1)$.}  \qed

\subsection{Reduced dynamics and boosted dynamics}
\label{ssec.reddyn}

Some parts of our proof are just obtained from the inspection of
certain given (finite) patterns. How can such a simple ingredient be
compatible with a classification of classes of arbitrary size? A
crucial notion is the definition of `reduced' dynamics and `boosted' dynamics,
i.e., the analysis of the behaviour of \emph{patterns} in
configurations under the dynamics.\footnote{Here `pattern' is said in
  the sense of Permutation Pattern Theory, and is clarified in the
  following.}

Given a permutation $\s$, we can partition the set of edges into two
colors: black and gray. Say that 
$c: [n] \to \{ \textrm{black}, \textrm{gray} \}$ describes this
colouration.  The permutation $\tau$ corresponding to the restriction
of $\s$ to the black edges is called the \emph{reduced permutation}
for $(\s,c)$. Conversely, call $\hat{c}$ a data structure required to
reconstruct $\s$ from the pair $(\tau, \hat{c})$. A choice for
$\hat{c}$ is as an unordered list of quadruples, $(i,k|j,h)$. Such a
quadruple means that there is a gray edge with bottom endpoint being
the $k$-th gray point at the right of the $i$-th black point, and top
endpoint being the $h$-th gray point at the left of the $j$-th black
point. These quadruples evolve in a simple way under the action of $L$
and $R$ (i.e., they evolve consistently with the base values $i$ and
$j$, while $k$ and $h$ stay put). As a corollary, if we consider the
dynamics on $\tau$ as an edge-labeled permutation, we can reconstruct
the evolution of $\hat{c}$ from the evolution of $\tau$.

Call \emph{pivots} of $\s$ the two edges $(1,\s(1))$ and
$(\s^{-1}(n),n)$, if we work in the $\perms_n$ dynamics. [If we are in
  $\permsex_n$, we call pivot also $(n,\s(n))$ and $(\s^{-1}(1),1)$]

For a pair $(\s,c)$, we say that $\s$ is \emph{proper} if no gray edge
of $\s$ is a pivot.  In this case the dynamics on $\tau$ extends to
what we call the \emph{boosted dynamics} on $\s$, defined as follows:
for every operator $H$ (i.e., $H\in\{L,R\}$ for $\perms_n$ and $H\in\{L,L',R,R'\}$ for $\permsex_n$), we define
$\alpha_H(\s,c)$ as the smallest positive integer such that
$H^{\alpha_H(\s,c)}(\s)$ is proper, and, for a sequence 
$S=H_k \cdots H_2 H_1$ acting on $\tau$, the sequence $B(S)$, the
\emph{boosted sequence} of $S$, acting on $\s$ is
$B(S)=H_k^{\alpha_k} \cdots H_2^{\alpha_2} H_1^{\alpha_1}$ for the
appropriate set of $\alpha_j$'s.

In other words, the boosted dynamics is better visualised as the
appropriate notion such that the following diagram makes sense:
for $(\tau, \hat{c})$ an edge-labeled permutation with colouring data,
as described above, and calling $\Phi$ the operator that reconstructs
$(\s,c)$ from it, we have
\[
\begin{CD}
(\tau, \hat{c})
@> S >>  (\tau', \hat{c}')  \\        
@VV{\Phi}V @VV{\Phi}V\\    
(\s,c)
@> B(S) >>   (\s',c')
\end{CD}
\]

\noindent
Working in the reduced dynamics gives concise certificates of
connectedness: we can prove that $\s \sim \s'$ by showing the
existence of a triple $(\tau,\hat{c},S)$ that allows to reconstruct
the full diagram above. The idea is that one can often show the
existence of a given pair $(\tau,S)$ (of finite size), and a family
$\{\hat{c}_n\}_{n \in \bN}$, which allows to prove the connectedness
of families $\{\s_n \sim \s'_n\}_{n \in \bN}$.  This is the method
used in section \ref{ssec.T_op} for the proof of the main theorem
\ref{thm_T_surjectif} on the $T$ operator.

In a slightly more general form, instead of having a finite sequence
$S$, we could have three finite sequences $S_-$, $S_0$ and $S_+$ such
that the connectedness of $\s_n$ and $\s'_n$ is proven through the
sequence $S_- S_0^n S_+$.  We use this generalised pattern only once,
in order to prove explicitly the connectedness of two permutations
inside a given class, in a case which is not covered by any of the
`big theorems' on the surgery operators (see Section \ref{ssec.tech2},
Lemma~\ref{lem.q2.maxrank})

\subsection{Square constructors for permutations}
\label{sec.sqconstr}

\begin{figure}[b!!]
\begin{center}
\includegraphics{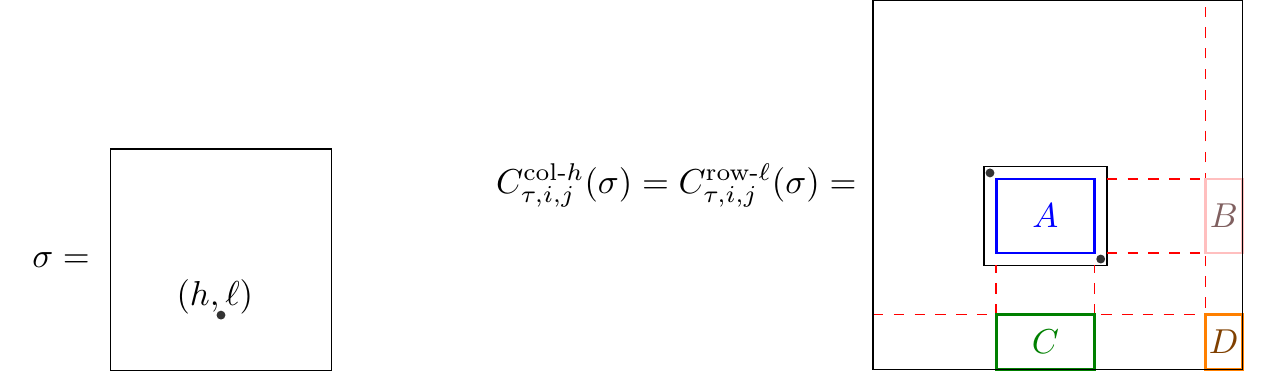}
\\
\setlength{\unitlength}{10pt}
\begin{picture}(27,13.5)(3,0)
\put(5.5,7.5){\makebox[0pt][c]{$(\tau,i,j)$}}
\put(3,3.8){\makebox[0pt][r]{$\scriptstyle{i}\left\{\rule{0pt}{13pt}\right.$}}
\put(3.1,5.1){\rotatebox{-90}{\makebox[0pt][r]{$\rotatebox{90}{$\scriptstyle{j}$}\left\{\rule{0pt}{21.5pt}\right.$}}}
\put(13,7.5){\makebox[0pt][c]{$\s$}}
\put(12.5,6.2){\makebox[0pt][c]{$\scriptstyle{h}$}}
\put(9.8,4.3){\makebox[0pt][r]{$\scriptstyle{\ell}$}}
\put(17.6,10){\makebox[0pt][r]{$C_{\tau,i,j}^{{\rm row-}\ell}(\s) =
C_{\tau,i,j}^{{\rm col-}h}(\s)$}}
\put(5,3.72){\makebox[0pt][c]{$A$}}
\put(7.5,3.72){\makebox[0pt][c]{$B$}}
\put(5,1.82){\makebox[0pt][c]{$C$}}
\put(7.5,1.82){\makebox[0pt][c]{$D$}}
\put(11,5.22){\makebox[0pt][c]{$X$}}
\put(14,5.22){\makebox[0pt][c]{$Y$}}
\put(11,2.22){\makebox[0pt][c]{$Z$}}
\put(14,2.22){\makebox[0pt][c]{$W$}}
\put(23,8.72){\makebox[0pt][c]{$A$}}
\put(29.5,8.72){\makebox[0pt][c]{$B$}}
\put(23,1.82){\makebox[0pt][c]{$C$}}
\put(29.5,1.82){\makebox[0pt][c]{$D$}}
\put(19,11.22){\makebox[0pt][c]{$X$}}
\put(27,11.22){\makebox[0pt][c]{$Y$}}
\put(19,5.22){\makebox[0pt][c]{$Z$}}
\put(27,5.22){\makebox[0pt][c]{$W$}}
\put(23,11.22){\makebox[0pt][c]{$\varnothing$}}
\put(19,8.72){\makebox[0pt][c]{$\varnothing$}}
\put(29.5,11.22){\makebox[0pt][c]{$\varnothing$}}
\put(27,8.72){\makebox[0pt][c]{$\varnothing$}}
\put(23,5.22){\makebox[0pt][c]{$\varnothing$}}
\put(19,1.82){\makebox[0pt][c]{$\varnothing$}}
\put(29.5,5.22){\makebox[0pt][c]{$\varnothing$}}
\put(27,1.82){\makebox[0pt][c]{$\varnothing$}}
\put(3,0){\begin{tikzpicture}[scale=0.35]
\permutationtwo{3,4,1,5,2}{4}{5};
\end{tikzpicture}}
\put(13,4.5){\oval(5.8,.8)}
\put(12.5,3){\oval(.8,5.8)}
\put(10,0){\begin{tikzpicture}[scale=0.35]
\permutation{4,2,5,1,6,3};
\end{tikzpicture}}
\put(24,7.5){\oval(11.8,.8)}
\put(20.5,6){\oval(.8,11.8)}
\put(24,10.5){\oval(11.8,.8)}
\put(25.5,6){\oval(.8,11.8)}
\put(18,0){\begin{tikzpicture}[scale=0.35]
\permutationtwo{7,5,
11,
3,9,1,10,
8,
4,12,6,
2
}{4}{12};
\end{tikzpicture}}
\end{picture}
\end{center}
\caption{Action of the square constructors on a permutation $\s$, with
  $\tau$ being a matrix with blocks $A$, $B$, $C$ and $D$. Top: the
  general structure. Bottom: an example.\label{figure.sqrcst}}
\end{figure}
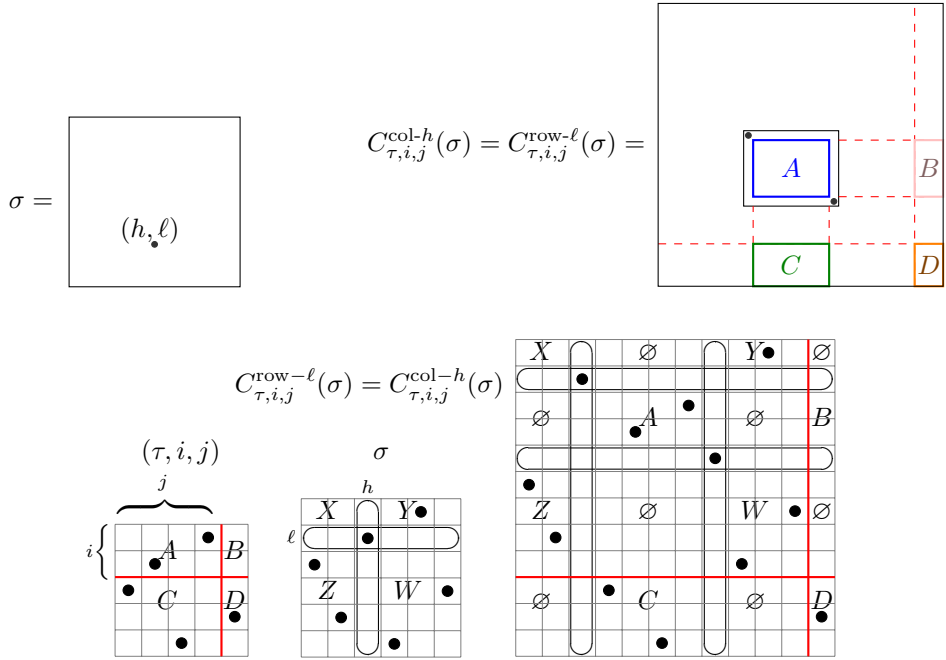

In this section we define a rather general surgery operation on
permutations, which behaves in a simple way under the Rauzy dynamics:
\begin{definition}[Square constructor]
Let $\tau$ be a permutation matrix of size $k$, and let $0 \leq i,j
\leq k$, with $i+j>0$. The pair $i,j$ describes a way to decompose
$\tau$ into four (possibly empty) rectangular blocks,
$\tau=\left( \begin{smallmatrix}A&B\\C&D\end{smallmatrix} \right)$,
so that $A$ and $B$ have $i$ rows, and $A$ and $C$ have $j$
columns. For $h,\ell \in \{1,\ldots,n\}$, and
$\s$ a permutation of size $n$ with $\ell=\s(h))$,
the \emph{square constructors} $C_{\tau,i,j}^{{\rm col-}h}$ and
$C_{\tau,i,j}^{{\rm row-}\ell}$, acting on $\s$, produce the same
configuration $\s'$, of size $n+k+1$, 
as described in Figure~\ref{figure.sqrcst}.
\end{definition}

\begin{lemma}[reduced dynamics for constructors]
\label{lem.reddynconstr}
Let $(\tau,i,j)$ as above.  Let $\s$ and $\s'$ be two permutations of
the same size, with one coloured edge (say, red and the rest is
black). Let $(h,\s(h))$ be the only entry of $\s$ marked in red and
$(h',\s'(h'))$ be the only entry of $\s'$ marked in red.
If $\s \sim \s'$ for the dynamic $\perms$, with colours as above, then
$C_{\tau,i,j}^{{\rm col-}h}(\s) \sim 
C_{\tau,i,j}^{{\rm col-}h'}(\s')$.
\end{lemma}

\pf 
This is a case of reduced dynamics, where, in 
$C_{\tau,i,j}^{{\rm col-}h}(\s)$, the entries of $\tau$ plus the extra
entry at the bottom-right corner of block $A$ are gray. The red point
is the point at the top-left of block $A$, this reduced permutation is
$(\s,h)$ where $h$ is the index of the red entry. Thus we can define a
bijection $\phi$ between $C_{\tau,i,j}^{{\rm col-}h}(\s)$ and
$(\s,h)$. Now suppose $\s$ and $\s'$ are connected by the sequence of
operators $S$, i.e.\ $S(\s)=\s'$. Then we need to show that
$B(S)(\phi^{-1}(\s,h))=\phi^{-1}(\s',h')$.

By transitivity of connectedness, it is not necessary to consider
arbitrary sequences, it is enough to consider $\s'=L\s$ and $\s'=R\s$.
The symmetry of the definition of the square constructor allows to
consider just one case, and we choose $L$.

\begin{figure}[t!!]
\begin{center}
\includegraphics{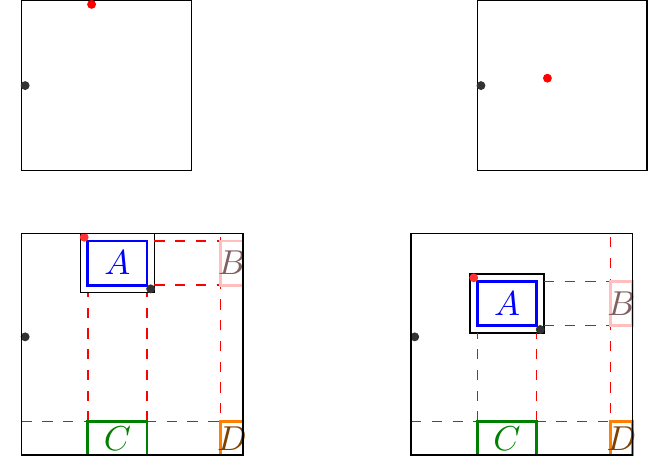}
\caption{\label{fig.red_pivot}Top left: $(\s,h)$. Top right:
  $(L(\s),h)$. Bottom left: $C_{\tau,i,j}^{{\rm col-}h}(\s)$. Bottom
  right: $B(L)(C_{\tau,i,j}^{{\rm col-}h}(\s))$.}
\end{center}
\end{figure}

If the red point is not a pivot, then neither are the gray edges, and
$B(L)=L$, so $L(\phi^{-1}(\s,h))=\phi^{-1}(\s',h')$ as wanted. Thus we
only need to consider the case when the red point is a pivot. In such
a case, (i.e.\ the red point is $(h,n)$), then $B(L)=L^{i+1}$ and
$B(L)(\phi^{-1}(\s,h))=\phi^{-1}(L(\s),h)$ as wanted (see
figure~\ref{fig.red_pivot}).  \qed

\begin{figure}[b!!]
\begin{center}
\includegraphics{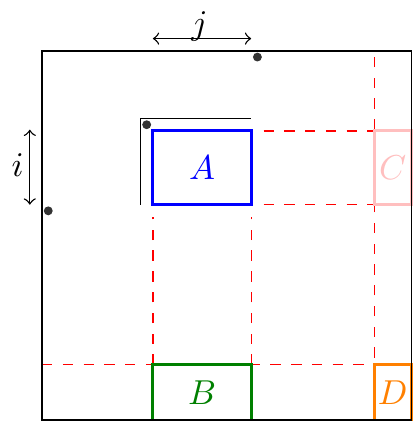}
\end{center}
\caption{Structure of a configuration $\s'$, such that $R^{-j-1} \s' =
  C_{\tau,i,j}^{{\rm col-}1}(\s)$ and $L^{-i-1} \s' =
  C_{\tau,i,j}^{{\rm row-}n}(\s)$.\label{figure.sqrLR}}
\end{figure}

\begin{lemma}[square transportation]\label{lem.squaretransportation}
Let $(\tau,i,j)$ as above. Then, the action of the constructor at
different locations gives equivalent configurations, i.e.\ for all
$\s$ of size $n$ and all $h, h' \in [n]$ we have
$C_{\tau,i,j}^{{\rm col-}h}(\s) \sim C_{\tau,i,j}^{{\rm col-}h'}(\s)$.
\end{lemma}
\pf
We start by proving that the lemma holds for
$(h,h')=(1,\s^{-1}(n))$. This is not hard to see. In fact, 
$R^{j+1} C_{\tau,i,j}^{{\rm col-}1}(\s)$ and $L^{i+1}
C_{\tau,i,j}^{{\rm row-}n}(\s)$ are the same
configuration, namely the one illustrated in
Figure~\ref{figure.sqrLR}.

\begin{figure}
\[
\raisebox{-6pt}{\includegraphics[scale=.8]{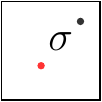}}
\xrightarrow{\mathrm{Std}_1}  
\raisebox{-6pt}{\includegraphics[scale=.8]{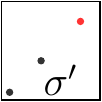}}
\xrightarrow{L^iEL^{-i}}
\raisebox{-6pt}{\includegraphics[scale=.8]{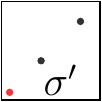}}
\xrightarrow{L^jEL^{-j}}
\raisebox{-6pt}{\includegraphics[scale=.8]{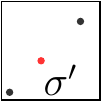}}
\xrightarrow{\mathrm{Std}_2^{-1}}
\raisebox{-6pt}{\includegraphics[scale=.8]{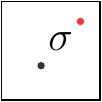}}
\]
\caption{\label{fig_seq_square_lemma}The sequence of operations which
  allow to transport the red mark from one position to another. Here
  $\mathrm{Std}$ is the standardasing sequence of the (uncoloured)
  configuration $\s$.}
\end{figure}

Now consider the reduced dynamics as in Lemma \ref{lem.reddynconstr}.
Let $(\s,c_L)$ the colouring of $\s$ with the red entry on the
$L$-pivot, and $(\s,c_R)$ the one with the red entry on the
$R$-pivot. Let $\phi$ be the map of the reduced dynamics.  We have
just proven that $\phi^{-1}(\s,c_L) \sim \phi^{-1}(\s,c_R)$.  Let $E$
be the operator that exchanges $(\s,c_L)$ and $(\s,c_R)$, i.e.
\[
\includegraphics{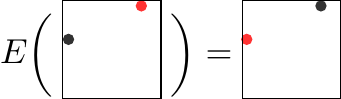}
\]
By Lemma \ref{lem.stdPerms} we know that, in the reduced dynamics, we
have a sequence $S$ in the dynamics such that $S \s = \s'$, and $\s'$ is
standard. In this configuration, either the
$L$-pivot is red, or exactly one configuration in the associated
standard family has a red $R$-pivot. In the first case, $E$ can be
applied to any configuration in the standard family. This shows that,
for $\s$ standard, all configurations $(\s,c_j)$ for the red point on
the $j$-th column are connected to $(\s,c_L)$, by a sequence of the form
$L^{-k}EL^{k}$, and thus are connected among themselves by a sequence of
the form $L^{-h}EL^{h-k}EL^{k}$. Lemma \ref{lem.reddynconstr} allows
us to conclude, for $\s'$ a standard permutation, and, by conjugating
with the sequence $S$, for $\s$ a generic permutation in the class
(see Figure~\ref{fig_seq_square_lemma}).

\qed

\begin{corollary}
\label{cor.construtor}
Let $(\tau,i,j)$ as above.  Let $\s$ and $\s'$ such that $\s \sim \s'$
for the dynamic $\perms$, then, for every $1 \leq h, h' \leq n$,
$C_{\tau,i,j}^{{\rm col-}h}(\s) \sim C_{\tau,i,j}^{{\rm col-}h'}(\s')$.
\end{corollary}

\noindent
A further corollary of this lemma is the statement of Proposition
\ref{prop.primhomoPre}, for which it suffices to take $A$ a diagonal
matrix and $B$, $C$, $D$ empty blocks.

\section{The sign invariant}
\label{sec.signinv}

\subsection{Arf functions for permutations}
\label{ssec.arf_inv}

For $\s$ a permutation in $\kS_n$, let 
\be
\chi(\s)
= \# \{ 1\leq i<j \leq n \;|\; \s(i)<\s(j) \}
\ee
i.e.\ $\chi(\s)$ is the number of pairs of non-crossing edges in the
diagram representation of~$\s$. 
 
Let $E=E(\s)$ be the subset of $n$ edges in $\cK_{n,n}$ described by
$\s$.  For any $I \subseteq E$ of cardinality $k$, the permutation
$\s|_{I} \in \kS_k$ is defined in the obvious way, as the one
associated to the subgraph of $\cK_{n,n}$ with edge-set $I$, with
singletons dropped out, and the inherited total ordering of the two
vertex-sets.

Define the two functions
\begin{align}
A(\s)
&:=
\sum_{I \subseteq E(\s)} (-1)^{\chi(\s|_I)}
\ef;
&
\Abar(\s)
&:=
\sum_{I \subseteq E(\s)} (-1)^{|I|+\chi(\s|_I)}
\ef.
\end{align}
When $\s$ is understood, we will just write $\chi_I$ for
$\chi(\s|_I)$. The quantity $A$ is accessory in the forthcoming
analysis, while the crucial fact for our purpose is that the quantity
$\Abar$ is invariant both in the $\perms_n$ and in the $\permsex_n$
dynamics. We will prove this in the following subsection.

As a result of a well-known property of Arf functions, rederived in
Section \ref{ssec.signindu} (Lemma \ref{lem.induSign}), the quantity
$\Abar(\s)$ is either zero or $\pm$ a power of 2.
For now, we content ourselves with the definition
\begin{definition}
For $\s$ a permutation,
define $s(\s)$, the \emph{sign} of $\s$, as the quantity
\be
s(\s) := \mathrm{Sign}(\Abar(\s)) \in \{0,\pm1\}
\ef.
\ee
\end{definition}
\noindent
However, we will see (in lemma \ref{lem.induSign}) that $s(\s)= 2^{-\frac{n+\ell}{2}} \Abar(\s)$,
where $\ell$ is the number of cycles (the rank path is not counted) .

\subsection{Calculating with Arf functions}
\label{ssec.arfcalcseasy}

Evaluating the functions $A$ and $\Abar$ on a generic `large'
permutation, starting from the definition, seems a difficult task, as
we have to sum an exponentially large number of terms. However, as the
name of Arf function suggests, these quantities are combinatorial
counterparts of the classical Arf invariant for manifolds (in this
case, the translation surfaces where the Rauzy dynamics is defined,
see Section \ref{ssec.geo_inv} for more details), and inherits from
them the covariance associated to the change of basis in the
corresponding quadratic form, and its useful consequences. We will
exploit this to some extent, in a simple and combinatorial way, that
we now explicitate.

The point is that we will \emph{not} try to evaluate Arf functions of
large configurations starting from scratch. We will rather compare the
Arf functions of two (or more) configurations, which differ by a
finite number of edges, and establish linear relations among their Arf
functions.  Yet another tool is proving that the Arf function of a
given configuration is zero, by showing that it contains some finite
\emph{pattern} that implies this property.

\begin{figure}[tb!]
\begin{center}
\includegraphics{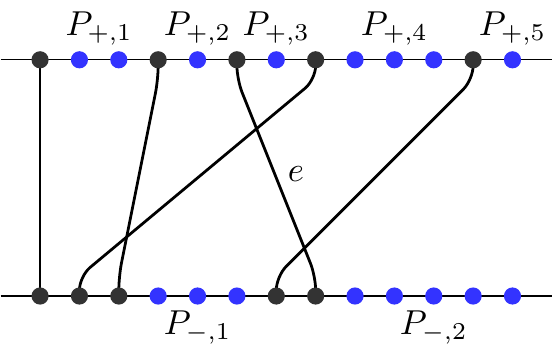}
\end{center}
\caption{\label{fig.arf_ex_def}A subset of edges in a permutation, as
  in Definition \ref{def.637647}. The represented edge set is $E'$,
  the blue vertices at the bottom constitute the set $V''_-$, while
  those at the top are $V''_+$. The two partitions denoted by the
  symbols $P_{\pm,i}$ are the coarsest partitions which are compatible
  with $E'$. We cannot show the full matrix $Q$ for such a big
  example, but we can give one row, for the edge which has the label
  $e$ in the drawing. The row $Q_e$ reads
$(Q_e)_{11, 12, \ldots, 15, 21, \ldots, 25}=(1,1,0,0,0,\,0,0,1,1,1)$.}
\end{figure}

In order to have the appropriate terminology for expressing this
strategy, let us define the following:


\begin{definition}
\label{def.637647}
Call $V_{\pm}$ the two vertex-sets of $\cK_{n,n}$.  For $\s$ a
permutation with edge set $E$, and $E' \subseteq E$ of cardinality
$m$, call $V''_{\pm}$ the subset of $V_{\pm}$ not containing the
endpoints of $E'$. A pair of partitions $P_+=(P_{+,1},\ldots,P_{+,h})
\in \cP(V''_+)$ and $P_-=(P_{-,1},\ldots,P_{-,k}) \in \cP(V''_-)$ are
said to be \emph{compatible with $E'$} if each block contains indices
which are consecutive in $[n]$.

Define the $m \times (hk)$ matrix valued in $\gf_2$
\be
Q_{e,ij} :=
\left\{
\begin{array}{ll}
1 & \textrm{edge $e \in E'$ does not cross the segment connecting
  $P_{-,i}$ to $P_{+,j}$,}\\
0 & \textrm{otherwise.}
\end{array}
\right.
\ee
For $v \in \gf_2^d$, let $|v|$ be the number of entries equal to
$1$. Similarly, identify $v$ with the corresponding subset of $[d]$.
Given such a construction, introduce the following functions on
$(\gf_2)^{hk}$
\begin{align}
A_{\s,E',P}(v)
&:=
\sum_{u \in (\gf_2)^{E'}}
(-1)^{\chi_u + (u,Qv)}
\ef;
&
\Abar_{\s,E',P}(v)
&:=
\sum_{u \in (\gf_2)^{E'}}
(-1)^{|u|+\chi_u + (u,Qv)}
\ef.
\end{align}
\end{definition}
The construction is illustrated in
Figure~\ref{fig.arf_ex_def}.

Let us comment on the reasons for introducing such a definition. The
quantities $A_{\s,E',P}(v)$ allows to sum together many contributions
to the function $A$, which behave all in the same way once the subset
restricted to $E'$ is specified. In a way similar to the philosophy
behind the reduced dynamics of Section \ref{ssec.reddyn}, our goal is
to have $E'$ of fixed size, while
$E \setminus E'$ is arbitrary and of unbounded size, so that the
verification of our properties, as it is confined to the matrix $Q$,
involves a finite data structure.

Indeed, let us split in the natural way the sum over subsets $I$ that
defines $A$ and $\Abar$, namely
\[
\sum_{I \subseteq E} f(I)
=
\sum_{I' \subseteq E'} 
\sum_{I'' \subseteq E \setminus E'} f(I' \cup I'')
\]
For $I$ and $J$ two disjoint sets of edges, call $\chi_{I,J}$ the number
of pairs $(i,j) \in I \times J$ which do not cross. Then clearly
\[
\chi_{I \cup J} = \chi_{I} + \chi_{J} + \chi_{I,J}
\]
Now let $u(I') \in \{0,1\}^{E'}$ be the vector with entries $u_e=1$ if
$e \in I'$ and $0$ otherwise. Let $m(I'')=\{m_{ij}(I'')\}$ be the 
$k \times h$ matrix describing the number of edges connecting the
intervals $P_{-,i}$ to $P_{+,j}$ in $\s$, and
$v(I'')=\{v_{ij}(I'')\}$, $v_{ij} \in \{0,1\}$, as the parities of the
$m_{ij}$'s. Call $I''_{ij}$ the restriction of $I''$ to edges
connecting $P_{-,i}$ and $P_{+,j}$.  Clearly,
$\chi_{I',I''}=\sum_{ij} \chi_{I',I''_{ij}} 
= \sum_{e,ij} u_e Q_{e,ij} m_{ij}$, 
which has the same parity as the analogous expression with $v$'s
instead of $m$'s. Now, while the $m$'s are in $\bN$, the vector $v$ is
in a linear space of finite cardinality, which is crucial for allowing
a finite analysis of our expressions.

As a consequence,
\begin{align}\label{eq.arf_explain_def}
A(\s)
&=
\sum_{I \subseteq E(\s)\setminus E'(\s)}(-1)^{\chi_I}A_{\s,E',P}(v(I))
\ef;
\\\label{eq.arf_explain_def_1}
\Abar(\s)
&=\sum_{I \subseteq E(\s)\setminus E'(\s)}(-1)^{|I|+\chi_I}\Abar_{\s,E',P}(v(I))
\ef.
\end{align}

\noindent
Then we have two criteria for establishing relations among Arf
functions
\begin{proposition}
\label{prop.Arf0}
Let $\s$, $E'$ and $P$ as above. If for all $v \in (\gf_2)^{hk}$ we
have $A_{\s,E',P}(v)=0$, then $A(\s)=0$.  The same holds for $\Abar$.
\end{proposition}
\begin{proposition}
\label{prop.ArfProp}
Let $\s$ and $\t$ be permutations (possibly of different size), with edge
sets $E \cup E'_{\s}$ and $E \cup E'_{\t}$ respectively.  Let 
$P=(P_-, P_+)$ be a pair of partitions of size $|E|$ compatible with both
$E'_{\s}$ and $E'_{\t}$. If there exists $K \in \bQ$ such that
for all $v \in (\gf_2)^{hk}$ we have $A_{\s,E'_{\s},P}(v)=K\,
A_{\t,E'_{\t},P}(v)$, then $A(\s)=K\, A(\t)$. The same holds with one
or both of the $A$'s replaced by $\Abar$.
Analogous statements hold for linear combinations of Arf functions
associated to more than two configurations.\footnote{We refer here to
  the obvious generalisation, of the form, for $\s_j$ having edge-set
  $E \cup E'_{\s_j}$, $P=(P_-, P_+)$ a pair of partitions of size
  $|E|$ compatible with all $E'_{\s_j}$'s, if there exist $K_j \in
  \bQ$ s.t.\ for all $v \in (\gf_2)^{hk}$ we have $\sum_j K_j
  A_{\s_j,E'_{\s_j},P}(v)=0$, then $\sum_j K_j A(\s_j)=0$, and
  similarly with some of the $A$'s replaced by $\Abar$'s.}
\end{proposition}
\noindent
The proof of these propositions is an immediate consequence of
equations~(\ref{eq.arf_explain_def}) and (\ref{eq.arf_explain_def_1}).

For this purpose of our main classification theorems, we need four facts
which are specialisations of the propositions above. One of them
establish the invariance of function $\Abar$ under the dynamics, the
other three relate the Arf functions on configurations obtained from
the ``surgery operators'' sketched in Section~\ref{ssec.TQSintro} and
discussed in Section~\ref{sec.operators}, and a few further, simpler
manipulations.
\begin{proposition}[Invariance of the sign]
\label{prop.signinvdyn}
\be
\put(52,-20){$P_{-,1}$}\put(88,-20){$P_{-,2}$}\put(36,20){$P_{+,1}$}\put(75,20){$P_{+,2}$}
\Abar\bigg(\,\t=\raisebox{-12pt}{\includegraphics[scale=1.8]{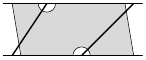}}\,\bigg)
\put(-65,0){$e_1$}\put(-44,-0){$e_2$}
=
\put(52,-20){$P_{-,1}$}\put(88,-20){$P_{-,2}$}\put(36,20){$P_{+,1}$}\put(75,20){$P_{+,2}$}
\Abar\bigg(\,\s=\raisebox{-12pt}{\includegraphics[scale=1.8]{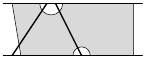}}\,\bigg)
\put(-65,0){$e_1$}\put(-44,-0){$e_2$}
\ee
\end{proposition}
\proof
We have in this case
\begin{align}
Q_{\t}&=
\begin{pmatrix}
0 & 1 & 0 & 1 \\
1 & 1 & 0 & 0
\end{pmatrix}
\ef;
&
Q_{\s}&=
\begin{pmatrix}
0 & 1 & 0 & 1 \\
1 & 0 & 0 & 1
\end{pmatrix}
\ef;
\end{align}
Checking that the conditions of Proposition \ref{prop.ArfProp} are
met, with $K=1$, is a straightforward calculation. The resulting
function of $v$ for $\tau$ is
\begin{align*}
\Abar_{\tau,\{e_1,e_2\},P}(v)
&=\sum_{u \in (\gf_2)^{2}}
(-1)^{|u|+\chi_u + (u,Qv)}\\
&=
(-1)^{0+0 + 0}+(-1)^{1+0 + \begin{pmatrixsm}0&1&0&1\end{pmatrixsm}\cdot v}+(-1)^{1+0 + \begin{pmatrixsm}1&1&0&0\end{pmatrixsm}\cdot v}
\\
&\qquad\qquad\quad  +(-1)^{0+1 +  \left(\begin{pmatrixsm}0&1&0&1\end{pmatrixsm}+ \begin{pmatrixsm}1&1&0&0\end{pmatrixsm}\right)\cdot v}\\
&=1
-(-1)^{\begin{pmatrixsm}0&1&0&1\end{pmatrixsm}\cdot v}
-(-1)^{\begin{pmatrixsm}1&1&0&0\end{pmatrixsm}\cdot v}
-(-1)^{\begin{pmatrixsm}1&0&0&1\end{pmatrixsm}\cdot v}
\end{align*}
and likewise the resulting function of $v$ for $\s$ is  

\begin{align*}
\Abar_{\s,\{e_1,e_2\},P}(v)
&=\sum_{u \in (\gf_2)^{2}}
(-1)^{|u|+\chi_u + (u,Qv)}\\
&=
(-1)^{0+0 + 0}+(-1)^{1+0 + \begin{pmatrixsm}0&1&0&1\end{pmatrixsm}\cdot v}+(-1)^{1+0 + \begin{pmatrixsm}1&0&0&1\end{pmatrixsm}\cdot v}\\
&\qquad\qquad\quad+ (-1)^{0+1 +  \left(\begin{pmatrixsm}0&1&0&1\end{pmatrixsm}+ \begin{pmatrixsm}1&0&0&1\end{pmatrixsm}\right)\cdot v}
\\
&=1
-(-1)^{\begin{pmatrixsm}0&1&0&1\end{pmatrixsm}\cdot v}
-(-1)^{\begin{pmatrixsm}1&0&0&1\end{pmatrixsm}\cdot v}
-(-1)^{\begin{pmatrixsm}1&1&0&0\end{pmatrixsm}\cdot v}
\end{align*}
thus we have
\be
\Abar_{\tau,\{e_1,e_2\},P}(v)=\Abar_{\s,\{e_1,e_2\},P}(v) 
\text{\qquad for all $v \in \{0,1\}^4$.}
\ee
\noindent
This proposition implies the invariance of $\Abar$ under the operation
$L$, as $\t=L \s$. The invariance under $R$, $L'$ and $R'$ is deduced
from the symmetry of the definition of $A$ and $\Abar$ under the
dihedral group on permutation diagrams.

It is convenient to introduce the notation
$\vec{A}(\s)=\begin{pmatrixsm} \Abar(\s) \\ A(\s) \end{pmatrixsm}$.
We have
\begin{proposition}
\label{prop.fingerred}
\begin{align}
\label{eq.546455a}
\vec{A}\bigg(\,
\raisebox{-12pt}{\includegraphics[scale=1.8]{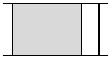}}\,\bigg)
&=
\begin{pmatrix}
1 & -1 \\ 1 & 1
\end{pmatrix}
\vec{A}\bigg(\,
\raisebox{-12pt}{\includegraphics[scale=1.8]{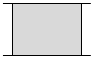}}\,\bigg)
\ef;
\\
\label{eq.546455b}
\vec{A}\bigg(\,
\raisebox{-12pt}{\includegraphics[scale=1.8]{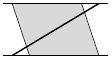}}\,\bigg)
&=
\begin{pmatrix}
0 & 0 \\ 0 & 2
\end{pmatrix}
\vec{A}\bigg(\,
\raisebox{-12pt}{\includegraphics[scale=1.8]{Figure4_fig_arfproof_s3.pdf}}\,\bigg)
\ef;
\\
\label{eq.546455c}
\vec{A}\bigg(\,
\raisebox{-12pt}{\includegraphics[scale=1.8]{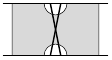}}\,\bigg)
&=
2
\vec{A}\bigg(\,
\raisebox{-12pt}{\includegraphics[scale=1.8]{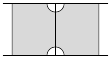}}\,\bigg)
\ef;
\\
\label{eq.546455d}
\vec{A}\bigg(\,
\raisebox{-12pt}{\includegraphics[scale=1.8]{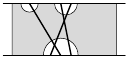}}\,\bigg)
&=
\vec{A}\bigg(\,
\raisebox{-12pt}{\includegraphics[scale=1.8]{Figure4_fig_arfproof_t8.pdf}}\,\bigg)
\ef.
\end{align}
\end{proposition}
\noindent
These relations are all straightforward applications of Proposition
\ref{prop.ArfProp}, with matrices $Q$ of rather small dimension.
Equation (\ref{eq.546455d}), involving a larger matrix $Q$, may be
derived with a shortcut: on the LHS, the sum over $u \in (\gf_2)^3$
contains the four terms contributing to the RHS, where the left-most
edge is absent, and four other terms, which cancel pairwise for
obvious reasons, even with no need for analysing the matrices $Q$ in
detail.  \qed

Equation (\ref{eq.546455a}) implies
\begin{corollary}\label{cor.id_sign}
\be
\begin{split}
\vec{A}(\id_n)
&=
\begin{pmatrix}
1 & -1 \\ 1 & 1
\end{pmatrix}^{n+1}
\begin{pmatrix}
1\\0
\end{pmatrix}
=
2^{\frac{n+1}{2}}
\begin{pmatrix}
\cos \frac{(n+1)\pi}{4} 
\\
\sin \frac{(n+1)\pi}{4} 
\end{pmatrix}
\\
&=
(-1)^{\lfloor\frac{n+1}{4}\rfloor}
\left\{
\begin{array}{ll}
\begin{pmatrix}
2^{\frac{n+1}{2}}
\\0
\end{pmatrix}
&
n \equiv -1 \quad \mathrm{(mod\ 4);}
\\
\rule{0pt}{20pt}%
\begin{pmatrix}
2^{\frac{n}{2}}\\ 2^{\frac{n}{2}}
\end{pmatrix}
&
n \equiv 0 \quad \mathrm{(mod\ 4);}
\\
\rule{0pt}{22pt}%
\begin{pmatrix}
0 \\ 2^{\frac{n+1}{2}}
\end{pmatrix}
&
n \equiv 1 \quad \mathrm{(mod\ 4);}
\\
\rule{0pt}{20pt}%
\begin{pmatrix}
-2^{\frac{n}{2}}\\
2^{\frac{n}{2}}
\end{pmatrix}
&
n \equiv 2 \quad \mathrm{(mod\ 4);}
\end{array}
\right.
\end{split}
\ee
\end{corollary}

\noindent
In particular, $\Abar(\id_4)=-4<0$ and $\Abar(\id_6)=8>0$, this fact is used in section \ref{ssec.tech}.

Equations (\ref{eq.546455c}) and (\ref{eq.546455d}) 
imply
\begin{corollary}
\label{cor.Tsign}
\be
\Abar\bigg(\,\s=\raisebox{-12pt}{\includegraphics[scale=1.8]{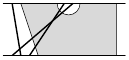}}\,\bigg)
=
2
\Abar\bigg(\,\t=\raisebox{-12pt}{\includegraphics[scale=1.8]{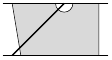}}\,\bigg)
\ee
\end{corollary}

\noindent
This corollary states that, for the operator $T$ defined in Section
\ref{ssec.T_op}, $\Abar(T \t)=2\, \Abar(\t)$.

Equations (\ref{eq.546455a}) and (\ref{eq.546455b}) 
imply
\begin{corollary}
\label{cor.q1sign}
\be
\Abar\bigg(\,\s=\raisebox{-12pt}{\includegraphics[scale=1.8]{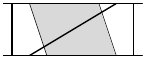}}\,\bigg)
=
2 \,
\Abar\bigg(\,\t=\raisebox{-12pt}{\includegraphics[scale=1.8]{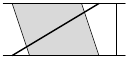}}\,\bigg)
\ef.
\ee
\end{corollary}
\noindent
This corollary states that, for the operator $q_1$ defined in Section
\ref{ssec.q_op}, $\Abar(q_1 \t)=2\, \Abar(\t)$, because in fact
$q_1 \s = R \t$.

More generally,
we have
\begin{proposition}
\label{prop.854423343}
\be
\vec{A}\bigg(\,\s=\raisebox{-12pt}{\includegraphics[scale=1.8]{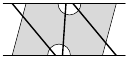}}\,\bigg)
=
\begin{pmatrix}
2 & 0 \\ 0 & 0
\end{pmatrix}
\vec{A}\bigg(\,\t=\raisebox{-12pt}{\includegraphics[scale=1.8]{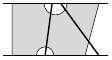}}\,\bigg)
\ef.
\ee
\end{proposition}
\noindent
that gives, using twice equation (\ref{eq.546455a}) in
Proposition \ref{prop.fingerred},
\begin{corollary}
\label{cor.q2sign}
\be
\Abar\bigg(\,\raisebox{-12pt}{\includegraphics[scale=1.8]{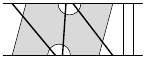}}\,\bigg)
=0
\ef.
\ee
\end{corollary}
\noindent
This corollary states that, for the operator $q_2$ defined in Section
\ref{ssec.q_op}, $\Abar(q_2 \t)=0$.

The proof of Proposition \ref{prop.854423343} is done again by applying
Proposition \ref{prop.ArfProp}.
We have in this case
\begin{align}
Q_{\s}&=
\begin{pmatrix}
0 &0 &1 &1 \\
1 &0 &0 &1 \\
1 &0 &1 &0 
\end{pmatrix}
\ef,
&
Q_{\t}&=
\begin{pmatrix}
1 &0 &0 &1 \\
1 &0 &1 &0
\end{pmatrix}
\ef,
\end{align}
and in particular
$Q_{\s} = \begin{pmatrixsm}1 & 1 \\ 1 & 0 \\ 0 & 1\end{pmatrixsm}
  Q_{\t}$.
Checking that the conditions of Proposition \ref{prop.ArfProp} are
met, with $K=2$ and $0$ in the two cases, is a straightforward
calculation. The resulting function $\Abar$ of $v$, as calculated for
the RHS, is
\be
\label{eq.6543564b}
1
+(-1)^{1 + \begin{pmatrixsm}1&0&0&1\end{pmatrixsm}\cdot v}
+(-1)^{1 + \begin{pmatrixsm}1&0&1&0\end{pmatrixsm}\cdot v}
+(-1)^{3 + \begin{pmatrixsm}0&0&1&1\end{pmatrixsm}\cdot v}
\ee
while for the LHS we have
\begin{multline}
1
+(-1)^{1 + \begin{pmatrixsm}0&0&1&1\end{pmatrixsm}\cdot v}
+(-1)^{1 + \begin{pmatrixsm}1&0&0&1\end{pmatrixsm}\cdot v}
+(-1)^{1 + \begin{pmatrixsm}1&0&1&0\end{pmatrixsm}\cdot v}
\\
+(-1)^{3 + \begin{pmatrixsm}1&0&1&0\end{pmatrixsm}\cdot v}
+(-1)^{3 + \begin{pmatrixsm}1&0&0&1\end{pmatrixsm}\cdot v}
+(-1)^{3 + \begin{pmatrixsm}0&0&1&1\end{pmatrixsm}\cdot v}
+(-1)^{6}
\end{multline}
which, after some simplifications, reduces to twice the expression in
(\ref{eq.6543564b}).

The calculation for function $A$ is analogous. For example,
%
for the LHS we have
\begin{multline}
\label{eq.6543564d}
1
+(-1)^{\begin{pmatrixsm}0&0&1&1\end{pmatrixsm}\cdot v}
+(-1)^{\begin{pmatrixsm}1&0&0&1\end{pmatrixsm}\cdot v}
+(-1)^{\begin{pmatrixsm}1&0&1&0\end{pmatrixsm}\cdot v}
\\
+(-1)^{1 + \begin{pmatrixsm}1&0&1&0\end{pmatrixsm}\cdot v}
+(-1)^{1 + \begin{pmatrixsm}1&0&0&1\end{pmatrixsm}\cdot v}
+(-1)^{1 + \begin{pmatrixsm}0&0&1&1\end{pmatrixsm}\cdot v}
+(-1)^{3}
\ef,
\end{multline}
an expression which
is identically zero.
\qed

\section{Surgery operators}
\label{sec.operators}

In this section we introduce and study the `surgery operators' which
have been outlined in Section \ref{ssec.TQSintro}, and whose geometric
interpretation has been described in Section
\ref{ssec.surgery_geo_combi}. These operators have a crucial role in
the proof of the classification theorem. In particular, we introduce
the notion:
\begin{definition}[Pullback function]
Let $f: X_n \to X_{n+1}$ be a function, $\sim$ be an equivalence
relation on $X_n$, with classes $Y_n$, and compatible with $f$, and
let $\tilde{f}$ be the function which leads to the commuting diagram
\be
\begin{CD}
X_n
@> f >>
X_{n+1}
\\
@VV{\cdot/\sim}V @VV{\cdot/\sim}V\\
Y_n
@> \tilde{f} >>
Y_{n+1}
\end{CD}
\label{eq.commuPullb1}
\ee
Call $Y(x)$ the function which associates to each $x \in X_n$ its
class in $Y_n$.  We say that $f$ is a \emph{pullback function} if, for
all $x' \in X_{n+1}$, and for all $y \in Y_{n}$ such that
$\tilde{f}(y)=Y(x)$, there exists a $x \in X_n$ such that $Y(x)=y$ and
$f(x)=x'$.
\end{definition}

\noindent
In other words, each diagram on the left can be completed to the
diagram on the right:
\begin{align}
\forall \quad
&
\begin{CD}
@.
x'
\\
@. @VV{\cdot/\sim}V\\
y
@> \tilde{f} >>
y'
\end{CD}
&
\exists\; x \textrm{~s.t.} \quad
&
\begin{CD}
x
@> f >>
x'
\\
@VV{\cdot/\sim}V @VV{\cdot/\sim}V\\
y
@> \tilde{f} >>
y'
\end{CD}
\label{eq.commuPullb2}
\end{align}
Our operators $\bar{T}$, $\bar{q}_1$ and $\bar{q}_2$ are functions on
the sets of (non-empty) classes of a give size, and with rank in a
certain range. In Theorems \ref{thm_T_surjectif}, \ref{thm.q1} and
\ref{thm.q2} we will establish that they are pullback functions
w.r.t.\ the equivalence relation given by the invariant.

\subsection{Operator $T$}
\label{ssec.T_op}

We define a first operator, in terms of the `square constructors'
defined in Section~\ref{sec.sqconstr}.

\begin{figure}[t!!]
\begin{center}
\includegraphics{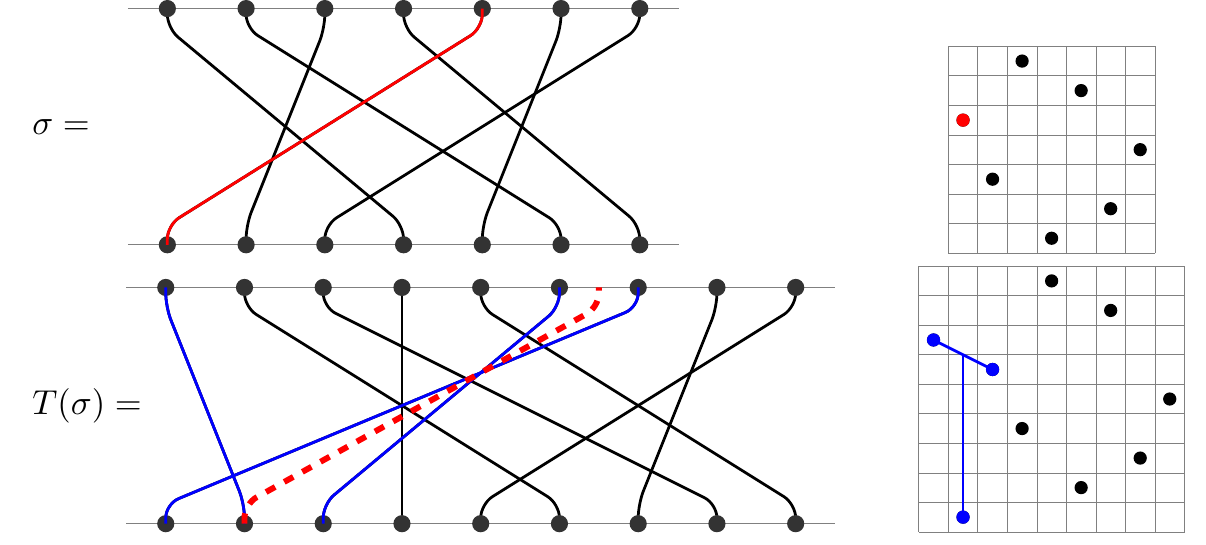}
\caption{\label{fig.T_example} An example of the $T$ operator}
\end{center}
\end{figure}

\begin{definition}
We define the $T$ operator $T: \kS_n \rightarrow \kS_{n+2}$ as 
$T= C_{\tau,0,1}^{{\rm col-}1}$ with 
$\tau=\left(\;\textrm{\raisebox{-2pt}{\makebox[0pt][l]{\raisebox{8pt}{\rule{11pt}{.5pt}}}$\bullet$\;\rule{.5pt}{11pt}\;}}\right)$.
\end{definition}
\noindent
Figure \ref{fig.T_example} illustrates this with an example.

\begin{remark}
\label{lem.T_homomorphism}
From Corollary \ref{cor.construtor}, we know that
if $\s \sim \s'$, then $T(\s)\sim T(\s')$.
\end{remark}
\noindent
As a consequence, we can define $\bar{T}$ as the map from classes at
size $n$ to classes at size $n+2$, such that $\bar{T}(C)=C'$ if there
exists $\s \in C$ with $T(\s) \in C'$.

\begin{figure}[b!!]
\begin{center}
\includegraphics{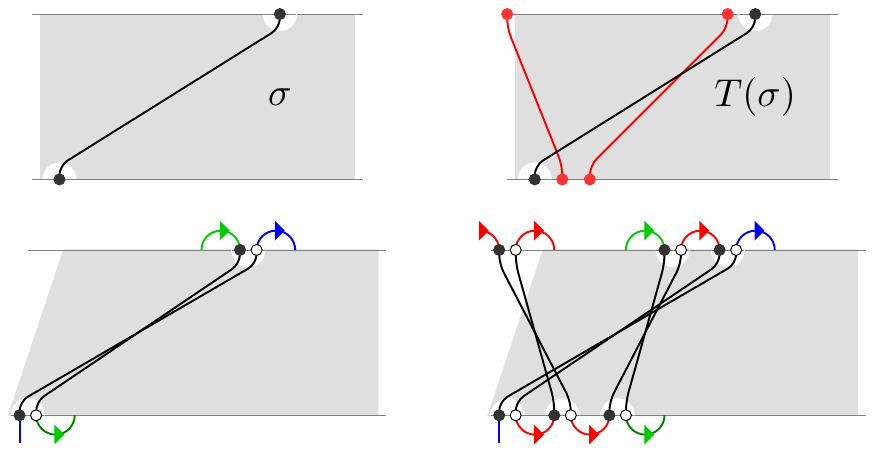}
\caption{\label{fig.T_lambda_r_chg} Left: the permutation $\s$. Right:
  the permutation $T \s$.
  It is clear that the only change in the cycle invariant is the
  addition of two to the rank length. On the bottom, the edges are
  drawn `doubled', as in the diagram construction of the cycle
  invariant.}
\end{center}
\end{figure}

\begin{lemma}\label{lem_T_operator}
Let $C$ be a class with invariant $(\lambda,r,s)$. Then $\bar{T}(C)$
has invariant $(\lambda,r+2,s)$.
\end{lemma}
\proof In light of Corollary \ref{cor.Tsign}, we know that the sign
invariant does not change.  For what concerns the cycle invariant,
Figure \ref{fig.T_lambda_r_chg} illustrates our claim, which holds
more generally also for the constructor
$C_{\tau,0,1}^{{\rm col-}\ell}$, with arbitrary $\ell$.
\qed

Finally, as outlined in Section \ref{ssec.TQSintro}, we need to
establish the following crucial property:
\begin{theorem}[Pullback of $T$]
\label{thm_T_surjectif}
For every non-exceptional class $C$ at size $n+2$ with invariant
$(\lambda,r,s)$, and $r>2$, there exists an irreducible $\s$ at size
$n$ such that $T(\s) \in C$. Or equivalently, $\bar{T}$ is surjective
onto the set of (primitive irreducible) non-exceptional classes with
rank larger than~2, and more specifically it is a pullback function.
\end{theorem}
\noindent
The proof of this theorem takes the largest part of this subsection.
Before starting this proof, we need to set some notation for dealing
with reduced permutations, as explained in Section \ref{ssec.reddyn},
which are of arbitrary size, but with a finite number of `black' edges
(and an arbitrarily large number of `gray' ones).

In the in-line notation for permutations -- i.e., the string 
$(\s(1), \ldots, \s(n))$ -- we will use the notation $*$ for a sublist
of arbitrary length, possibly zero, and the symbol $n$ to denote the
maximal element, while indices such as $h$, $k$, \ldots\ will be used
for generic elements. For example, the permutation
$\s=(1,4,2,5,6,8,9,3,7)$ is one of those of the form
$(1,*,k,k+1,*,n,*)$, because it starts with $1$, the consecutive $5$
and $6$ are candidate for $k$ and $k+1$, and $n=9$ comes after
them. Note that, when wild characters like $k$ are used, a permutation
can be of a given form in more than one way (e.g., this would have
been the case for $\s=(1,4,2,5,6,7,9,3,8)$, as both 5 and 6 are
possible choices for $k$).
In some specially taylored expressions this never happens. This would
have been the case, for example, for the pattern $(1,*,2,k,k+1,*,n,*)$
as the explicit $2$ forces $k$ to admit at most one realisation.

\begin{notation}
We want to represent graphically permutations in reduced dynamics,
explicitating only the pattern $\tau$ of black edges, and, of the
auxiliary structure $\hat{c}$, which pairs of consecutive black-edge
endpoints have no gray endpoints in between them. In a matrix diagram,
we use red lines between consecutive rows and columns\footnote{And/or
  before the first/after the last row and column.} to denote the
certified absence of gray points (in absence of a red line, there may
or may not be gray points in between).  
\end{notation}
\noindent
Such a data structure, when ``not too big'', is also conveniently
encoded by an in-line expression of the form mentioned above.
For example, the in-line patterns $(1,2,*,3,k,*,(k+1),*)$ 
and $(k,1,(k-1),*,n,*,2)$ are represented
graphically as
\begin{align}
&
(1,2,*,3,k,*,(k+1),*)\;:
&& 
(k,1,(k-1),*,n,*,2)\;:
\\
& 
\begin{tikzpicture}[scale=0.4]
\permutationtwo{1,2,3,4,5}{1,2,3,5}{1,2,4};
\end{tikzpicture}
&&
\begin{tikzpicture}[scale=0.4]
\permutationtwo{4,1,3,5,2}{1,2,4,6}{1,2,3,6};
\end{tikzpicture}
\end{align}
We will use extensively this graphical convention in the following
proof of Theorem~\ref{thm_T_surjectif}.

In particular, in this section, all our patterns are of the form
$(1,2,\ldots)$, so we have two red lines at the bottom, and two at the
left. This implies that we know the position of the $L$-pivot, which
is the $\s(1)=1$ entry. If (and only if) we do not have a red line on
top, we do not know the position of the $R$-pivot. Nonetheless, up to
an operator $L^{a}$ for some $a$, we can move all the gray entries
above the top-most black entry, and in turns trade the horizontal red
line at position $2$ for a red line on top.
For example,
\begin{center} 
\begin{tikzpicture}[scale=0.4]
\permutationtwo{1,2,4,5,3}{1,2}{1,2};
\end{tikzpicture}
\quad
\raisebox{25pt}{$\sim$}
\quad
\begin{tikzpicture}[scale=0.4]
\permutationtwo{1,2,4,5,3}{1,6}{1,2};
\end{tikzpicture}
\end{center}
From this moment on, in the reduced dynamics we will always have a
red line on the left and top side of the diagram, and a full control
on the position of both pivots. I.e., in the reduced dynamics the
pivots are always both `black'.

Our aim is to prove that, for a collection of patterns covering all
possible cases of the theorem (of standard non-exceptional primitive
irreducible configurations of rank at least 3), in the reduced
dynamics we can reach a pattern that certifies the presence of a
$T$-structure, which, when removed, leaves with an irreducible
configuration.
The presence of a $T$-structure is easily verified (in our drawings,
we use blue construction lines to evidentiate it). Irreducibility
after $T^{-1}$ is certified by a zig-zag path (that we represent in
green), which reaches red lines at its endpoints (for certifying the
absence of an irreducible block constituted exclusively of gray
entries).  For primitivity, in principle there is nothing to check, as
the number of descents is preserved both by the dynamics, and by
removing the $T$.  Nonetheless, for reasons explained below, we will
keep track explicitly of possible descents involving black edges.  A
typical ``winning outcome'' of the dynamics, with the associated
construction lines, is as follows:
\[
\begin{tikzpicture}[scale=0.4]
\draw [blue,very thick] (1.5,6.5)--(3.5,5.5);
\draw [blue,very thick] (2.5,6)--(2.5,1.5);
\draw [green!50!black!100,very thick] (3.5,5.5)--(6.5,7.5)--(5.5,3.5)--(8.5,4.5)--(7.5,2.5);
\permutationtwo{6,1,5,8,3,7,2,4}{1,2,6,9}{1,2,3};
\end{tikzpicture}
\]
We call such a pattern a \emph{$T$-structure certificate}.

The irreducibility of the outcome is not for granted. It is not
uncommon that shorter and simpler candidate sequences have to be
rejected for lack of this property.  An example which is instructive
in retrospective is the pattern associated to case 4.1 in Table
\ref{tab.T_tab_all_case}. Na\"ively, one could have guessed that a
good choice of sequence is $\bar{R}\bar{L}^3$, just as for case
2. This indeed produces a $T$-structure, however the pattern resulting
after removing the structure is not irreducible (we have a block of
size 4, followed by one of size 2). This forced us to search for the
longer sequence appearing in the table.


Patterns are stable by inclusion. If we add points to a pattern $\tau$
in correspondence of crossing of non-red lines, producing a larger
pattern $\tau'$, and $\tau$ could reach a certificate as above, also
$\tau'$ can, this because the relevant features of both the
$T$-structure and the zig-zag path are related to red lines, which do
not interfere with the insertion.

The relation in the other direction is slightly more involved: if
$\tau'$ can reach a certificate, and a point $p$ in $\tau'$ is neither
used in the $T$-structure certificate, nor used as a pivot in the
sequence, then $\tau = \tau' \setminus p$ can reach a certificate.

We will use this relation mostly in the second form.  In some cases,
when we get a pattern $\tau'$ from our case analysis, instead of
producing a related certificate, we will make the ansatz that a
certificate exists also under a certain reduction, and produce the
associated (generally smaller) certificate. These reductions will also
allow to merge together a number of branches in the case analysis.

Now we can state the following
\begin{proposition}
\label{prop.fortab.T_tab_all_case}
Table \ref{tab.T_tab_all_case} lists triples $(\tau, \tau', S)$ such
that $S \tau = \tau'$.  (The ``names'' for the triples are for future
reference).

Gray bullets denote arbitrary blocks, possibly empty, while gray
bullets with a black bullet inside denote arbitrary non-empty blocks.

All $\tau'$ are $T$-structure certificates, with the exception of
pattern 7, which has a zig-zag path if and only if at least one of the
two gray bullets is a non-empty block.

All sequences $S$ have alternation length at most 6, all zig-zag paths
in the certificates have length at most 5.
\end{proposition}
\noindent
Curiously, as a corollary (not useful at our purposes), we get that
patterns 3.1 and 4.1 are connected, as well as 5.1, 5.2 and 6.1, fact
that was not obvious \emph{a priori}.

We are now ready for proving our theorem.

\begin{table}[p!]
\noindent
\begin{tabular}{c}
\rule{0pt}{12pt}%
A: $\bar{R}$
\\
\begin{tikzpicture}[scale=0.38]
\permutationtwo{1,2,3,4}{1,2,3,5}{1,2,4,5};
\end{tikzpicture}
\\
\begin{tikzpicture}[scale=0.38]
\draw [blue,very thick] (1.5,3.5)--(3.5,2.5);
\draw [blue,very thick] (2.5,3)--(2.5,1.5);
\draw (4.5,4.5) [fill,green!50!black!100] circle (.35);
\draw (4.5,4.5) [fill,white] circle (.25);
\permutationtwo{3,1,2,4}{1,2,3,5}{1,2,3,5};
\end{tikzpicture}
\\
\end{tabular}
\!\!\!
\begin{tabular}{c}
\rule{0pt}{12pt}%
B: $\bar{R}L$
\\
\begin{tikzpicture}[scale=0.38]
\permutationtwo{1,2,5,3,4}{1,3,5,6}{1,2,5,6};
\end{tikzpicture}
\\
\begin{tikzpicture}[scale=0.38]
\draw [blue,very thick] (1.5,4.5)--(3.5,3.5);
\draw [blue,very thick] (2.5,4)--(2.5,1.5);
\draw (5.5,5.5) [fill,green!50!black!100] circle (.35);
\draw (5.5,5.5) [fill,white] circle (.25);
\permutationtwo{4,1,3,2,5}{1,2,4,6}{1,2,3,6};
\end{tikzpicture}
\\
\end{tabular}
\!\!\!
\begin{tabular}{c}
\rule{0pt}{12pt}%
C: $\bar{R}L$
\\
\begin{tikzpicture}[scale=0.38]
\permutationtwo{1,2,3,5,6,4}{1,3,6,7}{1,2,4,7};
\end{tikzpicture}
\\
\begin{tikzpicture}[scale=0.38]
\draw [blue,very thick] (1.5,4.5)--(3.5,3.5);
\draw [blue,very thick] (2.5,4)--(2.5,1.5);
\draw [green!50!black!100,very thick] (3.5,3.5)--(6.5,5.5);
\permutationtwo{4,1,3,6,2,5}{1,2,4,7}{1,2,3,7};
\end{tikzpicture}
\\
\end{tabular}
\!\!\!
\begin{tabular}{c}
\rule{0pt}{12pt}%
1: $\bar{R}LR\bar{L}\bar{R}$
\\
\begin{tikzpicture}[scale=0.38]
\draw (7.5,6.5) [fill,black!40] circle (.35);
\permutationtwo{1,2,3,4,7,5,6}{1,3,5,8}{1,2,4};
\end{tikzpicture}
\\
\begin{tikzpicture}[scale=0.38]
\draw [blue,very thick] (1.5,4.5)--(3.5,3.5);
\draw [blue,very thick] (2.5,4)--(2.5,1.5);
\draw [green!50!black!100,very thick] (3.5,3.5)--(7.5,6.5)--(6.5,2.5);
\draw (7.5,6.5) [fill,black!40] circle (.35);
\permutationtwo{4,1,3,7,5,2,6}{1,2,4,8}{1,2,3};
\end{tikzpicture}
\\
\end{tabular}
\!\!\!
\begin{tabular}{c}
\rule{0pt}{12pt}%
2: $\bar{R}\bar{L}^3$
\\
\begin{tikzpicture}[scale=0.38]
\draw (7.5,7.5) [fill,black!40] circle (.35);
\permutationtwo{1,2,6,3,4,5,7}{1,3,5,8}{1,2,5};
\end{tikzpicture}
\\
\begin{tikzpicture}[scale=0.38]
\draw [blue,very thick] (1.5,6.5)--(3.5,5.5);
\draw [blue,very thick] (2.5,6)--(2.5,1.5);
\draw [green!50!black!100,very thick] (5.5,7.5)--(4.5,3.5)--(7.5,4.5)--(6.5,2.5);
\draw (7.5,4.5) [fill,black!40] circle (.35);
\permutationtwo{6,1,5,3,7,2,4}{1,2,6,8}{1,2,3};
\end{tikzpicture}
\\
\end{tabular}

\noindent
\hspace*{-4mm}
\begin{tabular}{c}
\rule{0pt}{12pt}%
3.1: $\bar{R}LR\bar{L}^2\bar{R}\bar{L}^2$
\\
\begin{tikzpicture}[scale=0.38]
\draw (8.5,8.5) [fill,black!40] circle (.35);
\permutationtwo{1,2,4,5,3,6,7,8}{1,3,7,9}{1,2,6};
\end{tikzpicture}
\\
\begin{tikzpicture}[scale=0.38]
\draw [blue,very thick] (1.5,7.5)--(3.5,6.5);
\draw [blue,very thick] (2.5,7)--(2.5,1.5);
\draw [green!50!black!100,very thick](5.5,8.5)--(4.5,3.5)--(8.5,5.5)--(7.5,2.5);
\draw (8.5,5.5) [fill,black!40] circle (.35);
\permutationtwo{7,1,6,3,8,4,2,5}{1,2,7,9}{1,2,3};
\end{tikzpicture}
\\
\end{tabular}
\!\!\!
\begin{tabular}{c}
\rule{0pt}{12pt}%
3.2: $\bar{R}L^2R\bar{L}\bar{R}\bar{L}^6$
\\
\begin{tikzpicture}[scale=0.38]
\draw (8.5,8.5) [fill,black!40] circle (.35);
\permutationtwo{1,2,5,3,6,4,7,8}{1,3,7,9}{1,2,5};
\end{tikzpicture}
\\
\begin{tikzpicture}[scale=0.38]
\draw [blue,very thick] (1.5,6.5)--(3.5,5.5);
\draw [blue,very thick] (2.5,6)--(2.5,1.5);
\draw [green!50!black!100,very thick] (6.5,8.5)--(5.5,3.5)--(8.5,4.5)--(7.5,2.5);
\draw (8.5,4.5) [fill,black!40] circle (.35);
\permutationtwo{6,1,5,7,3,8,2,4}{1,2,6,9}{1,2,3};
\end{tikzpicture}
\\
\end{tabular}
\!\!\!
\begin{tabular}{c}
\rule{0pt}{12pt}%
3.3: $\bar{R}LR\bar{L}^2\bar{R}$
\\
\begin{tikzpicture}[scale=0.38]
\draw (8.5,7.5) [fill,black!40] circle (.35);
\permutationtwo{1,2,4,8,3,5,6,7}{1,3,6,9}{1,2,6};
\end{tikzpicture}
\\
\begin{tikzpicture}[scale=0.38]
\draw [blue,very thick] (1.5,4.5)--(3.5,3.5);
\draw [blue,very thick] (2.5,4)--(2.5,1.5);
\draw [green!50!black!100,very thick] (3.5,3.5)--(8.5,6.5)--(7.5,2.5);
\draw (8.5,6.5) [fill,black!40] circle (.35);
\permutationtwo{4,1,3,7,8,5,2,6}{1,2,4,9}{1,2,3};
\end{tikzpicture}
\\
\end{tabular}
\!\!\!
\begin{tabular}{c}
\rule{0pt}{12pt}%
4.1: $\bar{R}LR\bar{L}^3RL$
\\
\begin{tikzpicture}[scale=0.38]
\draw (8.5,6.5) [fill,black!40] circle (.35);
\permutationtwo{1,2,7,8,3,4,5,6}{1,3,5,9}{1,2,6};
\end{tikzpicture}
\\
\begin{tikzpicture}[scale=0.38]
\draw [blue,very thick] (1.5,7.5)--(3.5,6.5);
\draw [blue,very thick] (2.5,7)--(2.5,1.5);
\draw [green!50!black!100,very thick] (5.5,8.5)--(4.5,3.5)--(8.5,5.5)--(7.5,2.5);
\draw (8.5,5.5) [fill,black!40] circle (.35);
\permutationtwo{7,1,6,3,8,4,2,5}{1,2,7,9}{1,2,3};
\end{tikzpicture}
\\
\end{tabular}

\noindent
\hspace*{-4mm}
\begin{tabular}{c}
\rule{0pt}{12pt}%
5.1: $\bar{R}L^2RL\bar{R}\bar{L}^6$
\\
\begin{tikzpicture}[scale=0.38]
\draw (8.5,8.5) [fill,black!40] circle (.35);
\permutationtwo{1,2,3,6,4,5,7,8}{1,3,7,9}{1,2,4};
\end{tikzpicture}
\\
\begin{tikzpicture}[scale=0.38]
\draw [blue,very thick] (1.5,6.5)--(3.5,5.5);
\draw [blue,very thick] (2.5,6)--(2.5,1.5);
\draw [green!50!black!100,very thick] (3.5,5.5)--(6.5,7.5)--(5.5,3.5)--(8.5,4.5)--(7.5,2.5);
\draw (8.5,4.5) [fill,black!40] circle (.35);
\permutationtwo{6,1,5,8,3,7,2,4}{1,2,6,9}{1,2,3};
\end{tikzpicture}
\\
\end{tabular}
\!\!\!
\begin{tabular}{c}
\rule{0pt}{12pt}%
5.2: $\bar{R}\bar{L}^4$
\\
\begin{tikzpicture}[scale=0.38]
\draw (8.5,8.5) [fill,black!40] circle (.35);
\permutationtwo{1,2,3,5,7,4,6,8}{1,3,6,9}{1,2,4};
\end{tikzpicture}
\\
\begin{tikzpicture}[scale=0.38]
\draw [blue,very thick] (1.5,6.5)--(3.5,5.5);
\draw [blue,very thick] (2.5,6)--(2.5,1.5);
\draw [green!50!black!100,very thick] (3.5,5.5)--(6.5,7.5)--(5.5,3.5)--(8.5,4.5)--(7.5,2.5);
\draw (8.5,4.5) [fill,black!40] circle (.35);
\permutationtwo{6,1,5,8,3,7,2,4}{1,2,6,9}{1,2,3};
\end{tikzpicture}
\\
\end{tabular}
\!\!\!
\begin{tabular}{c}
\rule{0pt}{12pt}%
6.1: $\bar{R}LRL\bar{R}\bar{L}^5$
\\
\begin{tikzpicture}[scale=0.38]
\draw (8.5,8.5) [fill,black!40] circle (.35);
\permutationtwo{1,2,3,4,6,7,5,8}{1,3,5,9}{1,2,4};
\end{tikzpicture}
\\
\begin{tikzpicture}[scale=0.38]
\draw [blue,very thick] (1.5,6.5)--(3.5,5.5);
\draw [blue,very thick] (2.5,6)--(2.5,1.5);
\draw [green!50!black!100,very thick] (3.5,5.5)--(6.5,7.5)--(5.5,3.5)--(8.5,4.5)--(7.5,2.5);
\draw (8.5,4.5) [fill,black!40] circle (.35);
\permutationtwo{6,1,5,8,3,7,2,4}{1,2,6,9}{1,2,3};
\end{tikzpicture}
\\
\end{tabular}
\!\!\!
\begin{tabular}{c}
\rule{0pt}{12pt}%
7:
$\bar{R}L^2\bar{R}^2\bar{L}^4R$
\\
\begin{tikzpicture}[scale=0.38]
\permutationtwo{1,2,3,4,5,8,9,6,7}{1,3,5,10}{1,2,4};
\draw (6.5,8.5) [fill,black!40] circle (.35);
\draw (9.5,7.5) [fill,black!40] circle (.35);
\end{tikzpicture}
\\
\begin{tikzpicture}[scale=0.38]
\draw [blue,very thick] (1.5,8.5)--(3.5,7.5);
\draw [blue,very thick] (2.5,8)--(2.5,1.5);
\draw [green!50!black!100,very thick](5.5,9.5)--(4.5,4.5)--(8.5,6)--(6.5,2.5);
\draw [green!50!black!100,very thick](7.5,6.5)--(9.5,5.5);
\permutationtwo{8,1,7,4,9,2,6,3,5}{1,2,8,10}{1,2,3};
\draw (7.5,6.5) [fill,black!40] circle (.35);
\draw (9.5,5.5) [fill,black!40] circle (.35);
\end{tikzpicture}
\\
\end{tabular}
\\
\phantom{empty line}
\caption{\label{tab.T_tab_all_case}Triples as is
  Proposition~\ref{prop.fortab.T_tab_all_case}.}
\end{table}

\bigskip

\pfof{Theorem \ref{thm_T_surjectif}}
Our aim is to prove that, for each irreducible primitive
non-exceptional class $C$ of rank at least 3, there exists a
configuration $\s'$ in the image of the operator $T$, such that its
preimage is also primitive and irreducible.

The proof will go as follows: we describe how to break the problem
into a finite number of cases, ameaneable to reduced dynamics, and in
fact each of a form considered in Proposition
\ref{prop.fortab.T_tab_all_case}.  To this end it is convenient to
assume that, to start with, we have some reference configuration $\s$
which is standard (in the sense of Section \ref{ssec.std_perm}).  In
particular, its in-line expression starts with $(1,2,\ldots)$.  This
is legitimate because, from the results of Section \ref{ssec.std_perm}
(see Lemma \ref{lem.stdPerms}), we know that each irreducible class has at
least one standard family.  As a matter of fact (and for what we know
from Appendix \ref{sec.excp_class}) the only exceptional configuration of this form of
rank at least 3 is $\id_n$. Our case analysis is tree-like, and we
will see that this case emerge from \emph{a single branch} of our
tree, though this was not granted in advance.  As it was obvious that
it had to come from at least one branch, the relevant feature is that
the tree is finite, and in particular that there is a finite number of
branches leading to $\id_n$ configurations (instead, e.g., of one
branch per value of $n$), which would have made the reasoning more
cumbersome, if not impractical.

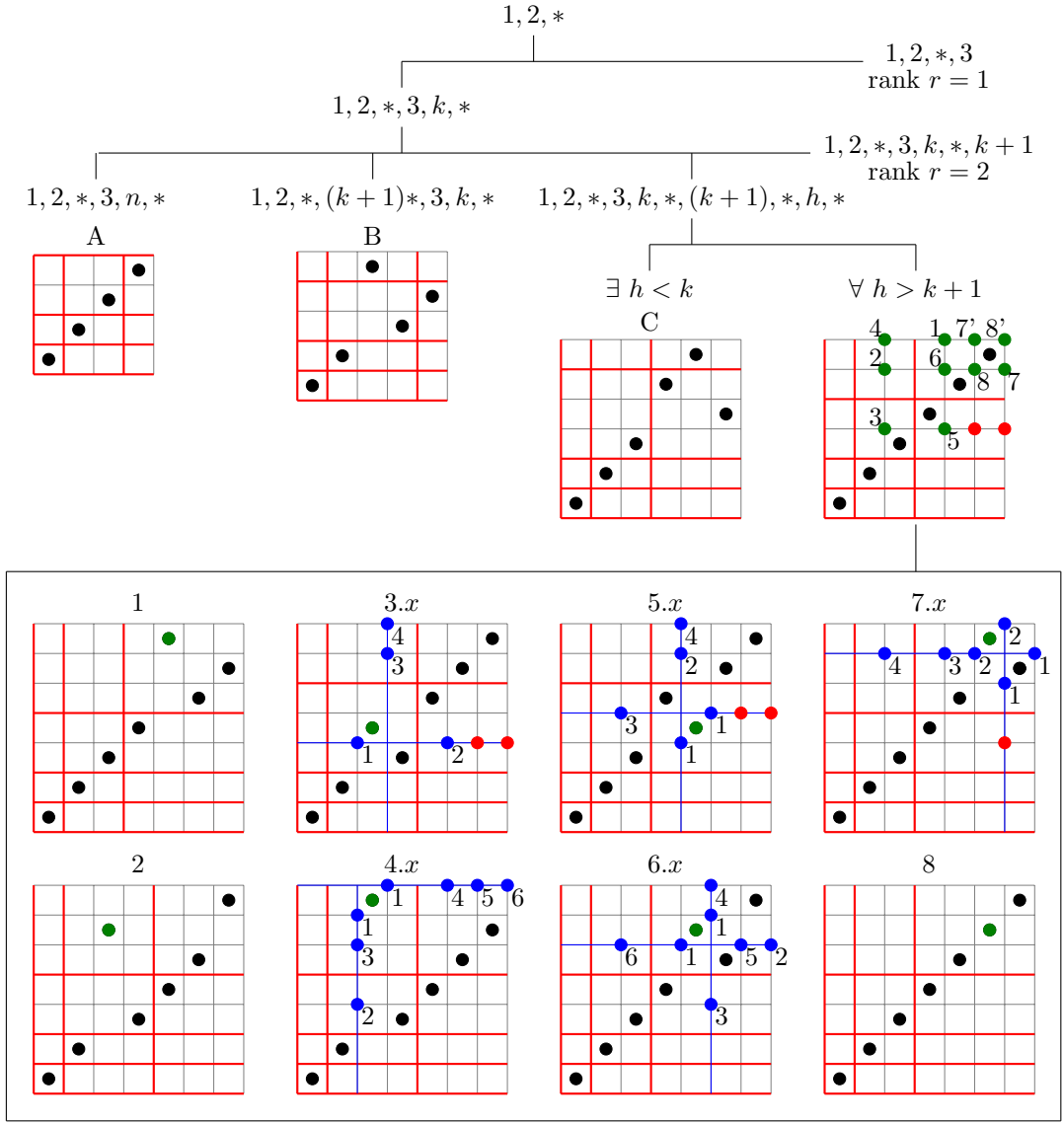
\begin{figure}[tb!]
\[
\setlength{\unitlength}{100pt}
\begin{picture}(3.8,4.2)
\put(-.1,-.1){\line(1,0){4}}
\put(-.1,-.1){\line(0,1){2.1}}
\put(3.9,2.){\line(-1,0){4}}
\put(3.9,2.){\line(0,-1){2.1}}
\put(3.35,2.){\line(0,1){.18}}

\put(0.4,1.85){\gostrC{1}}
\put(0.4,0.85){\gostrC{2}}
\put(1.4,1.85){\gostrC{3$.x$}}
\put(1.4,0.85){\gostrC{4$.x$}}
\put(2.4,1.85){\gostrC{5$.x$}}
\put(2.4,0.85){\gostrC{6$.x$}}
\put(3.4,1.85){\gostrC{7$.x$}}
\put(3.4,0.85){\gostrC{8}}

\put(1.9,4.1){\gostrC{$1,2,*$}}

\put(3.4,3.95){\gostrC{$1,2,*,3$}}
\put(3.4,3.85){\gostrC{rank $r=1$}}

\put(3.4,3.6){\gostrC{$1,2,*,3,k,*,k+1$}}
\put(3.4,3.5){\gostrC{rank $r=2$}}

\put(1.4,3.75){\gostrC{$1,2,*,3,k,*$}}
\put(0.25,3.6){\line(1,0){2.7}}
\put(1.4,3.6){\line(0,1){.1}}
\put(1.4,3.85){\line(0,1){.1}}
\put(1.4,3.95){\line(1,0){1.75}}
\put(1.9,3.95){\line(0,1){.1}}

\put(0.24,3.4){\gostrC{$1,2,*,3,n,*$}}
\put(0.24,3.5){\line(0,1){.1}}
\put(0.24,3.25){\gostrC{A}}
\put(0,2.75){%
\begin{tikzpicture}[scale=0.4]
\permutationtwo{1,2,3,4}{1,2,3,5}{1,2,4};
\end{tikzpicture}
}

\put(1.29,3.4){\gostrC{$1,2,*,(k+1)*,3,k,*$}}
\put(1.29,3.5){\line(0,1){.1}}
\put(1.29,3.25){\gostrC{B}}
\put(1,2.65){%
\begin{tikzpicture}[scale=0.4]
\permutationtwo{1,2,5,3,4}{1,2,3,5}{1,2,5};
\end{tikzpicture}
}

\put(2.5,3.4){\gostrC{$1,2,*,3,k,*,(k+1),*,h,*$}}
\put(2.5,3.5){\line(0,1){.1}}
\put(2.5,3.25){\line(0,1){.1}}
\put(2.34,3.25){\line(1,0){1.01}}
\put(2.34,3.15){\line(0,1){.1}}
\put(3.35,3.15){\line(0,1){.1}}

\put(2.34,3.05){\gostrC{$\exists\ h<k$}}
\put(2.34,2.92){\gostrC{C}}
\put(2,2.2){%
\begin{tikzpicture}[scale=0.4]
\permutationtwo{1,2,3,5,6,4}{1,2,3,6}{1,2,4};
\end{tikzpicture}
}

\put(3.35,3.05){\gostrC{$\forall\ h>k+1$}}
\put(3,2.2){%
\begin{tikzpicture}[scale=0.4]
\permutationtwo{1,2,3,4,5,6}{1,2,3,5}{1,2,4};
\draw (6,4) [fill,red] circle (.2);
\draw (7,4) [fill,red] circle (.2);
\draw (5,4) [fill,green!50!black!100] circle (.2);
\draw (5,7) [fill,green!50!black!100] circle (.2);
\draw (5,6) [fill,green!50!black!100] circle (.2);
\draw (3,4) [fill,green!50!black!100] circle (.2);
\draw (3,6) [fill,green!50!black!100] circle (.2);
\draw (3,7) [fill,green!50!black!100] circle (.2);
\draw (6,7) [fill,green!50!black!100] circle (.2);
\draw (6,6) [fill,green!50!black!100] circle (.2);
\draw (7,7) [fill,green!50!black!100] circle (.2);
\draw (7,6) [fill,green!50!black!100] circle (.2);
\node at (4.7,7.4) {1};
\node at (2.7,6.4) {2};
\node at (2.7,4.4) {3};
\node at (2.7,7.4) {4};
\node at (4.7,6.4) {6};
\node at (5.7,7.4) {7'};
\node at (6.7,7.4) {8'};
\node at (4.7+.6,4.4-.8) {5};
\node at (6.7+.6,6.4-.8) {7};
\node at (5.7+.6,6.4-.8) {8};
\end{tikzpicture}
}
\put(0,1){%
\begin{tikzpicture}[scale=0.4]
\permutationtwo{1,2,3,4,7,5,6}{1,2,3,5}{1,2,4};
\draw (5+.5,7+.5) [fill,green!50!black!100] circle (.2);
\end{tikzpicture}
}
\put(0,0){%
\begin{tikzpicture}[scale=0.4]
\permutationtwo{1,2,6,3,4,5,7}{1,2,3,5}{1,2,5};
\draw (3+.5,6+.5) [fill,green!50!black!100] circle (.2);
\end{tikzpicture}
}
\put(1,1){%
\begin{tikzpicture}[scale=0.4]
\permutationtwo{1,2,4,3,5,6,7}{1,2,3,6}{1,2,5};
\draw (3+.5,4+.5) [fill,green!50!black!100] circle (.2);
\draw [blue] (4,1)--(4,8);
\draw [blue] (1,4)--(8,4);
\draw (3,4) [fill,blue] circle (.2);
\draw (6,4) [fill,blue] circle (.2);
\draw (4,7) [fill,blue] circle (.2);
\draw (4,8) [fill,blue] circle (.2);
\draw (7,4) [fill,red] circle (.2);
\draw (8,4) [fill,red] circle (.2);
\node at (3.35,3.6) {1};
\node at (6.35,3.6) {2};
\node at (4.35,6.6) {3};
\node at (4.35,7.6) {4};
\end{tikzpicture}
}
\put(1,0){%
\begin{tikzpicture}[scale=0.4]
\permutationtwo{1,2,7,3,4,5,6}{1,2,3,5}{1,2,5};
\draw (3+.5,7+.5) [fill,green!50!black!100] circle (.2);
\draw [blue] (3,1)--(3,8);
\draw [blue] (1,8)--(8,8);
\draw (4,8) [fill,blue] circle (.2);
\draw (3,7) [fill,blue] circle (.2);
\draw (6,8) [fill,blue] circle (.2);
\draw (7,8) [fill,blue] circle (.2);
\draw (3,4) [fill,blue] circle (.2);
\draw (3,6) [fill,blue] circle (.2);
\draw (8,8) [fill,blue] circle (.2);
\node at (3.35,3.6) {2};
\node at (3.35,5.6) {3};
\node at (3.35,6.6) {1};
\node at (4.35,7.6) {1};
\node at (6.35,7.6) {4};
\node at (7.35,7.6) {5};
\node at (8.35,7.6) {6};
\end{tikzpicture}
}
\put(2,1){%
\begin{tikzpicture}[scale=0.4]
\permutationtwo{1,2,3,5,4,6,7}{1,2,3,6}{1,2,4};
\draw (5+.5,4+.5) [fill,green!50!black!100] circle (.2);
\draw [blue] (5,1)--(5,8);
\draw [blue] (1,5)--(8,5);
\draw (5,4) [fill,blue] circle (.2);
\draw (5,7) [fill,blue] circle (.2);
\draw (5,8) [fill,blue] circle (.2);
\draw (3,5) [fill,blue] circle (.2);
\draw (6,5) [fill,blue] circle (.2);
\draw (7,5) [fill,red] circle (.2);
\draw (8,5) [fill,red] circle (.2);
\node at (3.35,4.6) {3};
\node at (6.35,4.6) {1};
\node at (5.35,3.6) {1};
\node at (5.35,6.6) {2};
\node at (5.35,7.6) {4};
\end{tikzpicture}
}
\put(2,0){%
\begin{tikzpicture}[scale=0.4]
\permutationtwo{1,2,3,4,6,5,7}{1,2,3,5}{1,2,4};
\draw (5+.5,6+.5) [fill,green!50!black!100] circle (.2);
\draw [blue] (6,1)--(6,8);
\draw [blue] (1,6)--(8,6);
\draw (6,4) [fill,blue] circle (.2);
\draw (6,7) [fill,blue] circle (.2);
\draw (6,8) [fill,blue] circle (.2);
\draw (3,6) [fill,blue] circle (.2);
\draw (5,6) [fill,blue] circle (.2);
\draw (7,6) [fill,blue] circle (.2);
\draw (8,6) [fill,blue] circle (.2);
\node at (3.35,5.6) {6};
\node at (5.35,5.6) {1};
\node at (7.35,5.6) {5};
\node at (8.35,5.6) {2};
\node at (6.35,7.6) {4};
\node at (6.35,6.6) {1};
\node at (6.35,3.6) {3};
\end{tikzpicture}
}
\put(3,1){%
\begin{tikzpicture}[scale=0.4]
\permutationtwo{1,2,3,4,5,7,6}{1,2,3,5}{1,2,4};
\draw (6+.5,7+.5) [fill,green!50!black!100] circle (.2);
\draw (7+.5,6+.5) [fill,black] circle (.2);
\draw [blue] (7,1)--(7,8);
\draw [blue] (1,7)--(8,7);
\draw (7,4) [fill,red] circle (.2);
\draw (7,6) [fill,blue] circle (.2);
\draw (7,8) [fill,blue] circle (.2);
\draw (3,7) [fill,blue] circle (.2);
\draw (5,7) [fill,blue] circle (.2);
\draw (6,7) [fill,blue] circle (.2);
\draw (8,7) [fill,blue] circle (.2);
\node at (3.35,6.6) {4};
\node at (5.35,6.6) {3};
\node at (6.35,6.6) {2};
\node at (8.35,6.6) {1};
\node at (7.35,7.6) {2};
\node at (7.35,5.6) {1};
\end{tikzpicture}
}
\put(3,0){%
\begin{tikzpicture}[scale=0.4]
\permutationtwo{1,2,3,4,5,6,7}{1,2,3,5}{1,2,4};
\draw (6+.5,6+.5) [fill,green!50!black!100] circle (.2);
\end{tikzpicture}
}
\end{picture}
\]
\caption{\label{figure.Tthm}Decomposition of cases in the proof of
  Theorem~\ref{thm_T_surjectif}. Numbers on blue bullets correspond to
  different values of $x$ in case label.}
\end{figure}

Figure \ref{figure.Tthm} summarises the main steps of the case
decomposition,
which we also describe in words in the following paragraphs.

As we said, by Lemma \ref{lem.stdPerms}, we can suppose that $\s$ is a
standard permutation with $\s(2)=2$, i.e., in in-line notation, with
the pattern $\s=(1,2,*)$. We now investigate the possible preimages of
$3$.
\begin{enumerate}
 \item If 3 is at the end i.e. $\s=(1,2,*,3)$, then $\s$ has rank 1 and the theorem does not apply.
 \item Otherwise $\s=1,2,*,3,k,*$ for some $k>3$.
 \begin{enumerate}
  \item If $k=n$ i.e. $\s=(1,2,*,3,n,*)$, we are in case `A' of
    Figure~\ref{figure.Tthm}, which is solved by the homonymous case in
    Proposition~\ref{prop.fortab.T_tab_all_case}.
  \item Otherwise $k<n$ and we can investigate the possible positions
    of the image of $k+1$.
  \begin{enumerate}
   \item If $\s=(1,2,*,(k+1),*,3,k,*)$, we are in case `B', which yet
     again is solved by the homonymous case in
     Proposition~\ref{prop.fortab.T_tab_all_case}.
   \item If $\s=(1,2,*,3,k,*,(k+1))$ then $\s$ has rank 2 and the theorem does not apply.
   \item Otherwise $\s=(1,2,*,3,k,*,(k+1),*,h,*)$ for one or more $h$:
   \begin{enumerate}
    \item If there exists such a $h$ with $3<h<k$, we are in case `C'.
    \item Otherwise, we shall analyse with care the case
      $\s=(1,2,*,3,k,*,(k+1),*,h,*)$ with all entries after $k+1$
      being larger than $k+1$.
   \end{enumerate}
  \end{enumerate}
 \end{enumerate}
\end{enumerate}
We have reached a case that requires a better control on the extra
data structure $\hat{c}$ of the reduced dynamics. At all intersections
of non-red lines, we may or may not have a non-empty block of
points. We will take some order on these points, and continue a case
analysis of the form ``given that the first $a-1$ blocks are empty,
and the $a$-th one has at least one point,\ldots''. 

In some of these cases, adding one such block produces a pattern which
may have a descent involving the black points, thus we need to further
refine the analysis with a further entry certifying the absence of
such descents.  Although not necessary by itself (as we control the
number of descents in certificates), this is important for the
following reason: differently from the primitive case, at every given
size there are several non-primitive exceptional standard
configurations with rank at least 3, which thus would proliferate in
all sorts of branches in the decomposition tree, making the analysis
inconclusive. It is restriction to primitive configurations that
allows to concentrate the identity configurations of all sizes into a
single branch (indeed, a result of Appendix
\ref{app.ididp} is that $\id_n$ is the only standard permutation in
$\Id_n$ starting with $(1,2,\ast)$).

The following picture, which is the pertinent crucial node in the tree
of Figure \ref{figure.Tthm}, represents the position of the blocks, by
green bullets. The red bullets correspond to blocks which are
certified to be empty, by the fact that, as we said, all entries at
the right of $k+1$ are larger than $k+1$.
\begin{center}
    \includegraphics{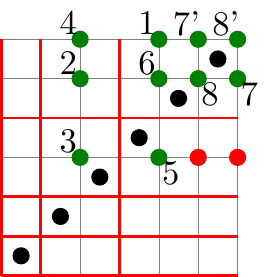}
\end{center}

\noindent
We have ten candidate blocks, labeled $\{1,2,\ldots,7,8;7',8'\}$, that
we now analyse. First of all, it is easily seen that cases $7'$ and
$8'$ coincide, up to relabeling, to cases $7$ and $8$, respectively.
Then, cases $1$ and $2$ can be analysed in one stroke, by adding one
entry in the corresponding block, while cases
from 3 to 7 require a finer analysis in order to ensure primitivity,
and thus involve the insertion of two entries, one in the
corresponding block, the other one in an available block among those
splitting the descent.  We have thus several cases, labeled from $3.1$
to $7.4$ (4 cases for blocks 3, 5, 7 and 6 cases for blocks 4 and 6).
The case of block 8 is the one that reduces to the sole
configuration $\id_n$.

All cases from $1$ to $7.4$ correspond to some case appearing in
Proposition~\ref{prop.fortab.T_tab_all_case} (which explains the
naming of cases in the proposition).  Sometimes this occurs after
dropping one entry. The following scheme summarizes this
correspondence:
\[
\begin{array}{rl}
1 & 
\casb{1}
\\[2mm]
2 & 
\casb{2}
\\[2mm]
3 & 
\left\{
\begin{array}{rcl}
3.1 && \casb{3.1}
\\
3.2 && \casb{3.2}
\\
3.3 && \casb{3.3}
\\
3.4 & \dro{3} & \cas{2}
\end{array}
\right.
\\[9mm]
4 & 
\left\{
\begin{array}{rcl}
4.1 && \casb{4.1}
\\
4.2 && \cas{3.3}
\\
4.3 & \dro{4} & \cas{2}
\\
4.4 & \dro{3} & \cas{1}
\\
4.5 & \dro{8} & \cas{2}
\\
4.6 & \dro{7} & \cas{2}
\end{array}
\right.
\end{array}
\quad
\begin{array}{rl}
5 & 
\left\{
\begin{array}{rcl}
5.1 && \casb{5.1}
\\
5.2 && \casb{5.2}
\\
5.3 && \cas{3.2}
\\
5.4 & \dro{6} & \cas{1}
\end{array}
\right.
\\[9mm]
6 & 
\left\{
\begin{array}{rcl}
6.1 && \casb{6.1}
\\
6.2 & \dro{7} & \cas{1}
\\
6.3 && \cas{5.2}
\\
6.4 & \dro{5} & \cas{1}
\\
6.5 & \dro{8} & \cas{1}
\\
6.6 & \dro{6} & \cas{2}
\end{array}
\right.
\\[13mm]
7 & 
\left\{
\begin{array}{rcl}
7.1 && \casb{7}
\\
7.2 && \casb{7}
\\
7.3 && \cas{6.4}
\\
7.4 & \dro{8} & \cas{2}
\end{array}
\right.
\end{array}
\]
On the left column, we have the case in the
decomposition. Possibly, the annotation ``$\dro{\cdot\,}$'' specifies
that we have to drop one entry, in the given column. In
``$\cas{\,\cdot\,}$'' we put the label of
Table~\ref{tab.T_tab_all_case} (in boldface, if it is the first
appearence in the list).

The two cases $7.1$ and $7.2$ correspond to pattern $7$ in which one
or the other gray block is certified to be non-empty, which is a
sufficient condition for having a zig-zag path certificate.

The lower part of Figure \ref{figure.Tthm} specifies which is which
among all subcases 3.1, \ldots, 7.4 associated to blocks 3 to~7 (the
same information is also carried in a more detailed way in 
the tables at the end of this subsection:
we remind the configuration associated to the block, analyse the
possible positions for a further entry splitting the descent, using
blue construction-lines, and then list the corresponding cases).

As we anticipated, if all blocks $\{1,2,\ldots,7;7'\}$ are empty,
and only blocks 8 (and $8'$) are possibly non-empty,
then the permutation is $\id_n$ for $n \geq 6$,
and the theorem does not apply, since the class $\Id_n$ is
exceptional.
\qed

\vfill

\noindent
\begin{tabular}{|c|cc|}
\hline
\begin{tikzpicture}[scale=0.35]
\useasboundingbox (0,0) rectangle (8,7);
\permutationtwo{1,2,3,4,5,6}{1,2,3,5}{1,2,4};
\draw (3,4) [fill,green!50!black!100] circle (.2);
\node at (4,7.7) {block~3};
\end{tikzpicture}
&
\begin{tikzpicture}[scale=0.35]
\permutationtwo{1,2,4,5,3,6,7,8}{1,2,3,7}{1,2,6};
\draw (3+.5,4+.5) [fill,green!50!black!100] circle (.2);
\draw (4+.5,5+.5) [fill,green!50!black!100] circle (.2);
\node at (5,9.7) {3.1};
\end{tikzpicture}
&
\begin{tikzpicture}[scale=0.35]
\permutationtwo{1,2,5,3,6,4,7,8}{1,2,3,7}{1,2,5};
\draw (3+.5,5+.5) [fill,green!50!black!100] circle (.2);
\draw (6+.5,4+.5) [fill,green!50!black!100] circle (.2);
\node at (5,9.7) {3.2};
\end{tikzpicture}
\\
\begin{tikzpicture}[scale=0.35]
\useasboundingbox (2,.5) rectangle (7,8);
\permutationtwo{1,2,4,3,5,6,7}{1,2,3,6}{1,2,5};
\draw (3+.5,4+.5) [fill,green!50!black!100] circle (.2);
\draw [blue] (4,1)--(4,8);
\draw [blue] (1,4)--(8,4);
\draw (3,4) [fill,blue] circle (.2);
\draw (6,4) [fill,blue] circle (.2);
\draw (4,7) [fill,blue] circle (.2);
\draw (4,8) [fill,blue] circle (.2);
\draw (7,4) [fill,red] circle (.2);
\draw (8,4) [fill,red] circle (.2);
\end{tikzpicture}
&
\begin{tikzpicture}[scale=0.35]
\permutationtwo{1,2,4,8,3,5,6,7}{1,2,3,6}{1,2,6};
\draw (3+.5,4+.5) [fill,green!50!black!100] circle (.2);
\draw (4+.5,8+.5) [fill,green!50!black!100] circle (.2);
\node at (5,9.7) {3.3};
\end{tikzpicture}
&
\begin{tikzpicture}[scale=0.35]
\permutationtwo{1,2,4,7,3,5,6,8}{1,2,3,6}{1,2,6};
\draw (3+.5,4+.5) [fill,green!50!black!100] circle (.2);
\draw (4+.5,7+.5) [fill,green!50!black!100] circle (.2);
\node at (5,9.7) {3.4};
\end{tikzpicture}
\\
\hline
\end{tabular}
\\
\rule{0pt}{10pt}
\\
\begin{tabular}{|c|ccc|}
\hline
\begin{tikzpicture}[scale=0.35]
\useasboundingbox (0,0) rectangle (8,7);
\permutationtwo{1,2,3,4,5,6}{1,2,3,5}{1,2,4};
\draw (3,7) [fill,green!50!black!100] circle (.2);
\node at (4,7.7) {block~4};
\end{tikzpicture}
&
\begin{tikzpicture}[scale=0.35]
\permutationtwo{1,2,7,8,3,4,5,6}{1,2,3,5}{1,2,6};
\draw (3+.5,7+.5) [fill,green!50!black!100] circle (.2);
\draw (4+.5,8+.5) [fill,green!50!black!100] circle (.2);
\node at (5,9.7) {4.1};
\end{tikzpicture}
&
\begin{tikzpicture}[scale=0.35]
\permutationtwo{1,2,4,8,3,5,6,7}{1,2,3,6}{1,2,6};
\draw (3+.5,4+.5) [fill,green!50!black!100] circle (.2);
\draw (4+.5,8+.5) [fill,green!50!black!100] circle (.2);
\node at (5,9.7) {4.2};
\end{tikzpicture}
&
\begin{tikzpicture}[scale=0.35]
\permutationtwo{1,2,6,8,3,4,5,7}{1,2,3,5}{1,2,6};
\draw (3+.5,6+.5) [fill,green!50!black!100] circle (.2);
\draw (4+.5,8+.5) [fill,green!50!black!100] circle (.2);
\node at (5,9.7) {4.3};
\end{tikzpicture}
\\
\begin{tikzpicture}[scale=0.35]
\useasboundingbox (2,.5) rectangle (7,8);
\permutationtwo{1,2,7,3,4,5,6}{1,2,3,5}{1,2,5};
\draw (3+.5,7+.5) [fill,green!50!black!100] circle (.2);
\draw [blue] (3,1)--(3,8);
\draw [blue] (1,8)--(8,8);
\draw (4,8) [fill,blue] circle (.2);
\draw (3,7) [fill,blue] circle (.2);
\draw (6,8) [fill,blue] circle (.2);
\draw (7,8) [fill,blue] circle (.2);
\draw (3,4) [fill,blue] circle (.2);
\draw (3,6) [fill,blue] circle (.2);
\draw (8,8) [fill,blue] circle (.2);
\end{tikzpicture}
&
\begin{tikzpicture}[scale=0.35]
\permutationtwo{1,2,7,3,4,8,5,6}{1,2,3,5}{1,2,5};
\draw (3+.5,7+.5) [fill,green!50!black!100] circle (.2);
\draw (6+.5,8+.5) [fill,green!50!black!100] circle (.2);
\node at (5,9.7) {4.4};
\end{tikzpicture} 
&
\begin{tikzpicture}[scale=0.35]
\permutationtwo{1,2,7,3,4,5,8,6}{1,2,3,5}{1,2,5};
\draw (3+.5,7+.5) [fill,green!50!black!100] circle (.2);
\draw (7+.5,8+.5) [fill,green!50!black!100] circle (.2);
\node at (5,9.7) {4.5};
\end{tikzpicture} 
&
\begin{tikzpicture}[scale=0.35]
\permutationtwo{1,2,7,3,4,5,6,8}{1,2,3,5}{1,2,5};
\draw (3+.5,7+.5) [fill,green!50!black!100] circle (.2);
\draw (8+.5,8+.5) [fill,green!50!black!100] circle (.2);
\node at (5,9.7) {4.6};
\end{tikzpicture}
\\
\hline
\end{tabular}
\\
\rule{0pt}{10pt}
\\
\begin{tabular}{|c|cc|}
\hline
\begin{tikzpicture}[scale=0.35]
\useasboundingbox (0,0) rectangle (8,7);
\permutationtwo{1,2,3,4,5,6}{1,2,3,5}{1,2,4};
\draw (5,4) [fill,green!50!black!100] circle (.2);
\node at (4,7.7) {block~5};
\end{tikzpicture}
&
\begin{tikzpicture}[scale=0.35]
\permutationtwo{1,2,3,6,4,5,7,8}{1,2,3,7}{1,2,4};
\draw (5+.5,4+.5) [fill,green!50!black!100] circle (.2);
\draw (6+.5,5+.5) [fill,green!50!black!100] circle (.2);
\node at (5,9.7) {5.1};
\end{tikzpicture}
&
\begin{tikzpicture}[scale=0.35]
\permutationtwo{1,2,3,5,7,4,6,8}{1,2,3,6}{1,2,4};
\draw (5+.5,7+.5) [fill,green!50!black!100] circle (.2);
\draw (6+.5,4+.5) [fill,green!50!black!100] circle (.2);
\node at (5,9.7) {5.2};
\end{tikzpicture}
\\
\begin{tikzpicture}[scale=0.35]
\useasboundingbox (2,.5) rectangle (7,8);
\permutationtwo{1,2,3,5,4,6,7}{1,2,3,6}{1,2,4};
\draw (5+.5,4+.5) [fill,green!50!black!100] circle (.2);
\draw [blue] (5,1)--(5,8);
\draw [blue] (1,5)--(8,5);
\draw (5,4) [fill,blue] circle (.2);
\draw (5,7) [fill,blue] circle (.2);
\draw (5,8) [fill,blue] circle (.2);
\draw (3,5) [fill,blue] circle (.2);
\draw (6,5) [fill,blue] circle (.2);
\draw (7,5) [fill,red] circle (.2);
\draw (8,5) [fill,red] circle (.2);
\end{tikzpicture}
&
\begin{tikzpicture}[scale=0.35]
\permutationtwo{1,2,5,3,6,4,7,8}{1,2,3,7}{1,2,5};
\draw (3+.5,5+.5) [fill,green!50!black!100] circle (.2);
\draw (6+.5,4+.5) [fill,green!50!black!100] circle (.2);
\node at (5,9.7) {5.3};
\end{tikzpicture}
&
\begin{tikzpicture}[scale=0.35]
\permutationtwo{1,2,3,5,8,4,6,7}{1,2,3,6}{1,2,4};
\draw (5+.5,8+.5) [fill,green!50!black!100] circle (.2);
\draw (6+.5,4+.5) [fill,green!50!black!100] circle (.2);
\node at (5,9.7) {5.4};
\end{tikzpicture}
\\
\hline
\end{tabular}
\\
\rule{0pt}{10pt}
\\
\begin{tabular}{|c|ccc|}
\hline
\begin{tikzpicture}[scale=0.35]
\useasboundingbox (0,0) rectangle (8,7);
\permutationtwo{1,2,3,4,5,6}{1,2,3,5}{1,2,4};
\draw (5,6) [fill,green!50!black!100] circle (.2);
\node at (4,7.7) {block~6};
\end{tikzpicture}
&
\begin{tikzpicture}[scale=0.35]
\permutationtwo{1,2,3,4,6,7,5,8}{1,2,3,5}{1,2,4};
\draw (5+.5,6+.5) [fill,green!50!black!100] circle (.2);
\draw (6+.5,7+.5) [fill,green!50!black!100] circle (.2);
\node at (5,9.7) {6.1};
\end{tikzpicture}
&
\begin{tikzpicture}[scale=0.35]
\permutationtwo{1,2,3,4,7,5,8,6}{1,2,3,5}{1,2,4};
\draw (5+.5,7+.5) [fill,green!50!black!100] circle (.2);
\draw (8+.5,6+.5) [fill,green!50!black!100] circle (.2);
\node at (5,9.7) {6.2};
\end{tikzpicture}
&
\begin{tikzpicture}[scale=0.35]
\permutationtwo{1,2,3,5,7,4,6,8}{1,2,3,6}{1,2,4};
\draw (5+.5,7+.5) [fill,green!50!black!100] circle (.2);
\draw (6+.5,4+.5) [fill,green!50!black!100] circle (.2);
\node at (5,9.7) {6.3};
\end{tikzpicture}
\\
\begin{tikzpicture}[scale=0.35]
\useasboundingbox (2,.5) rectangle (7,8);
\permutationtwo{1,2,3,4,6,5,7}{1,2,3,5}{1,2,4};
\draw (5+.5,6+.5) [fill,green!50!black!100] circle (.2);
\draw [blue] (6,1)--(6,8);
\draw [blue] (1,6)--(8,6);
\draw (6,4) [fill,blue] circle (.2);
\draw (6,7) [fill,blue] circle (.2);
\draw (6,8) [fill,blue] circle (.2);
\draw (3,6) [fill,blue] circle (.2);
\draw (5,6) [fill,blue] circle (.2);
\draw (7,6) [fill,blue] circle (.2);
\draw (8,6) [fill,blue] circle (.2);
\end{tikzpicture}
&
\begin{tikzpicture}[scale=0.35]
\permutationtwo{1,2,3,4,6,8,5,7}{1,2,3,5}{1,2,4};
\draw (5+.5,6+.5) [fill,green!50!black!100] circle (.2);
\draw (6+.5,8+.5) [fill,green!50!black!100] circle (.2);
\node at (5,9.7) {6.4};
\end{tikzpicture}
&
\begin{tikzpicture}[scale=0.35]
\permutationtwo{1,2,3,4,7,5,6,8}{1,2,3,5}{1,2,4};
\draw (5+.5,7+.5) [fill,green!50!black!100] circle (.2);
\draw (7+.5,6+.5) [fill,green!50!black!100] circle (.2);
\node at (5,9.7) {6.5};
\end{tikzpicture}
&
\begin{tikzpicture}[scale=0.35]
\permutationtwo{1,2,6,3,4,7,5,8}{1,2,3,5}{1,2,5};
\draw (3+.5,6+.5) [fill,green!50!black!100] circle (.2);
\draw (6+.5,7+.5) [fill,green!50!black!100] circle (.2);
\node at (5,9.7) {6.6};
\end{tikzpicture}
\\
\hline
\end{tabular}
\\
\rule{0pt}{10pt}
\\
\begin{tabular}{|c|cc|}
\hline
\begin{tikzpicture}[scale=0.35]
\useasboundingbox (0,0) rectangle (8,7);
\permutationtwo{1,2,3,4,5,6}{1,2,3,5}{1,2,4};
\draw (6,7) [fill,green!50!black!100] circle (.2);
\node at (4,7.7) {block~7};
\end{tikzpicture}
&
\begin{tikzpicture}[scale=0.35]
\permutationtwo{1,2,3,4,5,8,6,7}{1,2,3,5}{1,2,4,9};
\draw (7+.5,6+.5) [fill,green!50!black!100] circle (.2);
\draw (6+.5,8+.5) [fill,green!50!black!100] circle (.2);
\node at (5,9.7) {7.1};
\end{tikzpicture}
&
\begin{tikzpicture}[scale=0.35]
\permutationtwo{1,2,3,4,5,7,8,6}{1,2,3,5}{1,2,4,9};
\draw (6+.5,7+.5) [fill,green!50!black!100] circle (.2);
\draw (7+.5,8+.5) [fill,green!50!black!100] circle (.2);
\node at (5,9.7) {7.2};
\end{tikzpicture}
\\
\begin{tikzpicture}[scale=0.35]
\useasboundingbox (2,.5) rectangle (7,8);
\permutationtwo{1,2,3,4,5,7,6}{1,2,3,5}{1,2,4};
\draw (6+.5,7+.5) [fill,green!50!black!100] circle (.2);
\draw (7+.5,6+.5) [fill,black] circle (.2);
\draw [blue] (7,1)--(7,8);
\draw [blue] (1,7)--(8,7);
\draw (7,4) [fill,red] circle (.2);
\draw (7,6) [fill,blue] circle (.2);
\draw (7,8) [fill,blue] circle (.2);
\draw (3,7) [fill,blue] circle (.2);
\draw (5,7) [fill,blue] circle (.2);
\draw (6,7) [fill,blue] circle (.2);
\draw (8,7) [fill,blue] circle (.2);
\end{tikzpicture}
&
\begin{tikzpicture}[scale=0.35]
\permutationtwo{1,2,3,4,7,5,8,6}{1,2,3,5}{1,2,4,9};
\draw (5+.5,7+.5) [fill,green!50!black!100] circle (.2);
\draw (7+.5,8+.5) [fill,green!50!black!100] circle (.2);
\node at (5,9.7) {7.3};
\end{tikzpicture}
&
\begin{tikzpicture}[scale=0.35]
\permutationtwo{1,2,7,3,4,5,8,6}{1,2,3,5}{1,2,5,9};
\draw (3+.5,7+.5) [fill,green!50!black!100] circle (.2);
\draw (7+.5,8+.5) [fill,green!50!black!100] circle (.2);
\node at (5,9.7) {7.4};
\end{tikzpicture}
\\
\hline
\end{tabular}


\subsection{Operators $q_1$ and $q_2$}
\label{ssec.q_op}

We have seen that 
$\bar{T}\left(\kS_{n-2}^{\rm (prim)}\right) \rightarrow \kS_{n}^{\rm (prim)}|_{\rank>2}$,
it is surjective on this set, and is a pullback function.
In order to produce a complete decomposition, we need to introduce two
more operators, $q_1$ and $q_2$, such that 
\begin{align}
\bar{q_1}\left(\kS_{n-1}^{\rm (prim)}\right)
&\rightarrow \kS_{n}^{\rm (prim)}|_{\rank=1}
\ef;
&
\bar{q_2}\left(\kS_{n-1}^{\rm (prim)}\right)
&\rightarrow \kS_{n}^{\rm (prim)}|_{\rank=2}
\ef.
\end{align}
We do this in this section.

\begin{definition}
Let $\s$ be a permutation of size $n$. We define $\add_i(\s)$ for
$i=1,\ldots, n$ to be the permutation of size $n+1$
$\add_i(\s):=\s(1)+1,\ldots,\s(i-1)+1,1,\s(i)+1,\ldots,\s(n)+1$.
\end{definition}
\noindent 
I.e., $\add_i(\s)$ is described by the diagrammatic manipulation below:
\[
\mraisebox{10pt}{\add_i}
\;
\mleftbox{%
\rule{5pt}{0pt}%
\underbrace{\rule{25pt}{0pt}}_{i-1} }%
\includegraphics[scale=1.8]{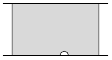}
\mraisebox{10pt}{=}
\mleftbox{%
\rule{5pt}{0pt}%
\underbrace{\rule{25pt}{0pt}}_{i-1} \; \mraisebox{-8pt}{i} }%
\includegraphics[scale=1.8]{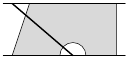}
\]
The two operators $q_1$ and $q_2$ are defined in terms of $\add_i$'s
for suitable $i$'s. In order to endow them with the appropriate
properties, we shall investigate the difference between the cycle
invariants of $\s$ and $\add_i(\s)$. This requires a case analysis
based on the `type' of permutations (among $H$- and $X$-type, see
Definition \ref{def.XHtype}). This is presented in Table
\ref{table.addi}, and illustrated in the following paragraph.

\begin{proposition}\label{pro.add1}
Let $i\in\{1,2\}$. Let $\s$ be a primitive permutation with invariant
\Lr\ with $r>i$.
\begin{itemize}
\item If $\s$ has type $X(r,j)$ then there exists exactly one index
  $\ell$ such that $\rank(\add_{\ell}(\s))=i$. 
\item If $\s$ has type $H(r_1,r_2)$ with $r_2 \neq i$ then there
  exists exactly one index $\ell$ such that \mbox{$\rank(\add_{\ell}(\s))=i$.}
\item If $\s$ has type $H(r_1,r_2)$ with $r_2=i$ then there exist
  exactly two indices $\{\ell,m\}=\{1,\s^{-1}(n)\}$ such that
  $\rank(\add_{\ell}(\s))=\rank(\add_{m}(\s))=i$; moreover
  $\add_{\ell}(\s)$ and $\add_{m}(\s)$ are in the same class, as
  $\add_{\ell}(\s)=R^{-1}\left(\add_{m}(\s)\right)$.
\end{itemize}
\end{proposition}

\begin{table}[p!]
\setlength{\unitlength}{72pt}
\begin{picture}(5.45,8.33)
\put(0,0){\includegraphics[scale=1.2]{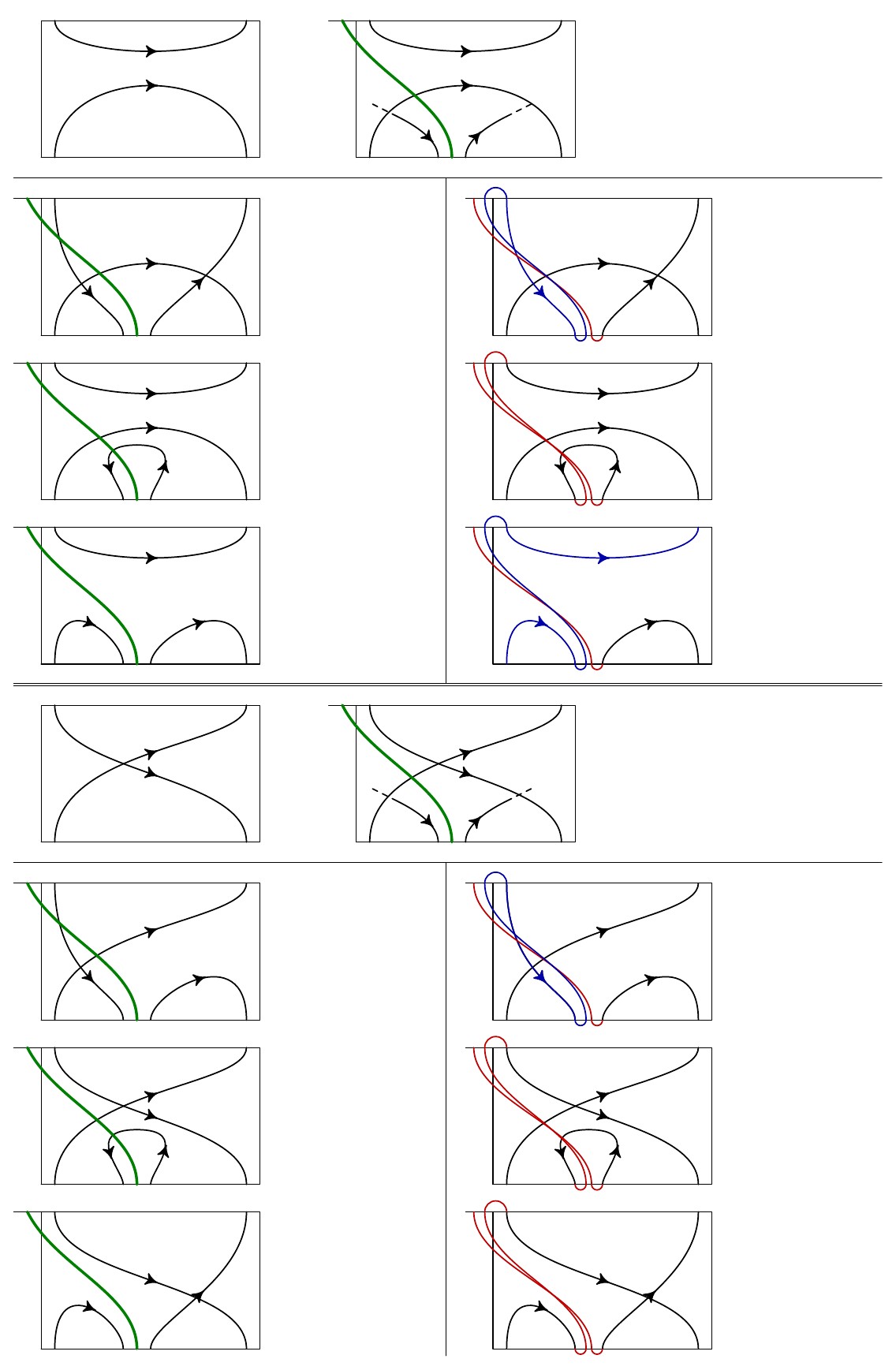}}
\put(0.24,7.4){\begin{picture}(1.34,.83)
\put(.6,.72){$2r_1$}
\put(.6,.3){$2r_2$}
\end{picture}}
\put(1,7.5){\begin{picture}(1.34,.83)
\put(1.25,0.27){?}
\put(2.25,0.27){?}
\put(2.8,.6){$H$-type\quad $H(r_1,r_2)$}
\end{picture}}
\put(0.24,6.315){\begin{picture}(1.34,.83)
\put(.13,.7){$2s-1$}
\put(.42,.54){$2r_1-2s+1$}
\put(.6,.3){$2r_2$}
\put(2.25,.75){\ding{172}}
\put(1.42,.45){\begin{minipage}{2.5cm}
\begin{center}
Once per\\
$1 \leq s \leq r_1$
\end{center}
\end{minipage}}
\end{picture}}
\put(3,6.315){\begin{picture}(1.34,.83)
\put(.13,.7){$2s-1$}
\put(.42,.54){$2r_1-2s+1$}
\put(.6,.3){$2r_2$}
\put(1.42,.45){\begin{minipage}{2.5cm}
\begin{align*}
&H(r'_1,r'_2)\\
r_1'&=r_1\!-\!s\!+\!1 \\
r_2'&=r_2 \\
\lambda'&=\lambda \!\cup\! \{s\}
\end{align*}
\end{minipage}}
\end{picture}}
\put(0.24,5.315){\begin{picture}(1.34,.83)
\put(.6,.72){$2r_1$}
\put(.92,.43){$2r_2$}
\put(.82,.19){$2\ell$}
\put(2.25,.75){\ding{173}}
\put(1.42,.4){\begin{minipage}{2.5cm}
\begin{center}
$\ell$ times per\\
each cycle of $\lam$\\
of length $\ell$
\end{center}
\end{minipage}}
\end{picture}}
\put(3,5.315){\begin{picture}(1.34,.83)
\put(.6,.72){$2r_1$}
\put(.92,.43){$2r_2$}
\put(.82,.19){$2\ell$}
\put(1.42,.45){\begin{minipage}{2.5cm}
\begin{align*}
&H(r'_1,r'_2)\\
r_1'&=r_1\!+\!\ell+1 \\
r_2'&=r_2 \\
\lambda'&=\lambda \!\setminus\! \{\ell\}
\end{align*}
\end{minipage}}
\end{picture}}
\put(0.24,4.315){\begin{picture}(1.34,.83)
\put(.6,.72){$2r_1$}
\put(.17,.12){$2s$}
\put(.77,.31){$2r_2-2s$}
\put(2.25,.75){\ding{174}}
\put(1.42,.45){\begin{minipage}{2.5cm}
\begin{center}
Once per\\
$0 \leq s \leq r_2$
\end{center}
\end{minipage}}
\end{picture}}
\put(3,4.315){\begin{picture}(1.34,.83)
\put(.6,.72){$2r_1$}
\put(.17,.12){$2s$}
\put(.77,.31){$2r_2-2s$}
\put(1.42,.45){\begin{minipage}{2.5cm}
\begin{align*}
&X(r',i')\\
r'&=r_2\!-\!s \\
i'&=r_1\!+\!s \\
\lambda'&=\lambda \!\cup\! \{i'\}
\end{align*}
\end{minipage}}
\end{picture}}
\put(0.24,3.21){\begin{picture}(1.34,.83)
\put(.73,.71){$2i+1$}
\put(.9,.36){$2r+1$}
\end{picture}}
\put(1,3.31){\begin{picture}(1.34,.83)
\put(1.25,0.27){?}
\put(2.25,0.27){?}
\put(2.8,.6){$X$-type\quad $X(r,i)$}
\end{picture}}
\put(0.24,2.125){\begin{picture}(1.34,.83)
\put(.73,.71){$2i+1$}
\put(.85,.32){$2r-2s$}
\put(.13,.7){$2s-1$}
\put(2.25,.75){\ding{175}}
\put(1.42,.45){\begin{minipage}{2.5cm}
\begin{center}
Once per\\
$1 \leq s \leq r$
\end{center}
\end{minipage}}
\end{picture}}
\put(3,2.125){\begin{picture}(1.34,.83)
\put(.73,.71){$2i+1$}
\put(.85,.32){$2r-2s$}
\put(.13,.7){$2s-1$}
\put(1.42,.45){\begin{minipage}{2.5cm}
\begin{align*}
&X(r',i')\\
r'&=r\!-\!s \\
i'&=i \\
\lambda'&=\lambda \!\cup\! \{s\}
\end{align*}
\end{minipage}}
\end{picture}}
\put(0.24,1.125){\begin{picture}(1.34,.83)
\put(.73,.71){$2i+1$}
\put(.9,.36){$2r+1$}
\put(.81,.13){$2\ell$}
\put(2.25,.75){\ding{176}}
\put(1.42,.4){\begin{minipage}{2.5cm}
\begin{center}
$\ell$ times per\\
each cycle of $\lam$\\
of length $\ell$
\end{center}
\end{minipage}}
\end{picture}}
\put(3,1.125){\begin{picture}(1.34,.83)
\put(.73,.71){$2i+1$}
\put(.9,.36){$2r+1$}
\put(.81,.13){$2\ell$}
\put(1.42,.45){\begin{minipage}{2.5cm}
\begin{align*}
&X(r',i')\\
r'&=r\!+\!\ell\!+\!1 \\
i'&=i \\
\lambda'&=\lambda \!\setminus\! \{\ell\}
\end{align*}
\end{minipage}}
\end{picture}}
\put(0.24,0.125){\begin{picture}(1.34,.83)
\put(.58,.71){$2i-2s+1$}
\put(.4,.54){$2r+1$}
\put(.17,.12){$2s$}
\put(2.25,.75){\ding{177}}
\put(1.42,.45){\begin{minipage}{2.5cm}
\begin{center}
Once per\\
$0 \leq s \leq i$
\end{center}
\end{minipage}}
\end{picture}}
\put(3,0.125){\begin{picture}(1.34,.83)
\put(.58,.71){$2i-2s+1$}
\put(.4,.54){$2r+1$}
\put(.17,.12){$2s$}
\put(1.42,.45){\begin{minipage}{2.5cm}
\begin{align*}
&H(r'_1,r'_2)\\
r'_1&=i\!-\!s\!+\!1 \\
r'_2&=r\!+\!s\!+\!1 \\
\lambda'&=\lambda \!\setminus\! \{i\}
\end{align*}
\end{minipage}}
\end{picture}}
\end{picture}
\caption{\label{table.addi}Modification to the cycle invariant of a
  configuration produced by $\add_i$. In green, the newly-added edge.
  In red, parts which get added to the rank path. In blue, parts which
  are singled out to form a new cycle.}
\end{table}

\proof This emerges from the analysis of Table \ref{table.addi}. We
discuss the three propositions one by one, with reference to the 6
rows of the table.
\begin{itemize}
\item If $X(r,j)$, it can be checked that the only possibility for
  having rank $i$ is case \dingD$\, \!$ with $s=r-i$. The reason why we need
  $r>i$ is because, e.g.\ if $r=1$, either $s=0$ or 1, but both
  eventualities are impossible, since cycles of length 0 do not exist
  and cycles of length 1 are not allowed in primitive permutations.
\item If $H(r_1,r_2)$ with $r_2 \neq i$, it can be checked that the
  only possibility is case \dingC$\, \!$ with $s=r_2-i$, if $r_2>i$, and case
  \dingA$\, \!$ with $s=r_1$, if $r_2=1$. (In fact, the only case with $r_2
  <i$ is when $i=2$ and $r_2=1$).
\item If $H(r_1,r_2)$ with $r_2=i$, then there are exactly two
  possibilities: case \dingC $\, \!$ of the table with $s=0$, and case \dingA$\, \!$
  with $s=r_1$. Clearly we go from one target configuration to another
  with the operator $R^{\pm 1}$, since these only differ by the edge
  $(\alpha,1)$ which is $(1,1)$ in case \dingC$\, \!$ and $(\s^{-1}(n),1)$ in
  case \dingA$\, \!$.
\end{itemize}
\qed

\noindent
The previous proposition shows that for configurations $\s$ with rank
$i$, there exists one and only one configuration $\tau$ and index
$\ell$ such that $\s=\add_\ell(\tau)$, with an exceptional case in
which there are two such configurations, which are however easily
shown to be in the same class. We can thus just ``break the tie'' for
the third situation in Proposition \ref{pro.add1} (see Figure
\ref{fig.duecasiQ1} for an illustration), and give the following
definition:

\begin{figure}[tb!]
\begin{center}
\setlength{\unitlength}{60pt}
\begin{picture}(6,2.8)
\put(0,0){\includegraphics{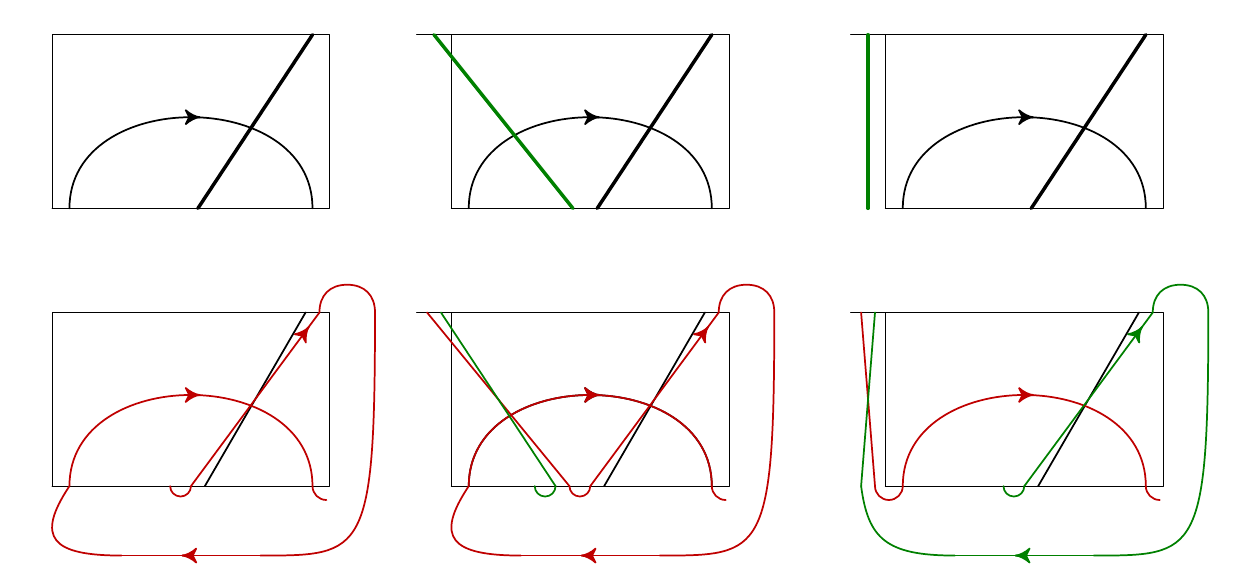}}
\put(.8,2.23){$2r$}
\put(0.1,2.65){$\s$}
\put(2.1,2.65){$\tau$}
\put(4.1,2.65){$\tau'$}
\end{picture}
\end{center}
\caption{The two possible consistent definitions of $q_r$ on a
  configuration $\s$ of type $H(r',r)$, mentioned in the third case of
  Proposition \ref{pro.add1}.  The two configurations $\tau$ and
  $\tau'$ are such that $\tau= R \tau'$.  Our choice is to define $q_r
  \s = \tau$.
\label{fig.duecasiQ1}}
\end{figure}

\begin{definition}
Let $i\in\{1,2\}$.  Let $\tau$ be a primitive permutation with rank
$r>i$.  Define $q_i(\tau)$ to be the unique $\add_{\ell}(\tau)$ such
that $\rank(\add_{\ell}(\tau))=i$, if we are in one of the first two
situations of Proposition \ref{pro.add1}, and
$q_i(\tau)=\add_{\tau^{-1}(n)}(\tau)$ if we are in the third
situation.
\end{definition}

Note that the edge $(\ell,1)$ in $q_i(\s)$ is never $(1,1)$. This is
true by a further investigation of the table, for the first two
situations, and holds at sight for the third situation, exactly
because of our choice of convention (if it were, this would mean that
$\s^{-1}(n)=1$, which implies that $\s$ is reducible). Thus the edge
$(i,1)$ is never a pivot, a fact that will be useful when we pass to a
reduced dynamics after the action of a surgery operator $q_i$.

Now that we have defined our $q_i$'s, we need to show that they satisfy
the properties outlined in Section \ref{ssec.TQSintro}, i.e.\ the
statement of the following lemma:
\begin{lemma}
Let $\s$ and $\s'$ be configurations in the same class. Then 
$q_i(\s) \sim q_i(\s')$ for both $i=1,2$.
\end{lemma}
\proof
Clearly both $q_i$'s are injective. They are actually `almost
bijective', as the number of preimages is at most 2 (and, for what it
matters, in large size there is a unique preimage in the majority of
configurations).

Let $\tau$ and $\tau'$ be primitive permutations in the same class
$C$, and let $S$ be a sequence such that $S(\tau)=\tau'$. Call
$\s=q_i(\tau)$. The reduced permutation of $\s$ where the edge
$(\ell,1)$ is the only gray edge is $\tau$. Since this edge is not a
pivot, we can use the boosted dynamics and define $\s'=B(S)(\s)$, and
we have that $\s'=q_i(\tau')$.

Indeed, by definition of $B(S)$, removing the edge $(\ell',1)$ in
$\s'$ gives $\tau'$, and $(\ell',1)$ is not a pivot, which coincides
exactly with our characterisation of $q_i(\tau')$ (since $\s$ has rank
$i$, so has $\s'$). Thus $q_i(\s) \sim q_i(\s')$.  \qed

\noindent
As a consequence of this lemma, analogously to what we have done for
the operator $T$, we can define the two maps $\bar{q}_1$ and $\bar{q}_2$
on classes. 

First of all, we observe the triviality
\begin{lemma}\label{lem.q1.id}
$\bar{q}_1(\Id_n)=\tree_{n+1}$.
\end{lemma}
\noindent
(this is seen by direct inspection of the canonical representatives, as
$q_1(\id_n) = R^{-2}\, \idp_{n+1}$).

Then, yet again, these maps have a definite behaviour for what
concerns the invariants.  Let us introduce the useful notation
\begin{notation}
Let $\lambda$ be a cycle invariant, and let $j \in \lambda$. We use
$\lambda(j)$ as a shortcut for $\lambda \setminus \{j\}$.
\end{notation}
Then:
\begin{lemma}\label{lem_q1_operator}
Let $C$ be a class with invariant $(\lambda,r,s)$ with $r>1$. Then
$C' = \bar{q}_1(C)$ has invariant $(\lambda\cup\{r\},1,s)$.
\end{lemma}

\proof
The claim on the sign comes as an application of Corollary
\ref{cor.q1sign}. Indeed, as seen in Proposition
\ref{prop.signinvdyn}, the sign is a class invariant, so we can take
$\s \in C'$ to be standard with $\s(n)=n$. Since $\s$ has rank 1, we
know that $\s(2)=n-1$. This is the permutation $\s$ on the LHS of
Corollary \ref{cor.q1sign},
and removing the edge $(1,1)$, i.e.\ taking the preimage of
$R^{-1} \s$ under $q_1$, we obtain the permutation $\t$ on the RHS of
Corollary \ref{cor.q1sign} (indeed, we have $q_1(\tau)=R \s$ because
we are in the third case of the definition of $q_1$).

For what concerns cycle and rank invariants, this comes from a
straighforward inspection (cf.\ again Table~\ref{table.addi}).
\qed

\begin{corollary}
\label{cor.goback123a}
Let $\s$ be a permutation with invariant $(\lambda,1,s)$ then either
$\tau := q_1^{-1}(\s)$ or $\tau' := q_1^{-1}(R\s)$ or both are
defined. Each of these configurations have invariant
$(\lambda(i),i,s)$ for some $i\in \lambda$.
\end{corollary}

\begin{lemma}\label{lem_q2_operator}
Let $C$ be a class with invariant $(\lambda,r,s)$ with $r>2$. Then
$\bar{q}_2(C)$ has invariant $(\lambda\cup\{r-1\},2,0)$.
\end{lemma}
\proof The claim on the sign comes as an application of Corollary
\ref{cor.q2sign}. Again we can choose a standard representative with
$\s(n)=n$. Then there must exist $j$ such that $\s'=R^j(\s)$ has
$\s'(n)=n$ and $\s'(n-1)=n-1$. Since $\s'$ has rank 2, this
configuration is of the form of the LHS permutation of Corollary
\ref{cor.q2sign}, and removing the edge $(\s'^{-1}(1),1)$, which
corresponds to take the preimage of the surgery operator $q_2$,
i.e.\ produces a permutation of smaller size (call it $\tau$) such
that $q_2(\tau)=\s'$ (or $q_2(\tau)=R(\s')$ if $\s'^{-1}(1)=1$). This
proves that, regardless from the sign of $\tau$, the sign of $\s'$ is
zero.

Again, for what concerns cycle and rank invariants, this comes from a
straighforward inspection.
\qed

\begin{corollary}
\label{cor.goback123b}
Let $\s$ be a permutation with invariant $(\lambda,2,0)$ then either
$q_2^{-1}(\s)$ or $q_2^{-1}(R\s)$ are defined.  This configuration has
invariant $(\lambda(i-1),i,s)$ for some $i\in \lambda$ and
$s\in\{-1,0,+1\}$, more precisely,
\begin{itemize}
\item if $\lambda(i-1)$ has a positive number of even cycles, then $s=0$;
\item if $\lambda(i-1)$ has no even cycles, then $s=\pm 1$.
\end{itemize}
\end{corollary}

\bigskip
\noindent
At this point, we are left with the investigation of the pullback
properties of our operators.
\begin{lemma}\label{lem.inv.st.one}
Let $\s=(1,\s(2),\ldots,\s(n))$ be a permutation of type $X(i,j)$ with
cycle and rank invariant $(\lambda,i)$. Then
$\tau=(\s(2)-1,\ldots,\s(n)-1)$ has type $H(j,i)$ and invariant
$(\lambda(j),r')$, with $r'=i+j-1$.
\end{lemma}
\begin{align*}
\s &=
\raisebox{-30pt}{\setlength{\unitlength}{60pt}
\begin{picture}(2,.4)
\put(0,0){\includegraphics{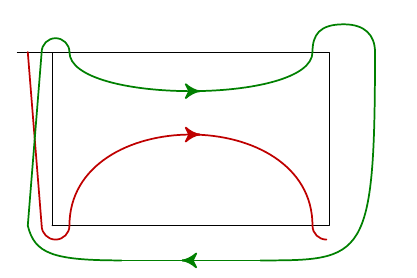}}
\put(.85,0.56){$2i$}
\put(.85,0.94){$2j$}
\end{picture}}
&
\tau
&=
\raisebox{-30pt}{\setlength{\unitlength}{60pt}
\begin{picture}(2,1.05)
\put(0,0){\includegraphics{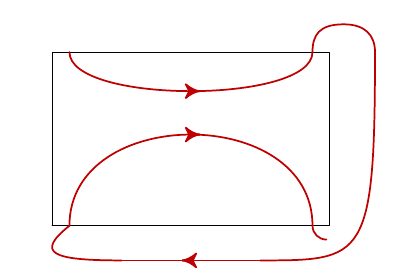}}
\put(.85,0.56){$2i$}
\put(.85,0.94){$2j$}
\end{picture}}
\end{align*}
\proof
This is the reverse implication of case \dingC\ in Table
\ref{table.addi}, specialised to $s=0$.
\qed
\begin{theorem}[Pullback of $\bar{q}_1$]
\label{thm.q1}
Let $C$ be a class with invariant $(\lambda,1,s)$ then for every
$j\in\lambda$ there exists a class $B_j$ with invariant
$(\lambda(j),j,s)$ such that $\bar{q}_1(B_j)=C$. In other words,
$\bar{q}_1$ is a pullback function w.r.t.\ the equivalence relation
given by the invariant, i.e.\ $C \xrightarrow{\cdot/\sim} (\lam,r,s)$.
\end{theorem}

\proof First of all, let us remark that, by Corollary
\ref{cor.goback123a}, the set of invariants $(\lambda(j),j,s)$
coincides with the set of preimages of $(\lambda,1,s)$ under
$\widetilde{\bar{q}_1}$.\footnote{Where $\widetilde{\cdot}$ is the
  notion introduced in the definition of pullback function.}

Now, we take a standard family of $C$. Then, by Lemma
\ref{trivialbutimp}, for every $j\in \lambda$ there exists
$\s_j=(1,\s_j(2)\ldots)$ of type $X(1,j)$ such that
$\tau_j=(\s_j(2)-1,\ldots,\s_j(n)-1)$ is irreducible.

Moreover, $\tau_j$ has invariant $(\lambda(j),j,s)$. We know this, for
cycle and rank, from Lemma \ref{lem.inv.st.one} (in the case $i=1$),
and, for the sign, from Lemma \ref{lem_q1_operator}, the fact that
$q_1(\tau_j)=R(\s_j)$, and that $\s_j$ has sign $s$. Thus for every
$j$ the class $B_j\ni \tau_j$ satisfies $\bar{q}_1(B_j)=C$.  
\qed

\noindent
For operator $q_2$, the facts we can establish at this point are
a bit weaker, the trouble coming from the ``loss of memory'' on the
sign, which becomes 0 regardless of what was its value on the
preimage. As a result, we will need some extra work in the induction
step of the next section, while at this point we establish the
following:

\begin{theorem}[Pullback of $\bar{q}_2$]
\label{thm.q2}
Let $C$ be a class with invariant $(\lambda,2,0)$. Then for every
$j\in\lambda$ there exists a class $B_j$ with invariant
$(\lambda(j),j+1,s)$, for some $s$, such that $\bar{q}_2(B_j)=C$.  In
other words, $\bar{q}_2$ is a pullback function w.r.t.\ the
equivalence relation given by the cycle invariant, i.e.\ $C
\xrightarrow{\cdot/\sim} (\lam,r)$.
\end{theorem}

\proof The proof is analogous to the one above (with the discussion of
the sign left out).  

By Corollary \ref{cor.goback123b}, the set of cycle invariants
$(\lambda(j),j+1)$ coincides with the set of preimages of $(\lambda,2)$
under $\widetilde{\bar{q}_2}$.

We take a standard family of $C$, then by Lemma
\ref{trivialbutimp} for every $j\in \lambda$ there exists
$\s_j=(1,\s_j(2)\ldots)$ with type $X(2,j)$ such that
$\tau_j=(\s_j(2)-1,\ldots,\s_j(n)-1)$ is irreducible.

Moreover, $\tau_j$ has cycle and rank invariant $(\lambda(j-1),j)$ by
Lemma \ref{lem.inv.st.one} (in the case $i=2$).  Since
$q_2(\tau_j)=R(\s_j)$, for every $j$ the class $B_j\ni \s_j'$
satisfies $\bar{q}_2(B_j)=C$.
\qed

\section{The induction}
\label{sec.induction}

This section is devoted to the main induction of the paper, that
implies the classification theorem. The induction mainly works at the
level of non.exceptional classes, so, before doing this, we need to
exclude the proliferation of distint classes with the same invariant
due to the action of our surgery operators on exceptional
classes. This requires to establish `fusion lemmas', i.e.\ lemmas of
the form, for a given exceptional class $I=\{I_n\}$, for $n \geq n_0$,
there exists non-exceptional classes $C=\{C_n\}$, such that a certain
surgery operator $X$, acting on $I$, gives the same class as the
action on $C$ (in formulas, $\exists\; C \;:\; X(C) \sim X(I)$).  In
the paragraph above, in principle, $I$ may be $\Id$ or $\Id'$, while
$X$ may be $T$, $q_1$ or $q_2$. However, we will see that not all of
these cases need to be considered.

In proving these lemmas, we will make the ansatz that the classes
$C_n$ at different $n \geq n_0$ have one representative with a `nice
structure' in $n$, e.g.\ a certain permutation of size $n_0$ and an
identity block of size $n-n_0$ in a specific position, and that the
sequence connecting this representative to a canonical representative
of the exceptional class, once written in terms of dynamics operators
and their inverses (so that the variable-size identity block is never
broken up along the dynamics), is the same for all sizes.

In particular, we need to establish that the forementioned
representative is in fact in a class $C$ which is non-exceptional.
This can be done in two ways, either by using our knowledge of the
structure of \emph{all} configurations in the special classes $\Id$
and $\tree$, presented in Appendix~\ref{sec.excp_class}, or by using a
result, also proven in Appendix~\ref{sec.excp_class}, stating that
both $\Id$ and $\tree$ have a unique standard family (in particular,
in each of these classes there is only one standard permutation $\pi$
with $\pi(1)=1$ and $\pi(n)=n$, namely $\pi=\id_n$ for $\Id_n$ and
$\pi=\id'_n$ for $\tree_n$), so that, if the representative is
standard, the check is straightforward. The second method is more
compact, so we will adopt it in this section.

The very last lemma of this form will be based on a more general
ansatz (see the pattern in the following
Definition~\ref{def.Vblock}). 
The sequence\footnote{Recall that we use the shortcuts $\bar{L}$,
  $\bar{R}$ for $L^{-1}$ and $R^{-1}$.}  
$w_n \in \{L,R,\bar{L},\bar{R}\}^*$ is not the same for
all $n \geq n_0$, but rather has the form
$w_n = w_{\rm end} w^{n-n_0} w_{\rm start}$. Indeed, instead of having
a `large' identity block of size $n-n_0$, we have two `large' simple
blocks, of size $n-n_0-k$ and $k$, and the iterated sequence $S$ in
the middle of $S_n$ is used to change the value of $k$ by one.

Such intelligible patterns allow to give unified proofs of these
fusion lemmas which hold for all $n$ large enough, as desired.  We
admit that these patterns could have hardly been guessed without the
aid of computer search. We have been lucky in the respect that not
only our ansatz holds, but also the sequences implementing our ansatz
are the shortest ones connecting the two guessed representatives, so
that, once the `good' representatives have been found by the computer
at finite size, we could invent a proof just by analysing the
structure of the numerically-found geodesic path in the Cayley graph.

\subsection{Simple properties of the invariant}
\label{ssec.signindu}

Here we prove a few facts on the invariants, that we have previously
stated without proof. A stronger characterisation will emerge from the
main induction, however it is instructive to deduce the following
lemma, at the light of the previous results, and with a small extra
effort.

\begin{lemma}
\label{lem.induSign}
Let $C$ be a class with invariant $(\lam,r,s)$ and let $\ell$ be the
number of cycles (not including the rank). The following statements
hold:
\begin{enumerate}
\item The list $\{r,\lam_1,\ldots,\lam_\ell\}$ contains an even number of even entries.
\item The sign of the class, $s(C)=\mathrm{Sign}(\Abar(C))$ can be
  written as
  $s(C)=2^{-\frac{n+\ell}{2}} \Abar(C)$.
\item The sign $s(C)$ is zero if and only if some of the entries $r$
  and $\lam_j$ of the cycle invariant $(\lam,r)$ are even.
\end{enumerate}
\end{lemma}

\proof First of all, the statements can be verified for the
(exceptional) class of the identity, since computing (by induction)
the cycle invariant $(\lam,r)$ of $id_n$ is an easy task, while the
calculation of the Arf invariant was done in
Corollary~\ref{cor.id_sign}. The results are summarised in
Table~\ref{tab.guardaId}, top.

Then, the statements follow for the other exceptional class,
$\tree_n$, by using the fact that $\bar{q_1}(\Id_n)=\tree_{n+1}$
(Lemma \ref{lem.q1.id}) and the action of $\bar{q}_1$ on the Arf
invariant of a class (Corollary \ref{cor.q1sign}). The results are
summarised in Table~\ref{tab.guardaId}, bottom.


Finally, the statements follow inductively for all non-exceptional
classes, for the three cases of rank equal to 1, 2, or at least 3,
using the pullback results for the operators $\bar{q}_1$,
$\bar{q}_2$ and $\bar{T}$ (Theorems \ref{thm.q1}, \ref{thm.q2} and
\ref{thm_T_surjectif}, respectively) and the action of the operators
on the cycle invariant of a class (Lemma \ref{lem_q1_operator},
\ref{lem_q2_operator} and \ref{lem_T_operator} respectively),
established in Section \ref{sec.operators}, and on the Arf invariant
of a class (Corollary \ref{cor.q1sign}, \ref{cor.q2sign} and
\ref{cor.Tsign} respectively), established in
Section~\ref{sec.signinv}.

In order to illustrate how the argument works, let us analyse it in
full detail in the case of operator $T$. Let $C$ be a class of size
$n$ with invariant $(\lam,r,s)$, with $r \geq 3$.  Then, by the
pullback theorem for $T$ (Theorem \ref{thm_T_surjectif}), there
exists a class $B$ such that $\bar{T}(B)=C$, and $B$ shall have
invariant $(\lam,r-2,s)$ in agreement with Lemma \ref{lem_T_operator}.
Clearly $T$ preserves the parity of the number of even entries in the
cycle invariant (only $r$ has changed, and by 2), so, as by induction
$B$ has an even number of even cycles, this also holds for~$C$.

By induction we know that $s(B)=2^{-\frac{n-2+\ell}{2}}\Abar(B)$ where
$\ell=\ell(B)$ is the number of cycles in $B$. However, we have just
noticed that $\ell(C)=\ell(B)=\ell$.  By Corollary \ref{cor.Tsign},
$\Abar(C)=2\Abar(B)$, thus 
$2^{-\frac{n+\ell}{2}}\Abar(C)=2 \times 2^{-\frac{n+\ell}{2}}\Abar(B)=
2^{-\frac{n-2+\ell}{2}}\Abar(B)=s(B)$.
Hence, since we have $s(C)=s(B)$ by Lemma \ref{lem_T_operator}, we
have consistently
$s(C)=\mathrm{Sign}(\Abar(C))=2^{-\frac{n+\ell}{2}}\Abar(C)$.

Finally, by induction $s(B)$ is zero if and only if some of the
entries among $r-2$ and $\lam_j$ are even. As we have just established
that $s(C)=s(B)$ and that $C$ has as many even cycles as $B$,
the statement holds also for $C$.  \qed

\begin{table}[tb]
\[
\begin{array}{r|cccccccc}
n & 2 & 3 & 4 & 5 & 6 & 7 & 8 & 9 
\\
\hline
(\lam,r) \text{ of } \Id_n &
(\emptyset,1) & (\{1\},1) &
(\emptyset,3) & (\{2\},2) &
(\emptyset,5) & (\{3\},3) &
(\emptyset,7) & (\{4\},4) 
\\
\Abar(\Id_n) &
-2 & -4 & -4 & 0 & 8 & 16 & 16 & 0
\\
\hline
(\lam,r) \text{ of } \Id'_n &
&&
(\{1,1\},1) & (\{3\},1) &
(\{2,2\},1) & (\{5\},1) &
(\{3,3\},1) & (\{7\},1) 
\\
\Abar(\Id'_n) &
&& -8 & -8 & 0 & 16 & 32 & 32
\end{array}
\]
\caption{\label{tab.guardaId}The cycle and sign invariants of $\Id_n$
  and $\Id'_n$. The pattern of $(\lam,r)$ is repeated with period 2,
  while the pattern of $\Abar$ is repeated with period 8. We include
  also the non-primitive classes $\Id_2$ and $\tree_4$, in order to
  illustrate the mechanism on smaller values of~$n$.}
\end{table}

The results in Table \ref{tab.guardaId} imply the following simple
fact
\begin{lemma}
\label{lem.TididpDiffer}
For all $n \geq 4$, the cycle invariants $(\lam,r)$ of
$\bar{T}(\Id_n)$ and of $\bar{T}(\Id'_n)$ are distinct.
\end{lemma}
\proof Just combine Table \ref{tab.guardaId} and
Lemma~\ref{lem_T_operator}.
\qed

\subsection{Fusion lemmas for $T$}
\label{ssec.tech}

\begin{lemma}[Fusion lemma for $(T,\Id)$]
\label{lem.T.id}
Let $n \geq 8$. There exists some non-exceptional class $A_n$ such that
$\bar{T}(A_n)=\bar{T}(\Id_n)$.
\end{lemma}

\proof
We just searched (by computer) a sequence such that $T(\id_n)$ is transformed
into a $\s=T(\tau)$ where $\tau\notin \Id_n$.
For $n\geq 7$ the configuration
$\bar{R}L^3RL\bar{R}\bar{L}R\bar{L}RL^2\bar{R}\bar{L}\, T(\id_n)$
is shown in Figure \ref{fig.lem.T.id}, left, where the red square is
an identity of size $n-6$.  This configuration has a $T$-structure and
its pre-image w.r.t.\ $T$ is the configuration shown in Figure
\ref{fig.lem.T.id}, middle. Finally, by applying the algorithm of
standardization to $\tau$, we see that
$\tau_s=\bar{R}^2\bar{L}R^3L(\tau)$ is the standard permutation shown
in Figure \ref{fig.lem.T.id}, right, where the red square is an
identity of size~$n-7$. When this block is non-empty, this
configuration has $\tau_s(1)=1$ and $\tau_s(n)=n$, but it is not
$\id_n$, from which we conclude.
\qed

\begin{figure}[t!!]
\[
\id_n
\rightarrow
\text{\begin{tikzpicture}[scale=0.38]
\permutation{7,1,6,10,4,8,9,2,3,5};
\draw[red,thick] (6,8) rectangle (8,10);
\draw[blue,thick] (1+.5,7+.5) -- (3+.5,6+.5);
\draw[blue,thick] (2+.5,1+.5) -- (2+.5,7);
\end{tikzpicture}}
\rightarrow
\text{\begin{tikzpicture}[scale=0.38]
\permutation{5,8,3,6,7,1,2,4};
\draw[red,thick] (4,6) rectangle (6,8);
\draw [green!50!black!100, very thick] (1.5,5.5)--(5.5,7.5)--(3.5,3.5)--(8.5,4.5)--(6.5,1.5);
\end{tikzpicture}}
\rightarrow
\text{\begin{tikzpicture}[scale=0.38]
\permutation{1,2,3,6,4,5,7,8};
\draw[red,thick] (4,6) rectangle (5,7);
\end{tikzpicture}}
\neq \id_n
\]
\makebox[0pt][l]{\setlength{\unitlength}{5pt}
\begin{picture}(1,0)(-5,0)
\put(0,0){\vector(1,3){1}}
\put(1,0){$\bar{R}L^3RL\bar{R}\bar{L}R\bar{L}RL^2\bar{R}\bar{L}\, T$}
\end{picture}}
\makebox[0pt][l]{\setlength{\unitlength}{5pt}
\begin{picture}(1,0)(-30.5,0)
\put(0,0){\vector(0,1){3}}
\put(1,0){$T^{-1}$}
\end{picture}}
\makebox[0pt][l]{\setlength{\unitlength}{5pt}
\begin{picture}(1,0)(-50.5,0)
\put(0,0){\vector(0,1){3}}
\put(1,0){$\bar{R}^2\bar{L}R^3L$}
\end{picture}}
\makebox[0pt][r]{\setlength{\unitlength}{5pt}
\begin{picture}(1,0)(-69.5,0)
\put(0,1.5){if $n\geq 8$}
\end{picture}}
\caption{\label{fig.lem.T.id}The three configurations involved in the
  proof of Lemma~\ref{lem.T.id}.}
\end{figure}

\begin{lemma}[Fusion lemma for $(T,\tree)$]
\label{lem.T.tree}
Let $n \geq 9$. There exists some non-exceptional class $C_n$ such that $\bar{T}(C_n)=\bar{T}(\tree_n)$.
\end{lemma}

\begin{figure}[t!!]
\[
\idp_n \rightarrow
\text{\begin{tikzpicture}[scale=0.38]
\permutation{10,1,9,6,7,11,3,8,2,4,5};
\draw[red,thick] (4,6) rectangle (6,8);
\draw[blue,thick] (1+.5,10+.5) -- (3+.5,9+.5);
\draw[blue,thick] (2+.5,1+.5) -- (2+.5,10);
\end{tikzpicture}}
\rightarrow
\text{\begin{tikzpicture}[scale=0.38]
\permutation{8,5,6,9,2,7,1,3,4};
\draw [green!50!black!100, very thick] (4.5,9.5)--(2.5,5.5) --(6.5,7.5)--(5.5,2.5)--(9.5,4.5)--(7.5,1.5);
\draw[red,thick] (2,5) rectangle (4,7);
\end{tikzpicture}}
\rightarrow
\text{\begin{tikzpicture}[scale=0.38]
\permutation{1,8,2,3,6,4,5,7,9};
\draw[red,thick] (5,6) rectangle (6,7);
\end{tikzpicture}}
\neq \id'
\]
\makebox[0pt][l]{\setlength{\unitlength}{5pt}
\begin{picture}(1,0)(-2.5,0)
\put(0,0){\vector(1,3){1}}
\put(1,0){$\bar{R}LR\bar{L}^2\bar{R}\bar{L}^3RL^2RL\bar{R}\bar{L} \, T$}
\end{picture}}
\makebox[0pt][l]{\setlength{\unitlength}{5pt}
\begin{picture}(1,0)(-30,0)
\put(0,0){\vector(0,1){3}}
\put(1,0){$T^{-1}$}
\end{picture}}
\makebox[0pt][l]{\setlength{\unitlength}{5pt}
\begin{picture}(1,0)(-52,0)
\put(0,0){\vector(0,1){3}}
\put(1,0){$\bar{R}^2\bar{L}^2\bar{R}L^2R$}
\end{picture}}
\makebox[0pt][r]{\setlength{\unitlength}{5pt}
\begin{picture}(1,0)(-73,0)
\put(0,1.5){if $n\geq 9$}
\end{picture}}
\caption{\label{fig.lem.T.tree}The three configurations involved in the
  proof of Lemma~\ref{lem.T.tree}.}
\end{figure}
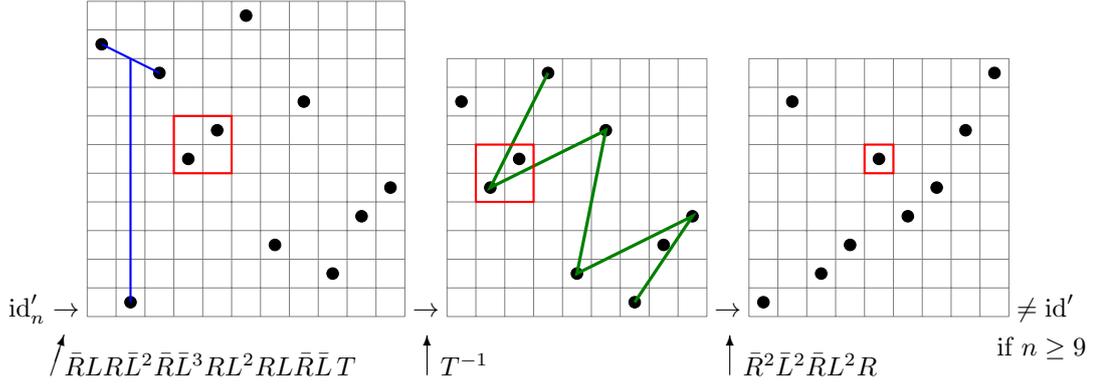

\proof
We searched a sequence such that $T(\idp_n)$ is transformed into a
$\s=T(\tau)$ where $\tau\notin \tree_n$. 
The configuration
$\bar{R}LR\bar{L}^2\bar{R}\bar{L}^3RL^2RL\bar{R}\bar{L} \, T(\idp_n)$
is shown in Figure \ref{fig.lem.T.tree}, left, where the red square is
an identity of size $n-7$.  This configuration has a $T$-structure and
its pre-image w.r.t.\ $T$ is the configuration shown in Figure
\ref{fig.lem.T.tree}, middle. Finally, by applying the algorithm of
standardization to $\tau$ we see that
$\tau_s=\bar{R}^2\bar{L}^2\bar{R}L^2R(\tau)$ is the standard
permutation shown in Figure \ref{fig.lem.T.tree}, right, where the red
square is an identity of size $n-8$. When this block is non-empty,
this configuration has $\tau_s(1)=1$ and $\tau_s(n)=n$, but it is not
$\id'_n$, from which we conclude.  \qed

\subsection{Fusion lemmas for $q_1$ and $q_2$}
\label{ssec.tech2}

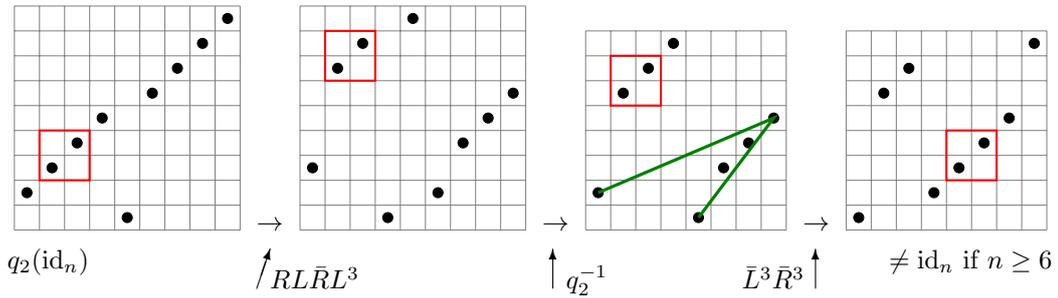
\begin{figure}
\[
\text{
\begin{tikzpicture}[scale=0.33]
\permutation{2,3,4,5,1,6,7,8,9}
\draw[red,thick] (2,3) rectangle (4,5);
\end{tikzpicture}
}
\rightarrow
\text{
\begin{tikzpicture}[scale=0.33]
\permutation{3,7,8,1,9,2,4,5,6}
\draw[red,thick] (2,7) rectangle (4,9);
\end{tikzpicture}
}
\rightarrow
\text{
\begin{tikzpicture}[scale=0.33]
\permutation{2,6,7,8,1,3,4,5}
\draw [green!50!black!100, very thick] (1.5,2.5)--(8.5,5.5)--(5.5,1.5);
\draw[red,thick] (2,6) rectangle (4,8);
\end{tikzpicture}
}
\rightarrow
\text{
\begin{tikzpicture}[scale=0.33]
\permutation{1,6,7,2,3,4,5,8}
\draw[red,thick] (2+3,6-3) rectangle (4+3,8-3);
\end{tikzpicture}
}
\]
\makebox[0pt][l]{\setlength{\unitlength}{5pt}
\begin{picture}(1,0)(0,0)
\put(0,1.5){$q_2 (\id_n)$}
\end{picture}}
\makebox[0pt][l]{\setlength{\unitlength}{5pt}
\begin{picture}(1,0)(-18,0)
\put(0,0){\vector(1,3){1}}
\put(1,0){$RL\bar{R}L^3$}
\end{picture}}
\makebox[0pt][l]{\setlength{\unitlength}{5pt}
\begin{picture}(1,0)(-39.5,0)
\put(0,0){\vector(0,1){3}}
\put(1,0){$q_2^{-1}$}
\end{picture}}
\makebox[0pt][l]{\setlength{\unitlength}{5pt}
\begin{picture}(1,0)(-57,0)
\put(1.5,0){\vector(0,1){3}}
\put(-4,0){$\bar{L}^3\bar{R}^3$}
\end{picture}}
\makebox[0pt][r]{\setlength{\unitlength}{5pt}
\begin{picture}(1,0)(-65,0)
\put(0,1.5){$\neq \id_n$ if $n\geq 6$}
\end{picture}}
\caption{\label{fig.lem.q2.id}The four configurations involved in the
  proof of Lemma~\ref{lem.q2.id}. Recall that, whenever 
  $\s(1) \neq 1$, $q_2^{-1} \s$ is obtained by removing the entry of
  $\s$ in the bottom-most row.}
\end{figure}

In this section we shall see that there is no proliferation of classes
due to the action of $q_1$ and $q_2$ on exceptional classes of large
size $\Id_n$ and $\Id'_n$. This makes, in principle, $2 \times 2=4$
cases. However, we know from Lemma \ref{lem.q1.id} that
$\bar{q}_1(\Id_n)=\tree_{n+1}$, while neither $q_1$ nor $q_2$ apply to
$\Id'_n$, because these classes have rank $1$, so this rules out three
cases out of four, and we just need the following:

\begin{lemma}[Fusion lemma for $(q_2, \Id)$]
\label{lem.q2.id}
Let $n \geq 6$. There exists some non-exceptional class $B_n$ such that
$\bar{q}_2(B_n)=\bar{q}_2(\Id_n)$.
\end{lemma}
\proof We searched a sequence such that $q_2(\id_n)$ is transformed
into a $\s=q_2(\tau)$ where $\tau\notin \Id_n$.  The configuration
$RL\bar{R}L^3\,q_2(\id_n)$ is shown in Figure \ref{fig.lem.q2.id}
(second image), where the red square is an identity of size $n-6$ (so
it is possibly empty).  Its pre-image w.r.t.\ $q_2$ is the
configuration shown in Figure \ref{fig.lem.q2.id} (third). Then, it
can be easily verified that for $n\geq 6$ the standard permutation
$\tau_s=\bar{L}^3\bar{R}^3(\tau)$ shown in Figure \ref{fig.lem.q2.id}
(fourth), is neither $\id_n$ nor $\id'_n$.\qed

However, we are not done still. We need to deal with the problem,
anticipated in Section \ref{ssec.q_op}, that the theorem on the
pullback of $\bar{q}_2$, Thm.~\ref{thm.q2}, is less precise than
what we need, as it does not rule out the possibility that classes of
different sign remain non-connected after the application of $q_2$
(which sets the sign to zero), in which case we would have two classes
with the same invariant $(\lambda,2,0)$ one coming from a sign $s=+1$,
the other from $s=-1$. For this we have a trivial preparatory
observation:

\begin{lemma}
Let $\lambda \cup \{2\}$ be an integer partition, with an even number of even
parts. The two following facts are equivalent:
\begin{enumerate}
\item There exists a value of $j$ such that 
  $(\lambda \setminus j) \cup \{j+1\}$ has a positive even number of even
  parts.
\item $\lambda$ does not consist of a single (even) part.
\end{enumerate}
\end{lemma}
\proof Clearly $\lambda$ has at least one even part. If it has also
some odd part, its value is a valid candidate for $j$. If it has three
or more even parts, all parts of $\lambda$ are valid candidates. On
the other side, if $\lambda=\{2k\}$, then $j=2k$ is the only
candidate, and results in
$(\lambda \setminus j) \cup \{j+1\} = \{2k+1\}$.
\qed

This is used in the following (again, we use $\lambda(j) \equiv
\lambda \setminus j$)
\begin{lemma}[Fusion of sign for $q_2$, ordinary case]
\label{lem.q2_sign_0}
Let $C$ be a class with invariant $(\lambda,2,0)$. Then, if
$\lambda$ does not consists of a single cycle of even length,
there exists a value of $j$, and a class $B_j$ with invariant
$(\lambda(j),j+1,0)$, such that $\bar{q}_2(B_j)=C$.
\end{lemma}

\proof By Theorem \ref{thm.q2} we know that, given a class with invariant
$(\lambda,2,0)$, then, for all $j \in \lam$, there exists a class
$B_j$ with cycle invariant $(\lambda(j),j+1)$, such that
$\bar{q}_2(B_j) = C$. By the previous lemma, if $\lambda \neq \{2k\}$,
there exists one such $j$ with $(\lambda \setminus j) \cup \{j+1\}$
containing even parts, thus, by Lemma \ref{lem.induSign}, the
corresponding class must also have $s(B_j)=0$, as claimed.  \qed

\noindent
In light of this lemma, the only possible instance of proliferation of
disjoint classes which is still open, for which we need a separate
argument, corresponds to the very special partitions $\lambda$
escaping the characterisation above, i.e., referring to notations as
in Lemma \ref{lem.q2_sign_0}, the $\lambda$'s consisting of a single
cycle of even length.  These special cases are analysed in the
following Lemma \ref{lem.q2.maxrank}. We observe the trivial fact

\begin{lemma}\label{lem.two_max_classes}
Define $C'_k=\bar{T}^{k-1}(\Id_4)$ and
$C''_k=\bar{T}^{k-2}(\Id_6)$. Then $C'_k$ has invariant
$(\emptyset,r=2k+1,-1)$ and $C''_k$ has invariant
$(\emptyset,r=2k+1,+1)$.
\end{lemma}
\proof
By corollary \ref{cor.id_sign} $\Id_4$ has invariant
$(\emptyset,3,+1)$ while $\Id_6$ has invariant
$(\emptyset,5,-1)$. Then, the statement follows from the properties of
$\bar{T}$ implied by Corollary~\ref{cor.Tsign}.
\qed

Then we have
\begin{lemma}[Fusion of sign for $q_2$, special case]
\label{lem.q2.maxrank}
For $k \geq 3$, define $C'_{k}=\bar{T}^{k-1}(\Id_4)$ and
$C''_{k}=\bar{T}^{k-2}(\Id_6)$, and similarly $c'_{k}=T^{k-1}(\id_4)$
and $c''_{k}=T^{k-2}(\id_6)$.  We have $q_2(c'_{k}) \sim
q_2(c''_{k})$, i.e.\ $\bar{q}_2(C'_{k})=\bar{q}_2(C''_{k})$.
\end{lemma}
\noindent
We split the proof in a few lemmas. We start with a definition.
\begin{definition}
\label{def.Vblock}
For $k \geq 3$ and $0 \leq h \leq k-3$, call $C_{k,h}$ the following
configuration
\[
\text{\begin{tikzpicture}[scale=0.3]
\permutation{4,12,1,7,2,6,3,5,13,21,14,20,15,19,16,18,8,10,17,9,11};
\draw [green!50!black!100] (4.2,8.2)--(4.4,8.6)--(9.4,6.1)--(9.2,5.7);
\draw [green!50!black!100] (4.2+6,8.2+14)--(4.4+6,8.6+14)--(9.4+8,6.1+13)--(9.2+8,5.7+13);
\draw (3+.5,1+.5) [fill,green!50!black!100] circle (.2);
\draw (5+.5,2+.5) [fill,green!50!black!100] circle (.2);
\draw (7+.5,3+.5) [fill,green!50!black!100] circle (.2);
\draw (8+.5,5+.5) [fill,green!50!black!100] circle (.2);
\draw (6+.5,6+.5) [fill,green!50!black!100] circle (.2);
\draw (4+.5,7+.5) [fill,green!50!black!100] circle (.2);
\draw (9+.5,13+.5) [fill,green!50!black!100] circle (.2);
\draw (11+.5,14+.5) [fill,green!50!black!100] circle (.2);
\draw (13+.5,15+.5) [fill,green!50!black!100] circle (.2);
\draw (15+.5,16+.5) [fill,green!50!black!100] circle (.2);
\draw (16+.5,18+.5) [fill,green!50!black!100] circle (.2);
\draw (14+.5,19+.5) [fill,green!50!black!100] circle (.2);
\draw (12+.5,20+.5) [fill,green!50!black!100] circle (.2);
\draw (10+.5,21+.5) [fill,green!50!black!100] circle (.2);
\node at (8,8.7) {$\bm{k-h-2}$};
\node at (16.5,21.4) {$\bm{h-1}$};
\end{tikzpicture}}
\]
\end{definition}
\noindent
We have
\begin{proposition}
\label{prop.RLRLL}
\be
\label{eq.inprop.RLRLL}
\text{\raisebox{-30pt}{
\begin{tikzpicture}[scale=0.4]
\permutationtwo{2,5,4,1,3}{2,3,4}{1,2,4,5};
\draw (3+.5,4+.5) [fill,green!50!black!100] circle (.2);
\draw (4+.5,1+.5) [fill,green!50!black!100] circle (.2);
\draw (5+.5,3+.5) [fill,green!50!black!100] circle (.2);
\node at (2+1,5+1) {$\bm{a}$};
\node at (5+1,5+1) {$\bm{b}$};
\node at (2+1,4+1) {$\bm{c}$};
\node at (5+1,4+1) {$\bm{d}$};
\node at (2+1,0+1) {$\bm{e}$};
\node at (5+1,0+1) {$\bm{f}$};
\end{tikzpicture}
}}
\xrightarrow{\bar{R}LR\bar{L}\bar{L}}
\text{\raisebox{-30pt}{
\begin{tikzpicture}[scale=0.4]
\permutationtwo{1,3,2,4,5}{2,4,6}{1,2,4,5};
\draw (3+.5,2+.5) [fill,green!50!black!100] circle (.2);
\draw (4+.5,4+.5) [fill,green!50!black!100] circle (.2);
\draw (5+.5,5+.5) [fill,green!50!black!100] circle (.2);
\node at (2+1,4+1) {$\bm{a}$};
\node at (5+1,4+1) {$\bm{b}$};
\node at (2+1,2+1) {$\bm{c}$};
\node at (5+1,2+1) {$\bm{d}$};
\node at (2+1,0+1) {$\bm{e}$};
\node at (5+1,0+1) {$\bm{f}$};
\end{tikzpicture}
}}
\ee
\end{proposition}
\noindent
Note, the sequence above does \emph{not} require to be boosted. In
other words, the exponents in the boosted sequence are $\alpha_j \in
\{\pm 1\}$.\footnote{Also note that the green bullets on the
  right-hand-side are not the image of those on the left, they are
  just coloured for use in the following statement.}
\begin{corollary}
\label{cor.ofprop.RLRLL}
For $k \geq 4$ and $0 \leq h < k-3$,
\[
C_{k,h}
\xrightarrow{\bar{R}LR\bar{L}\bar{L}}
C_{k,h+1}
\]
\end{corollary}
\noindent
It suffices to identify the blocks $a$, \ldots, $f$ in the boosted
dynamics of Proposition \ref{prop.RLRLL} as follows:
\[
\text{\begin{tikzpicture}[scale=0.3]
\permutation{4,12,1,7,2,6,3,5,13,21,14,20,15,19,16,18,8,10,17,9,11};
\draw (3+.5,1+.5) [fill,green!50!black!100] circle (.2);
\draw (5+.5,2+.5) [fill,green!50!black!100] circle (.2);
\draw (7+.5,3+.5) [fill,green!50!black!100] circle (.2);
\draw (8+.5,5+.5) [fill,green!50!black!100] circle (.2);
\draw (6+.5,6+.5) [fill,green!50!black!100] circle (.2);
\draw (4+.5,7+.5) [fill,green!50!black!100] circle (.2);
\draw (9+.5,13+.5) [fill,green!50!black!100] circle (.2);
\draw (11+.5,14+.5) [fill,green!50!black!100] circle (.2);
\draw (13+.5,15+.5) [fill,green!50!black!100] circle (.2);
\draw (15+.5,16+.5) [fill,green!50!black!100] circle (.2);
\draw (16+.5,18+.5) [fill,green!50!black!100] circle (.2);
\draw (14+.5,19+.5) [fill,green!50!black!100] circle (.2);
\draw (12+.5,20+.5) [fill,green!50!black!100] circle (.2);
\draw (10+.5,21+.5) [fill,green!50!black!100] circle (.2);
\draw [red,thick] (3,1) rectangle (6,3);
\draw [red,thick] (3,7) rectangle (6,12);
\draw [red,thick] (3,13) rectangle (6,22);
\draw [red,thick] (9,1) rectangle (22,3);
\draw [red,thick] (9,7) rectangle (22,12);
\draw [red,thick] (9,13) rectangle (22,22);
\end{tikzpicture}}
\]
\qed

\begin{figure}[tb!]
\begin{center}
\setlength{\unitlength}{180pt}
\begin{picture}(2.2,2.3)(-.1,-.2)
\put(-.1,1){%
\begin{tikzpicture}[scale=0.3]
\permutation{16,2,15,3,14,4,13,5,12,6,11,7,10,8,9,17,1,18,19,20,21};
\end{tikzpicture}
}
\put(1.1,1){%
\begin{tikzpicture}[scale=0.3]
\qquad
\permutation{7,18,1,13,2,12,3,11,4,10,5,9,6,8,19,21,14,16,20,15,17};
\draw (3+.5,1+.5) [fill,green!50!black!100] circle (.2);
\draw (5+.5,2+.5) [fill,green!50!black!100] circle (.2);
\draw (7+.5,3+.5) [fill,green!50!black!100] circle (.2);
\draw (9+.5,4+.5) [fill,green!50!black!100] circle (.2);
\draw (11+.5,5+.5) [fill,green!50!black!100] circle (.2);
\draw (13+.5,6+.5) [fill,green!50!black!100] circle (.2);
\draw (14+.5,8+.5) [fill,green!50!black!100] circle (.2);
\draw (12+.5,9+.5) [fill,green!50!black!100] circle (.2);
\draw (10+.5,10+.5) [fill,green!50!black!100] circle (.2);
\draw (8+.5,11+.5) [fill,green!50!black!100] circle (.2);
\draw (6+.5,12+.5) [fill,green!50!black!100] circle (.2);
\draw (4+.5,13+.5) [fill,green!50!black!100] circle (.2);
\draw (15+.5,19+.5) [fill,green!50!black!100] circle (.2);
\draw (16+.5,21+.5) [fill,green!50!black!100] circle (.2);
\end{tikzpicture}
}
\put(1,-.23){\makebox[0pt][c]{$\xrightarrow{
\bar{L}
R^2
L
\bar{R}^2
\bar{L}^4
\bar{R}
\bar{L}
}$}}
\put(1,2){\makebox[0pt][c]{$\xrightarrow{
\bar{R}
\bar{L}
R^2
L
\bar{R}^3
\bar{L}
\bar{R}
\bar{L}
R
\bar{L}
\bar{R}
\bar{L}^2
R^2}$}}
\put(1,.92){\makebox[0pt][c]{\raisebox{-4pt}{$\swarrow$}
$(\bar{R}LR\bar{L}\bar{L})^{n-2}$, by
iterated application of Corollary \ref{cor.ofprop.RLRLL}
\raisebox{+4pt}{$\swarrow$}}}
\put(-.1,-.15){%
\begin{tikzpicture}[scale=0.3]
\permutation{2,8,1,3,9,21,10,20,11,19,12,18,13,17,14,16,4,6,15,5,7};
\draw (3+2.5,1+8.5) [fill,green!50!black!100] circle (.2);
\draw (5+2.5,2+8.5) [fill,green!50!black!100] circle (.2);
\draw (7+2.5,3+8.5) [fill,green!50!black!100] circle (.2);
\draw (9+2.5,4+8.5) [fill,green!50!black!100] circle (.2);
\draw (11+2.5,5+8.5) [fill,green!50!black!100] circle (.2);
\draw (13+2.5,6+8.5) [fill,green!50!black!100] circle (.2);
\draw (14+2.5,8+8.5) [fill,green!50!black!100] circle (.2);
\draw (12+2.5,9+8.5) [fill,green!50!black!100] circle (.2);
\draw (10+2.5,10+8.5) [fill,green!50!black!100] circle (.2);
\draw (8+2.5,11+8.5) [fill,green!50!black!100] circle (.2);
\draw (6+2.5,12+8.5) [fill,green!50!black!100] circle (.2);
\draw (4+2.5,13+8.5) [fill,green!50!black!100] circle (.2);
\draw (3+.5,1+.5) [fill,green!50!black!100] circle (.2);
\draw (4+.5,3+.5) [fill,green!50!black!100] circle (.2);
\end{tikzpicture}
}
\put(1.1,-.15){%
\begin{tikzpicture}[scale=0.3]
\permutation{1,18,2,17,3,16,4,15,5,14,6,13,7,12,8,11,9,10,19,20,21};
\end{tikzpicture}
}
\end{picture}
\end{center}
\caption{\label{figure.T7toT8_full} Illustration of
  Lemma~\ref{lem.q2.maxrank}.  The configuration on the top left is
  $q_2 T^n(id_6)$, and the one on the bottom right is $\bar{L} q_2
  T^{n+1}(id_4)$, for $n=7$.}
\end{figure}
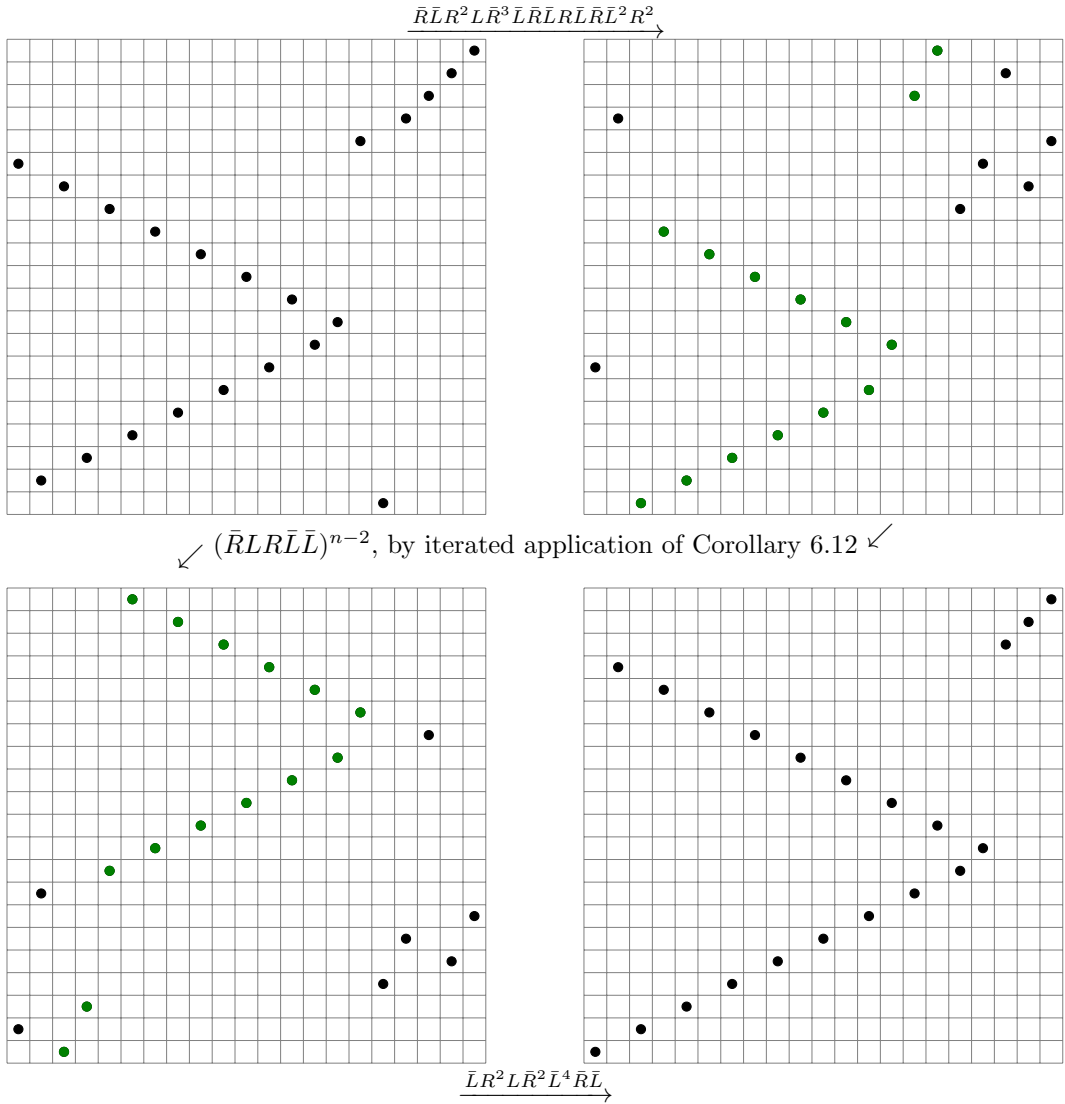

\noindent
Then we have the technical verifications
\begin{proposition}
\label{prop.1inq2tid}
\[
\bar{R}
\bar{L}
R^2
L
\bar{R}^3
\bar{L}
\bar{R}
\bar{L}
R
\bar{L}
\bar{R}
\bar{L}^2
R^2
q_2(c'_k) 
=
C_{k,0}
\ef.
\]
\end{proposition}
\begin{proposition}
\label{prop.2inq2tid}
\[
\bar{L}
R^2
L
\bar{R}^2
\bar{L}^4
\bar{R}
\bar{L}
\,
C_{k,k-3}
=
q_2(c''_k) 
\ef.
\]
\end{proposition}
\noindent
(the configurations involved in these two propositions are shown in
Figure \ref{figure.T7toT8_full}),
which, together with Corollary \ref{cor.ofprop.RLRLL},
imply Lemma~\ref{lem.q2.maxrank}.

\subsection{Classification of Rauzy classes}
\label{ssec.mainthproof}

We are now ready to collect the large number of lemmas established so
far into a proof of Theorem~\ref{thm.Main_theorem}. 

\begin{figure}
\begin{center}
\begin{tabular}{cc}
\begin{tikzpicture}[scale=0.4]
\useasboundingbox (0,0) rectangle (12,11);
\node at (4.5,10) {Operator $T$ ($n \to n+2$)};
\draw[->,blue,thick] (6,8) -- (6,4);
\draw[blue,thick] (6,4) -- (6,0);
\draw[->,blue,thick] (0,8) -- (0,4);
\draw[blue,thick] (0,4) -- (0,0);
\draw[->,blue,thick] (3,4) -- (3,2);
\draw[blue,thick] (3,2) -- (3,0);
\draw[blue,thick] plot (3,8) .. controls (4.6,5.5) and (4.6,2.5) .. (3,0);
\draw[->,blue,thick] (4.2,4.1) -- (4.2,4);
\draw[->,blue,thick,dashed] (9,8) -- (9,5);
\draw[->,blue,thick,dashed] (9,0) -- (9,-3);
\node at (12,5.5) {\begin{minipage}{3cm}
\begin{center}
not important\\
by Lemma\\
\ref{lem.T.tree}
\end{center}
\end{minipage}};
\node at (12,-2.5) {\begin{minipage}{3cm}
\begin{center}
not important\\
by Lemma\\
\ref{lem.T.id}
\end{center}
\end{minipage}};
\draw[blue,thick] plot (0,0) .. controls (0-2.4,0-3.6) and (0+2.4,0-3.6) .. (0,0);
\draw[->,blue,thick,dashed] (0+0,-2.7) -- (0+.1,-2.7);
\draw[blue,thick] plot (3,0) .. controls (3-2.4,0-3.6) and (3+2.4,0-3.6) .. (3,0);
\draw[->,blue,thick,dashed] (3+0,-2.7) -- (3+.1,-2.7);
\draw[blue,thick] plot (6,0) .. controls (6-2.4,0-3.6) and (6+2.4,0-3.6) .. (6,0);
\draw[->,blue,thick,dashed] (6+0,-2.7) -- (6+.1,-2.7);
\draw[fill,blue!15!white!100] (0,0) circle (.9);
\draw[blue,thick] (0,0) circle (.9);
\node at          (0,0) {$o^-_{\geq 3}$};
\draw[fill,blue!15!white!100] (3,0) circle (.9);
\draw[blue,thick] (3,0) circle (.9);
\node at          (3,0) {$e^{\phantom{i}}_{\geq 3}$};
\draw[fill,blue!15!white!100] (6,0) circle (.9);
\draw[blue,thick] (6,0) circle (.9);
\node at          (6,0) {$o^+_{\geq 3}$};
\draw[fill,blue!15!white!100] (3,4) circle (.9);
\draw[blue,thick] (3,4) circle (.9);
\node at          (3,4) {$e^{\phantom{i}}_{2}$};
\draw[fill,blue!15!white!100] (0,8) circle (.9);
\draw[blue,thick] (0,8) circle (.9);
\node at          (0,8) {$o^-_{1}$};
\draw[fill,blue!15!white!100] (3,8) circle (.9);
\draw[blue,thick] (3,8) circle (.9);
\node at          (3,8) {$e^{\phantom{i}}_{1}$};
\draw[fill,blue!15!white!100] (6,8) circle (.9);
\draw[blue,thick] (6,8) circle (.9);
\node at          (6,8) {$o^+_{1}$};
\draw[fill,blue!15!white!100] (9,8) circle (.9);
\draw[blue,thick] (9,8) circle (.9);
\node at          (9,8) {$\tree$};
\draw[fill,blue!15!white!100] (9,0) circle (.9);
\draw[blue,thick] (9,0) circle (.9);
\node at          (9,0) {$\Id$};
\end{tikzpicture}
&
\begin{tikzpicture}[scale=0.4]
\useasboundingbox (0,0) rectangle (12,11);
\node at (4.5,10) {Operator $q_1$ ($n \to n+1$)};
\draw[->,blue,thick] (6,0) -- (6,4);
\draw[blue,thick] (6,4) -- (6,8);
\draw[->,blue,thick] (0,0) -- (0,4);
\draw[blue,thick] (0,4) -- (0,8);
\draw[->,blue,thick] (3,4) -- (3,6);
\draw[blue,thick] (3,6) -- (3,8);
\draw[blue,thick] plot (3,8) .. controls (4.6,5.5) and (4.6,2.5) .. (3,0);
\draw[->,blue,thick] (4.2,3.9) -- (4.2,4);
\draw[->,blue,thick] (9,0) -- (9,4);
\draw[blue,thick] (9,4) -- (9,8);
\node at (10.8,3.5) {\begin{minipage}{3cm}
\begin{center}
by\\
Lemma\\
\ref{lem.q1.id}
\end{center}
\end{minipage}};
\draw[fill,blue!15!white!100] (0,0) circle (.9);
\draw[blue,thick] (0,0) circle (.9);
\node at          (0,0) {$o^-_{\geq 3}$};
\draw[fill,blue!15!white!100] (3,0) circle (.9);
\draw[blue,thick] (3,0) circle (.9);
\node at          (3,0) {$e^{\phantom{i}}_{\geq 3}$};
\draw[fill,blue!15!white!100] (6,0) circle (.9);
\draw[blue,thick] (6,0) circle (.9);
\node at          (6,0) {$o^+_{\geq 3}$};
\draw[fill,blue!15!white!100] (3,4) circle (.9);
\draw[blue,thick] (3,4) circle (.9);
\node at          (3,4) {$e^{\phantom{i}}_{2}$};
\draw[fill,blue!15!white!100] (0,8) circle (.9);
\draw[blue,thick] (0,8) circle (.9);
\node at          (0,8) {$o^-_{1}$};
\draw[fill,blue!15!white!100] (3,8) circle (.9);
\draw[blue,thick] (3,8) circle (.9);
\node at          (3,8) {$e^{\phantom{i}}_{1}$};
\draw[fill,blue!15!white!100] (6,8) circle (.9);
\draw[blue,thick] (6,8) circle (.9);
\node at          (6,8) {$o^+_{1}$};
\draw[fill,blue!15!white!100] (9,8) circle (.9);
\draw[blue,thick] (9,8) circle (.9);
\node at          (9,8) {$\tree$};
\draw[fill,blue!15!white!100] (9,0) circle (.9);
\draw[blue,thick] (9,0) circle (.9);
\node at          (9,0) {$\Id$};
\end{tikzpicture}
\\
\raisebox{50pt}{\begin{tabular}{r|l}
$o^+_r$ & classes with $r$, $\lambda_i$ all odd, and $s=+1$
\\
$o^-_r$ & classes with $r$, $\lambda_i$ all odd, and $s=-1$
\\
$e_r$ & classes with some $r$, $\lambda_i$ even, and $s=0$
\\
$\Id$, $\tree$ & special classes
\end{tabular}}
&
\begin{tikzpicture}[scale=0.4]
\useasboundingbox (0,-1) rectangle (12,13);
\node at (4.5,10) {Operator $q_2$ ($n \to n+1$)};
\draw[->,blue,thick] (6,0) -- (4.5,2);
\draw[blue,thick] (3,4) -- (4.5,2);
\draw[->,blue,thick] (0,0) -- (1.5,2);
\draw[blue,thick] (3,4) -- (1.5,2);
\draw[->,blue,thick] (3,0) -- (3,2);
\draw[blue,thick] (3,4) -- (3,2);
\draw[->,blue,thick,dashed] (9,0) -- (9,3);
\node at (9.5,5) {\begin{minipage}{3cm}
\begin{center}
not important\\
by Lemma\\
\ref{lem.q2.id}
\end{center}
\end{minipage}};
\draw[fill,blue!15!white!100] (0,0) circle (.9);
\draw[blue,thick] (0,0) circle (.9);
\node at          (0,0) {$o^-_{\geq 3}$};
\draw[fill,blue!15!white!100] (3,0) circle (.9);
\draw[blue,thick] (3,0) circle (.9);
\node at          (3,0) {$e^{\phantom{i}}_{\geq 3}$};
\draw[fill,blue!15!white!100] (6,0) circle (.9);
\draw[blue,thick] (6,0) circle (.9);
\node at          (6,0) {$o^+_{\geq 3}$};
\draw[fill,blue!15!white!100] (3,4) circle (.9);
\draw[blue,thick] (3,4) circle (.9);
\node at          (3,4) {$e^{\phantom{i}}_{2}$};
\draw[fill,blue!15!white!100] (0,8) circle (.9);
\draw[blue,thick] (0,8) circle (.9);
\node at          (0,8) {$o^-_{1}$};
\draw[fill,blue!15!white!100] (3,8) circle (.9);
\draw[blue,thick] (3,8) circle (.9);
\node at          (3,8) {$e^{\phantom{i}}_{1}$};
\draw[fill,blue!15!white!100] (6,8) circle (.9);
\draw[blue,thick] (6,8) circle (.9);
\node at          (6,8) {$o^+_{1}$};
\draw[fill,blue!15!white!100] (9,8) circle (.9);
\draw[blue,thick] (9,8) circle (.9);
\node at          (9,8) {$\tree$};
\draw[fill,blue!15!white!100] (9,0) circle (.9);
\draw[blue,thick] (9,0) circle (.9);
\node at          (9,0) {$\Id$};
\end{tikzpicture}
\end{tabular}
\end{center}
\caption{\label{fig.Induction_proof}Inductive
  scheme for the application of operators $T$, $q_1$ and $q_2$.}
\end{figure}
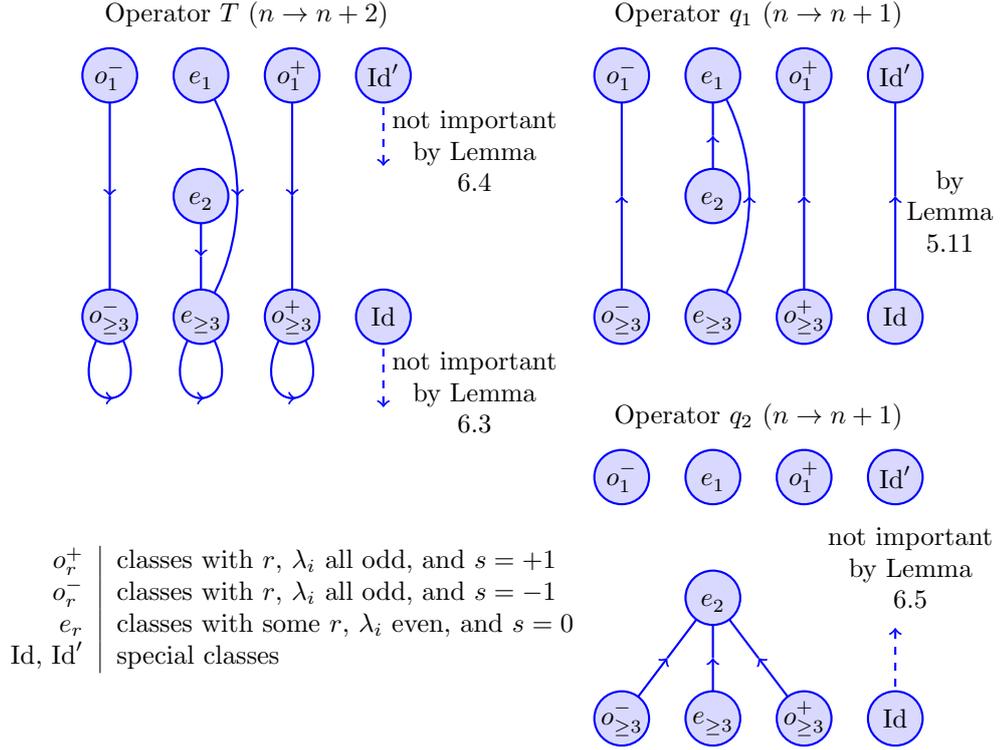

\pfof{Theorem~\ref{thm.Main_theorem}}
By the results of Appendix \ref{sec.excp_class}, we have a full understanding of the
two exceptional classes $\Id_n$ and $\tree_n$ (and, of course, we know
their invariant), so that we can concentrate solely on (primitive)
non-exceptional classes.

We proceed by induction. The theorem is established by explicit
investigation for classes up to $n=8$. Then, for the inductive step,
we suppose that non-exceptional classes of size $n-1$ and $n$ are
completely characterized by $(\lambda,r,s)$, in agreement with the
statement of Theorem~\ref{thm.Main_theorem}, and investigate classes
at size $n+1$.

Figure \ref{fig.Induction_proof}, in three copies for the three
surgery operators $\bar{T}$, $\bar{q}_1$ and $\bar{q}_2$, illustrates
a partition of the set of classes according to certain properties of
the cycle invariant.  This partition, despite containing `only' 5
blocks for non-exceptional cases, is fine enough so that classes
within the same block have a consistent behaviour w.r.t.\ all of our
three surgery operators.  These behaviours are represented through
arrows.  Recall that the operator $\bar{T}$ increases the size by $2$,
while operators $\bar{q}_{1,2}$ increase it by 1.  The pullback
theorems (Theorem \ref{thm.q1}, \ref{thm.q2} and
\ref{thm_T_surjectif}) of Section \ref{ssec.T_op} and \ref{ssec.q_op}
justify the fact that, in order to prove that there exists \emph{at
  least} one non-exceptional class per invariant, for the invariants
listed in the theorem statement, it suffices to observe that, in the
layer of largest size, all blocks of the partition have positive
in-degree (when the three copies are considered altogether). This is
steadily verified on the image.

For blocks which have in-degree higher than 1, we need our
fusion lemmas to conclude that there exists \emph{exactly} one
non-exceptional class per invariant, for the invariants listed in the
theorem statement.

In the following paragraph we recall which lemma justifies which arrow
of the diagram, and how the fusion lemmas are used.  We analyse the
arrows in some order, according to the rank of the image of the
operator.

We start with rank at least 3. By Theorem \ref{thm_T_surjectif},
$\bar{T}$ is a surjection from classes of size $n-1$ to classes of
size $n+1$ with rank at least 3 other than $\Id_{n+1}$. Notice however
that the classes of size $n-1$ include $\Id_{n-1}$ and $\tree_{n-1}$.
Nonetheless, by the `fusion' Lemmas \ref{lem.T.id} and
\ref{lem.T.tree}, there exist two non-exceptional classes, $C_1$ and
$C_2$, of size $n-1$, such that $\bar{T}(\Id_{n-1})=\bar{T}(C_1)$ and
$\bar{T}(\tree_{n-1})=\bar{T}(C_2)$, thus $\bar{T}$ provides a natural
bijection from the set of non-exceptional classes of size $n-1$ to the
set of non-exceptional classes of size $n+1$ with rank at least~3.

Now we pass to rank 1.  By (the easy) Lemma \ref{lem.q1.id}
$\bar{q}_1(\Id_n)=\tree_{n+1}$, furthermore $\bar{q}_1$ cannot be
applied to $\tree_n$ since $\tree_n$ has rank 1. Thus the case of
exceptional classes has been dealt with. By Theorem \ref{thm.q1}, if
there exist two classes $C_1$ and $C_2$ with invariant $(\lambda,1,s)$
of size $n+1$, then this would imply that, for some $i$ a part of
$\lambda$, there exist also two classes of size $n$ with invariant
$(\lambda(i),i,s)$.  So the induction step is proved for classes of
rank~1.

Finally, we analyse classes of rank 2.  Again the exceptional classes
are ruled out, since $\tree_n$ has rank 1, and, by the fusion Lemma
\ref{lem.q2.id}, there exists $C_1$ such that
$\bar{q}_2(\Id_n)=\bar{q}_2(C_1)$. So we are left only with the issue
of the pullback of $\bar{q}_2$, among non-exceptional classes.

Theorem \ref{thm.q2} tells us that for every $C$ with invariant
$(\lambda,2,0)$, for every $j\in\lambda$ there exists some class
$B_j$ with invariant $(\lambda(j),j+1)$ such that $\bar{q}_1(B_j)=C$,
so there are two cases:
\begin{itemize}
\item One of these $B_j$ has invariant $(\lambda(j),j+1,0)$. In this
  case, as was above with $q_1$, $C$ is the only class with invariant
  $(\lambda,2,0)$.
\item None of these $B_j$ has invariant $(\lambda(j),j+1,0)$. In this
  case it could be that all the classes $B_j^{+}$ with invariant
  $(\lambda(j),j+1,1)$ give a class $C^{+}$ with invariant
  $(\lambda,2,0)$ and all the classes $B_j^{-}$ with invariant
  $(\lambda(j),j+1,-1)$ give a class $C^{-}$, distinct from $C^+$,
  with invariant $(\lambda,2,0)$.  We need to rule out this
  possibility now, because this was not still excluded by the `weak'
  Lemma \ref{thm.q2}.

  Indeed, at this point we have all the elements to conclude that this
  does not happen: by Lemma \ref{lem.q2_sign_0}, if no $B_j$ has
  invariant $(\lambda(j),j+1,0)$ then $\lambda$ consists of a single
  cycle of even length (say, of length $2k$), so the two candidate
  classes $B^{+}$ and $B^{-}$ must have invariant $(\emptyset,2k+1,+1)$
  and $(\emptyset,2k+1,-1)$ respectively. By the induction and Lemma
\ref{lem.two_max_classes}
we know that
  $B^{+}=\bar{T}^{k-2}(\Id_6)$ and $B^{-}=\bar{T}^{k-1}(\Id_4)$, and
  by Lemma~\ref{lem.q2.maxrank} we know that for $n$ large enough
  $\bar{q}_2(B^{+})=\bar{q}_2(B^{-})$.
\end{itemize}
This proves the induction step for classes of rank 2, and allows to
conclude.  \qed

\noindent
Table \ref{tab.exinduct} illustrates a typical step of the induction,
at a size sufficiently large that all fusion lemmas are already in
place, and all typical situations do occur.

\begin{figure}
\begin{center}
\setlength{\unitlength}{50pt}
\begin{picture}(8,7.3)(-.4,-.63)
\put(-.4,-.33){\includegraphics[scale=1.24]{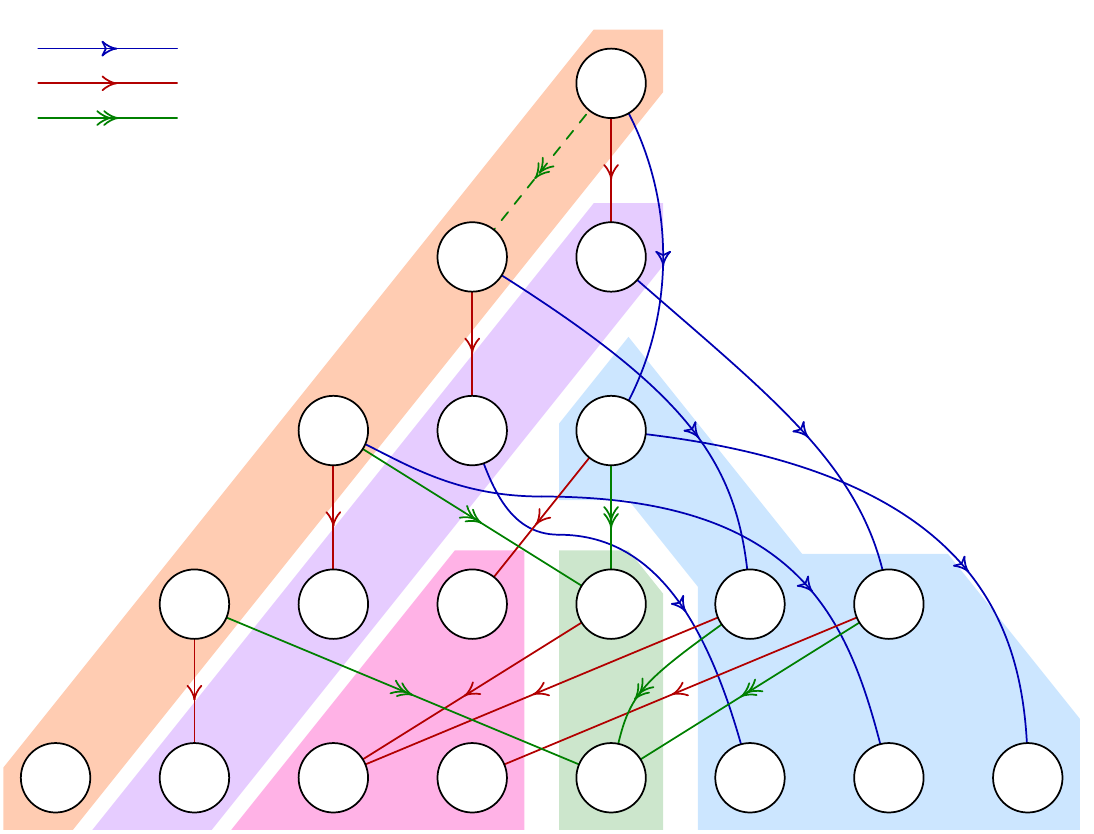}}
\put(1,5.25){$T$}
\put(1,5){$q_1$}
\put(1,4.75){$q_2$}
\put(2.5,-.55){\gostrC{rank 1}}
\put(4,-.55){\gostrC{rank 2}}
\put(6,-.55){\gostrC{rank $\geq 3$}}
\put(0,0.4){\gostrC{$\Id_8$}}
\put(1,1.65){\gostrC{$\Id_7$}}
\put(2,2.9){\gostrC{$\Id_6$}}
\put(3,4.15){\gostrC{$\Id_5$}}
\put(4,5.40){\gostrC{$\Id_4$}}
\put(1.2,0.4){\gostrC{$\Idp_8$}}
\put(2.2,1.65){\gostrC{$\Idp_7$}}
\put(3.2,2.9){\gostrC{$\Idp_6$}}
\put(4.2,4.15){\gostrC{$\Idp_5$}}
\put(0,0){\gostrC{$\emptyset|7+$}}
\put(1,0){\gostrC{$33|1+$}}
\put(2,0){\gostrC{$42|1$}}
\put(3,0){\gostrC{$33|1-$}}
\put(4,0){\gostrC{$32|2$}}
\put(5,0){\gostrC{$22|3$}}
\put(6,0){\gostrC{$\emptyset|7+$}}
\put(7,0){\gostrC{$\emptyset|7-$}}
\put(1,1.25){\gostrC{$3|3+$}}
\put(2,1.25){\gostrC{$5|1+$}}
\put(3,1.25){\gostrC{$5|1-$}}
\put(4,1.25){\gostrC{$4|2$}}
\put(5,1.25){\gostrC{$2|4$}}
\put(6,1.25){\gostrC{$3|3-$}}
\put(2,2.5){\gostrC{$\emptyset|5+$}}
\put(3,2.5){\gostrC{$22|1$}}
\put(4,2.5){\gostrC{$\emptyset|5-$}}
\put(3,3.75){\gostrC{$2|2$}}
\put(4,3.75){\gostrC{$3|1-$}}
\put(4,5){\gostrC{$\emptyset|3-$}}
\end{picture}
\end{center}
\caption{\label{fig.fulltree48}Full decomposition tree up to $n=8$.
Each balloon denotes a class, and contains its cycle and sign
invariant, in the form $\lam_1 \lam_2 \cdots \lam_{\ell} | r s$, with
$s \in \{+,-\}$ and omitted if zero.
}
\end{figure}

Figure \ref{fig.fulltree48} shows the top-most part 
($4 \leq n \leq 8$) of the full decomposition tree associated to the
action of $T$, $q_1$ and $q_2$ operators on the classes. At difference
with Table \ref{tab.exinduct}, this part is quite lacunary w.r.t.\ the
general pattern. The relevant fact is that, for $n>8$, there are no
more exceptions (cf.~the statements in Section \ref{ssec.tech} and
\ref{ssec.tech2}, which hold under the hypotheses $n \geq n_0$ for
certain values of $n_0$ all at most $8$), so that this picture
encompasses all of the small-size exceptional behaviour.

\begin{table}[p]
\noindent
\begin{tabular}{r|rcl|}
\multicolumn{1}{c}{}&
\multicolumn{1}{c}{$n=9,10$}&
\multicolumn{1}{c}{}&
\multicolumn{1}{c}{$n=11$}
\\
\cline{2-4}
&& new
& \classdrawlitX{1,2,3,4,5,6,7,8,9,10,11}{5}{5}{$\Id\ \oddclassX$}{-}{1023}{class_N11_Id_o_lam00000020000r5z.eps} \\
\cline{2-4}
&\classdrawlitX{1,2,3,4,5,7,8,9,6}{2}{6}{}{0}{31031}{class_N09__lam000100010r6z.eps}
& $\bar{T} \rightarrow$ & \classdrawlitX{1,2,3,4,5,6,7,9,10,11,8}{2}{8}{}{0}{2950101}{class_N11__lam00010000010r8z.eps}\\
&\classdrawlitX{1,2,3,4,5,6,9,7,8}{3}{5}{$\oddclassX$}{+}{6614}{class_N09_o_lam000010100r5z.eps}
& $\bar{T} \rightarrow$ & \classdrawlitX{1,2,3,4,5,6,7,11,9,10,8}{3}{7}{$\oddclassX$}{+}{785421}{class_N11_o_lam00001000100r7z.eps}\\
&\classdrawlitX{1,2,3,4,7,8,9,5,6}{3}{5}{$\oddclassX$}{-}{14709}{class_N09_o_lam000010100r5m.eps}
& $\bar{T} \rightarrow$ & \classdrawlitX{1,2,3,4,5,6,11,7,8,9,10}{3}{7}{$\oddclassX$}{-}{1220504}{class_N11_o_lam00001000100r7m.eps}\\
by Lemma \ref{lem.T.id}
&$
\left.
\begin{array}{r}
\classdrawlitX{1,2,3,4,5,6,7,8,9}{4}{4}{$\Id$}{0}{255}{class_N09_Id_lam000002000r4z.eps}\\
\classdrawlitX{1,2,3,4,5,8,9,6,7}{4}{4}{}{0}{15568}{class_N09__lam000002000r4z.eps}
\end{array}
\right\}
$
& $\bar{T} \rightarrow$ & \classdrawlitX{1,2,3,4,5,6,7,8,11,9,10}{4}{6}{}{0}{1495901}{class_N11__lam00000101000r6z.eps}\\
&\classdrawlitX{1,2,3,4,5,6,8,9,7}{5}{3}{$\oddclassX$}{+}{3954}{class_N09_o_lam000010100r3p.eps}
& $\bar{T} \rightarrow$ & \classdrawlitX{1,2,3,4,5,6,7,10,11,8,9}{5}{5}{$\oddclassX$}{+}{469943}{class_N11_o_lam00000020000r5p.eps}\\
&\classdrawlitX{1,2,3,4,8,9,5,6,7}{5}{3}{$\oddclassX$}{-}{8797}{class_N09_o_lam000010100r3m.eps}
& $\bar{T} \rightarrow$ & \classdrawlitX{1,2,3,4,5,8,9,10,11,6,7}{5}{5}{$\oddclassX$}{-}{729495}{class_N11_o_lam00000020000r5m.eps}\\
&\classdrawlitX{1,2,3,4,6,7,8,9,5}{6}{2}{}{0}{10543}{class_N09__lam000100010r2z.eps}
& $\bar{T} \rightarrow$ & \classdrawlitX{1,2,3,4,5,6,7,8,10,11,9}{6}{4}{}{0}{998333}{class_N11__lam00000101000r4z.eps}\\
&\classdrawlitX{1,2,3,4,9,6,7,8,5}{2,2,2}{2}{}{0}{1255}{class_N09__lam000000040r2z.eps}
& $\bar{T} \rightarrow$ & \classdrawlitX{1,2,3,4,11,5,6,10,7,8,9}{2,2,2}{4}{}{0}{119814}{class_N11__lam00000001030r4z.eps}\\
by Lemma \ref{lem.T.tree}
&$
\left.
\begin{array}{r}
\classdrawlitX{1,2,4,5,6,8,9,7,3}{7}{1}{$\oddclassX$}{+}{2679}{class_N09_o_lam001000001r1p.eps}\\
\classdrawlitX{1,2,4,5,6,7,8,9,3}{7}{1}{$\tree\ \oddclassX$}{+}{135}{class_N09_Idp_o_lam001000001r1p.eps}
\end{array}
\right\}
$
& $\bar{T} \rightarrow$ & \classdrawlitX{1,2,3,4,5,6,8,10,11,9,7}{7}{3}{$\oddclassX$}{+}{335525}{class_N11_o_lam00001000100r3p.eps}\\
&\classdrawlitX{1,2,4,5,8,9,6,7,3}{7}{1}{$\oddclassX$}{-}{6289}{class_N09_o_lam001000001r1m.eps}
& $\bar{T} \rightarrow$ & \classdrawlitX{1,2,3,4,5,6,8,9,10,11,7}{7}{3}{$\oddclassX$}{-}{521688}{class_N11_o_lam00001000100r3m.eps}\\
&\classdrawlitX{1,2,4,5,7,8,9,6,3}{3,2,2}{1}{}{0}{2569}{class_N09__lam000000121r1z.eps}
& $\bar{T} \rightarrow$ & \classdrawlitX{1,2,3,4,5,6,11,8,9,10,7}{3,2,2}{3}{}{0}{243536}{class_N11__lam00000000220r3z.eps}\\
\cline{2-4}
&$
\left.
\begin{array}{r}
\classdrawlitX{1,2,3,4,5,6,10,7,8,9}{3,2}{4}{}{0}{96434}{class_N10__lam0000001110r4z.eps}\\
\classdrawlitX{1,2,3,4,5,6,10,8,9,7}{3,3}{3}{$\oddclassX$}{+}{9876}{class_N10_o_lam0000000300r3z.eps}\\
\classdrawlitX{1,2,3,4,8,9,10,5,6,7}{3,3}{3}{$\oddclassX$}{-}{23167}{class_N10_o_lam0000000300r3m.eps}
\end{array}
\right\}
$
& $\bar{q}_2 \rightarrow$ & \classdrawlitX{1,2,3,4,6,7,11,8,9,10,5}{3,3,2}{2}{}{0}{164997}{class_N11__lam00000000220r2z.eps}\\
&$
\left.
\begin{array}{r}
\classdrawlitX{1,2,3,4,5,10,7,8,9,6}{2,2}{5}{}{0}{70886}{class_N10__lam0000010020r5z.eps}\\
\classdrawlitX{1,2,3,4,5,6,8,9,10,7}{4,2}{3}{}{0}{72006}{class_N10__lam0000001110r3z.eps}
\end{array}
\right\}
$
& $\bar{q}_2 \rightarrow$ & \classdrawlitX{1,2,3,4,6,7,9,10,11,8,5}{4,2,2}{2}{}{0}{181772}{class_N11__lam00000001030r2z.eps}\\
\begin{tabular}{r}
\rule{0pt}{20pt}\raisebox{-15pt}{\rule{0pt}{11pt}}%
by Lemma \ref{lem.q2.maxrank}\\
by Lemma \ref{lem.q2.id}
\end{tabular}\!\!\!
&$
\left.
\begin{array}{r}
\classdrawlitX{1,2,3,4,5,6,7,9,10,8}{}{9}{$\oddclassX$}{+}{233285}{class_N10_o_lam0100000000r9p.eps}\\
\classdrawlitX{1,2,3,4,5,7,8,9,10,6}{}{9}{$\oddclassX$}{-}{352697}{class_N10_o_lam0100000000r9m.eps}\\
\classdrawlitX{1,2,3,4,5,6,7,8,9,10}{}{9}{$\Id\ \oddclassX$}{-}{511}{class_N10_Id_o_lam0100000000r9m.eps}
\end{array}
\right\}
$
& $\bar{q}_2 \rightarrow$ & \classdrawlitX{1,2,3,4,6,7,8,9,10,11,5}{8}{2}{}{0}{747773}{class_N11__lam00010000010r2z.eps}\\
\cline{2-4}
&\classdrawlitX{1,2,3,4,5,6,7,9,10,8}{}{9}{$\oddclassX$}{+}{233285}{class_N10_o_lam0100000000r9p.eps}
& $\bar{q}_1 \rightarrow$ & \classdrawlitX{1,2,4,5,6,7,8,10,11,9,3}{9}{1}{$\oddclassX$}{+}{261005}{class_N11_o_lam00100000001r1p.eps}\\
&\classdrawlitX{1,2,3,4,5,7,8,9,10,6}{}{9}{$\oddclassX$}{-}{352697}{class_N10_o_lam0100000000r9m.eps}
& $\bar{q}_1 \rightarrow$ & \classdrawlitX{1,2,4,5,6,8,9,10,11,7,3}{9}{1}{$\oddclassX$}{-}{405607}{class_N11_o_lam00100000001r1m.eps}\\
by Lemma \ref{lem.q1.id}
&\classdrawlitX{1,2,3,4,5,6,7,8,9,10}{}{9}{$\Id\ \oddclassX$}{-}{511}{class_N10_Id_o_lam0100000000r9m.eps}
& $\bar{q}_1 \rightarrow$ & \classdrawlitX{1,2,4,5,6,7,8,9,10,11,3}{9}{1}{$\tree\ \oddclassX$}{-}{521}{class_N11_Idp_o_lam00100000001r1m.eps}\\
&$
\left.
\begin{array}{r}
\classdrawlitX{1,2,3,4,6,7,9,10,8,5}{5,2}{2}{}{0}{57606}{class_N10__lam0000010020r2z.eps}\\
\classdrawlitX{1,2,3,4,5,10,7,8,9,6}{2,2}{5}{}{0}{70886}{class_N10__lam0000010020r5z.eps}
\end{array}
\right\}
$
& $\bar{q}_1 \rightarrow$ & \classdrawlitX{1,2,4,5,6,11,8,9,10,7,3}{5,2,2}{1}{}{0}{145772}{class_N11__lam00000010021r1z.eps}\\
&$
\left.
\begin{array}{r}
\classdrawlitX{1,2,3,4,6,7,8,9,10,5}{4,3}{2}{}{0}{48954}{class_N10__lam0000001110r2z.eps}\\
\classdrawlitX{1,2,3,4,5,6,8,9,10,7}{4,2}{3}{}{0}{72006}{class_N10__lam0000001110r3z.eps}\\
\classdrawlitX{1,2,3,4,5,6,10,7,8,9}{3,2}{4}{}{0}{96434}{class_N10__lam0000001110r4z.eps}
\end{array}
\right\}
$
& $\bar{q}_1 \rightarrow$ & \classdrawlitX{1,2,4,5,6,7,9,10,11,8,3}{4,3,2}{1}{}{0}{246914}{class_N11__lam00000001111r1z.eps}\\
&\classdrawlitX{1,2,3,4,5,6,10,8,9,7}{3,3}{3}{$\oddclassX$}{+}{9876}{class_N10_o_lam0000000300r3z.eps}
& $\bar{q}_1 \rightarrow$ & \classdrawlitX{1,2,4,5,6,7,11,9,10,8,3}{3,3,3}{1}{$\oddclassX$}{+}{10796}{class_N11_o_lam00000000301r1z.eps}\\
&\classdrawlitX{1,2,3,4,8,9,10,5,6,7}{3,3}{3}{$\oddclassX$}{-}{23167}{class_N10_o_lam0000000300r3m.eps}
& $\bar{q}_1 \rightarrow$ &
\classdrawlitX{1,2,4,5,9,10,11,6,7,8,3}{3,3,3}{1}{$\oddclassX$}{-}{26887}{class_N11_o_lam00000000301r1m.eps}\\
\cline{2-4}
\end{tabular}
\caption{\label{tab.exinduct}A typical step of the induction, at a
  size sufficiently large that all fusion lemmas are already in place,
  and all typical situations do occur. All classes at size $n=11$
  (right column), except for $\Id_{11}$, are obtained from classes at
  $n=9$, under the action of $\bar{T}$, and classes at $n=10$ and rank
  large enough, under the action of $\bar{q}_1$ (on classes with $r
  \geq 2$) and $\bar{q}_2$ (on classes with $r \geq 3$). For
  comparison, the list of classes for $n\leq 8$, and their
  construction through surgery operators, is shown in
  figure~\ref{fig.fulltree48}.}
\end{table}

\subsection{Classification of extended Rauzy classes}
\label{ssec.permsextheo}

In this section we show how, by a crucial use of
Lemma~\ref{trivialbutimp}, the classification theorem for
$\permsex_n$, Theorem \ref{thm.Main_theorem_2} descends from the one
for $\perms_n$, Theorem \ref{thm.Main_theorem}, as a rather
straightforward corollary.

The idea is that, by including more operators for the dynamics on the
same set of configurations, we may only join classes. As, in fact,
most of the invariants for $\perms_n$ survive in $\permsex_n$
(essentially, only the rank is lost), as a matter of fact, and as we
shall prove, there is no residual lack of connectivity: all classes of
$\perms_n$ with the same structure of invariants except for the rank
ultimately join together in~$\permsex_n$.

Let us define two involutions on $\kS_n$, $S$ and $\bar{S}$.  The
first one is defined as $S \s = \s^{-1}$, i.e., in matrix
representation, it acts as a reflection along the main diagonal, and,
in diagram representation, as a reflection along the horizontal
axis. The operator $S$ intertwines between the four operators in the
$\permsex_n$ dynamics, i.e.\ $SLS=L'$ and $SRS=R'$.  Analogously, we
have $\bar{S} \s = \tau$ for $s(i)=j \Leftrightarrow \tau(n+1-i) =
n+1-j$.  Then, $\bar{S}$ acts on permutations in matrix representation
as a reflection along the main anti-diagonal, and on permutations in
diagram representation as a reflection along the vertical axis, and
intertwines between the four operators in the $\permsex_n$ dynamics as
$\bar{S}L\bar{S}=R'$ and $\bar{S}R\bar{S}=L'$.

Recall that, in Section \ref{ssec.typeXH}, we defined permutations of
type $X$ and $H$, and the \emph{principal cycle} in a type-$X$
permutation as the cycle passing through the `$-1$ mark'. It is easy to
establish
\begin{proposition}
\label{prop.6543654}
Let $\s$ be an irreducible permutation of type $H$, with invariant
$(\lam,r,s)$. Then $\s^{-1} = S \s$ is an irreducible
permutation of type $H$, with invariant $(\lam,r,s)$.

Let $\s$ be an irreducible permutation of type $X$, with invariant
$(\lam,r,s)$, and principal cycle of length $\bar{r}$. Let
$\bar{\lam}=\lam \setminus \{ \bar{r} \} \cup \{ r \}$. Then 
$\s^{-1} = S \s$ is an irreducible permutation of type $X$, with
invariant $(\bar{\lam},\bar{r},s)$ and principal cycle of length~$r$.
\end{proposition}
\noindent
This follows from the diagrammatic construction of the invariants,
which is symmetric w.r.t.\ the involution $S$, up to exchanging the
lengths of the rank and of the principal cycle in type-$X$
permutations. \qed

\noindent
Recall that we define $\lambda'=\lambda\cup \{r\}$, and that we have
determined that $\lam'$ is invariant in the $\permsex_n$ dynamics.

For the proof Theorem \ref{thm.Main_theorem_2} we need two
ingredients. On one side, we need to deal with exceptional classes.
For $\Id'_n$, it is especially easy: this class has rank 1 at all $n$,
thus it is not primitive in the $\permsex$ dynamics, and shall be
discarded. So we are left to prove that the class $\Id_n$ is not
connected to other classes by the two further operators $L'$ and $R'$
of the dynamics. This is easily evinced from the results of Appendix
\ref{ssec.appid}, in particular the organisation of the configurations
in $\Id_n$ in a complete binary tree. The involution $\bar{S}$ acts
just as the vertical symmetry on this tree. As a result, the class
$\Id_n$ is closed under the action of $L$, $R$ and $\bar{S}$, which,
as $L'=\bar{S}R\bar{S}$ and $R'=\bar{S}L\bar{S}$, implies closure
under the action of $L$, $R$, $L'$ and~$R'$.

On the other side, we shall prove that all the non-exceptional classes
of $\perms_n$ with same sign invariant, and cycle invariant 
$(\lam' \setminus \{ \lam'_i \}, \lam'_i)$, for $i$ running over the
different cycles of $\lam'$, are connected by the extended dynamics. A
first example occurs at size $7$, where we shall show that the classes
of $\perms_7$ with invariant $(\{2\},4,0)$ and $(\{4\},2,0)$ merge
into a unique class with invariant $(\{2,4\},0)$.

In fact, we shall prove an even stronger statement:
\begin{theorem}
\label{thm.perSex}
Any two configurations $\s$, $\s'$ in the same class for the dynamics
$\permsex_n$
are connected by a word of the form $w=w_1 (L')^k w_2$, with $k$ an
integer, and $w_1$, $w_2$ words in the $\perms_n$ dynamics.
\end{theorem}
\proof In other words, let $\lam'=\lam \cup \{r_1, r_2\}$, with 
$r_1 \neq r_2$. We want to prove that the classes $C_1$ and $C_2$ of
$\perms_n$, with invariants $(\lam \cup \{r_2\},r_1,s)$ and 
$(\lam \cup \{r_1\},r_2,s)$, respectively, are connected in
$\permsex_n$, and through a word $w$ of the type above.

Let $\s \in C_1$ be a configuration in a standard family, and of type
$H(r'_1,r''_1)$ (thus with $r'_1 + r''_1 = r_1$). From our
classification theorem for $\perms_n$ we can restrict to such
configurations with no loss of generality.  It is easily seen that
$S\s$ is also standard, and has the same invariant in $\perms_n$ (by
Proposition~\ref{prop.6543654}), thus $\s \sim S \s$ in the $\perms_n$
dynamics.

By Lemma \ref{trivialbutimp}, for all $r_2 \in \lam$ there exists $i$
such that $L^i \s$ has form $X(r_1,r_2)$, and thus (by
Proposition~\ref{prop.6543654}) $S L^i \s = (L')^i S \s$ has form
$X(r_2,r_1)$, and is thus in class $C_2$.
This completes our proof.  
\qed

\noindent
As remarked above, the stronger Theorem \ref{thm.perSex} implies
Theorem \ref{thm.Main_theorem_2} as a corollary.

\appendix

\section{Representation of Cayley graphs}
\label{app.repres}

This short appendix deals with a graphical convention on the
representation of the dynamics. 

Given a group dynamics with $h+k$ generators $\{\iota_1, \ldots,
\iota_h, \kappa_1, \ldots, \kappa_k\}$, of which $h$ involutive and
$k$ non-involutive, a \emph{Cayley (di)graph} can be constructed, with
one connected component per equivalence class. Vertices $v_x$ are
associated to configurations $x$. The edges are labeled by a
generator. They are unoriented in the case of involutive generators,
and oriented otherwise. We have the unoriented edge $(v_x,v_y)$ with
label $\iota_a$ if $\iota_a x = y$, and the oriented edge $(v_x,v_y)$
with label $\kappa_a$ if $\kappa_a x = y$. As a result, each vertex
has overall degree $h+2k$.

For dynamics with two (non-involutive) generators, such as $\perms_n$ and
$\matchs_n$, the Cayley graph has thus degree 4.
Even for relatively small classes, this makes a too complex structure
for being visualised in an effective way. In this section we introduce
a notation that, through an embedding on a two-dimensional grid,
allows to omit to draw a fraction of the edges with no loss of
information. As a result of this pruning, some of the Cayley graphs of
classes are reduced to regular intelligible structures, this being in
particular the case for exceptional classes, discussed in
Appendix~\ref{app.ididp}. From this point on, we concentrate on the
case of non-involutive generators.

Clearly, such an enhancement cannot be made for graphs in general, 
and we shall use the specialties coming from the fact that
these graphs are Cayley graphs of a group action. 

The iterated application of a single generator partitions a Cayley
graph into cyclic orbits. Let us embed our graph in dimension 
$d\geq 2$, and choose a direction $\theta_a$, in the sphere of
dimension $d-1$, for each generator $\kappa_a$.  If we take care of
representing a directed edge with label $\kappa_a$ as oriented in the
direction $\theta_a$, then one edge per orbit can be omitted without
loss of information: an orbit of length $\ell$ is turned into a linear
chain with $\ell$ vertices and $\ell-1$ collinear oriented edges, and
it will be understood without ambiguities that the image of the last
vertex of the chain under the operator is the first vertex of the
chain.  In fact, edges can be omitted completely, and only the angles
$\theta_a$ be specified, at the condition that no other vertices,
besides those in the orbit, lay on the line in $\mathbb{R}^d$
associated to the orbit. If $d=k$, the directions $\kappa_a$ can be
taken to constitute the canonical basis without loss of generality.

The most useful application of this strategy is when there are two
generators, and thus the graph is conveniently represented on a plane.
In this special case, we will use edges with a single- or
double-arrow, and in red or blue, for the two generators, and add a
small tag on the last vertex of any chain, for further enhancing the
visualisation.

\begin{figure}
\begin{center}
\includegraphics[scale=.75]{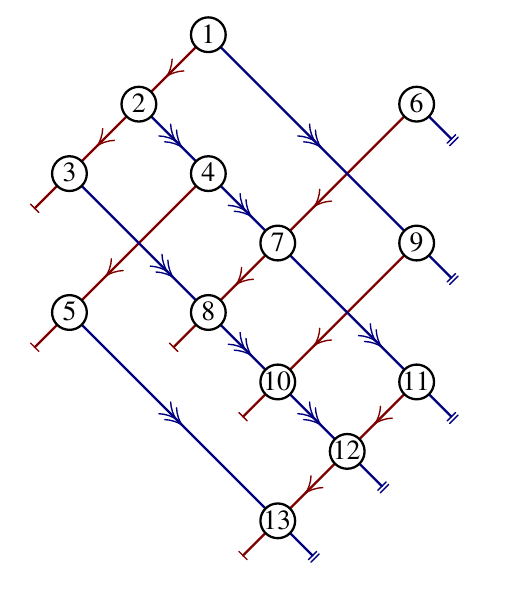}
\includegraphics[scale=.75]{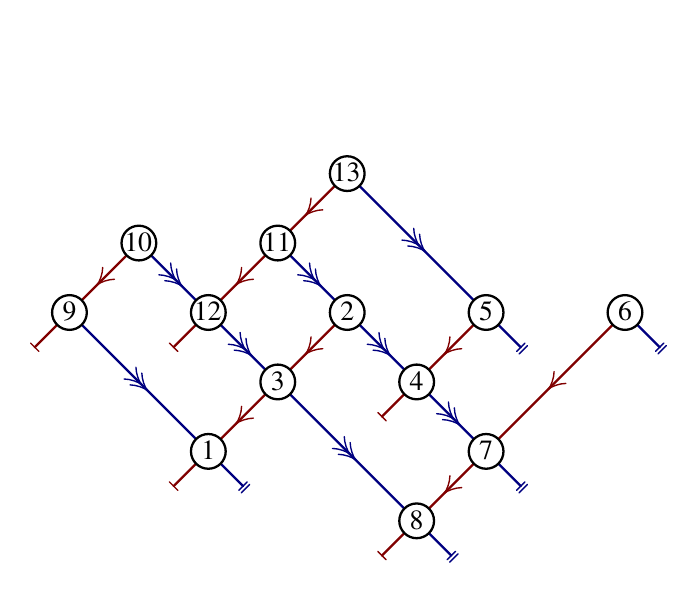}
\includegraphics[scale=.75]{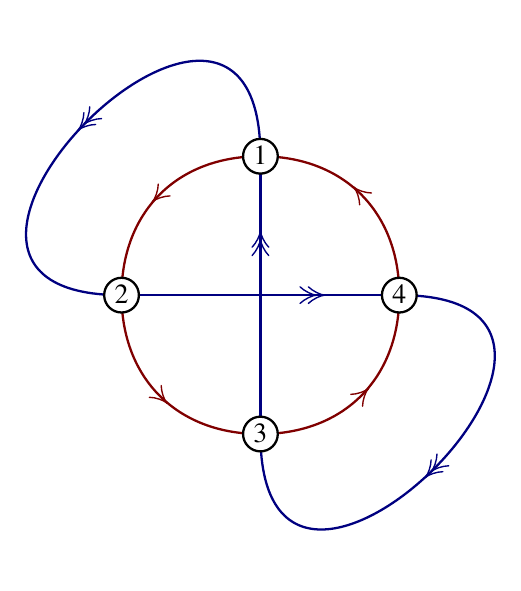}
\end{center}
\caption{\label{fig.straightgraph}Left: Straight representation of a
  digraph. This representation is \emph{not} planar.  Middle: a planar
  straight representation of the same digraph.  Right: a digraph which
  is not straight: the four vertices are in the same red and blue
  orbit, and cannot be simultaneously collinear along $x$ and along
  $y$, on two distinct directions and in the given cyclic orderings,
  while staying distinct.}
\end{figure}

See Figure \ref{fig.straightgraph}, left, for an example of
representation of such a digraph.

We say that a Cayley graph admitting a representation with these
properties is \emph{straight}, and if it admits a representation with
these properties, and such that the edges do not cross, it is
\emph{straight planar}.

It is easy to see that if a graph is straight, it also allows a
straight representation in which all vertices are distinct vertices of
$\mathbb{Z}^2$, and different orbits are in different rows/columns of
the grid, thus straight digraphs are also representable as finite
subsets of $\mathbb{Z}^2$, and a list of rational numbers (one per
generator) associated to the slopes.

It is also easy to see that not all digraphs are straight. See Figure
\ref{fig.straightgraph}, right, for a simple counter-example.

\section{Non-primitive classes}
\label{app.primitheo}

In this appendix we discuss how the classification theorem on
primitive classes implies, with a small amount of further reasoning, a
classification theorem for all classes, including non-primitive
ones. More generally, we prove how even the full structure of the
Cayley graph of a non-primitive class can be evinced from the Cayley
graph of the associated primitive class.  This result is announced in
Section~\ref{ssec.primi}.

Recall that we announced, in Corollary \ref{cor.primhomo}, that the
map $\prim{\,\cdot\,}$ is a homomorphism for the dynamics, so that, as
a result of the analysis performed below, if $\s \in C$, we can
naturally define the \emph{primitive} $C'=\prim{C}$ as the class
containing $\prim{\s}$ for any $\s \in C$.

We start by discussing the most basic example of irreducible
non-primitive class, namely the set of classes $T_n$, of size $n \geq
2$, such that $\prim{T_n}=\Id_2$.\;\footnote{This is the smallest
  example among irreducible configurations, as if $\prim{\s}=\id_1$
  and $\s \neq \id_1$, then $\s$ is reducible.}

\begin{figure}
\begin{center}
\setlength{\unitlength}{15pt}
\begin{picture}(16,9)
\put(0,.3){\includegraphics[scale=.75]{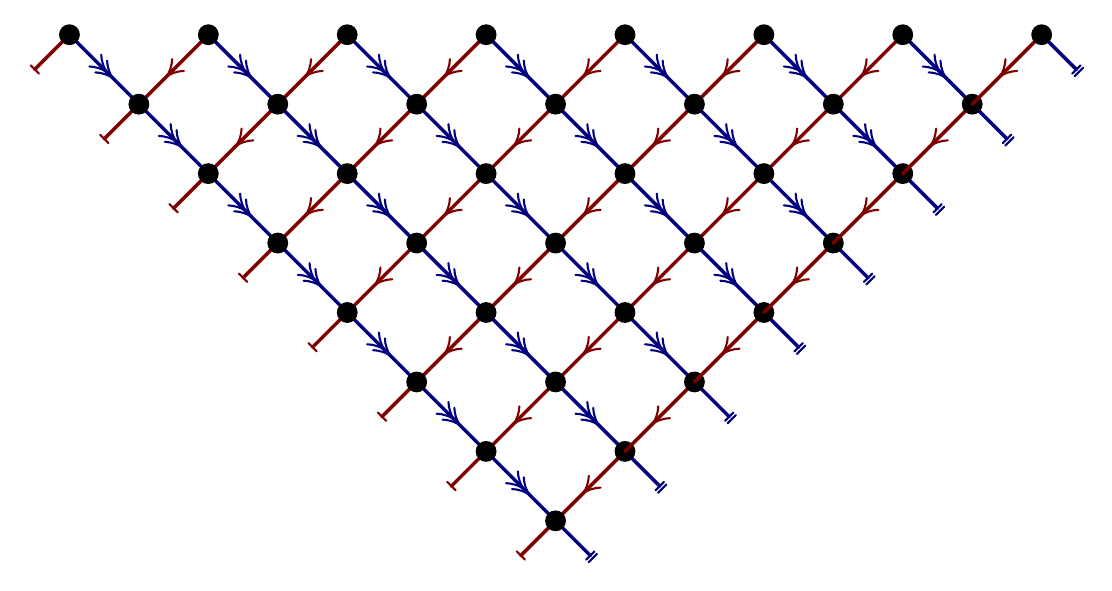}}
\put(2.5,2.5){\makebox[0pt][c]{$a$}}
\put(13.5,2.5){\makebox[0pt][c]{$b$}}
\put(7,0){\makebox[0pt][c]{${1}$}}
\put(6,1){\makebox[0pt][c]{${2}$}}
\put(5,2){\makebox[0pt][c]{${3}$}}
\put(4,3){\makebox[0pt][c]{${4}$}}
\put(3,4){\makebox[0pt][c]{${5}$}}
\put(2,5){\makebox[0pt][c]{${6}$}}
\put(1,6){\makebox[0pt][c]{${7}$}}
\put(0,7){\makebox[0pt][c]{${8}$}}
\put(9,0){\makebox[0pt][c]{${1}$}}
\put(10,1){\makebox[0pt][c]{${2}$}}
\put(11,2){\makebox[0pt][c]{${3}$}}
\put(12,3){\makebox[0pt][c]{${4}$}}
\put(13,4){\makebox[0pt][c]{${5}$}}
\put(14,5){\makebox[0pt][c]{${6}$}}
\put(15,6){\makebox[0pt][c]{${7}$}}
\put(16,7){\makebox[0pt][c]{${8}$}}
\end{picture}
\end{center}
\caption{\label{fig.classtri}The straight representation of the Cayley
  graph of the class $T_n$. Here $n=9$.}
\end{figure}
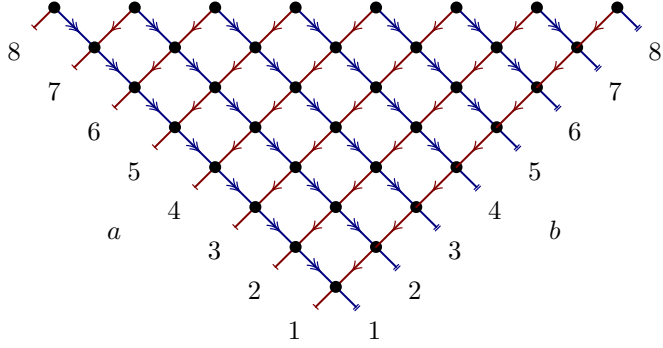

We claim that the class $T_n$ has cardinality $\binom{n}{2}$, and that
the Cayley graph admits an especially simple straight planar
representation, illustrated in Figure~\ref{fig.classtri}, namely a
triangular portion of the square grid.

\begin{figure}
\begin{align*}
&
\setlength{\unitlength}{15pt}
\begin{picture}(9,5.3)(0,-1.3)
\put(0,0){\includegraphics[scale=3]{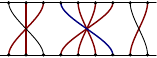}}
\put(0.1,-.2){$\underbrace{\rule{42pt}{0pt}}_{\displaystyle{a}}$}
\put(7.1,-.2){$\underbrace{\rule{27pt}{0pt}}_{\displaystyle{b}}$}
\end{picture}
&&
\setlength{\unitlength}{15pt}
\begin{picture}(9,5.3)(0,-1.3)
\put(0,0){\includegraphics[scale=3]{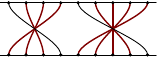}}
\put(0.1,-.2){$\underbrace{\rule{57pt}{0pt}}_{\displaystyle{a}}$}
\put(4.1,-.2){$\underbrace{\rule{72pt}{0pt}}_{\displaystyle{b}}$}
\end{picture}
\end{align*}
\caption{\label{fig.classtri23}Left, the configuration with coordinate
  $(a,b)=(3,2)$ in $T_9$, composed of three bundles of crossing edges,
  and having a special descent.  Right, the configuration with
  coordinate $(a,b)=(4,5)$ in $T_9$, composed of two bundles of
  crossing edges, and having no special descent.}
\end{figure}
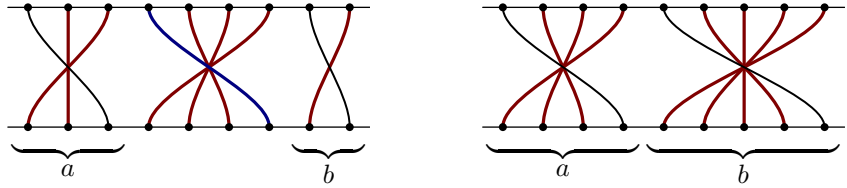

More precisely, w.r.t.\ the coordinates $a$ and $b$ shown in figure,
we claim that the configurations with coordinate $(a,b)$ are composed
of three bundles, and have a special descent, if $a+b<n$, and are
composed of two bundles, and have no special descent, if $a+b=n$, as
shown in Figure~\ref{fig.classtri23}. These facts are easily verified.

\noindent
For a permutation $\s$ of size $n$, with no special descent, call
$\tilde{\s}$ the permutation of size $n+1$ in which the special
descent has been added:
\begin{align*}
\s&=\raisebox{-12pt}{\includegraphics[scale=1.8]{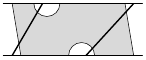}}
&
\tilde{\s}&=\raisebox{-12pt}{\includegraphics[scale=1.8]{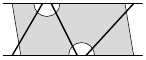}}
\end{align*}
Now we want to illustrate how the structure of classes $T_n$ is
sufficiently general to describe the local structure of non-primitive
classes $C$ in terms of the one of the class $C'=\prim{C}$.  Let $n$
and $n'$ be the sizes of classes $C$ and $C'$. Let $\s \in C$, without
a special descent, and $\t=\prim{\s}$. The configuration $\s$ thus
consists of $n'$ bundles of crossing edges. Let
$\bm{m}=(m_1,\ldots,m_{n'})$ be the cardinalities of these bundles (so
that $m_i \geq 1$ and $m_1+\cdots+m_{n'}=n$). In particular, we single
out the multiplicities of the two pivots, i.e.,
$\bm{m}=(m_1,\ldots,m_{\tau^{-1}(n')},\ldots)=:
(m_L,\ldots,m_R,\ldots)$.
The datum of the pair $(\t,\bm{m})$ is of course equivalent to
$\s$. The following properties are easily verified:
\begin{itemize}
\item
For $i < m_L$,  $L^i \s$, gives the permutation
associated to $(\tilde{\t},(m_L-i,\ldots,i,m_R,\ldots))$.  Similarly,
for $i < m_R$, $R^i\s$ gives $(\tilde{\t},(m_L,\ldots,i,m_R-i,\ldots))$.
\item
$L^{m_L} \s$ gives $(L\t,(\ldots,m_L,m_R,\ldots))$.
\item
Similarly, $R^{m_R} \s$ gives $(R\t,(m_L,\ldots,m_R,\ldots))$.
\end{itemize}
As a corollary, if $m_L \geq 2$ we have that $R^{-1}L \s$ gives
$(\t,(m_L-1,\ldots,m_R+1,\ldots))$, and if $m_R \geq 2$ then $L^{-1}R
\s$ gives $(\t,(m_L+1,\ldots,m_R-1,\ldots))$, which implies
\be
\label{eq.movedesc}
(\t,(m_L,\ldots,m_R,\ldots)) \sim
(\t,(m_L+c,\ldots,m_R-c,\ldots))
\qquad
\forall \quad
-m_L+1 \leq c \leq m_R-1
\ef.
\ee
More generally, combining the remarks above, we have
\be
\label{eq.movedesc2}
\begin{array}{rl}
(\t,(m_L,\ldots,m_R,\ldots)) 
& \!\!\!\!\sim
(\tilde{\t},(m_L-a,\ldots,a+b,m_R-b,\ldots))
\\
&
\qquad
\forall \ a \leq m_L-1, \ b \leq m_R-1, \ a+b \geq 0
\\
& \!\!\!\!\sim
(L \t,(\ldots,m_L+c,m_R-c,\ldots))
\\
&
\qquad
\forall \ -m_L+1 \leq c \leq m_R-1
\\
& \!\!\!\!\sim
(R \t,(m_L+c,\ldots,m_R-c,\ldots))
\\
&
\qquad
\forall \ -m_L+1 \leq c \leq m_R-1
\ef.
\end{array}
\ee
In other words, a portion of the Cayley graph for $C$, containing
$\s=(\t,(m_L,\ldots,m_R,\ldots))$, is isomorphic to a large patch of
the class $T_{m_L+m_R}$ represented in Figure \ref{fig.classtri}, with
the identification of parameters $(a,b)=(m_L,m_R)$, up to an important
difference. For the class $T_k$ in itself, according to the notation
introduced in Appendix \ref{app.repres}, when we exit the
straight-planar representation of $T_k$ from one boundary, we come
back to a vertex on another boundary of the same class. The portion of
$T_n$ within $C$, instead, is such that when we exit from the boundary
of its planar representation, we enter into a different copy $T_{k'}$,
and associated to a different set of omitted parameters $(\tau',
\bm{m'} \smallsetminus \{m'_L,m'_R\})$.

It is a fortunate coincidence that the simple straight representation
of the class $T_k$ presented above has the orbits cut exactly at the
position at which the different copies of $T_k$'s are patched
together.

\vspace{2mm}

\pfof{Proposition \ref{prop.primhomoPre}}
%
As another application of Lemma \ref{lem.reddynconstr}, 
we see that
\[
(\t,(m_1,m_2,\ldots,m_{n})) \sim (\t,(m'_1,m'_2,\ldots,m'_{n}))
\]
whenever $\sum_i m_i = \sum_i m'_i$, and all $m_i$, $m'_i$ are
strictly positive.
This implies the statement for $\s$ having no special descent. A
configuration $\s'$ with a special descent is always at alternating
distance 1 from a $\s$ with no special descent, this allowing to
complete the reasoning.  \qed

The clarification of the Cayley graph of $C$ in terms of the Cayley
graph of $C'$ has a corollary on the size of these graphs. Let $C$ be
a class of size $n+k$, with $k$ descents. Then $C'$ is a class of size
$n$. Let $\t \in C'$. Consider the number of configurations $\s \in C$
which can be written as $\s=(\t,(m_L,\ldots,m_R,\ldots))$
(if $\s$ has no special descent) or as
$\s=(\tilde{\t},(m_L,\ldots,m_{\rm pivot},m_R,\ldots))$ (if $\s$ has a
special descent).  Identify $(\t,(m_L,\ldots,m_R,\ldots)) \equiv
(\tilde{\t},(m_L,\ldots,0,m_R,\ldots))$, and consider the list
$(m_L-1,m_2-1,\ldots,m_{\rm pivot},m_R-1,\ldots,m_n-1)=:(k_1,\ldots,k_{n+1})$.  
At fixed $\t \in C'$, each list with $k_i \geq 0$ and $\sum_i k_i = k$
is obtained from a single $\s \in C$, and lists not satisfying the
properties above are never obtained. As a result, we have
\begin{proposition}
With notations as above,
\be
|C| = \binom{n+k}{n} |C'|
\ef.
\ee
\end{proposition}
A version of this proposition, for the case of the extended dynamcis,
has been first established by Delecroix in \cite[Thm.\ 2.4]{Del13}.
His derivation is quite different in spirit, as it results from the
analysis of a complicated general formula for $|C|$, instead of
deriving directly the weaker (and easier) result on $|C|/|C'|$.

\section{Exceptional classes\label{sec.excp_class}}
\label{app.ididp}

In this appendix we describe the structure of the Cayley graph of the
classes $\Id_n$ and $\tree_n$, which we have called 
\emph{exceptional classes} in the main text.  This structure is so
rigid that, for any permutation $\s$, it is `easily' tested if $\s
\Id_n$, $\s \in \tree_n$ or none of the two (by `easily' we mean, in
particular, that this can be tested in a time $\sim n$, instead of the
na\"ive upper bound $\sim 2^n$ based on the cardinality of these
classes).  This result is of crucial importance for our classification
theorem, as the latter is based on a result of the form ``a list of
invariants for the dynamics is complete''. Our main invariants are the
cycle (with the rank) and the sign. We have painfully proven that
these invariants distinguish non-exceptional classes, however at
generic $n$ each of the two exceptional classes have the same cycle
and sign invariant of one non-exceptional class, so that the complete
list of invariants must include the boolean `exceptionality'
invariant, i.e.\ the outcome of the forementioned membership testing
algorithm.

In this appendix, when using a matrix representation of
configurations, it is useful to adopt the following notation:
The symbol $\epsilon$ denotes the $0 \times 0$ empty matrix.
The symbol
\raisebox{-2.5pt}{\includegraphics[width=4.8mm]{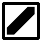}}
denotes a square block in a matrix (of any size $\geq 0$), filled with
an identity matrix.  A diagram, containing these special symbols and
the ordinary bullets used through the rest of the paper, describes the
set of all configurations that could be obtained by specifying the
sizes of the identity blocks.
In such a syntax, we can write
equations of the like
\begin{align}
\id 
&:=
\raisebox{-1mm}{\includegraphics[width=4.8mm]{FigureA2_fig_matr_id.pdf}}
=
\epsilon \cup
\raisebox{-3mm}{\includegraphics[width=8.8mm]{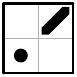}}
=
\epsilon \cup
\raisebox{5.8mm}{\includegraphics[width=8.8mm,
    angle=180]{FigureA2_fig_matr_Bid.pdf}}
\ef;
&
\id'
&:=
\raisebox{-7mm}{\includegraphics[width=16.8mm]{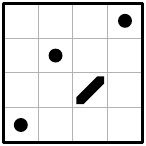}}
\ef.
\end{align}
The sets $\id$ and $\id'$ contain one element per size, $\id_n$ and
$\id'_n$, for $n \geq 0$ and $n \geq 3$ respectively.

The two exceptional classes $\Id_n$ and $\tree_n$ contain the
configurations $\id_n$ and $\id'_n$, respectively (and are primitive
for $n \geq 4$ and $n \geq 5$, respectively).

\subsection{Classes $\Id_n$}
\label{ssec.appid}

\begin{figure}[tb!]
\begin{center}
\begin{align*}
X_{LL}^{(k)}
&=
\quad
\makebox[0pt][l]{\raisebox{-16mm}{\includegraphics[width=33mm]{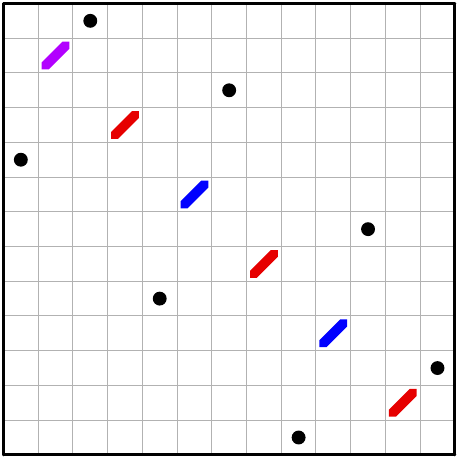}}}
\hspace{-16pt}
\raisebox{-30pt}{$\rotatebox{45}{$k \left\{\rule{0pt}{50pt}\right. $}$}
&
X_{RR}^{(k)}
&=
\makebox[0pt][l]{\raisebox{-16mm}{\includegraphics[width=33mm]{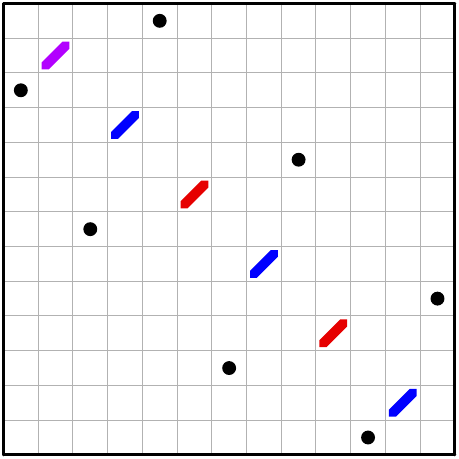}}}
\hspace{28pt}\raisebox{14pt}{$\rotatebox{45}{$\left. \rule{0pt}{50pt}\right\} k $}$}
\\
X_{RL}^{(k)}
&=
\makebox[0pt][l]{\raisebox{-16mm}{\includegraphics[width=38mm]{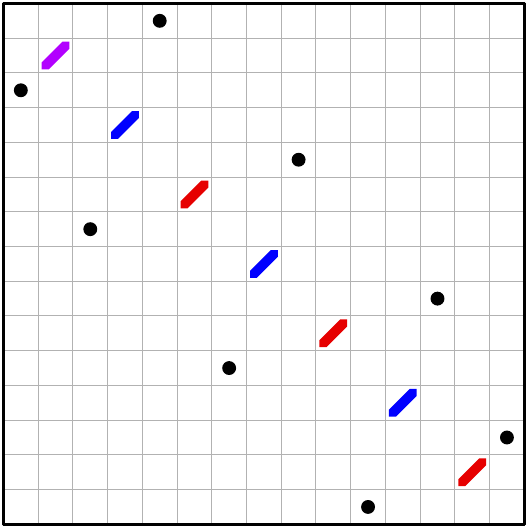}}}
\hspace{-2pt}
\raisebox{-30pt}{$\rotatebox{45}{$k \left\{\rule{0pt}{50pt}\right. $}$}
&
X_{LR}^{(k)}
&=
\makebox[0pt][l]{\raisebox{-16mm}{\includegraphics[width=38mm]{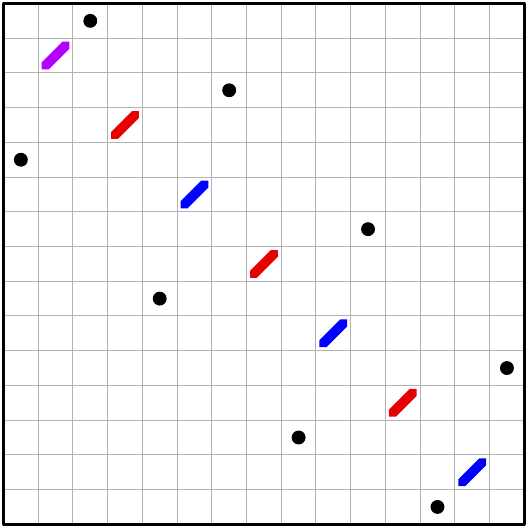}}}
\hspace{42pt}\raisebox{14pt}{$\rotatebox{45}{$\left. \rule{0pt}{50pt}\right\} k $}$}
\end{align*}
\end{center}
\caption{Description of the structure of configurations in
the identity class $\Id_n$, besides the configuration $\id_n$.
\label{fig.struct_id}}
\end{figure}

The structure of the classes $\Id_n$ is summarised by the following
relation:
\be
\Id := \bigcup_n \, \Id_n
=
\bigg(
\bigcup_{k \geq 1}
(X_{RL}^{k} \cup X_{LR}^{k} \cup X_{LL}^{k} \cup X_{RR}^{k})
\bigg)
\cup
\id
\ee
where the configurations $X_{\cdot \cdot}^{k}$ are defined as in
figure~\ref{fig.struct_id} (discard colours for the moment).

In other words, we claim that the configurations in classes $\Id_n$
are partitioned into five disjoint sets: the identity configurations
$\id_n$, and those contained in the four sets $X_{\cdot \cdot}^{k}$,
with $\cdot=R,L$ and total size $n$. More in detail, the configuration
$\s_{RL}^{n;(i_1,j_1,\ldots,i_k,j_k)} := R^{j_k+1}\cdots
L^{i_2+1}R^{j_1+1}L^{i_1+1} \id_n$ is in $X_{RL}^{k}$ whenever
$i_1+j_1+i_2+\cdots+j_k = n-2k-2+\delta$, with $\delta \geq 0$, and it
is represented exactly as in figure \ref{fig.struct_id}, bottom-left,
with red blocks having size $i_1$, $i_2$, \ldots, $i_k$, (from
bottom-right to top-left), blue blocks having size $j_1$, $j_2$,
\ldots, $j_k$ (still from bottom-right to top-left), and the violet
box having size $\delta$.  The other three sets have similar
definitions.

\begin{figure}
\[
\includegraphics[scale=.52, angle=90]{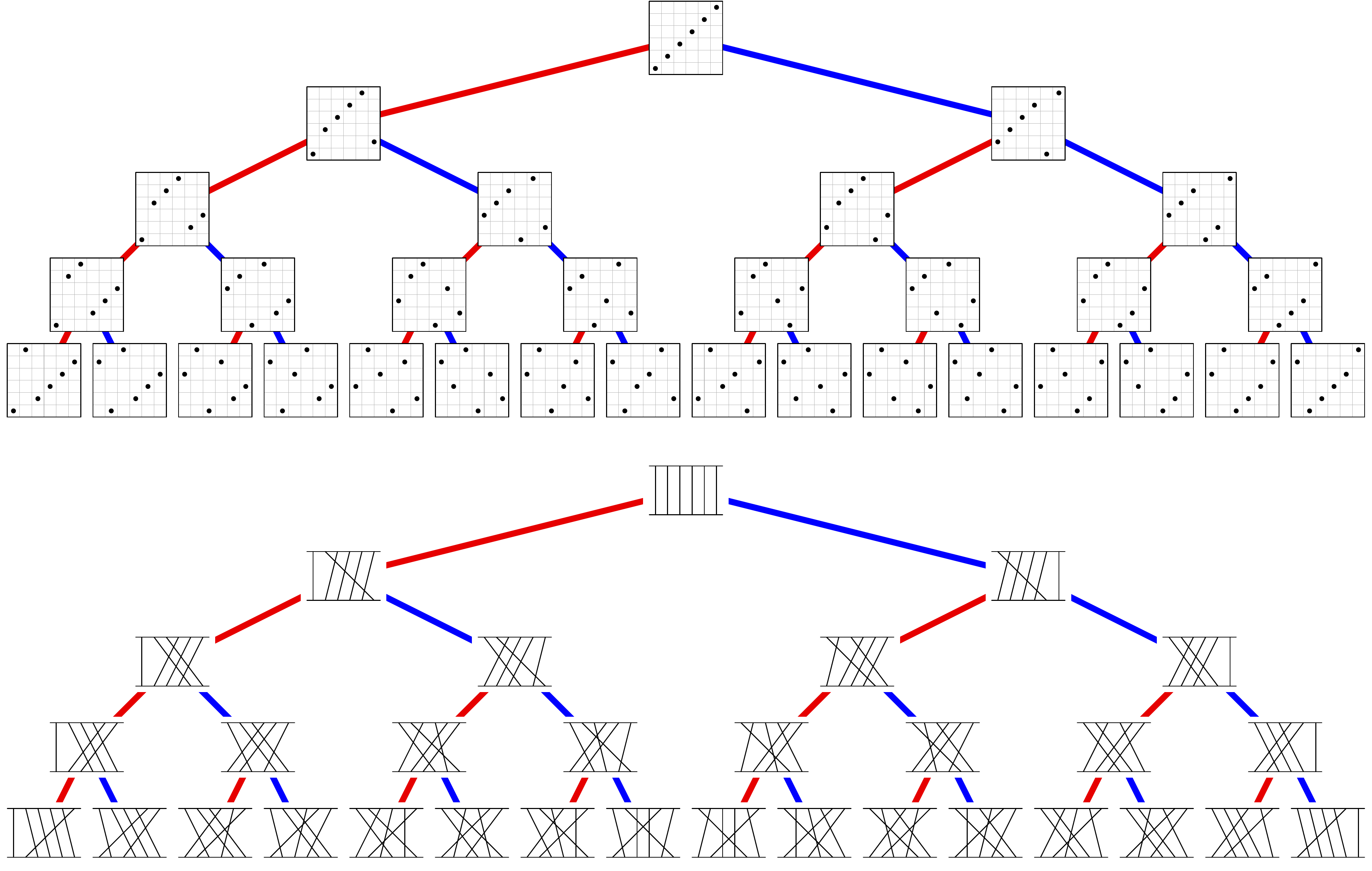}
\]
\caption{\label{fig.id6}The Cayley graph of $Id_6$, with
  configurations in matrix and
  diagram representation. The image shall be observed rotated by 90
  degrees (i.e., the root of the tree is on top). Red and blue edges
  correspond to actions of $L$ and $R$ respectively, and are omitted
  with the same rules of straight representations.}
\end{figure}

We also claim that the class $\Id_n$ has a straight planar
representation (in the sense of Appendix \ref{app.repres}) consisting
of a complete binary tree of height $n-2$ (and, in particular,
$|\Id_n| = 2^{n-1}-1$).  An illustration of this fact for $\Id_6$ is
presented in Figure \ref{fig.id6}.  The root configuration of the tree
is $\id_n$. All the vertices of the tree can be described in terms of
the unique path that reaches them starting from the root. When we have
a vertex such that the corresponding path starts with $i_1+1$
left-steps, followed by $j_1+1$ right steps, followed by $i_2+1$
left-steps, and so on, and finally terminating with $j_k+1$ right
step, the corresponding configuration is the forementioned
$\s_{RL}^{n;(i_1,j_1,\ldots,i_k,j_k)}$.  The other three cases are
described analogously.

\begin{figure}[b!!]
\[
\begin{array}{rcc}
&
\delta \geq 1
&
\delta = 0
\\
\s_{RL}^{n;(i_1,j_1,\ldots,i_k,j_k)}
&
\makebox[0pt][l]{\raisebox{-16mm}{\includegraphics[width=35.5mm]{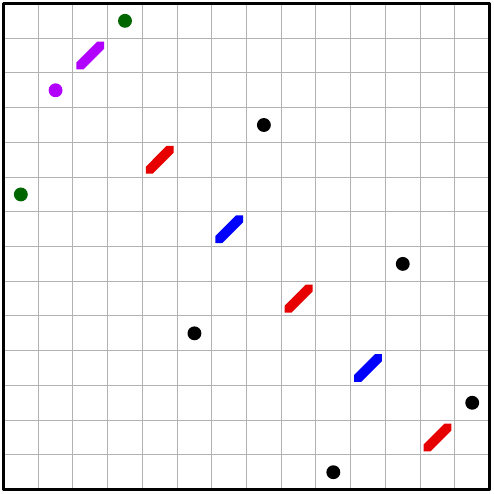}}}
\setlength{\unitlength}{2.5mm}
\begin{picture}(14,8)(0,7)
\put(2.8,12.2){$\delta-1$}
\put(5.2,10){$i_k$}
\end{picture}
&
\makebox[0pt][l]{\raisebox{-16mm}{\includegraphics[width=30.5mm]{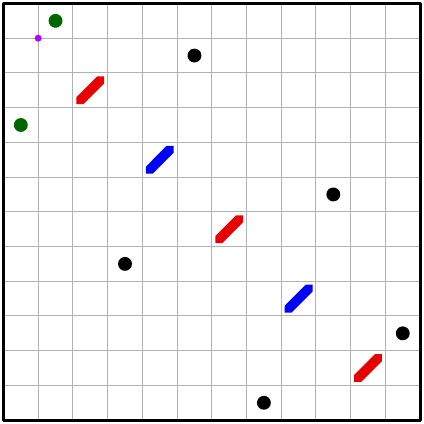}}}
\setlength{\unitlength}{2.5mm}
\begin{picture}(12,7)(0,7)
\put(3.2,10){$i_k$}
\end{picture}
\\
L \, \s_{RL}^{n;(i_1,j_1,\ldots,i_k,j_k)}
&
\makebox[0pt][l]{\raisebox{-16mm}{\includegraphics[width=35.5mm]{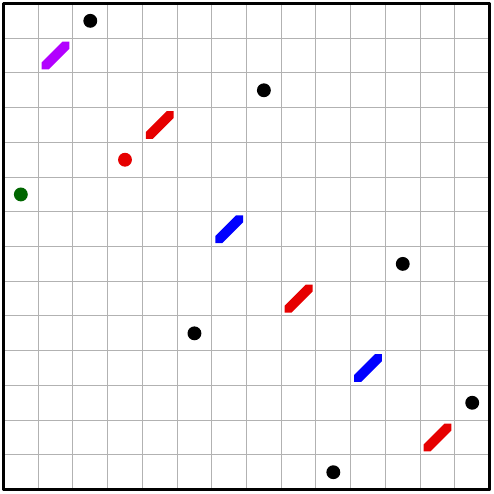}}}
\setlength{\unitlength}{2.5mm}
\begin{picture}(14,8.5)(0,7)
\put(1.8,12.2){$\delta-1$}
\put(5.2,11){$i_k$}
\end{picture}
&
\makebox[0pt][l]{\raisebox{-16mm}{\includegraphics[width=30.5mm]{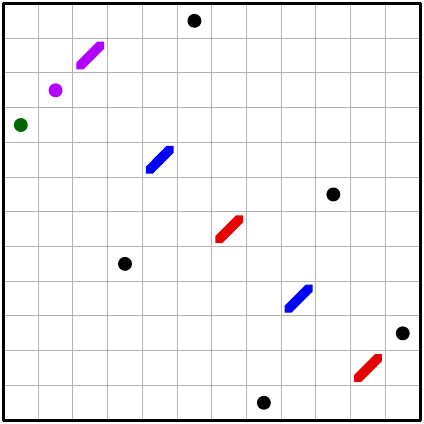}}}
\setlength{\unitlength}{2.5mm}
\begin{picture}(12,7)(0,7)
\put(3.2,11){$i_k$}
\end{picture}
\\
R \, \s_{RL}^{n;(i_1,j_1,\ldots,i_k,j_k)}
&
\makebox[0pt][l]{\raisebox{-16mm}{\includegraphics[width=35.5mm]{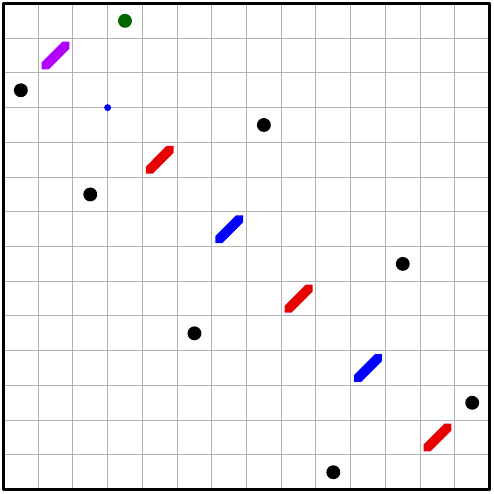}}}
\setlength{\unitlength}{2.5mm}
\begin{picture}(14,8.5)(0,7)
\put(2.3,12.7){$\delta-1$}
\put(3.7,11.5){$j_k=0$}
\put(5.2,10){$i_k$}
\end{picture}
&
\makebox[0pt][l]{\raisebox{-16mm}{\includegraphics[width=30.5mm]{FigureA2_fig_matr_CidGen1_del0.pdf}}}
\setlength{\unitlength}{2.5mm}
\begin{picture}(12,7)(0,7)
\put(3.2,10){$i_k$}
\end{picture}
\end{array}
\]
\caption{Action of the dynamics on a configuration in $\Id_n$. We
  analyse separately the cases of $\delta\geq 1$ (in this case we can
  single out one point in the top-right identity block) and $\delta=0$
  (in this case we can omit the top-right identity block). In all the
  four cases we obtain a configuration in one of the four forms
  considered in figure \ref{fig.struct_id}, and all the sizes of the
  identity blocks, except for a finite number of blocks on the
  top-left corner, are left unchanged, this allowing for the analysis
  of a generic configuration.\label{fig.RLactId}}
\end{figure}

In order to see why this correspondence holds, we have to verify that
the action of $R$ and $L$ on these configurations is consistent with
this straight representation.  The relations implied by the straight
representation as a binary tree are
\begin{subequations}
\begin{align}
L\, \s_{RL}^{n;(i_1,j_1,\ldots,i_k,j_k)}
&=
\left\{
\begin{array}{ll}
\s_{LL}^{n;(i_1,j_1,\ldots,i_k,j_k,1)}
&
\delta \geq 1
\\
\rule{0pt}{14pt}%
\s_{RL}^{n;(i_1,j_1,\ldots,i_k,j_k)}
&
\delta =0
\end{array}
\right.
\\
R\, \s_{RL}^{n;(i_1,j_1,\ldots,i_k,j_k)}
&=
\left\{
\begin{array}{ll}
\s_{RL}^{n;(i_1,j_1,\ldots,i_k,j_k+1)}
&
\delta \geq 1
\\
\rule{0pt}{14pt}%
\s_{LL}^{n;(i_1,j_1,\ldots,i_k)}
&
\delta =0
\end{array}
\right.
\end{align}
\end{subequations}
Similarly we have
\begin{subequations}
\label{eqs.verifths665475}
\begin{align}
L\, \s_{LL}^{n;(i_1,j_1,\ldots,i_k)}
&=
\left\{
\begin{array}{ll}
\s_{LL}^{n;(i_1,j_1,\ldots,i_k+1)}
&
\delta \geq 1
\\
\rule{0pt}{14pt}%
\s_{RL}^{n;(i_1,j_1,\ldots,j_{k-1})}
&
\delta =0
\end{array}
\right.
\\
R\, \s_{LL}^{n;(i_1,j_1,\ldots,i_k)}
&=
\left\{
\begin{array}{ll}
\s_{RL}^{n;(i_1,j_1,\ldots,i_k,1)}
&
\delta \geq 1
\\
\rule{0pt}{14pt}%
\s_{LL}^{n;(i_1,j_1,\ldots,i_k)}
&
\delta =0
\end{array}
\right.
\end{align}
\end{subequations}
(the other two cases are treated analogously).  Recalling that the
graphical description of the dynamics, in matrix representation, is
given by figure \ref{fig.defDyn3}, we find that these equations are
verified, as is illustrated in Figure~\ref{fig.RLactId} for equations
(\ref{eqs.verifths665475}) (the other cases are treated analogously)

\begin{figure}[b!!]
\[
\includegraphics[scale=.52, angle=0]{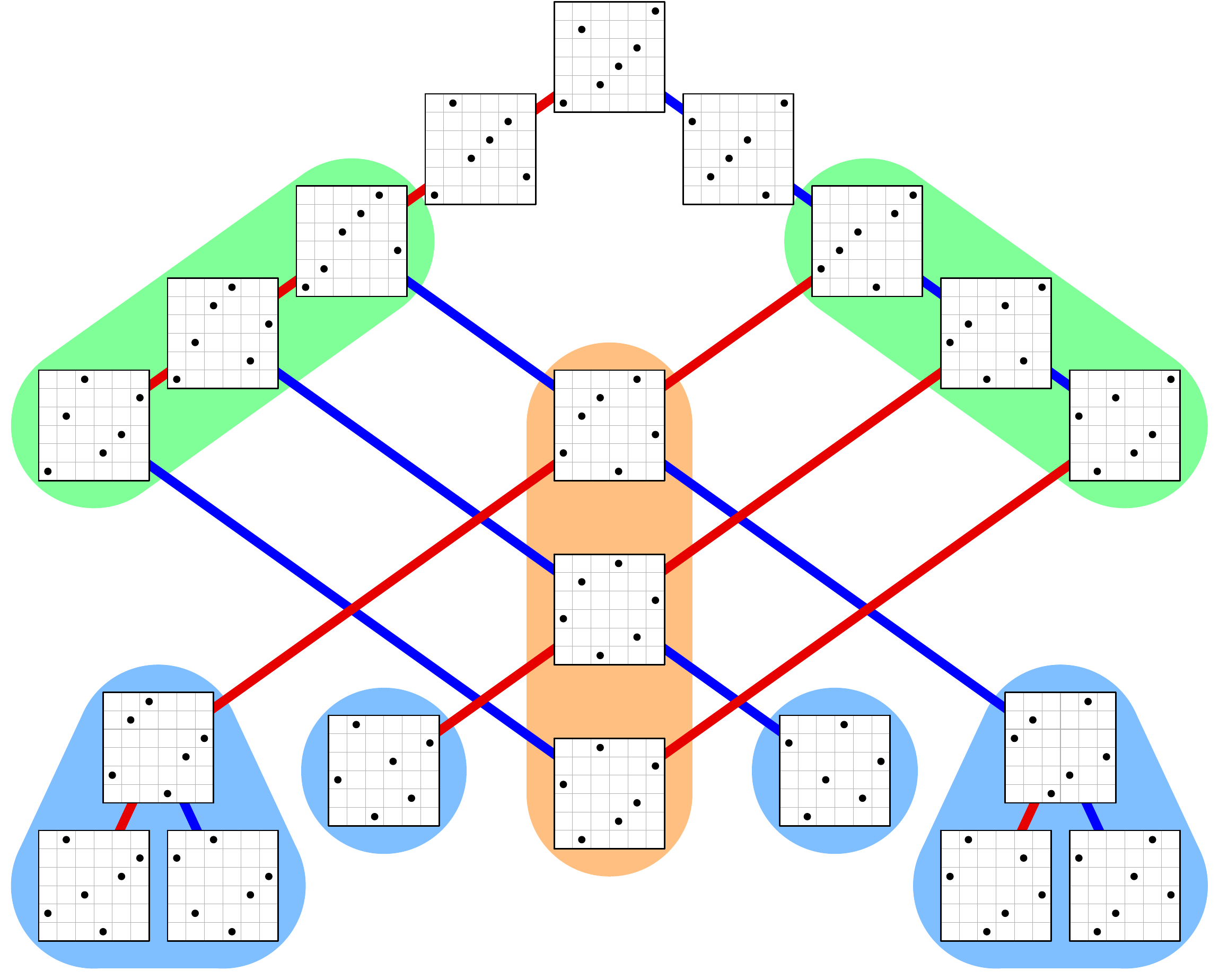}
\]
\caption{\label{fig.idp6}The Cayley graph of $\Id'_6$, with
  configurations in matrix representation. Red and blue edges
  correspond to actions of $L$ and $R$ respectively, and are omitted
  with the same rules of straight representations. Green, orange and
  cyan blocks of configurations correspond to configurations with
  different structure, illustrated in the text. In particular, the
  cyan blocks are isomorphic to classes $\Id_k$, except for the action
  of $L^{-1}$ (or $R^{-1}$, depending on the position w.r.t.\ the
  vertical axis) on the root of the binary tree, and the action of $L$
  (resp.\ $R$) on the left-most leaf of the tree (resp.\ right-most).}
\end{figure}

The characterisation of $\Id_n$ has a couple of interesting
consequences. A first one is the fact that $\Id_n$ has a unique
standard family, containing $\id_n$ (this comes from the direct
inspection of figure \ref{fig.struct_id}). The second one is the
following

\begin{corollary}
\label{cor.lineIde}
We can decide in linear time if $\s \in \Id_n$.
\end{corollary}

\begin{pf}
It is easily seen that we can check in linear time if $\s$ has one of
the four structures of figure \ref{fig.struct_id}. In fact, for each
of the two sets $\id \cup X_{RL} \cup X_{LL}$ and $\id \cup X_{LR}
\cup X_{RR}$, the test takes $\sim 2n$ accesses to the matrix, if
successful, and less than this if unsuccessful.
\end{pf}

Another way of seeing this is to use a linear-time standardisation
algorithm, and the uniqueness of the standard family. However, the
na\"ive standardisation algorithm is quadratic in time. A linear
algorithm exists, but we postpone its description to future work.

\subsection{Classes $\tree_n$}
\label{ssec.appidp}

This class is somewhat more complicated than $\Id_n$, however, within
its exponential cardinality, the vast majority of configurations and
transitions follow the same basic mechanism of $\Id_n$ (and are
arranged into complete binary trees), the new ingredients being
confined to a number, linear in $n$, of special configurations, with a
simple structure (i.e., at generic $n$, described by a finite number
of bullets and identity blocks).

We use again a straight representation of the Cayley graph. This could
be represented in a planar way (with some edges being very stretched),
but we will instead use a non-planar representation, illustrated in
Figure \ref{fig.idp6} for the case of $\Id'_6$.

\begin{figure}[tb!]
\[
\setlength{\unitlength}{5.2pt}
\begin{picture}(66.2,53.2)(-0.1,-0.1)
\put(0,0){\includegraphics[scale=.52,angle=0]{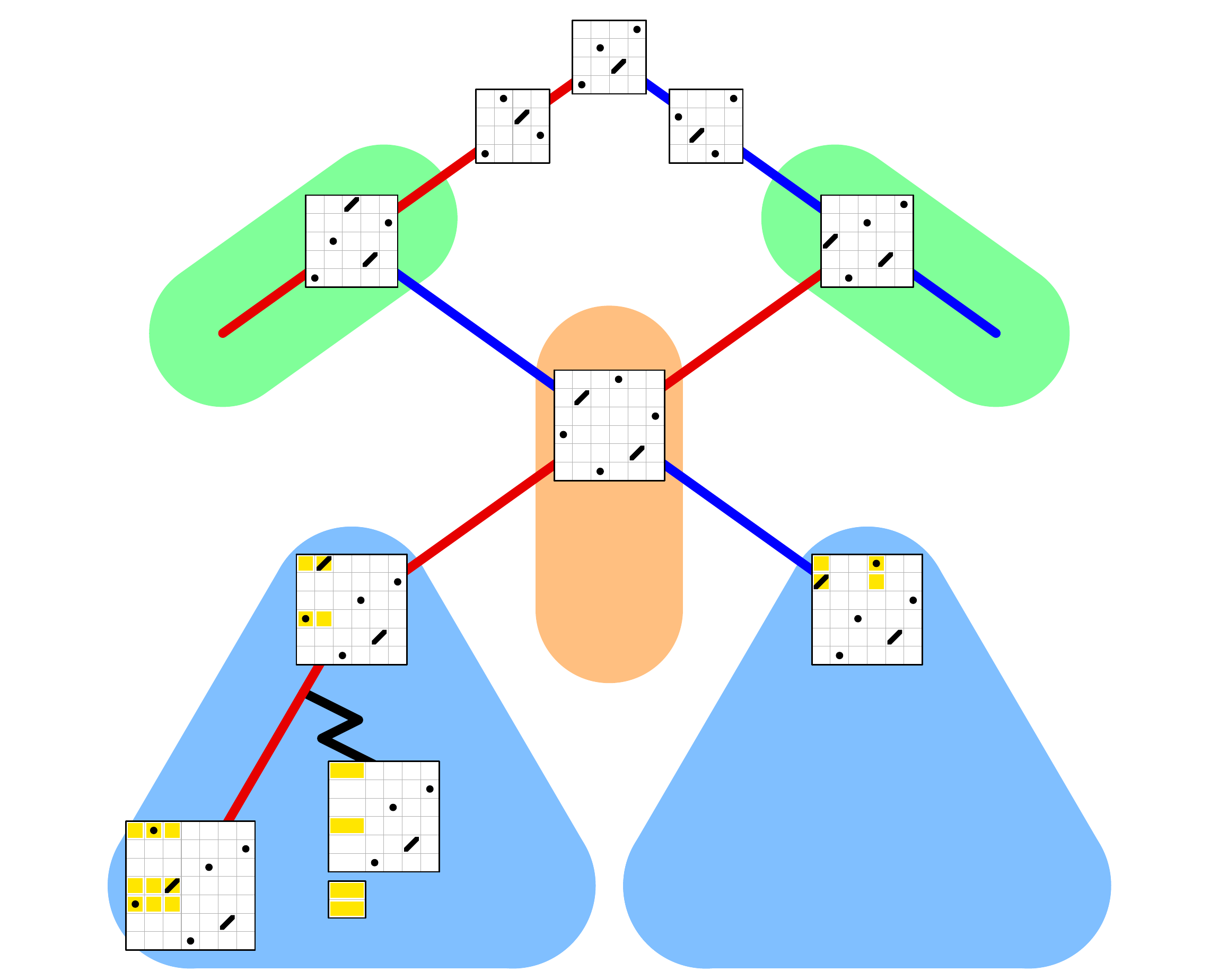}}

\put(36,52){$\id'_n$}
\put(34.5,49.5){{\scriptsize $\bm{n\!-\!3}$}}
\put(29.3,46.7){{\scriptsize $\bm{n\!-\!3}$}}
\put(15,45){{\scriptsize $i:$}}
\put(19,45){{\scriptsize $1$}}
\put(15.5,42.5){{\scriptsize $2$}}
\put(7.5,37.5){{\scriptsize $n\!-\!3$}}
\put(11,40){{\scriptsize \rotatebox{35.5}{$\cdots$}}}
\put(21,39){{\scriptsize $\bm{i\!-\!1}$}}
\put(20,42){{\scriptsize $\bm{n\!-\!i\!-\!2}$}}
\put(35.5,28.5){{\scriptsize $\bm{i\!-\!1}$}}
\put(32.5,31.5){{\scriptsize $\bm{n\!-\!i\!-\!3}$}}
\put(21.5,18.5){{\scriptsize $\bm{i\!-\!1}$}}
\put(18.25,22.25){{\scriptsize $\bm{n\!-\!i\!-\!3}$}}

\put(17.5,16){\rotatebox{-35.5}{$\updownarrow$}\raisebox{-2pt}{$j$}}
\put(20,13){$S$}
\put(23.3,6.5){{\scriptsize $\bm{i\!-\!1}$}}
\put(24.3,8.2){{\scriptsize $\bm{\}\; j\!+\!1}$}}

\put(20.7,3.7){$=S L^j \id_{n-i-2}$}


\end{picture}
\]
\caption{\label{fig.idpN}A schematisation of the Cayley graph of $\Id'_n$, with
  configurations in matrix representation. Only a few branches are
  represented (the full structure is evinced from comparison with
  Figure \ref{fig.idp6}). The configuration on the bottom-left corner
  has been added for clarity, in order to illustrate its `anomalous
  behaviour' under the application of $L$, w.r.t.\ the copy of
  $\Id_{n-i-2}$ in which it is contained.}
\end{figure}

Figure \ref{fig.idpN} illustrates the general structure of the class
$\Id'_n$. We observe four families of configurations:
\begin{itemize}
\item
The three configurations $\id'$, $L\, \id'$ and $R\, \id'$, which are
on the top part of our straight representation of the Cayley graph.
\item
Two linear families of configurations, related one another by the symmetry of
reflection along the diagonal, and denoted in green in the figures
\ref{fig.idp6} and \ref{fig.idpN}. An index $1 \leq i \leq n-3$ is
associated to the size of one of the two diagonal blocks.
\item
One linear family of configurations, which are symmetric,
and denoted in orange in the figures.
Again, an index $1 \leq i \leq n-3$ is
associated to the size of one of the two diagonal blocks.
\item
Two linear families, related one another by the symmetry, of subgraphs
of the Cayley graph, isomorphic to the graph of $\Id_{n-i-2}$ except
for one transition (interpret $\Id_0$ as an empty graph).  These are
denoted by the cyan triangles in the figures.  In figure
\ref{fig.idpN}, the restriction to the yellow blocks of these
configurations coincides with the corresponding configurations in
$\Id_{n-i-2}$.
\end{itemize}
At the light of the results of the previous section on the structure
of $\Id$, it is easy to verify that this Cayley graph indeed describes
the dynamics on this class.




\begin{corollary}
$\tree_n$ has size $2^{n-2}+n-2$.
\end{corollary}

\begin{pf}
We have the three configurations $\id'$, $L\, \id'$ and $R\, \id'$,
then $2(n-3)$ `green' configurations (w.r.t.\ the colours of Figure
\ref{fig.idp6}), $n-3$ `orange' ones, and $2\sum_{i=1}^{n-3}
(2^{i-1}-1)=2^{n-2}-2n+4$ for the `cyan' ones. Collecting all the
summands gives the statement.
\end{pf}

\begin{corollary}
\label{cor.lineIdeP}
We can decide in linear time if $\s \in \tree_n$.
\end{corollary}

\begin{pf}
All configurations in $\Id'_n$ have either a structure described by a
finite number of bullets and identity blocks, within the finite list
of Figure \ref{fig.idpN}, or a structure of this form, plus a block
which is in the class $\Id_k$ for some $k$. With an argument similar
to the one of Corollary \ref{cor.lineIde}, and in light of the latter,
we can thus conclude.
\end{pf}


\subsection{Synthetic presentation of classes $\Id_n$ and $\tree_n$}

Instead of using the explicit descriptions of figures
\ref{fig.struct_id} and \ref{fig.idpN}, we can describe the sets $\Id$
and $\Id'$ in a synthetic recursive way. To this aim we need a few
definitions.

\begin{definition}
A configuration $\s \in \Id_n$ is \emph{of type $(L,j)$}
if $\s(n)=j>\s^{-1}(1)$, and it is \emph{of type $(R,j)$} if
$\s(n)<j=\s^{-1}(1)$. We denote by $\Id_n^{(L,j)}$ and $\Id_n^{(R,j)}$
the corresponding sets.
\end{definition}

\begin{proposition}
The configurations in $\Id_n$ are the disjoint union of $\id_n$,
configurations of type $(L,j)$ for $2 \leq j \leq n-1$, and
configurations of type $(R,j)$ for $2 \leq j \leq n-1$.  The set
$\Id_n^{(R,j)}$ consists of the transpose of matrices in set
$\Id_n^{(L,j)}$.
\end{proposition}
\noindent
For compactness of the following expressions, we do the identification $\id
\equiv \Id_n^{(L,1)} \equiv \Id_n^{(R,1)}$, and we introduce
$\Id_n(L)=\cup_{j \geq 1} \Id_n^{(L,1)}$, and similarly for
$\Id_n(R)$.


Define also
recursively the family of square matrices divided into two rectangular blocks
\raisebox{-9.5pt}{\includegraphics[width=8.8mm]{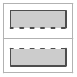}}
as
\be
\raisebox{-9.5pt}{\includegraphics[width=8.8mm]{FigureA2_fig_matr_MMl.pdf}}
=
\raisebox{-9.5pt}{\includegraphics[width=8.8mm]{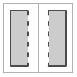}}
=
\epsilon
\cup
\raisebox{-20.5pt}{\includegraphics[width=18.4mm]{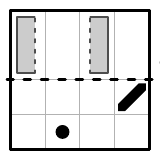}}
\ee
Then the set $\Id$ defined as
\begin{align}
\Id &= \Id(L) \cup \Id(R)
\ef;
&
\Id(L) \cap \Id(R)
&= \id
\ef;
\end{align}
\begin{align}
\Id(L) &:=
\raisebox{-5mm}{\includegraphics[width=12.8mm]{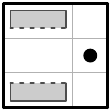}}
\ef;
&
\Id(R) &:=
\raisebox{-5mm}{\includegraphics[width=12.8mm]{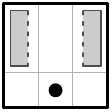}}
\ef;
\end{align}
consists of configurations of size $\leq 3$, which are not primitive,
and the union of $\Id_n$ for $n \geq 4$.

Similarly, the set $\tree$ defined as follows
\begin{align}
&\tree = \tree(L) \cup \tree(R) \cup \tree(T)
\ef;
&
&\tree(R) = \big( \tree(L) \big)^T
\ef;
\\
&(\tree(L) \cup \tree(R)) \cap \tree(T) = \emptyset
\ef;
&
&\tree(L) \cap \tree(R) = \id'
\ef;
\end{align}
\begin{align}
\tree(T) &=
\raisebox{-11mm}{\includegraphics[width=24.8mm]{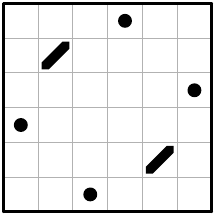}}
\ef;
&
\tree(L) &=
\raisebox{-7mm}{\includegraphics[width=16.8mm]{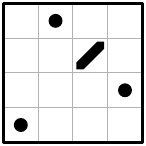}}
\;\cup\;
\raisebox{-9mm}{\includegraphics[width=20.8mm]{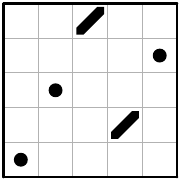}}
\;\cup\;
\raisebox{-11mm}{\includegraphics[width=24.8mm]{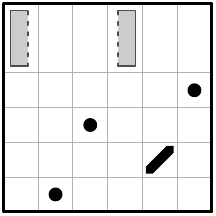}}
\ef;
\end{align}
consists of configurations of size $3$ and $4$, which are not
primitive, and the union of $\tree_n$ for $n \geq 5$.

Note that, with respect to the notations of Section \ref{ssec.appid},
we have $\Id(L) = \id \cup X_{RL} \cup X_{LL}$ and $\Id(R) = \id \cup
X_{LR} \cup X_{RR}$, and with respect to the notations of Section
\ref{ssec.appidp} we have that $\tree(T)$ contains the orange
configurations, and $\tree(L)$, $\tree(R)$ both contain $\id'$, and
for the rest contain each half of the configurations which are not
symmetric w.r.t.\ reflection along the diagonal.


\end{document}